\documentclass[12pt]{report}
\usepackage{fullpage}
\usepackage{longtable}
\usepackage{graphicx}

\usepackage{amsmath, amssymb, amscd,eucal}
\input epsf.tex
\newtheorem{thm}{Theorem}[section]

\newtheorem{cor}[thm]{Corollary}

\newtheorem{lem}[thm]{Lemma}
\newtheorem{prop}[thm]{Proposition}
\newtheorem{con}[thm]{Conjecture}

\newenvironment{proof}[1][Proof]{\begin{trivlist}
\item[\hskip \labelsep {\bfseries #1}]}{\end{trivlist}}
\def\qed{\quad \vrule height7.5pt width4.17pt depth0pt} \def\ex{\text{ex}}

\def\[{\left(} \def\]{\right)}

\input epsf.tex

\begin{document}

\title{
\vspace*{-24pt} Random Surfaces}
\author{Scott Sheffield}
\date{}
\maketitle

\begin{abstract}
We study the statistical physical properties of (discretized) ``random surfaces,'' which are random
functions from $\mathbb Z^d$ (or large subsets of $\mathbb Z^d$) to
$E$, where $E$ is $\mathbb Z$ or $\mathbb R$.  Their laws are
determined by convex, nearest-neighbor, gradient Gibbs potentials
that are invariant under translation by a full-rank sublattice
$\mathcal L$ of $\mathbb Z^d$; they include many discrete and
continuous height function models (e.g., domino tilings, square ice,
the harmonic crystal, the Ginzburg-Landau $\nabla \phi$ interface
model, the linear solid-on-solid model) as special cases.

We prove a {\it variational principle}---characterizing gradient
phases of a given slope as minimizers of the specific free
energy---and an empirical measure large deviations principle (with a
unique rate function minimizer) for random surfaces on mesh
approximations of bounded domains. We also prove that the surface
tension is strictly convex and that if $u$ is in the interior of the
space of finite-surface-tension slopes, then there exists a minimal
energy gradient phase $\mu_u$ of slope $u$.

Using a new geometric technique called {\em cluster swapping} (a
variant of the Swendsen-Wang update for Fortuin-Kasteleyn clusters),
we show that $\mu_u$ is unique if at least one of the following
holds: $E = \mathbb R$, $d \in \{1, 2 \}$, there exists a {\it rough
gradient phase} of slope $u$, or $u$ is irrational. When $d=2$ and
$E = \mathbb Z$, we show that the slopes of all {\it smooth phases}
(a.k.a. {\it crystal facets}) lie in the dual lattice of $\mathcal L$.

In the case $E = \mathbb Z$ and $d=2$, our results resolve and
greatly generalize a number of conjectures of Cohn, Elkies, and
Propp---one of which is that there is a unique ergodic Gibbs measure
on domino tilings for each non-extremal slope.  We also prove
several theorems cited by Kenyon, Okounkov, and Sheffield in their
recent exact solution of the dimer model on general planar lattices.
In the case $E = \mathbb R$, our results generalize and extend many
of the results in the literature on Ginzurg-Landau $\nabla
\phi$-interface models.
\end{abstract}

\begin{center} {\bf Acknowledgements} \end{center}

This work is a revised version of a Ph.D. thesis advised by Amir
Dembo, whom I thank for tremendous assistance.  We spent many
hours---sometimes entire days---fleshing out ideas, wading through
chapter drafts, and talking through conceptual challenges.  He is
devoted to his students, has a wonderful sense of humor, and is
truly a pleasure to work with.

Thanks also to Persi Diaconis, Michael Harrison, Jon Mattingly, Rafe
Mazzeo, and George Papanicolaou for service on oral exam, defense,
and reading committees.  Particular thanks to Persi Diaconis for
extensive advice about the problems in this text and for supervising
my first readings in Gibbs measures and statistical physics.

For summer and post-doctoral support and many valuable discussions,
I thank the members of the Theory Group at Microsoft Research,
particularly Henry Cohn, David Wilson, Oded Schramm, Jennifer
Chayes, and Christian Borgs.

Thanks to Richard Kenyon and Andrei Okounkov for their advice and
collaboration on the subject of random surfaces derived from perfect
matchings in the plane (the {\it dimer model}).  Although the work
with Kenyon and Okounkov cites my thesis, the two projects were
actually completed in tandem, and the dimer model intuition inspired
many of the results presented here.

Thanks also to Jean-Dominique Deuschel for his course on
Ginzurg-Landau $\nabla \phi$-interface models at Stanford, and to
both Stefan Adams and Jean-Dominique Deuschel for their kind
hospitality and extensive, extremely helpful discussions during my
visit to Berlin.

Thanks to Jamal Najim for finding and introducing me to Cianchi's
modern treatment of Orlicz-Sobolev spaces.  Applying these beautiful
results to random surfaces was truly a pleasure.

Thanks to Alan Hammond for many relevant conversations, and to a
kind friend, Michael Shirts, for managing my thesis submission at
Stanford while I was in Paris.

Finally, thanks to my family, especially my wife, Julie, for
unwavering patience, love, and support.

\tableofcontents

\chapter{Introduction}\label{introchapter}

The following is a fundamental problem of variational calculus:
given a bounded open subset $D$ of $\mathbb R^d$ and a free energy
function $\sigma: D \times \mathbb R^d \times \mathbb R^{d \times m}
\mapsto \mathbb R$, find the differentiable function $f:D \mapsto
\mathbb R^m$ that (possibly subject to boundary conditions)
minimizes the free energy integral: $$\int_{D} \sigma(x,f(x), \nabla
f(x))dx.$$

Since the seventeenth century, these free-energy-minimizing
functions have been popular models for determining (among other
things) the shapes assumed by solid objects in the presence of
outside forces: ropes suspended between poles, elastic sheets
stretched to boundary conditions, and twisted or otherwise strained
three-dimensional solids.  They are also useful in modeling surfaces
of water droplets and other phase interfaces.  Rigorous formulations
and solutions to these problems rank among the great achievements of
classical analysis (including work by Fermat, Newton, Leibniz, the
Bernoullis, Euler, Lagrange, Legendre, Jacobi, Hamilton,
Weierstrass, etc. \cite{Go}) and play important roles in physics and
engineering.

All of these models assume that matter is continuous and distributes
force in a continuous way. One of the goals of statistical physics
has become not merely to solve variational problems but to
understand and, in some sense, to {\it justify} them in light of the
fact that matter is comprised of individual, randomly behaving
atoms.  To this end, one begins by postulating a simple form for the
local particle interactions: one approach---the one we will study in
this work---is to represent the ``atoms'' of the solid crystal by
points in a subset $\Lambda$ of $\mathbb Z^d$, each of which has a
``spatial position'' given by a function $\phi: \Lambda \mapsto
\mathbb R^m$, and to specify the interaction between the particles
by a {\it Gibbs potential} $\Phi$ that possesses certain natural
symmetries. The next step is to show that---at least in some
``thermodynamic limit''---a random Gibbs configuration will
approximate a free-energy-minimizing function like the ones
described above.

Another problem, which has no analog in the deterministic, non-atomic classical theory, is the
investigation of local statistics of a physical system.  How likely are particular microscopic
configurations of atoms to occur as sub-configurations of a larger system? How are these
occurrences distributed?  To what extent is matter {\it homogenous} throughout small but
non-microscopic regions?  Our solutions to these problems will involve {\it large deviations
principles}, which we precisely define later on.

Finally, we want to investigate more directly the connections between the Gibbs potential $\Phi$
and the kinds of behavior that can occur in these small but non-microscopic regions.  This will
require us to ask, given $\Phi$, what are the ``gradient phases'' (i.e., the ergodic gradient Gibbs
measures with finite specific energy) $\mu$ of a given slope?  Does the $\mu$-variance of the
height difference of points $n$ units apart remain bounded independently of $n$ or does it tend to
infinity with $n$?  When is the {\it surface tension} function $\sigma$ (defined precisely in
Chapter \ref{surfacetensionchapter}) strictly convex?

Before we state our results precisely and describe some of the
previous work in this area, we will need several definitions. While
we attempt to make our exposition relatively self-contained---and
define the terms we use precisely---we will also draw heavily from
the results in some standard texts: {\it Sobolev Spaces} by Adams
\cite{A} and recent extensions by Cianchi (\cite{C}, \cite{C1},
\cite{C2}); {\it Large Deviations Techniques and Applications} by
Dembo and Zeitouni \cite{DZ}; {\it Large Deviations} by Deuschel and
Stroock \cite{DS}; and {\it Gibbs Measures and Phase Transitions} by
Georgii \cite{G}. We will carefully state, if not prove, the outside
theorems we use.

\section{Random surfaces and gradient Gibbs measures} \label{intrograd}
\subsection{Gradient potentials}
The study of random functions $\phi$ from the lattice $\mathbb Z^d$
to a measure space $(E, \mathcal E)$ is a central component of
ergodic theory and statistical physics.  In many classical models
from physics (e.g., the Ising model, the Potts model, Shlosman's
plane rotor model), $E$ is a space with a finite underlying measure
$\lambda$, $\mathcal E$ is the Borel $\sigma$-field of $E$, and
$\phi(x)$ has a physical interpretation as the {\it spin} (or some
other internal property) of a particle at location $x$ in a crystal
lattice. (See e.g., \cite{G}.) In the models of interest to us, $(E,
\mathcal E)$ is a space with an infinite underlying measure
$\lambda$---either $\mathbb R^m$ with Lebesgue measure or $\mathbb
Z^m$ with counting measure---where $\mathcal E$ is the Borel
$\sigma$-algebra of $E$ and $\phi(x)$ usually has a physical
interpretation as the {\it spatial position} of a particle (or the
vertical height of a phase interface) at location $x$ in a lattice.
For example, if $m=d=3$, $\phi$ could describe the spatial positions
of the components of an elastic crystal; if $m=1$ and $d=2$, $\phi$
could describe the solid-on-solid or Ginzburg-Landau approximations
of a phase interface \cite{FS}.

Throughout the exposition, we denote by $\Omega$ the set of
functions from $\mathbb Z^d$ to $E$ and by $\mathcal F$ the Borel
$\sigma$-algebra of the product topology on $\Omega$.  If $\Lambda
\subset \mathbb Z^d$, we denote by $\mathcal F_{\Lambda}$ the
smallest $\sigma$-algebra with respect to which $\phi(x)$ is
measurable for all $x \in \Lambda$.  We write $\mathcal T_{\Lambda}
= \mathcal F_{\mathbb Z^d - \Lambda}$.  We write $\Lambda \subset
\subset \mathbb Z^d$ if $\Lambda$ is a finite subset of $\mathbb
Z^d$. A subset of $\Omega$ is called a {\it cylinder set} if it
belongs to $\mathcal F_{\Lambda}$ for some $\Lambda \subset \subset
\mathbb Z^d$. Let $\mathcal F$ be the smallest $\sigma$-algebra on
$\Omega$ containing the cylinder sets.  We write $\mathcal T$ for
the intersection of $\mathcal T_{\Lambda}$ over all finite subsets
$\Lambda$ of $\mathbb Z^d$; the sets in $\mathcal T$ are called {\it
tail-measurable sets}.

We will also always assume that we are given a family $\Phi$ of
measurable {\it potential functions} $\Phi_{\Lambda}: \Omega \mapsto
\mathbb R \cup \{ \infty \}$ (one for each finite subset $\Lambda$
of $\mathbb Z^d$); each $\Phi_{\Lambda}$ is $\mathcal F_{\Lambda}$
measurable.  We will further assume that $\Phi$ is invariant under
the group $\Theta$ of translations of $\mathbb Z^d$ by members of
some rank-$d$ lattice $\mathcal L$ --- i.e., if $s \in \mathcal L$,
then $\Phi_{\Lambda+s}(\phi_s) = \Phi_{\Lambda}(\phi)$, where
$\phi_s$ is defined by $\phi_s(i) = \phi(i-s)$. (In many
applications, we can take $\mathcal L = \mathbb Z^d$.)  We also
assume that $\Phi$ is invariant under a group $\tau$ of
measure-preserving translations of $E$ --- i.e.,
$\Phi_{\Lambda}(\phi) = \Phi_{\Lambda}(\tau\phi)$, where $\tau \phi$
is simply defined by $(\tau \phi) (x) = \tau(\phi(x))$. Potentials
$\Phi$ satisfying the above requirements are called {\it $\Theta
\times \tau$-invariant potentials} or {\it $\mathcal L \times
\tau$-invariant potentials}. For all of our main results, we will
assume that $\tau$ is the full group of translations of $\mathbb
Z^m$ or $\mathbb R^m$; in this case, each $\Phi_{\Lambda}(\phi)$ is
a function of the {\it gradient} of $\phi$, written $\nabla \phi$
and defined by $$\nabla \phi(x) = (\phi(x+e_1) - \phi(x),
\phi(x+e_2) - \phi(x), \ldots, \phi(x + e_d) - \phi(x)),$$ where
$e_i$ are the standard basis vectors of $\mathbb Z^d$.  In this
setting, we will refer to $\mathcal L \times \tau$-invariant
potentials as {\it $\mathcal L$-periodic} or {\it $\mathcal
L$-invariant} gradient potentials.  We use the term {\it
shift-invariant} to mean $\mathcal L$-invariant when $\mathcal L =
\mathbb Z^d$.  In some of our applications, we also restrict our
attention to {\it nearest-neighbor potentials}, i.e., those
potentials $\Phi$ for which $\Phi_{\Lambda} = 0$ unless $\Lambda$ is
a single pair of adjacent vertices in $\mathbb Z^d$.  We say that
$\Phi$ has {\it finite range} if there exists an $r$ such that
$\Phi_{\Lambda} = 0$ whenever the diameter of $\Lambda$ is greater
than $r$. For each finite subset $\Lambda$ of $\mathbb Z^d$ we also
define a {\it Hamiltonian}: $H_{\Lambda} (\phi) = \sum_{\Delta \cup
\Lambda \not = \emptyset} \Phi_{\Delta} (\phi)$, where the sum is
taken over finite subsets $\Delta$ of $\mathbb Z^d$.

We define the {\it interior Hamiltonian} of $\Lambda$, written
$H^o_{\Lambda}(\phi)$, to be: $$H^o_{\Lambda}(\phi) = \sum_{\Delta
\subset \Lambda} \Phi_{\Delta}(\phi).$$ This is different from
$H_{\Lambda}$ because the former sum includes sets $\Delta$ that
intersect $\Lambda$ but are not strictly contained in $\Lambda$.  On
the other hand, $H^o_{\Lambda}$ is $\mathcal F_{\Lambda}^\tau$
measurable, which is not true of $H_{\Lambda}$.  (This
$H^o_{\Lambda}$ is sometimes called the {\it free boundary}
Hamiltonian for $\Lambda$.)

\subsection{Gibbs Measures} \label{gibbsdefinition}
To define Gibbs measures and gradient Gibbs measures, we will need
some additional notation.  Let $(X, \mathcal X)$ and $(Y, \mathcal
Y)$ be general measure spaces. A function $\pi: \mathcal X \times Y
\mapsto [0, \infty]$ is called a {\it probability kernel} from $(Y,
\mathcal Y)$ to $(X, \mathcal X)$ if
\begin{enumerate}
\item $\pi(\cdot|y)$ is a probability measure on $(X, \mathcal X)$
for each fixed $y \in Y$, and \item $\pi(A|\cdot)$ is $\mathcal
Y$-measurable for each fixed $A \in \mathcal X$. \end{enumerate}

Since a probability kernel maps each point in $Y$ to a probability
measure on $X$, we may interpret a probability kernel as giving the
law for a random transition from an arbitrary point in $Y$ to a
point in $X$.  A probability kernel maps each measure $\mu$ on $(Y,
\mathcal Y)$ to a measure $\mu \pi$ on $(X, \mathcal X)$ by
$$\mu\pi(A) = \int \pi(A|\cdot)d\mu.$$ The following is a probability
kernel from $(\Omega, \mathcal T_{\Lambda})$ to $(\Omega, \mathcal
F)$; in particular, for any fixed $A \in \mathcal F$, it is a
$\mathcal T_{\Lambda}$-measurable function of $\phi$:
$$\gamma^{\Phi}_{\Lambda}(A, \phi) = Z_{\Lambda}(\phi)^{-1} \int \prod_{x \in \Lambda} d \phi(x)
\exp\[-H_{\Lambda}(\phi)\] 1_A(\phi).$$ (When the choice of
potential $\Phi$ is clear from context, we write
$\gamma^{\Phi}_{\Lambda}$ as $\gamma_{\Lambda}$.) In this
expression, $Z_{\Lambda}(\phi)$ (which is also $\mathcal
T_{\Lambda}$ measurable) is defined as follows:
$$Z_{\Lambda}(\phi) = \int \prod_{x \in \Lambda} d \phi(x)
\exp\[-H_{\Lambda}(\phi)\],$$ where $d\phi(x)$ is the underlying
(Lebesgue or counting) measure on $E$.  Informally, $\gamma_\Lambda$
is a random transition from $\Omega$ to itself that takes a function
$\phi$ and then ``rerandomizes'' $\phi$ within the set $\Lambda$.

We say $\phi$ has {\it finite energy} if $\Phi_{\Lambda}(\phi) <
\infty$ for all $\Lambda \subset \subset \mathbb Z^d$.  We say
$\phi$ is $\Phi$-{\it admissible} if each $Z_{\Lambda}(\phi)$ is
finite and non-zero.  Given a measure $\mu$ on $(\Omega, \mathcal
F)$, we define a new measure $\mu \gamma_{\Lambda}$ by $$\mu
\gamma_{\Lambda}(A) = \int \gamma_{\Lambda}(A|\cdot) d \mu.$$  We
say a probability measure $\mu$ on $(\Omega, \mathcal F)$ is a {\it
Gibbs measure} if $\mu$ is supported on the set of $\Phi$-admissible
functions in $\Omega$ and for all finite subsets $\Lambda$ of
$\mathbb Z^d$, we have $\mu \gamma_{\Lambda} = \mu$.  (In other
words, $\mu$ is Gibbs if and only if $\gamma_{\Lambda}$ describes a
regular conditional probability distribution, where the
$\mu$-conditional distribution of the values of $\phi(x)$ for $x \in
\Lambda$ is given by $\gamma_{\Lambda}(\cdot| \phi)$.)

A fundamental result in Gibbs measure theory is that for any $\Phi$,
the set of $\Theta$-invariant Gibbs measures is convex and its
extreme points are $\Theta$-ergodic. (See, e.g., Chapters 14 of
\cite{G}.  More details also appear in Chapter
\ref{decompositionchapter} of this text.) Since $\Theta$ is
understood to be the group of translations by a sublattice $\mathcal
L$ of $\mathbb Z^d$, we will also use the terms {\it $\mathcal
L$-invariant} and {\it $\mathcal L$-ergodic}. In physics jargon, the
$\mathcal L$-ergodic measures are the {\it pure phases} and a {\it
phase transition} occurs at potentials $\Phi$ which admit more than
one $\mathcal L$-ergodic Gibbs measure.

\subsection{Gradient Gibbs Measures} Let $\tau$ be the group of translations of $E$, and let
$\mathcal F^{\tau}$ be the $\sigma$-algebra containing
$\tau$-invariant sets of $\mathcal F$; this is the smallest
$\sigma$-algebra containing the sets of the form $\{ \phi | \phi(y)
- \phi(x) \in \mathcal S \}$ where $x,y \in \mathbb Z^d$ and
$\mathcal S \in \mathcal E$.  In other words, $\mathcal F^{\tau}$ is
the subset of $\mathcal F$ containing those sets that are invariant
under translations $ \phi \mapsto \phi + z$ for $z \in E$.
(Similarly, we write $\mathcal T_{\Lambda}^{\tau} = \mathcal
T_{\Lambda} \cap \mathcal F^{\tau}$ and $\mathcal F_{\Lambda}^\tau =
\mathcal F_{\Lambda} \cap \mathcal F^{\tau}$.)  Let $\Phi$ be an
$\mathcal L$-invariant gradient potential.  Since, given any $A \in
\mathcal F^{\tau}$, the kernels $\gamma^{\Phi}_{\Lambda}(A|\phi)$
are $\mathcal F^{\tau}$-measurable functions of $\phi$, it follows
that the kernel sends a given measure $\mu$ on $(\Omega, \mathcal
F^{\tau})$ to another measure $\mu \gamma^{\Phi}_{\Lambda}$ on
$(\Omega, \mathcal F^{\tau})$. A measure $\mu$ on $(\Omega, \mathcal
F^{\tau})$ is called a {\it gradient Gibbs measure} if it is
supported on admissible functions and $\mu \gamma^{\Phi}_{\Lambda} =
\mu$ for every $\mu$.  Note that this is the same as the definition
of Gibbs measure except that in this case the $\sigma$-algebra that
is different.

Clearly, if $\mu$ is a Gibbs measure on $(\Omega, \mathcal F)$, then its restriction to $\mathcal
F^{\tau}$ is a gradient Gibbs measure.  A gradient Gibbs measure is said to be {\it localized} or
{\it smooth} if it arises as the restriction of a Gibbs measure in this way.  Otherwise, it is {\it
non-localized} or {\it rough}.  (Many natural Gibbs measures are rough when $d \in \{1,2\}$; for
example, all the ergodic gradient Gibbs measures of the continuous, nearest-neighbor Gaussian
models in these dimensions are rough---see, e.g., \cite{G}.)  Moreover, the restriction of $\mu$ to
$\mathcal F^{\tau}$ may be $\mathcal L$-invariant even when $\mu$ itself is not.

Denote by $\mathcal P_{\mathcal L}(\Omega, \mathcal F^{\tau})$ the
set of $\mathcal L$-invariant probability measures on $(\Omega,
\mathcal F^{\tau})$ and by $\mathcal G_{\mathcal L}(\Omega, \mathcal
F^{\tau})$ the set of $\mathcal L$-invariant gradient Gibbs
measures. We say that a $\mu \in P_{\mathcal L}(\Omega, \mathcal
F^{\tau})$ has {\it finite slope} if $\mu(\phi(y) - \phi(x))$ is
finite for all pairs $x,y \in \mathbb Z^d$.  (Throughout the text,
we use the notation $\mu(f) = \int_{\Omega} f(\phi)d\mu(\phi)$.) One
easily checks that there is a unique $m \times d$ matrix $u$ such
that $\mu(\phi(y) - \phi(x)) = u(y-x)$ (where $u(y-x)$ denotes the
matrix product of $u$ and $(y-x)$) whenever $y - x \in \mathcal L$.
In this case, we call $u$ the {\it slope} of $\mu$, which we write
as $S(\mu)$.  Analogously to the non-gradient case, the extreme
points of $\mu \in \mathcal P_{\mathcal L}(\Omega, \mathcal
F^{\tau})$ are called {\it $\mathcal L$-ergodic gradient measures}
and the extreme points of $\mathcal G(\Omega, \mathcal F^{\tau})$
are called {\it extremal gradient Gibbs measures}.  We discuss these
notions in more detail in Chapter \ref{decompositionchapter}.
Although the term ``phase'' has many definitions in the physics
literature, when a full rank sublattice $\mathcal L$ of $\mathbb
Z^d$ is given, we will always use the term {\it gradient phase} to
mean an $\mathcal L$-ergodic gradient Gibbs measure with finite
specific free energy (a term we define precisely in Chapter
\ref{SFEchapter}).  A {\it minimal gradient phase} is a gradient
phase of some slope $u$ for which the specific free energy is
minimal among the set of all slope $u$, $\mathcal L$-invariant
gradient measures.

\subsection{Classes of periodic gradient potentials}

When $m=1$, we say a potential $\Phi$ is {\it simply attractive} if
$\Phi$ is an $\mathcal L$-invariant nearest-neighbor gradient
potential such that for each adjacent pair of vertices $x$ and $y$,
with $x$ preceding $y$ in the lexicographic ordering of $\mathbb
Z^d$, we have $\Phi_{\{x,y\}}(\phi) = V_{x,y}(\phi(y) - \phi(x))$,
where the $V_{x,y}:\mathbb R \mapsto [0, \infty]$ are convex and
positive, and $\lim_{\eta \mapsto \infty}V_{x,y}(\eta) = \lim_{\eta
\mapsto -\infty}V_{x,y}(\eta) = \infty$.  As before, we assume here
that $\mathcal L$ is a full-rank sublattice of $\mathbb Z^d$.  For
convenience, we will always assume $V_{x,y} = 0$ if $x$ and $y$ are
not adjacent or $x$ does not precede $y$ in the lexicographic
ordering of $\mathbb Z^d$.  When we refer to the nearest neighbor
potential for ``an adjacent pair $x,y$'' we will assume implicitly
that $x$ precedes $y$ in lexicographic ordering.

Note that each $V_{x,y}$ has a minimum at at least one point $\eta_0
\in \mathbb R$. In many applications, we can assume $\eta_0 = 0$; in
this case, the requirement that $V_{x,y}$ be convex implies that the
model is ``attractive'' or ``ferromagnetic'' in the sense that the
energy is lower when neighboring heights are close to one another
than when they are far apart.

We chose to invent the term ``simply attractive potential'' because
the obvious alternatives were either too long (``convex
nearest-neighbor periodic difference potential'') or too overloaded
and/or imprecise (``ferromagnetic potential,'' ``elastic
potential,'' ``anharmonic crystal potential,'' ``solid-on-solid
potential'').  Elsewhere in the literature, the latter terms have
definitions that are more general or more specific than ours,
although they usually agree in spirit.

Also, when $m=1$, we say $\Phi$ is {\it isotropic} if for some
$V:\mathbb R \mapsto \mathbb R$ (which must be positive, convex, and
even---i.e., $V(\eta) = V(-\eta)$) we have $V_{x,y}(\eta) = V(\eta)$
for all adjacent pairs $x,y \in \mathbb Z^d$.  We say $\Phi$ is {\it
Lipschitz} if there exist $\eta_1, \eta_2 \in \mathbb R$ such that
for all adjacent $x,y \in \mathbb Z^d$, we have
$V^{\Phi}_{x,y}(\eta) = \infty$ whenever $\eta<\eta_1$ or
$\eta>\eta_2$.  We will frequently use the following abbreviations:
\begin{enumerate}
\item {\it SAP}: Simply attractive potential
\item {\it ISAP}: Isotropic simply attractive potential.  We write $\Phi_V$ to denote the isotropic
simply attractive potential in which each $V_{x,y} = V$
\item {\it LSAP}: Lipschitz simply attractive potential
\end{enumerate}

Most of the simply attractive models discussed in the statistical
physics literature are either Lipschitz and have $E=\mathbb Z$
(e.g., height function models for perfect matchings on lattice
graphs \cite{K} and square ice \cite{Ba}) or isotropic (e.g., linear
solid-on-solid, Gaussian, and Ginzburg-Landau models, \cite{DGI}).

We say that a potential $\Phi$ {\it strictly dominates} a potential
$\Psi$ if there exists a constant $0 < c < 1$ such that
$|H^{\Psi}_{\Lambda}(\phi)| < c|H^{\Phi}_{\Lambda}(\phi)|$ for all
$\Lambda \subset \subset \mathbb Z^d$ and $\phi \in \Omega$.  (If
$m>1$, we replace the absolute value signs in this definition by the
Euclidean norm.)

When $m>1$, we can write any $\phi \in \Omega$ as $\phi_1 e^1 +
\phi_2 e^2 + \ldots + \phi_m e^m$ where the $\phi_i$ are real valued
(or integer valued) and the $e^i$ are the standard basis vectors in
$E$. In this setting, we say that $\Phi$ is an {\it SAP} (resp.,
{\it ISAP}, {\it LSAP}) if it can be written as $\Phi(\phi) = \sum
\Phi_i (\phi_i)$, where each of the $\Phi_i$ is a one-dimensional
simply attractive potential.  For any $m$, a {\it perturbed SAP}
(resp., {\it perturbed LSAP}, {\it perturbed ISAP}) is an $\mathcal
L$-periodic gradient potential of the form $\Phi + \Psi$ where
$\Phi$ is an SAP (resp., LSAP, ISAP), $\Psi$ has finite range, and
$\Phi$ strictly dominates $\Psi$.

Note that when $m>1$, our class of simply attractive potentials is
rather restrictive; each one can be decomposed into a sum of $m$
simply attractive potentials, one in each coordinate direction. The
class of perturbed SAPs is much larger.  For example, if $\Phi$ is a
nearest neighbor gradient potential defined by $\Phi_{\{x,y\}}(\phi)
= V_{x,y}(\phi(x)-\phi(y))$ when $m=1$, then we can define a
radially symmetric higher dimensional potential $\overline{\Phi}$ by
$\overline{\Phi}_{x,y}(\psi) = V_{x,y}(||\psi(x) - \psi(y)||)$, for
$\psi: \mathbb Z^d \mapsto \mathbb R^m$.  If the $V_{x,y}$ are
increasing on $[0, \infty)$ and for some $b>0$ satisfy the condition
that $V_{x,y}(m\eta) \leq b V_{x,y}(\eta)$ for all $\eta$, then
$\overline{\Phi}$ is a perturbed simply attractive potential.
(Observe that $\Psi_{\{x,y\}}(\psi) = b \sum_{i = 1}^m
V_{x,y}(|\psi(x)_i - \psi(y)_i|)$ is simply attractive.  Note that
$m \sup_{1 \leq i \leq m} |\eta_i| \geq ||\eta|| \geq \sup_{1 \leq i
\leq m} |\eta_i|$. Thus, $b\sum V_{x,y}(|\eta_i|) \geq
V_{x,y}(||\eta||) \geq \sum \frac{1}{m} V_{x,y}(|\eta_i|) $; it
follows that $\Psi \geq \overline{\Phi} \geq \frac{1}{mb}\Psi$ and
$\Psi$ strictly dominates $\Psi - \overline{\Phi}$.)  It is also
easy to see that the sum of a perturbed simply attractive potential
and any bounded potential is (at least after adding a constant) a
perturbed simply attractive potential.

SAPs and perturbed SAPs are (respectively) the most general convex
and not-necessarily-convex potentials we consider. Most of the
constructions in Chapters \ref{SFEchapter},
\ref{decompositionchapter}, \ref{surfacetensionchapter} apply to all
perturbed SAPs and are valid for any $E = \mathbb R^m$ or $E =
\mathbb Z^m$. The results of Chapter \ref{orliczsobolevchapter} are
analytical results used in later chapters; most of them are stated
in terms of ISAPs with $m=1$. The variational and large deviations
principle results of Chapters \ref{LDPempiricalmeasurechapter} and
\ref{LDPchapter} apply to perturbed ISAPs and perturbed LSAPs and
are valid for any $E=\mathbb R^m$ or $E = \mathbb Z^m$.  (We will
actually prove the results first for ISAPs and LSAPs when $m=1$ and
then observe that extensions to perturbed versions and to general
$m$ are straightforward.) All of the surface tension strict
convexity and gradient phase classifications in Chapters
\ref{clusterswappingchapter} and \ref{discretegibbschapter} apply to
SAPs in the case $m=1$.

\section{Overview of remaining chapters} \subsection{Specific free energy and surface tension}
In Chapter \ref{SFEchapter} we define the {\it specific free energy}
of a measure $\mu \in \mathcal P_{\mathcal L}(\Omega, \mathcal
F^{\tau})$ (denoted $SFE(\mu)$) and prove several consequences of
that definition. In particular, we show that if $\mu$ has slope $u$
and has minimal specific free energy among measures of slope $u$,
then $\mu$ is necessarily a gradient Gibbs measure.  (This is the
first half of our {\it variational principle}.) We discuss ergodic
and extremal decompositions in Chapter \ref{decompositionchapter}
and prove that $SFE(\mu)$ can be written as the $\mu$-expectation of
a particular tail-measurable function that is independent of $\mu$.
(In particular, $SFE$ is affine.)  These definitions and results are
analogous to those of the standard reference text \cite{G}, although
the setting is different.  In Chapter \ref{surfacetensionchapter},
we define the surface tension $\sigma(u)$ to be the infimum of
$SFE(\mu)$ over all slope-$u$ measures $\mu \in \mathcal P_{\mathcal
L}(\Omega, \mathcal F^{\tau})$.  The {\em pressure} of a potential
$\Phi$ is the infimum of the values assumed by $\sigma$ and denoted
$P(\Phi)$. Let $U_{\Phi}$ be the interior of the set of slopes $u$
for which $\sigma(u) < \infty$. We will see that whenever $\Phi$ is
a perturbed SAP, the set $U_{\Phi}$ is either all of $\mathbb R^{d
\times m}$ or the intersection of finitely many half spaces. Several
equivalent definitions of specific free energy and surface tension
are contained in Chapter \ref{LDPempiricalmeasurechapter}.

\subsection{Orlicz-Sobolev spaces and other analytical results}
In Chapter \ref{orliczsobolevchapter}, we define {\it Orlicz-Sobolev
spaces} and cite a number of standard results about them
(compactness of embeddings, equivalence of spaces, miscellaneous
bounds, etc.) from \cite{C}, \cite{C1}, \cite{C2}, \cite{A}, and
\cite{M}. The Orlicz-Sobolev space theory will enable us to derive
(in some sense) the strongest possible topology on surface shapes
(usually a topology induced by the norm of an Orlicz-Sobolev space)
in which our large deviations principles on surface shapes apply.

For example, this will enable us to prove that our large deviations
principles for the two-dimensional Ginzburg-Landau models apply in
any $L^p$ topology with $p < \infty$, whereas these results were
only proved for the $L^2$ topology in \cite{FS} and \cite{DGI}. This
allows us in particular to produce stronger concentration
inequalities---to show that typical random surfaces are ``close'' to
free-energy minimizing surfaces in an $L^p$ sense instead of merely
an $L^2$ sense. One of the reasons that Orlicz-Sobolev space theory
was developed was to provide tight conditions for the existence of
bounded solutions to PDE's and to variational problems involving the
minimization of free energy integrals; so it is not too surprising
that these tools should be applicable to the discretized/randomized
versions of these problems as well.

\subsection{Large deviations principle}
In Chapter \ref{LDPempiricalmeasurechapter} we derive several
equivalent definitions of the specific free energy and surface
tension.  We also complete the proof of the variational principle
(for perturbed ISAP and discrete LSAP models), which states the
following: if $\mu \in \mathcal P_{\mathcal L}(\Omega, \mathcal
F^{\tau})$ is $\mathcal L$-ergodic and has finite specific free
energy and slope $u$, then $\mu$ is a gradient Gibbs measure if and
only if $SFE(\mu) = \sigma(u)$. In particular, every gradient phase
of slope $u$ is a minimal gradient phase. In Chapter
\ref{LDPchapter} we derive a large deviations principle for
normalized height function shapes and ``empirical measure
profiles.''  Following standard notation (see, e.g., Section 1.2 of
\cite{DZ}), we say that a sequence of measures $\mu_n$ on a
topological space $(X, \mathcal X)$ {\it satisfy a large deviations
principle with rate function $I$ and speed $n^d$} if $I: X
\rightarrow [0, \infty]$ is lower-semicontinuous and for all sets $B
\in \mathcal X$, $$-\inf_{x \in B^o} I(x) \leq \liminf_{n
\rightarrow \infty} n^{-d} \log \mu_n(B) \leq \limsup_{n\rightarrow
\infty} n^{-d} \log \mu_n(B) \leq -\inf_{x \in \overline{B}} I(x).$$
Here $B^o$ is the interior and $\overline{B}$ the closure of $B$.
Roughly speaking, we can think of $I(x)$ as describing the
exponential ``rate'' (in terms of $n^d$) at which $\mu_n(B_x)$
decreases when $B_x$ is a very small neighborhood of $x$.  Also,
note that if $I$ obtains its minimum at a unique value $x_0 \in X$
and $B$ is any neighborhood of $x_0$, then $\mu_n(X \backslash B)$
decays exponentially at rate $\inf_{x \in X \backslash B}I(x)$
(whenever this value is non-zero). We refer to bounds of this form
as {\it concentration inequalities}, since they bound the rate at
which $\mu_n$ tends to be concentrated in small neighborhoods of
$x_0$.

By choosing the topological spaces appropriately, we will produce a
large deviations principle on random surface measures
which---although its formulation is rather technical---encodes a
great deal of information about both the typical local statistics
and global ``shapes'' of the surfaces. Though we defer a complete
formulation until Chapter \ref{LDPchapter}, a rough but almost
complete statement of our main large deviations result is the
following. Let $D$ be a bounded domain in $\mathbb R^d$ (satisfying
a suitable isoperimetric inequality), and write $D_n = n D \cap
\mathbb Z^d$. Let $\Phi$ be a perturbed ISAP or LSAP, and use
$H^o_{D_n}$ to define a Gibbs measure $\mu_n$ on gradient
configurations on $D_n$. Given $\phi_n: D_n \mapsto E$, we define an
{\it empirical profile measure} $R_{\phi_n, n} \in \mathcal P(D
\times \Omega)$ as follows:
$$R_{\phi_n, n} = \int_{D} \delta_{(x, \theta_{\lfloor nx \rfloor}
\phi_n)}dx,$$ (where $(\theta_y \phi)(x) := \phi(x+y)$). Informally,
to sample a point $(x, a)$ from $R_{\phi_n, n}$, we choose $x$
uniformly from $D$ and then set $a = \theta_{\lfloor nx \rfloor}
\phi_n$ (where $\phi_n$ is defined to be zero or some other
arbitrary value outside of $D_n$), i.e., $a$ is $\phi_n$ shifted so
that the origin is near $x$.  Also, using $\phi_n$, we will define a
function $\tilde \phi_n:D \mapsto \mathbb R^m$ by interpolating the
function $\frac{1}{n}\phi_n(nx)$ to a continuous, piecewise linear
functions on $D$; each such $\tilde \phi_n$ will be a member of an
appropriate {\it Orlicz-Sobolev space} $L^A$ (actually, a slight
extension $L^A_0$ of $L^A$ to include functions defined on most but
not all of $D$) where $A = \overline{V}^*$, a function we define
later.

Let $\rho_n$ be the measure on $\mathcal P(D \times \Omega) \times
L^A_0$ induced by $\mu_n$ and the map $\phi_n \mapsto (R_{\phi_n,n},
\tilde\phi_n)$. We say a measure $\mu \in \mathcal P(D \times
\Omega)$ is {\it $\mathcal L$-invariant} if $\mu (\cdot, \Omega)$ is
Lebesgue measure on $D$ and for any $D' \subset D$ of positive
Lebesgue measure, $\mu (D', \cdot)$ is an $\mathcal L$-invariant
measure on $\Omega$. Given any subset $D'$ of $D$ with positive
Lebesgue measure, we can write $S(\mu(D',\cdot))$ for the slope of
the measure $\mu(D', \cdot)/\mu(D' \times \Omega)$ (we have
normalized to make this a probability measure) times $\mu(D' \times
\Omega)$. The map $D' \rightarrow S(\mu(D', \cdot))$ is a signed,
vector-valued measure on $D$, and in particular, when $\psi$ is
smooth, we can define integrals $\int \psi(x) S(\mu(x, \cdot))dx$,
which we expect to be the same as the integral of the gradient of
the limiting surface shape $f$, i.e., $\int \psi(x) \nabla f (x)
dx$.  In fact, we show that the $\rho_n$ satisfy a large deviations
principle with speed $n^d$ and rate function
    $$I(\mu,f) = \begin{cases}
        SFE (\mu(D, \cdot)) - P(\Phi)  &
        \text{$\mu$ is $\mathcal L$-invariant and $S(\mu(x,\cdot)) = \nabla f(x)$}\\
        & \text{  as a distribution} \\
        \infty & \text{otherwise.} \\
        \end{cases} $$ in an appropriate topology on the space $\mathcal P(D \times \Omega) \times L^A_0$.
Contraction to the first coordinate yields an ``empirical profile'' large deviations principle;
contraction to the second coordinate yields a ``surface shape'' large deviations principle.
Analogous results apply in the presence of boundary restrictions on the $\phi_n$.

In Section \ref{externalfields}, we will see that the introduction of ``gravity'' or other
``external fields'' to our models alters the rate function $I$ in a predictable way; by computing
the rate function minimizer of the modified systems, we can describe the way ``typical surface
behavior'' changes in the presence of external fields. In fact, the ease of making changes of this
form is one of the main appeals of the large deviations formalism in statistical physics in
general: the rate function tells not only the ``typical'' macroscopic behavior but also the
relative free energies of all of the ``less likely'' behaviors which may become typical when the
system is modified.

\subsection{Surface tension strict convexity and Gibbs measure classifications}
The results in Chapters \ref{clusterswappingchapter} and Chapter \ref{discretegibbschapter} pertain
only to the case that $m=1$ and $\Phi$ is a simply attractive potential. In Chapter
\ref{clusterswappingchapter}, we introduce a geometric construction called {\it cluster swapping}
that we use to prove that the surface tension $\sigma$ is strictly convex and to classify gradient
Gibbs measures.  In some cases, these results will allow us to prove the uniqueness of the minimum
of the rate function of the LDP derived in Chapter \ref{LDPchapter}---and hence, also some
corresponding concentration inequalities.  For every $u \in U_{\Phi}$, there exists at least one
minimal gradient phase $\mu_u$ of slope $u$. We prove that each of the following is sufficient to
guarantee that this minimal gradient phase is unique:

\begin{enumerate}
\item $E = \mathbb R$
\item There exists a rough minimal gradient phase of slope $u$.
\item One of the $d$ components of $u$ is irrational.
\end{enumerate}

Each of the first two conditions also implies that $\mu_u$ is {\it
extremal}.  Whenever a minimal gradient phase of slope $u$ fails to
be extremal, it is necessarily smooth. We show that the extremal
components of $\mu_u$ can be characterized by their asymptotic
``average heights'' modulo $1$, and that every smooth minimal
gradient phase is characterized by its slope and its ``height offset
spectrum''---which is a measure on $[0,1)$ that is ergodic under
translations (modulo one) by the inner products $(u,x)$, for $x \in
\mathcal L$.  We give examples of models with non-trivial height
offset spectra and minimal gradient phase multiplicity---a kind of
phase transition---that occur when $d \geq 3$, $E =\mathbb Z$, and
$u$ is rational.

In Chapter \ref{discretegibbschapter}, we specialize to models in
which $\Phi$ is simply attractive, $d=2$, and $E = \mathbb Z$; many
classical models (e.g., perfect matchings of periodic weighted
graphs, square ice, certain six-vertex models) belong to this
category.  These models are sometimes used to describe the surface
of a crystal at equilibrium.  We show that in this setting, the
height offset spectra of smooth phases are always point masses in
$[0,1)$.  In this setting, $\mu_u$ is unique and extremal for every
$u \in U_{\phi}$ and the slopes of {\it all} smooth minimal phases
(also called {\it crystal facets}) lie in the dual lattice
$\tilde{\mathcal L}$ of $\mathcal L$.

\subsection{Differences from previous work}
Before reading on, the reader may wish to know which aspects of our
research we would expect a researcher with years of experience in
large deviations theory and statistical physics to find new or
surprising.

For readers who have studied the variational principle in the
context of, say, the Ising model, our random surface
formulation---that an ergodic gradient measure is Gibbs if and only
if it minimizes specific free energy among measures {\it of that
slope}---may not come as a huge surprise.  Indeed, it may surprise
the reader that nobody had formulated and proved this fundamental
result before.

The fact that the large deviations principle extends to empirical
measure profiles requires many technical advances, but the result
itself is also not shocking (in light of the many similar results
known for, say, the Ising model). Readers who learned Sobolev space
theory a couple of decades ago may be somewhat surprised to learn
how much stronger, simpler, and more intuitive the theory has
become---and how much of it seems to have been custom-made for our
research. Instead of imposing lots of ad hoc conditions, we can now
derive large deviations principles in the ``right'' topologies and
with the ``right'' boundary conditions while citing most of the
needed analytical results from other sources.

But in our view, the most surprising aspect of our large deviations
principles is the generality in which we prove uniqueness of the
rate function minimizer.  This uniqueness is a consequence of two
key results: the strict convexity of $\sigma$ and the uniqueness of
the gradient Gibbs measure of a given slope.  Both of these results
are proved in Chapters 8 and 9 using the variational principle and a
new geometric construction called ``cluster swapping.''

Before our work, some researchers suspected that if $V$ failed to be
strictly convex, then the surface tension $\sigma$ corresponding to
the ISAP $\Phi_V$ would also fail to be strictly convex.  In
particular, it was unknown whether the surface tension corresponding
to the linear solid-on-solid model $V(x)=|x|$ was everywhere
strictly convex in both the discrete and continuous height versions.

Also, although the uniqueness of the gradient Gibbs measure of a
given slope was known for Ginzburg-Landau $\nabla \phi$ models and
conjectured for some discrete models (see \cite{CEP}, \cite{CKP},
and the next section), our statement---particularly in the discrete
case $E = \mathbb Z$---is much more general than had been
conjectured.

Finally, our discrete model analysis of the smooth-phase/rough-phase
distinction in Chapters 8 and 9 is new.  The ``height offset
spectrum'' decomposition for general $d$, and the fact that when
$d=2$ all smooth phases have slopes in the dual lattice of $\mathcal
L$, were both, to our knowledge, unexpected.  Indeed, the dimer
model analog of the smooth phase classification theorem is one of
the more surprising qualitative results in \cite{KOS}.  In
additional to cluster swapping, our proofs of these results use, in
a new way, the FKG inequality and the homotopy theory of the
countably punctured plane.

\section{Two important special cases}
Special cases of what we call simply attractive potentials have been
very thoroughly studied in a variety of settings.  In this section,
we will briefly review relevant facts about two of the most
well-understood random surface models: domino tiling height
functions (here $E = \mathbb Z$) and the Ginzburg-Landau $\nabla
\phi$-interface models (here $E = \mathbb R$).  Each of these models
is the subject of a sizable literature, and each has features that
make it easier to work with than general simply attractive or
perturbed simply attractive potentials.

An exhaustive survey of the myriad physical, analytical,
probabilistic, and combinatorial results available for even these
two models --- let alone all simply attractive models --- is beyond
the scope of this work.  But we will mention a few of the papers and
conjectures that directly inspired our results and provide pointers
to the broader literature.  See the survey papers \cite{Ke4},
\cite{Gi}, \cite{Gi3}, and \cite{Gi4} for more details.

\subsection{Domino tilings} \label{dominosubsection}
Though we mentioned the modeling of solids and phase interfaces as
one motivation for our work, the models we describe have been used
for other purposes as well.  When $E= \mathbb Z$ and $\Phi$ is
chosen appropriately, the finite-energy surfaces $\phi \in \Omega$
correspond to the so-called {\it height functions} that are known to
be in one-to-one correspondence with spaces of domino tilings,
square ice configurations, and other discrete statistical physics
models.

We will now explicitly describe a well-known correspondence between
the set of domino tilings of a simply connected subset $R$ of the
squares of the $\mathbb Z^2$ lattice and the set of {\it height
functions} from the vertices of $R$ to $\mathbb Z$ that satisfy
certain boundary conditions and difference constraints.  Let
$\epsilon: \mathbb Z^2 \mapsto \{0,1,2,3\}$ be such that if
$i=(i_1,i_2) \in \mathbb Z^2$, then $\epsilon(i_1,i_2)$ assumes the
values $0$, $1$, $2$, and $3$ as the value of $i$ modulo $2$ is
respectively $(0,0)$, $(0,1)$, $(1,1)$, and $(1,0)$.  Given a
perfect matching of the squares of $\mathbb Z^2$ (which we can think
of as a ``domino tiling'', where a domino corresponds to an edge in
the matching), a height function $\psi$ on the vertices of $\mathbb
Z^2$ is determined, up to an additive constant, by the following two
requirements:
\begin{enumerate} \item $\psi(x) = \epsilon(x) \mod 4$
\item If $x$ and $y$ are neighboring vertices, then $|\psi(x) - \psi(y)| = 3$ if the edge between
them crosses a domino and $1$ otherwise. \end{enumerate} In the
height functions thus defined, the set of possible values of
$\psi(x)$ depends on the parity of $x$; in order to describe these
height functions as the finite-energy surfaces of a Gibbs potential,
we will instead use a slight modification: $\phi(x)=\frac{\psi(x) -
\epsilon(x)}{4}$. Now the set of possible heights at any vertex is
equal to $\mathbb Z$.  The set of height functions $\phi$ of this
form that arise from tilings are precisely those functions
$\phi:\mathbb Z^2 \mapsto \mathbb Z$ for which
$H^{\Phi}_{\Lambda}(\phi)$ is finite for all $\Lambda$, where $\Phi$
is the LSAP determined by the following nearest neighbor potentials:
$$V_{x,y}(\eta) = \begin{cases} 0 & \eta = 0 \\
    0 & \epsilon(x) > \epsilon(y) \text{ and } \eta = 1 \\
    0 & \epsilon(x) < \epsilon(y) \text { and } \eta = -1 \\
    \infty & \text{otherwise.} \\ \end{cases}$$ Because of this correspondence, we may think of a domino
tiling chosen uniformly from the set of all domino tilings of a
simply-connected subset $R$ of the squares of $\mathbb Z^2$ as a
discretized, incrementally varying {\it random surface}.  We can
also think of domino tilings as perfect matchings of a bipartite
graph.  It is not hard to compute the number of perfect matchings of
a bipartite graph using the permanent of an adjacency matrix.
Kastelyn observed in 1965 that by replacing the $1$'s in the
adjacency matrix with other roots of unity, one can convert the
(difficult) problem of permanent calculation into a (much easier)
determinant calculation.  The rich algebraic structure of
determinants has rendered tractable many problems that appear
difficult for more general families of random surfaces.

In one recent paper \cite{CKP}, Cohn, Kenyon, and Propp proved the
following: Suppose $R^n$ is a sequence of domino-tilable regions
such that the boundary of $\frac{1}{n}R^n$ converges to that of a
simple region $R$ in $\mathbb R^2$.  Let $\mu_n$ be the uniform
measure on the set of tilings of $R^n$.  Each $\mu_n$ also induces a
measure on the set of possible height functions $\phi_n$ on the
vertices of $R^n$; the values of $\phi_n$ on the boundary of $R^n$
are determined by the shape of $R^n$ independently of the tiling.
Suppose further that the boundary height functions
$\frac{1}{n}\phi_n(nx)$ converge (in a certain sense) to a function
$f_0$ defined on the boundary of $R$.  Then, \cite{CKP} shows that
as $n$ gets large, the normalized height functions of tilings chosen
from the $\mu_n$ approach the unique Lipschitz (with respect to an
appropriate norm) function $f$ that agrees with $f_0$ on the
boundary of $R$ and minimizes a surface tension integral
$$I(f) = \int_{R}\sigma(\nabla f(x)) dx.$$ In fact, their results imply that this surface tension
integral is a rate function for a {\it large deviations principle} (under the supremum topology)
that holds for a sequence of random surface measures $\nu_n$---derived from the $\mu_n$ by standard
interpolations---on the space of Lipschitz functions on $R$ that agree with $f_0$ on the boundary.

These authors also explicitly describe an ergodic gradient Gibbs
measures $\mu_u$ of each slope $u = (u_1, u_2)$ inside the set
$U_{\Phi} = \{u:|u_1|+|u_2| < \frac{1}{2}\}$; they show that only
zero-entropy gradient ergodic Gibbs measures exist with slopes on
the boundary of $U_{\Phi}$, and no Gibbs measures exist with slopes
outside the closure of $U_{\Phi}$.  Since every tiling determines a
height function up to an additive constant, a gradient Gibbs measure
in this context is equivalent to a Gibbs measure on tilings, where
in both cases we can take $\mathcal L = 2\mathbb Z^2$.  They
conjecture that for each $u\in U_{\Phi}$, $\mu_u$ is the {\it only}
gradient phase of slope $u$.  A similar conjecture appears in an
earlier paper by Cohn, Elkies, and Propp \cite{CEP}.

We will resolve this conjecture in Chapter \ref{discretegibbschapter}.  We also resolve another of
their conjectures (concerning the local probability densities of domino configurations in large
random tilings) as a consequence of our large deviations principle on profiles in Chapter
\ref{LDPchapter}. We extend these results, as well as the large deviations principle on random
surface shapes produced in \cite{CKP}, to more general families of simply attractive random
surfaces.

Using Kastelyn determinants, the authors of \cite{CKP} were able to
compute the surface tension $\sigma$ and the ergodic Gibbs measures
$\mu_u$ {\it exactly} in terms of special functions, and their large
deviations results rely on these exact computations. See \cite{KOS}
for a generalization of these computations to perfect matchings of
other planar, doubly periodic graphs by the author and two
co-authors.  These authors use algebraic geometric constructions
called {\em amoebae} to ``exactly solve'' the dimer model on general
weighted doubly periodic lattices by explicitly computing $\sigma$
and the local probabilities in all of the $\mu_u$.  The
characterization of ergodic Gibbs measures on perfect matchings
given in \cite{KOS} makes use of a few results from Chapters 8 and 9
of this text, including the uniqueness of a measure $\mu_u$ of a
given slope.  While we prove for a much more general class of
two-dimensional models that all smooth phases have slopes in the
dual lattice $\mathcal L'$ of the lattice $\mathcal L$ of
translation invariance, the authors in \cite{KOS} use the exact
solvability to determine precisely which of these slopes admit
smooth phases.  The smooth phases in this context are in
correspondence with cusps of the surface tension, and depending on
the way the edges of the doubly periodic planar graph are weighted,
some, all, or none of the Gibbs measures $\mu_u$ corresponding to $u
\in \mathcal L'$ will actually be smooth.

Currently, it seems unlikely that the techniques of \cite{KOS} can
be extended to exactly solve more general random surface models ---
particularly those in dimensions higher than two; but we will prove
enough qualitative results (such as the strict convexity of $\sigma$
and the gradient Gibbs measure classification) to show that the
large deviations theorems apply in general.

\subsection{Ginzburg-Landau $\nabla \phi$-interface models}
Recent papers by Funaki and Spohn \cite{FS} and Deuschel, Giacommin,
and Ioffe \cite{DGI} derive similar results for a continuous
generalization of the harmonic crystal called the {\it
Ginzburg-Landau $\nabla \phi$-interface model}.  These models use
ISAPs in which $E = \mathbb R$ and $V_{x,y} = V$ for all adjacent
$x,y$.  Here $V:\mathbb R \mapsto \mathbb R$ is convex, symmetric,
and $C^2$, with second derivatives bounded above and below by
positive constants. Such potentials $V$ are bounded above and below
by quadratic functions --- and we may think of them as
``approximately quadratic'' generalizations of the (Gaussian)
harmonic crystal, for which $V(\eta) = \beta \eta^2$.

Calculations for these models typically make use of the fact that
Gibbs measures in these models are stationary distributions of
infinite-dimensional elliptic stochastic differentiable equations
(see, e.g., \cite{NS}, for descriptions and more references).  Given
a configuration $\phi$, the ``force'' on any given ``particle'' (and
hence the stochastic drift of that particle's position) is within a
constant factor of what it would be if the potential were Gaussian;
and the rate at which a pair of Gibbs measures converges in certain
couplings is also within a constant of what it would be in a
Gaussian model. Although the calculations in, for example, \cite{FS}
or \cite{DGI}, are still rather complicated, they appear to be
simpler than they would be for general simply attractive models.
These authors also derive static Gibbs measure results as
corollaries of more general dynamic results.  For example, Funaki
and Spohn prove the uniqueness of gradient phases of a given slope
$u$ using a dynamic coupling \cite{FS}. Although the Gibbs measure
classifications and surface shape large deviations principles are
proved in these papers, our large deviations principle for profiles
and our variational principle are new results for $\nabla
\phi$-interface models.  Also, as mentioned earlier, we derive our
large deviations principles with respect to stronger topologies than
\cite{FS} and \cite{DGI}.

See, e.g., Giacomin's survey papers (\cite{Gi}, \cite{Gi3}, \cite{Gi4}) for many more references
about Ginzburg-Landau $\nabla \phi$-interface models, including wetting transitions, entropic
repulsion, roughening transitions, etc.

\chapter{Specific free energy and variational principle} \label{SFEchapter} The notion of {\it
specific free energy} lies at the heart of all of our main results.
Although definitions and applications of specific free energy are
well known for certain families of shift-invariant measures on
$(\Omega, \mathcal F)$ (see, e.g., Chapter 14 and 15 of \cite{G}),
we need to check that these notions also make sense for our
$\mathcal L$-invariant gradient measures on $(\Omega, \mathcal
F^{\tau})$.  In this chapter, we provide a definition of the
specific free energy of an $\mathcal L$-invariant gradient measure
and prove some straightforward consequences, including the first
(and easier) half of the variational principle.  We will cite lemmas
directly from reference texts \cite{G} and \cite{DZ} whenever
possible.  First, we review some standard facts about relative
entropy and free energy.

\section{Relative entropy review}
Throughout this section, we let $(X, \mathcal X)$ be any {\it
Polish} space (i.e., a complete, separable metric space endowed with
the metric topology and the Borel $\sigma$-algebra $\mathcal X$),
$\mu$ and $\nu$ any probability measures on $(X, \mathcal X)$, and
$\mathcal A$ a sub $\sigma$-algebra of $\mathcal X$.  Write $\mu <<
\nu$ if $\mu$ is absolutely continuous with respect to $\nu$.  The
{\it relative entropy} of $\mu$ with respect to $\nu$ on $\mathcal
A$, denoted $\mathcal H_{\mathcal A}(\mu, \nu)$, is defined as
follows: $$\mathcal H_{\mathcal A}(\mu, \nu) = \begin{cases}
\nu(f_{\mathcal A}\log
f_{\mathcal A}) & \text{ $\mu < < \nu$ on $\mathcal A$} \\ \infty & \text{otherwise,}\\
\end{cases}$$ where $f_{\mathcal A}$ is the Radon-Nikodym derivative of $\mu$ with respect to $\nu$
when both measures are restricted to $\mathcal A$. (We often write
$\mathcal H$ to mean $\mathcal H_{\mathcal X}$.) Note that this
definition still makes sense if $\nu$ is a finite (positive),
non-zero measure (not necessarily a probability measure). If $a>0$
and $\mu$ is as above with $\mu << \nu$ on $\mathcal A$, we have:
$$\mathcal H_{\mathcal A}(\mu, a\nu) = a\nu (\frac{1}{a}f_{\mathcal
A} \log (\frac{1}{a} f_{\mathcal A})) = \mathcal H_{\mathcal A}(\mu,
\nu) - \log a.$$ We now cite Proposition 15.5 of \cite{G}:

\begin{lem}  \label{relativeentropylemma} If $\mu$ and $\nu$ are probability measures and $a>0$,
then \begin{enumerate}
\item $\mathcal H_{\mathcal A}(\mu, \nu) \geq 0$ (and thus $\mathcal H_{\mathcal A}(\mu, a\nu) \geq -\log
a$)
\item $\mathcal H_{\mathcal A}(\mu, a\nu) = -\log a$ if and only if $\mu = \nu$ on $\mathcal A$ \item $\mathcal
H_{\mathcal A}(\mu, a\nu)$ is an increasing function of $\mathcal A$
\item $\mathcal H_{\mathcal A}(\mu, \nu)$ is a convex function of
the pair $(\mu, \nu)$ when $\mu$ ranges over probability measures
and $\nu$ over non-zero finite measures on $(X, \mathcal X)$
\end{enumerate} \end{lem} The following very important fact about
relative entropy (Proposition 15.6 of \cite{G}) will enable us to
approximate relative entropy with respect to a subalgebra $\mathcal
A$ by the relative entropy with respect to convergent subalgebras.
Throughout this work, we denote by $\mathcal P(X, \mathcal X)$ the
set of probability measures on $(X, \mathcal X)$.

\begin{lem} \label{increasingalgebraentropy} Let $\mu, \nu \in \mathcal P(X, \mathcal X)$ and let
$\mathcal A_n$ be an increasing sequence of subalgebras of $\mathcal
X$, and $\mathcal A$ the smallest $\sigma$-algebra containing
$\cup_{n=1}^{\infty} \mathcal A_n$. Then: $$\lim_{n \rightarrow
\infty} \mathcal H_{\mathcal A_n}(\mu, \nu)  = \mathcal H_{\mathcal
A}(\mu,\nu) = \sup_{n} \mathcal H_{\mathcal A_n}(\mu, \nu).$$
\end{lem} Regular conditional probability distributions do not exist
for general probability measures on general measure spaces. However,
the following lemma states that they do exist for all of the spaces
and measures that will interest us here.  It also enables us to
express the relative entropy of a measure on a product space $X_1
\times X_2$ as the relative entropy of the first component plus the
expected {\it conditional} entropy of the second component given the
first. Here, we let $\eta = (\eta_1, \eta_2)$ denote a generic point
in $X$. (See Theorem D.3 and D.13 of \cite{DZ}.)

\begin{lem} \label{conditionalentropy} Suppose $X = X_1 \times X_2$, where each $X_i$ is Polish
with Borel $\sigma$-algebra $\mathcal X_i$.  Let $\mu_1$ and $\nu_1$ denote the projections of
$\mu, \nu \in \mathcal P(X, \mathcal X)$ to $X_1$. Then there exist regular conditional probability
distributions $\mu^{\eta_1}(\cdot)$ and $\nu^{\eta_1}(\cdot)$ on $X_2$ corresponding to the
projection map $\pi: X \mapsto X_1$. Moreover, the map: $$\eta_1 \mapsto \mathcal
H(\mu^{\eta_1}(\cdot),\nu^{\eta_1}(\cdot)): X_1 \mapsto [0,\infty]$$ is $\mathcal X_1$ measurable
and $$\mathcal H_{\mathcal X}(\mu,\nu) = \mathcal H_{\mathcal X_1}(\mu_1, \nu_1) + \int_{X}
\mathcal H(\mu^{\eta_1}(\cdot),\nu^{\eta_1}(\cdot)) \mu_1(d\eta_1).$$ \end{lem} This result and the
following simple corollary are the key observations behind the proof of the easy half of our
variational principle (which states that if a slope $u$ gradient measure has minimal specific free
energy among measures of slope $u$, then it must be a gradient Gibbs measure).

\begin{lem} \label{productentropy} Let $\mu_1$ and $\nu_1$ be probability measures on $X_1$, and
$\mu_2$ and $\nu_2$ probability measures on $X_2$. Suppose that
$\mu$ is a probability measure on $X_1 \times X_2$ with marginal
distributions given by $\mu_1$ and $\mu_2$, respectively. Then:
$$\mathcal H(\mu, \nu_1 \otimes \nu_2) \geq \mathcal H(\mu_1 \otimes \mu_2, \nu_1 \otimes \nu_2) =
\mathcal H(\mu_1, \nu_1) + \mathcal H(\mu_2, \nu_2).$$ If we assume
that $\mathcal H(\mu_1, \nu_1)< \infty$ and $\mathcal H(\mu_2,
\nu_2) < \infty$, then equality holds if and only if $\mu = \mu_1
\otimes \mu_2$. \end{lem}

\begin{proof} By the previous lemma applied to $\nu = \nu_1 \otimes \nu_2$, we may assume that
$\mathcal H(\mu_1,\nu_1)<\infty$, and it is enough to show that
$$\int_{X}\mathcal H(\mu^{\eta_1}(\cdot),\nu_2(\cdot))\mu_1(d\eta_1)
\geq \mathcal H(\mu_2, \nu_2).$$ Since $\int_{X}\mu^{\eta_1}(\cdot)
\mu_1(d\eta_1) = \mu_2$, this follows from Jensen's inequality and
the convexity of $\mathcal H(\cdot,\nu_2)$, stated in Lemma
\ref{relativeentropylemma}. This function is strictly convex on its
level sets, hence the characterization of equality. \qed \end{proof}
Finally, we will also be interested in the convergence of sequences
of probability measures. Denote by $\mathcal P(X, \mathcal X)$ the
space of probability measures on $(X, \mathcal X)$.  The {\it weak
topology} on $\mathcal P(X, \mathcal X)$ is the smallest topology
with respect to which the map $\nu \mapsto \nu(f)$ is continuous for
all bounded continuous functions $f$.  The {\it $\tau$-topology} is
the smallest topology with respect to which $\nu \mapsto \nu(f)$ is
continuous for all bounded $\mathcal X$-measurable functions $f$ on
$X$. The reader may check that this is also the smallest topology
with respect to which $\nu \mapsto \nu(A)$ is continuous for every
$A \in \mathcal X$.  In general, the $\tau$-topology is stronger
than the weak topology. The two topologies coincide when $(X,
\mathcal X)$ is discrete. We now cite two lemmas (Lemma 6.2.12 with
its proof and Lemma D.8 of \cite{DZ}):

\begin{lem} \label{compact} Fix $C \in \mathbb R$ and a finite measure $\nu$ on $(X, \mathcal X)$;
then the {\it level set} $M_{C, \nu}$ of probability measures $\mu$ on $(X, \mathcal X)$ with
$\mathcal H(\mu, \nu) \leq C$ is compact in the $\tau$-topology.  \end{lem} \begin{lem}
\label{weakmetrizable} The weak topology is metrizable on $\mathcal P(X, \mathcal X)$ and makes
$\mathcal P(X, \mathcal X)$ into a Polish space (i.e., a complete, separable metric space).
\end{lem} From these lemmas, we deduce the following:

\begin{lem} \label{weaktauequivalent} The $\tau$-topology restricted to a level set $M_{C, \nu}$ is
equivalent to the weak topology and hence also metrizable.  In particular the compactness of the
level sets (Lemma \ref{compact}) implies sequential compactness of the level sets in both
topologies. \end{lem}

\begin{proof} This is a well-known result (stated, for example, in the proof of Theorem 3.2.21 and
Exercise 3.2.23 of \cite{DS}), but we sketch the proof here.  It is
enough to prove that the $\tau$-topology on $M_{C,\nu}$ is contained
in the weak topology on $M_{C, \nu}$. We can prove this by showing
that if a measure $\mu$ lies in a base set $\mathcal B$ of the
$\tau$-topology, then there exists a base set $\mathcal B'$ of the
weak topology with $\mu \in \mathcal B' \subset \mathcal B$.
Precisely, we must show that if $A \in X$, then for each
$\epsilon>0$ and measure $\mu$, there exists a bounded, continuous
function $f$ and an $\epsilon_0>0$ such that $|\mu(f) - \mu'(f)| <
\epsilon_0$ implies $|\mu(A) - \mu'(A)| < \epsilon$ whenever $\mu'
\in \mathcal M_{C, \nu}$.  Note that: $$|\mu(A) - \mu'(A)| \leq
|\mu(A) - \mu(f)| + |\mu(f) - \mu'(f)| + |\mu'(f) - \mu'(A)|.$$ We
would like to show that by choosing $\epsilon$ and $f$
appropriately, we can force each of the right hand terms to be as
small as we like.  The second term is obvious.  For the first term,
it is enough to note that for any $\delta$, we can find a positive
continuous function $f$ such that $\mu|f - 1_A| < \delta$ and $0
\leq f \leq 1$. (Simply let take a closed set $A'$ with $\nu|1_{A'}
- 1_A| < \delta/2$; define $f(x) = 1$ for $x \in A'$ and $f(x) =
\sup (0, 1- \alpha d(A',x))$ otherwise, where $\alpha>0$ is
sufficiently small and $d(A',x)$ is the distance from $A'$ to $x$.)
For the third term, note that by taking $\delta$ and $\alpha$ small
enough in the above construction, we can also assume that $f$ and
$1_A$ are equal outside a set of $\nu$-measure at most $\gamma$ for
any $\gamma > 0$. Using the definition of relative entropy, it is
not hard to see that for any $\gamma'$ we can find a $\gamma$ small
enough so that for each $\mu' \in M_{C, \nu}$ we have $\mu'(B) <
\gamma'$ whenever $B \in \mathcal X$ is such that $\nu(B) < \gamma$;
this puts a bound of $2 \gamma '$ on the third term. \qed
\end{proof}

In particular, this lemma implies that the level sets are closed in both topologies, which implies
the following:

\begin{cor} \label{lowersemicontinuity}  For fixed $\nu$, the function $\mathcal H(\mu, \nu)$ is a
lower semicontinuous function on $\mathcal P(X,\mathcal X)$, endowed with either the weak topology
or the $\tau$-topology. \end{cor}

Our motivation for the last few lemmas is that, using these results, we will later define a
topology on $(\Omega, \mathcal F^{\tau})$ (the {\it topology of local convergence}) with respect to
which ``specific relative entropy'' and specific free energy have compact level sets. This will
allow us to deduce, for example, that the specific free energy achieves its minimum on sets that
are closed in this topology.  And this will lead to proofs of the existence of gradient Gibbs
measures of particular slopes.

\section{Free energy}
Let $\lambda$ be an ``underlying'' probability measure on $X=\mathbb
R^m$ (or $\mathbb Z^m$) and $V$ an $\mathcal X$-measurable
Hamiltonian potential for which $Z = \int_X
e^{-V(\eta)}\lambda(d\eta)$ is finite. We define the {\it free
energy} of any measure $\mu$ on $(X, \mathcal X)$ as the following
relative entropy: $$FE^V(\mu) = \mathcal H(\mu, e^{-V} \lambda),$$
where we use the convention that $f \lambda$ is the measure whose
Radon-Nikodym derivative with respect to $\lambda$ is $f$  We will
write $FE^V(\mu) = FE(\mu)$ when the choice of potential is clear
from the context.  We can also write this expression as $\mu(V+ \log
f)$, where $f$ is the Radon-Nikodym derivative of $\mu$ with respect
to $\lambda$.  When they exist, we refer to $\mu(V)$ as the {\it
energy} of $\mu$ and to $-\mu (\log f)$ as the ${\it entropy}$ of
$\mu$ (which is $-\infty$ if $\mu$ is not absolutely continuous with
respect to Lebesgue measure).  The probability measure $\mu_V =
Z^{-1}e^{-V}\lambda$ is called the {\it Gibbs} measure for the
Hamiltonian $V$. The free energy of the Gibbs measure is simply
$-\log Z$, and the free energy of a general measure $\mu$ can also
be written $\mathcal H(\mu,\mu_V) - \log Z$. A trivial consequence
of this fact and Lemma \ref{relativeentropylemma} is the following
so-called {\it finite dimensional variational principle}:

\begin{lem} \label{freeenergybound} The free energy of any probability measure $\mu$ on $X$ is
equal to or greater than that of $\mu_V$; equality holds if and only if $\mu= \mu_V$.  In other
words, a measure is Gibbs if and only if it has minimal free energy, equal to $-\log Z$. \end{lem}
The following monotonicity is also an easy consequence of the definitions.

\begin{lem} \label{freeenergymonotonicity} If $V_1(\eta) \geq V_2(\eta)$ for all $x$, then
$FE^{V_1}(\mu) \geq FE^{V_2}(\mu)$ for all $\mu$. \end{lem}

It will sometimes be useful to know that an upper bound on the free energy $\mu$ allows us to put a
lower bound on the amount of mass of $\mu$ that lies outside of a particular compact set.
\begin{lem} \label{probabilitybound} If $X$ is $\mathbb Z^m$ or $\mathbb R^m$ and $Z = \int_X
e^{-V(\eta)}d\eta$ is finite, then for every $c>0$ and every $d$,
there exists a compact set $S \subset X$ such that $\mu(S) \leq 1-c$
implies $FE^V(\mu) \geq d$ whenever $\mu$ is a probability measure
on $X$.
\end{lem}

\begin{proof} From Lemma \ref{freeenergybound}, we know that if we require $\mu(X-S)=1$, then the
minimal free energy $\mu$ can have is given by $$-\log \[\int_{X -
S} e^{-V(\eta)}\lambda(d\eta)\].$$  Similarly, if $\mu(X-S) = a$,
the minimal free energy is at least $$a\log (a)+(1-a)\log (1-a) - a
\log \int_{X-S}e^{-V(\eta)}d\eta - (1-a) \log
\int_{S}e^{-V(\eta)}d\eta.$$ (See Lemma \ref{conditionalentropy}.)
If we choose $S$ to be a large enough closed ball containing a given
point such that $\int_{X-S} e^{-V(\eta)}d\eta \leq 1/n$, we can make
the latter expression at least $$-\log 2 + a \log n - (1-a)\log
(Z-1/n),$$ which tends to $\infty$ with $n$. \qed \end{proof} Also,
because free energy can be defined as a relative entropy, each of
the lemmas proved in the previous section applies to free energy as
well.

\section{Specific free energy: existence via superadditivity} \label{SFEdefsection}
We now return to our infinite dimensional setting.  That is,
$(\Omega, \mathcal F^{\tau})$ is the set of functions from $\mathbb
Z^d$ to $E$ endowed with the $\sigma$-algebra $\mathcal F^{\tau}$
described in the introduction, and $\Phi$ is an $\mathcal L \times
\tau$ invariant gradient potential.  In this section, we use limits
to give a definition of specific free energy for $\mathcal
L$-invariant gradient measures on $(\Omega, \mathcal F^{\tau})$. We
prove the existence of these limits using subadditivity arguments.

Let $\mu_{\Lambda}$ be the measure on $E^{|\Lambda|}$ obtained by
restricting $\mu$ to $\mathcal F_{\Lambda}$.

Let $\Lambda_n$ denote the box $[0,kn - 1]^d \subset \mathbb Z^d$,
where $k$ is chosen so that $k \mathbb Z^d \subset \mathcal L$. When
$\lambda$ is a finite measure on $E$, a standard definition of the
specific free energy of an ordinary (not gradient) shift-invariant
measure on $(\Omega, \mathcal F)$ is the following limit of
normalized relative entropies: $$\lim_{n \rightarrow \infty}
|\Lambda_n|^{-1} \mathcal H(\mu_{\Lambda_n},
e^{-H^o_{\Lambda_n}}\lambda^{|\Lambda_n|}).$$ We will use a similar
approach in our gradient setting except that we will only compute
relative entropies with respect to the subalgebra of gradient
measurable sets.  That is, if $\mu$ is a $\mathcal L$-invariant
measure on $\mu$, we write:
$$SFE(\mu) = \lim_{n \rightarrow \infty} |\Lambda_n|^{-1} \mathcal H (\mu_{\Lambda_n},
e^{-H^o_{\Lambda_n}}\lambda^{|\Lambda_n -1|}),$$ which we interpret
as follows: Fix a reference vertex $v_0 \in \Lambda_n$ and let $v_1,
\ldots, v_{|\Lambda_n| - 1}$ be an enumeration of the remaining
vertices. In this context, by $\lambda^{|\Lambda_n| - 1}$ we mean
the measure $\nu$ on $(\Omega, \mathcal F^{\tau}_{\Lambda_n})$ such
that for any measurable $A \subset E^{|\Lambda_n|-1}$, the value
$$\nu ( \{\phi | (\phi(v_1)-\phi(v_0), \phi(v_2)-\phi(v_0), \ldots,
\phi_{|\Lambda_n|-1}-\phi(v_0)) \in A \})$$ is equal to the measure
of $A$ in the product measure $\lambda^{|\Lambda_n|-1}$.  (The
reader may check that this definition is independent of the choice
of $v_0$.)

Also, when $\mu$ is a gradient measure --- only defined on $\mathcal F^{\tau}$ --- then we write
$\mu_{\Lambda}$ to mean the restriction of $\mu$ to $\mathcal F_{\Lambda}^\tau$. The latter is also
the smallest $\sigma$-algebra with respect to which $\phi(v)- \phi(v_0)$ is measurable for each $v
\in \Lambda$, so we can think of $\mu_{\Lambda}$ as a measure on this $|\Lambda|-1$ dimensional
space, $E^{|\Lambda|-1}$.  In this context, the expression $\mathcal H (\mu_{\Lambda_n},
e^{-H^o_{\Lambda_n}}\lambda^{|\Lambda_n -1|})$ makes sense.

As a convenient shorthand, we also write $FE_{\Lambda}^{\Phi}(\mu) =
\mathcal H (\mu_{\Lambda},e^{-H^o_{\Lambda}}\lambda^{|\Lambda|-1})$
and refer to this as the {\it free energy of $\mu$ restricted to
$\Lambda$}.  (We occasionally drop the $\Phi$ from
$FE^{\Phi}_{\Lambda}$ when the choice of potential is understood.)
Let $Z^o_{\Lambda}$ be the integral of $e^{-H^o_{\Lambda}}$ over
entire space $E^{|\Lambda|-1}$ of gradient functions, as described
above, and refer to this as the {\em free boundary partition
function} of $\Lambda$ with respect to $\Phi$. It is clear that ---
at least for perturbed simply attractive models
--- this value is always finite.

Moreover, from Lemma \ref{freeenergybound}, it follows that
$FE^{\Phi}_{\Lambda}(\mu) \geq -\text{log} Z^o_{\Lambda}$ for all
$\mu$. When the choice of $\Phi$ is understood, we write
$W(\Lambda)=-\text{log} Z^o_{\Lambda}$ for this minimal free energy.
We say that a potential $\Phi$ is {\it positive} if
$\Phi_{\Lambda}(\phi) \geq 0$ for all $\Lambda$ and all $\phi$.

\begin{lem}  \label{subadditive} Suppose that $\Phi_{\Lambda}$ is a positive potential and that
$\Lambda_1$ and $\Lambda_2$ are finite connected subsets of $\mathbb
Z^d$ (of finite weight with respect to $\Phi$) that have exactly one
vertex $w$ in common.  Then $W(\Lambda_1 \cup \Lambda_2) \geq
W(\Lambda_1) + W(\Lambda_2) $. Furthermore, for every measure $\mu$
on $(\Omega, \mathcal F^{\tau}_{\Lambda_1 \cup \Lambda_2})$ we have:
$$FE_{\Lambda_1 \cup \Lambda_2}(\mu) \geq FE_{\Lambda_1}(\mu) +
FE_{\Lambda_2}(\mu).$$ \end{lem}

\begin{proof} We take $w$ to be the reference vertex $v_0$ for $\Lambda_1$, $\Lambda_2$, and
$\Lambda_1 \cup \Lambda_2$.  Write $X = E^{|\Lambda_1|-1}$ and
$Y=E^{|\Lambda_2|-1}$ and view $E^{|\Lambda_1 \cup \Lambda_2|-1}$ as
$X \times Y$.  Then $H^o_{\Lambda_1 \cup \Lambda_2}$ is equal to
$H^o_{\Lambda_1} + H^o_{\Lambda_2}$ {\it plus} the sum of
$\Phi_{\Lambda'}$ over all finite subsets $\Lambda'$ of $\mathbb
Z^d$ such that $\Lambda' \subset \Lambda_1 \cup \Lambda_2$ but
$\Lambda' \not \subset \Lambda_1$ and $\Lambda' \not \subset
\Lambda_2$.  Since the latter sum is positive, we have
$$H^o_{\Lambda_1 \cup \Lambda_2} \geq H^o_{\Lambda_1} +
H^o_{\Lambda_2},$$ and the lemma is now an immediate consequence of
Lemma \ref{productentropy} and \ref{freeenergymonotonicity}: \qed
\end{proof}

In some cases, this gives us a useful lower bound on the free energy in a set $\Lambda$.
\begin{cor} \label{alphabound} Suppose that $W(e) \geq \alpha$ whenever $e$ is a pair of adjacent
edges in $\mathbb Z^d$. Then for any $\Lambda$, $W(\Lambda) \geq
(|\Lambda|-1) \alpha$.  In particular, $FE^{\Phi}_{\Lambda}(\mu)
\geq (|\Lambda|-1)\alpha$ for any finite, connected subset $\Lambda$
of $\mathbb Z^d$. \end{cor}

\begin{proof} Apply Lemma \ref{subadditive} to the $|\Lambda|-1$ edges in any spanning tree of
$\Lambda$. \qed \end{proof}

We denote by $\mathcal D_{\mathcal L}$ be the space of positive
$\mathcal L \times \tau$-invariant potentials $\Phi$ for which
$W(e)$ is finite for every edge $e$. This is a convenient class of
potentials in which to state a few lemmas; most importantly, for
this purposes of this paper, every $\mathcal L$-invariant perturbed
SAP is in $\mathcal D_{\mathcal L}$.  Since $\Phi$ is $\mathcal
L$-invariant, $W(e)$ is also $\mathcal L$ invariant and hence
assumes only finitely many values on edges $e$ in $\mathbb Z^d$.
Thus, the above corollary applies to all potentials in $\mathcal
D_{\mathcal L}$, which leads us to another corollary.  Let
$\epsilon(\Lambda)$ be the number of edges (i.e., pairs of adjacent
vertices) in the set $\Lambda$ and let $\alpha'=\inf \{ \alpha, 0
\}$, where $\alpha$ is as defined above.

\begin{cor} \label{supercor}  Fix $\Phi$ in $\mathcal D_{\mathcal L}$ and suppose $\mu$ is a
$\mathcal L$-invariant measure on $(\Omega, \mathcal F^\tau)$. Then the value $$FE'_{\Lambda}(\mu)
= FE^{\Phi}_{\Lambda}(\mu) - \epsilon(\Lambda) \alpha'$$ is superadditive in the sense that if
$\Lambda_1$ and $\Lambda_2$ are disjoint but at least one vertex of $\Lambda_1$ is adjacent to a
vertex of $\Lambda_2$, then $FE'_{\Lambda_1 \cup \Lambda_2}(\mu) \geq FE'_{\Lambda_1}(\mu) +
FE'_{\Lambda_2}(\mu)$ for any measure $\mu$. \end{cor}

\begin{proof} Let $e=(x,y)$ be the edge with $x \in \Lambda_1$ and $y \in \Lambda_2$ and apply
Lemma \ref{subadditive} twice, first to the pair $\Lambda_1$ and $\{x,y\}$, and then to the pair
$\Lambda_1 \cup \{x,y\}$ and $\Lambda_2$.  Then the above follows from the fact that
$FE^{\Phi}_{e}(\mu) \geq \alpha'$ and $\epsilon(\Lambda_1 \cup \Lambda_2) \geq \epsilon(\Lambda_1)
+ \epsilon(\Lambda_2) + 1$. \qed \end{proof}

Since $FE'_e(\mu)$ is positive for any edge, it is positive for any
connected set $\Lambda$. By Corollary \ref{supercor}, this implies
that $FE'_{\Lambda}$ is increasing as a function of the connected
set $\Lambda$: that is, $FE'_\Lambda(\mu) \geq FE'_{\Lambda'}(\mu)$
whenever $\Lambda' \subset \Lambda$ and both $\Lambda$ and
$\Lambda'$ are connected. Of course, the $FE'$ defined above is also
invariant under $\mathcal L$. It is not hard to see that there
exists an $N$, depending on $\mathcal L$, such that if $w$ is any
vector of integers in $\mathbb Z^d$, then $N w \in \mathcal L$. Now
consider, as a function of $w$, the value $FE'_{\Lambda_w}(\mu)$
where $\Lambda_w$ is the set of integer vectors $a \in \mathbb Z^d$
with $0 \leq a_i < w_i N$ for $1 \leq i \leq d$. The above corollary
implies that if $c_i = a_i + b_i$ and all other coordinates of $a$,
$b$, and $c$ are equal, then $FE'_{\Lambda_c}(\mu) \geq
FE'_{\Lambda_a}(\mu) + FE'_{\Lambda_b}(\mu)$ for every $\mathcal
L$-invariant $\mu$. This is the usual definition of {\it
superadditivity} for functions on $\mathbb Z^d$, and the following
lemma follows by standard methods (see, e.g. Lemma 15.11 of
\cite{G}):

\begin{lem} \label{SFEexist} If $\Phi \in \mathcal D_{\mathcal L}$, then the value
$SFE'_{\Lambda_w}(\mu) = |\Lambda_w|^{-1}FE'_{\Lambda_w}(\mu)$ with $FE'$ as defined above, tends
to a unique limit, $SFE'(\mu)$ in $[0,\infty]$ as the coordinates of $w$ tend to $\infty$.
\end{lem}

Now, we write $SFE(\mu) = SFE'(\mu) + d\alpha'$ and $SFE_{\Lambda}(\mu) =
\frac{1}{|\Lambda|}FE_{\Lambda}(\mu)$. It is not hard to see that $SFE_{\Lambda_w}(\mu)$ converges
to $SFE(\mu)$ as $w$ tends to $\infty$. We refer to $SFE(\mu)$ as the {\it specific free energy} of
$\mu$.

Note that the limits used in the definition of the specific free
energy assume that we have chosen a specific lattice $\mathcal L$:
we might write $SFE^{\mathcal L}(\mu)$ to denote the specific free
energy with respect to the lattice $\mathcal L$.  However, it is
clear that if $\mu$ is $\mathcal L$-invariant, then when $\mathcal
L$ is replaced by a full rank sublattice $\mathcal L'$, the limits
are not changed, so $SFE^{\mathcal L}(\mu) = SFE^{\mathcal
L'}(\mu)$.  Similarly, if $\mu$ is invariant with respect to any two
full rank lattices $\mathcal L$ and $\mathcal L'$, then we have
$SFE^{\mathcal L}(\mu) = SFE^{\mathcal L \cap \mathcal L'}(\mu) =
SFE^{\mathcal L'}(\mu)$. Next, using the $\mathcal L$-invariance of
$\mu$ and fact that $FE'_{\Delta}(\mu)$ is increasing in $\Delta$,
we have $FE'_{\Lambda_w}(\mu) \leq FE'_{\Lambda_{w+v}}(\theta_x \mu)
\leq FE'_{\Lambda_{w+2v}}(\mu)$, where $v \in \mathbb Z^d$ is the
vector with all of its coordinates equal to $1$. Taking limits and
using Lemma \ref{SFEexist}, it follows that $SFE^{\mathcal L}(\mu) =
SFE^{\mathcal L}(\theta_x \mu)$, where $\theta_x$ is any translation
of $\mathbb Z^d$ (and $x$ is not necessarily in $\mathcal L$).  We
state these facts as a lemma:

\begin{lem} \label{SFEdefequivalence} If $\mathcal L$ and $\mathcal L'$ are two full-rank
sublattice of $\mathbb Z^d$ and $\mu$ is both $\mathcal L$-invariant
and $\mathcal L'$-invariant, then $SFE^{\mathcal L}(\mu) =
SFE^{\mathcal L'}(\mu)$. Moreover, $SFE^{\mathcal L}(\mu) =
SFE^{\mathcal L}(\theta_x(\mu))$ for any $x \in \mathbb Z^d$ and any
$\mu \in \mathcal P(\Omega, \mathcal F^{\tau})$. \end{lem}

We have defined $SFE^{\mathcal L}(\mu)$ for all $\mu \in \mathcal P_{\mathcal L}(\Omega, \mathcal
F^{\tau})$; it is often convenient to extend the definition to all of $\mathcal P(\Omega, \mathcal
F^{\tau})$ by writing $SFE^{\mathcal L}(\mu) = \infty$ whenever $\mu \not \in \mathcal P_{\mathcal
L}(\Omega, \mathcal F^{\tau})$.  We will generally write $SFE(\mu)$ for $SFE^{\mathcal L}(\mu)$,
assuming the choice of lattice to be clear from the context.

Next, it will often be useful to us to have lower bounds on the specific free energy of $\mu$ in
terms of the free energies $FE_{\Lambda}(\mu)$.  Since $\alpha' \leq 0$, we have the following
bound for any $w$:

\begin{lem} \label{FESFEbound} $SFE(\mu) - d\alpha' \geq SFE'(\mu) \geq SFE'_{\Lambda_w}(\mu) \geq
SFE_{\Lambda_w}(\mu) = |\Lambda_w|^{-1} FE_{\Lambda_w}(\mu)$ \end{lem}

We can derive a similar bound involving $FE_{\Delta}$ for a non-rectangular set $\Delta$.  Let
$\Delta$ be a connected subset of $\Lambda_w$.  Then, using repeated applications of Lemma
\ref{subadditive}, we can show that $FE_{\Lambda_w}(\mu) \geq FE_{\Delta}(\mu) + \alpha\[
|\Lambda_w| - |\Delta|\]$. Thus, we have
$$SFE(\mu) \geq |\Lambda_w|^{-1}  \[ FE_{\Delta}(\mu) + \alpha\[|\Lambda_w| - \Delta\]\].$$

In particular, we can say the following:

\begin{lem} \label{SFEbound} For each $\Lambda \subset \mathbb Z^d$, there exist constants $C_1>0$ and
$C_2$ such that $SFE(\mu) \geq C_1 FE_{\Lambda}(\mu) + C_2$. \end{lem}

We can use this fact to check one more important result:

\begin{lem} \label{SEbound} For every constant $C \in \mathbb R$, there exist: \begin{enumerate}
\item A $C_1$ such that $SFE(\mu) \leq C$ implies $\mu\[V_{x,y}(\phi(y) - \phi(x))\] \leq C_1$ for
any adjacent pair $(x,y)$ in $\mathbb Z^d$.  (In fact, we can write
$C_1 = aC + b$ for some constants $a$ and $b$.)
\item A $C_2$ such that $SFE(\mu)\leq C$ implies $|S(\mu)| \leq C_2$ \end{enumerate} \end{lem}

\begin{proof} First, suppose $x$ and $y$ are fixed.  For an appropriate $C'$, we can use Lemma
\ref{SFEbound} to show that $SFE(\mu) \leq C$ implies
$FE_{\{x,y\}}(\mu) = \mathcal H(\mu_{\{x,y\}}, e^{-V_{x,y}}\lambda)
\leq C'$. Writing $\eta = \phi(y) - \phi(x)$ and letting $f$ be the
Radon-Nikodym derivative of $\mu_{x,y}$ with respect to $E^m$, we
can write the latter expression as $$\int f(\eta)\[\log f(\eta) +
V_{x,y}(\eta)\]d\eta \leq C'.$$  (Here $d\eta$ is understood to mean
$d\lambda(\eta)$.)  By Lemma \ref{relativeentropylemma}, there is a
$\beta$ such that $$\int f(\eta)\[\log f(\eta) + \frac{1}{2}
V_{x,y}(\eta)\]d\eta = \mathcal H(f\lambda,
e^{-\frac{1}{2}V_{x,y}}\lambda) \geq \beta$$ for all probability
densities $f$.  Taking the difference of the leftmost terms in the
preceding two equations, we conclude that $\int f V_{x,y}(\eta)
d\eta \leq 2(C' - \beta)$.  Finally, we can compute this last
expression for an edge $(x,y)$ in each of the equivalence classes of
edges modulo $\mathcal L$ and let $C_1$ be the infimum of these
values.

Next, if $SFE(\mu) \leq C$, then $\mu$ must be $\mathcal
L$-invariant; thus, to derive a uniform bound on $S(\mu)$, it is
enough to derive a uniform bound on $|\mu(\phi(y) - \phi(x))| \leq
\mu(|\phi(y) - \phi(x)|)$ for each adjacent pair $(x,y)$ in $\mathbb
Z^d$.  Since $V_{x,y}(\eta)$ increases at least linearly in
$|\eta|$, the existence of a uniform bound on $\mu(|\phi(y) -
\phi(x)|)$ follows immediately from the first part of this lemma.
\qed \end{proof}

\section{Specific free energy level set compactness}
Now that we have the specific free energy for measures in $\mathcal
D_{\mathcal L}$, we can begin to discuss its properties. One of the
most important concerns the {\it level sets} $M_C = \{ \mu| SFE(\mu)
\leq C \}$, as subsets of $\mathcal P_{\mathcal L}(\Omega, \mathcal
F^{\tau})$.  Define the {\it topology of local convergence} on
$\mathcal P(\Omega, \mathcal F^{\tau})$ to be the smallest topology
in which the maps $\mu \mapsto \mu(f)$ are continuous for every
bounded, gradient cylinder function $f$ (i.e., every bounded
function that is $\mathcal F_{\Lambda}^{\tau}$-measurable for some
$\Lambda \subset \subset \mathbb Z^d$) from $\Omega$ to $\mathbb R$.

\begin{lem} \label{levelsetmetric} Each level set $M_C$ of $\mathcal P(\Omega, \mathcal F^{\tau})$,
endowed with the restriction of the topology of local convergence to
that set, is a metric space (i.e., the topology of local convergence
restricted to $M_C$ can be induced by an appropriate metric).
\end{lem} \begin{proof} Let $\{ \Delta_j \}$ be an enumeration of
the connected finite subsets of $\mathbb Z^d$. Let $\delta_i(\mu,
\nu)$ be the distance between the restrictions $\mu_{\Delta_j}$ and
$\nu_{\Delta_j}$ in the metric for the weak topology on $\mathcal
P(\Omega, \mathcal F^{\tau}_{\Lambda})$.  Then $\delta(\mu, \nu) =
\sum_{i = 0}^{\infty} 2^{-i} \inf (1, \delta_i(\mu,\nu))$ is a
metric for the topology of local convergence on $M_C$.  It is clear
that $\mu^i$ converges to $\mu$ in this topology if and only if
$\mu^i_{\Delta_j}$ converges weakly to $\mu_{\Delta_j}$ for every
$j$. \qed \end{proof}

We next prove that $M_C$ is also compact (and hence Polish):

\begin{thm} \label{levelsetcompactness} The level sets $M_C$ are closed and sequentially compact in
the topology of local convergence on $\mathcal P(\Omega, \mathcal F^{\tau})$. Being metrizable,
they are thus compact, and hence Polish (i.e., complete and separable) metric spaces for this
topology. \end{thm}

\begin{proof} First of all, if we are given any sequence $\{ \mu^i \}$ of measures in $M_C$, then
by Lemma \ref{SFEbound}, Lemma \ref{subadditive}, and Lemma \ref{compact}, we can, for any fixed
$\Delta_j$, find a subsequence of $\mu^i$ on which the restrictions $\mu^i_{\Delta_j}$ converge in
the $\tau$-topology to a fixed probability measure $\nu_{\Delta_j}$ on $(\Omega, \mathcal
F_{\Delta_j}^\tau)$. By a standard diagonalization argument, we can take a subsequence such that
$\mu^i_{\Delta_j}$ converges in the $\tau$-topology for {\it each} of the countably many sets
$\Delta_j$ to some probability measure $\nu_{\Delta_j}$.

Kolmogorov's extension theorem for the Polish space $(\Omega,
\mathcal F^{\tau})$ then implies that there exists a measure $\nu
\in \mathcal P(\Omega, \mathcal F^{\tau})$ whose restrictions to the
$\Delta_j$ are in fact these measures. It is clear that this $\nu$
is a limit point of the $\mu_i$ in the topology of local
convergence; every $\Lambda \subset \subset \mathbb Z^d$ is
contained in some $\Delta_j$, and hence every cylinder set $A \in
\mathcal F_{\Lambda}$ is also contained in some $\mathcal
F_{\Delta_j}$ on which $\mu^i_{\Delta_j}$ converges to
$\nu_{\Delta_j}$. Since the restrictions $\nu_{\Delta_j}$ are
clearly $\mathcal L$-invariant---and the sets $\mathcal
F_{\Delta_j}$ generate $\mathcal F$---it follows that $\nu$ is
$\mathcal L$-invariant. Moreover, we must have $SFE(\nu) \leq C$.
If this were not the case, then by Lemmas \ref{FESFEbound} and
\ref{SFEexist} we would have $SFE'_{\Lambda_w}(\nu) > C-d\alpha'$
for some $w$. Now, Lemma \ref{compact} implies that there must be a
$\mu^i$ with $SFE'_{\Delta_j}(\mu^i) > C-d\alpha'$ (for $\Delta_j =
\Lambda_w$).  Applying Lemma \ref{FESFEbound}, this implies that
$SFE(\mu^i) > C$, a contradiction. \qed \end{proof}

\section{Minimizers of specific free energy are Gibbs measures}
\begin{lem} \label{minimizerexists} Whenever $\Phi \in \mathcal D_{\mathcal L}$, there exists an
$SFE^{\Phi}$-minimizing measure in $P_{\mathcal L}(\Omega, \mathcal F^{\tau})$, i.e., a measure
$\mu_0$ in $\mathcal P_{\mathcal L}(\Omega, \mathcal F^{\tau})$ such that, for any other measure
$\mu \in \mathcal P_{\mathcal L}(\Omega, \mathcal F^{\tau})$, we have $SFE^{\Phi}(\mu) \geq
SFE^{\Phi}(\mu_0)$. \end{lem} This minimal value is sometimes called the {\it pressure} of $\Phi$
and denoted $P(\Phi)$.

\begin{proof} Note that $\cap_{C > P(\Phi)} M_C$ is an intersection of non-empty, decreasing
compact sets; hence, it is nonempty. \qed \end{proof} The following is the easy half of our
variational principle. It is not hard to prove this result in more generality, but for simplicity
we will describe only the perturbed simply attractive case.

\begin{thm} \label{minimizersaregibbs} Let $\Phi$ be a perturbed simply attractive potential. If
$\mu$ has minimal specific free energy (with respect to $\Phi$) among $\mathcal L$-invariant
measures with slope $u$, then $\mu$ is a Gibbs measure. \end{thm}

\begin{proof} Suppose that $\mu$ is an $\mathcal L$-invariant measure with finite specific free
energy and slope $u$ that is {\it not} a Gibbs measure. We will show that in this case it is always
possible to modify $\mu$ to produce a measure $\overline{\mu}$ with slope equal to $u$ such that
$SFE^{\Phi}(\overline{\mu}) < SFE^{\Phi}(\mu)$.

If $\mu$ is not Gibbs, then for some $\Lambda$, we have $\mu
\gamma_{\Lambda} \not = \mu$, and hence $\mathcal H(\mu, \mu
\gamma_{\Lambda}) = D >0$. Since $\Phi$ has finite range, there is
an integer $r$ such that $\Phi_{\Delta} = 0$ whenever $\Delta
\subset \subset \mathbb Z^d$ contains two vertices of distance $r$
or more apart.

Now, let $\mathcal L'$ be a sublattice of $\mathcal L$ such that for
any non-zero $i \in \mathcal L'$, each vertex in $\Lambda$ is at
least distance $2r$ from each vertex of $\Lambda + i$.  Then we
define our modified measure: $$\overline{\mu} = \mu \prod_{i \in
\mathcal L'} \gamma_{\Lambda + i}.$$ Although the composition on the
righthand side is infinite, by choice of $\Lambda$, the kernels
$\gamma_{\Lambda+i}$ commute; hence, the order in which the kernels
are applied does not matter. Moreover, the infinite composition
converges in the topology of local convergence, since every
$\Lambda'$ intersects only finitely many $\Lambda+i$ sets.
Informally, it is easy to see why this measure has lower specific
free energy than $\mu$: applying of $\gamma_{\Lambda + i}$ increases
the free energy contained in supersets of $\Lambda+i$ by some fixed
amount.  So, naturally, applying the kernels at a positive fraction
of offsets in $\mathbb Z^2$ should increase the ``free energy per
site'' by a positive amount.  The formal proof that follows is not
very different from well-known proofs of standard (non-gradient)
analogs of this lemma.

Now, as before, take $\Lambda_n = [0,n-1]^d \subset \mathbb Z^d$.
Choose $n$ large enough so that each connected component of
$\Delta_n = \Lambda_n \cap (\mathcal L' + \Lambda)$ is completely
contained in $\Lambda_n$ and has all of its vertices at least $r$
units from the boundary of $\Lambda_n$. (If necessary, by Lemma
\ref{SFEdefequivalence}, we may replace $\mu$ with $\theta_v\mu$ and
$\Lambda$ with $\Lambda + v$ for some $v \in \mathbb Z^d$ in order
to make this possible.)  Fix a reference vertex $v_0 \in \Lambda_n
\backslash \Delta_n$.  We can decompose $\lambda^{|\Lambda_n|-1}$
into the product $\lambda^{|\Lambda_n| - |\Delta_n| - 1} \otimes
\lambda^{|\Delta_n|}$ by taking pairs into the product space to have
the form $(x,y)$ where the components of $x$ are the values $\phi(v)
- \phi(v_0)$ for $v \in \Lambda_n \backslash (\Delta_n \cup \{v_0
\})$ and the components of $y$ are $\phi(v) - \phi(v_0)$ for $v \in
\Delta_n$. For convenience in this proof only, we write $\mu_0 =
e^{-H^o_{\Lambda_n}}\lambda^{|\Lambda_n|-1}$. By Lemma
\ref{conditionalentropy}, there exist regular conditional
probability distributions $\mu_0^x$, $\overline{\mu}^x$, and $\mu^x$
on $\mathcal F_{\Delta_n}^{\tau}$ describing the distribution on $y$
{\it given} the value $x$, when $(x,y)$ has distribution $\mu_0$,
$\overline{\mu}$, and $\mu$, respectively.

Now, we claim the following: $$FE_{\Lambda_n}(\mu) -
FE_{\Lambda_n}(\overline{\mu}) = \mathcal H( \mu_{\Lambda_n}, \mu_0)
- \mathcal H( \overline{\mu}_{\Lambda_n}, \mu_0) = \mu \[\mathcal
H(\mu^x,\mu_0^x) - \mathcal H(\overline{\mu}^x,\mu_0^x)\].$$  The
first equality is true by definition. The second holds follows from
Lemma \ref{conditionalentropy} and the fact that
$\overline{\mu}_{\Lambda_n \backslash \Delta_n} = \mu_{\Lambda_n
\backslash \Delta_n}$.  By our choice of $\mathcal L'$, we know that
both $\mu_0^x$ and $\overline{\mu}^x$ are ($\mu$ almost surely)
products of their restrictions to the components $\Delta^i_n$ of
$\Delta_n$.  Thus, by Lemma \ref{productentropy}, we have that
$$\mathcal H(\mu^x,\mu_0^x) - \mathcal H(\overline{\mu}^x,\mu_0^x)
\geq \sum_i \[ \mathcal H_{\mathcal
F_{\Delta^i_n}^{\tau}}(\mu^x,\mu_0^x) - \mathcal H_{\mathcal
F_{\Delta^i_n}^{\tau}}(\overline{\mu}^x,\mu_0^x) \].$$ Note that
$\overline{\mu}^x_{\Lambda^i_n} = \mu_0^x|_{\Delta^i_n}$.  Hence,
$\mathcal H_{\mathcal F^{\tau}_{\Delta^i_n}}(\overline{\mu}^x,
\mu_0^x) = 0$. Moreover, the $x$ marginals of $\mu$ and
$\overline{\mu}$ coincide, hence
$$\mu(\mathcal H_{\mathcal F^{\tau}_{\Delta^i_n}}(\mu^x, \mu_0^x)) = \mu(\mathcal H_{\mathcal
F^{\tau}_{\Delta^i_n}}(\mu^x, \overline{\mu}^x)) = \mathcal
H_{\mathcal F^{\tau}_{\Delta^i_n}}(\mu, \overline{\mu}) = D,$$ where
the last equality uses the $\mathcal L$-invariance of $\mu$.  For
large $n$, the sum is $D$ times the number of components
$\Delta^i_n$ contained in $\Lambda_n$.  It follows that
$$SFE^{\mathcal L'}(\mu) - SFE^{\mathcal L'}(\overline{\mu}) =
\lim_{n \rightarrow \infty}
\frac{1}{|\Lambda_n|}(FE_{\Lambda_n}(\mu) -
FE_{\Lambda_n}(\overline{\mu})) \geq D/I,$$ where $I$ is the index
of $\mathcal L'$ in $\mathbb Z^d$.

While $\overline{\mu}$ is not necessarily $\mathcal L$ invariant, it is $\mathcal L'$-invariant. We
can make it $\mathcal L$-invariant by replacing it with an average $\overline{\mu}'$ over shifts by
elements of $\mathcal L$ modulo $\mathcal L'$.  By Lemma \ref{relativeentropylemma} and the
definition of the specific free energy, this averaging can only increase the specific free energy.
By Lemma \ref{SEbound}, $\overline{\mu}$ has finite slope; since $\mu$ and $\overline{\mu}$ have
the same laws on $\mathcal F^{\tau}_{\mathcal L' + \Lambda}$, it follows from the definition of
slope on $\mathcal L'$-invariant measures that that $S(\overline{\mu}) = S(\mu)$.  Since
$\overline{\mu}'$ is an $\mathcal L$-invariant measure that is an average of finitely many measures
of this slope, it is also clear that $S(\overline{\mu}') = S(\mu)$. \qed \end{proof}

\chapter{Ergodic/extremal decompositions and SFE} \label{decompositionchapter} In this chapter, we
will cite several standard results about ergodic and extremal
decompositions that we can apply to measures in $\mathcal P(\Omega,
\mathcal F^{\tau})$; in particular, we will see that every $\mathcal
L$-invariant gradient measure $\mu$ can be written, in a unique way,
as a weighted average of $\mathcal L$-ergodic gradient measures.
Moreover, we can compute $SFE(\mu)$ as the weighted average of the
specific free energy of the $\mathcal L$-ergodic components. The
latter result is well known (see Chapter $15$ of \cite{G}) for
ordinary (i.e., non-gradient) Gibbs measures on $\Omega = E^{\mathbb
Z^d}$ when $E$ has a finite underlying measure. However, we must
check that this result remains true for gradient Gibbs measures and
the specific free energy we have constructed in this context.
Throughout this chapter, we assume that a perturbed simply
attractive potential $\Phi$ is fixed; gradient Gibbs measure and
specific free energy are defined with respect to this $\Phi$.

\section{Funaki-Spohn gradient measures}
In this section, we describe an alternate (but equivalent)
formulation of gradient measures (which is also described in detail
in a work of Funaki and Spohn \cite{FS}).  The difference between
the formulation of \cite{FS} and our formulation is largely
cosmetic.  For the purposes of this text, ours is more convenient;
however, the results about $\mathcal L$-ergodic and extremal
decompositions described in \cite{G} and \cite{DS} apply more
directly to the formulation of \cite{FS} than to ours.

The main issue is that several of the basic facts that we will need
about extremal and ergodic decompositions of gradient Gibbs measures
(namely, Lemmas \ref{zeroone}, \ref{ergodiciffextreme},
\ref{extremalconvergence}, \ref{Gibbsdecomposition}, and
\ref{decompositionergodic}) have only been stated and proved in the
literature (for example, in the reference text \cite{G}) for
ordinary Gibbs measures.  Although these results are not terribly
difficult, reproving them individually in the gradient Gibbs measure
context would consume a good deal of space and provide little new
insight.

Instead, we will make a straightforward observation (following
\cite{FS}) that the laws of the gradients (defined below) of
functions sampled from gradient Gibbs measures {\em are} Gibbs
measures---not with respect to an ordinary Gibbs potential, but with
respect to a so-called {\em specification}, described below.
Focusing on these gradients will enable us to cite the above
mentioned lemmas directly from \cite{G} instead of proving them
ourselves.

For any function $\phi \in \Omega$, we define the {\it discrete
gradient} $\nabla \phi: \mathbb Z^d \mapsto E^d$ by
$$\nabla\phi(x) = (\phi(x+e_1) - \phi(x), \ldots, \phi(x+e_d)-\phi(x)),$$
where the $e_i$ are basis vectors of $\mathbb Z^d$.  Write
$\overline E = E^d$.  Denote by $\overline{\Omega}$ the set of
functions from $\mathbb Z^d$ to $\overline E$ and by
$\overline{\mathcal F}$ the Borel $\sigma$-algebra induced on
$\overline{\Omega}$ by the product topology.  Since $\nabla \phi$
only depends on the value of $\phi$ up to an additive constant, each
measure $\mu$ on $(\Omega, \mathcal F^{\tau})$ induces a measure
$\overline{\mu}$ on $(\overline{\Omega}, \overline{\mathcal F})$.

A function $\psi \in \overline{\Omega}$ is called a {\it gradient
function} if it can be written $\nabla \phi$ for some $\phi \in
\Omega$.  In \cite{FS}, the authors characterize the gradient
functions as those functions satisfying the ``plaquette condition''
$$\psi(x)_i + \psi(x+e_j)_i = \psi(x)_j + \psi(x+e_i)_j$$ for all $1
\leq i,j \leq d$ and $x \in \mathbb Z^d$.  (Here $\psi(x)_i$ denotes
the $i$th component of $\psi(x)$.)  Denote by $\overline{\Omega}_G
\subset \overline{\Omega}$ the set of gradient functions.

Instead of taking---as we do---the configuration space to be
$\Omega$ and using a $\sigma$-algebra $\mathcal F^\tau$ that only
measures properties of functions that are invariant under the
addition of a global constant, Funaki and Spohn use
$\overline{\Omega}$ as their configuration space and stipulate
further that all the measures they consider are supported on
$\overline{\Omega}_G$ (i.e., satisfy the plaquette condition almost
surely).

Define the {\em topology of local convergence on $\mathcal
P(\overline{\Omega}_G, \overline{\mathcal F})$} to be the smallest
in which $\mu \mapsto \mu(f)$ is continuous for every bounded
cylinder function $f: \overline{\Omega}_G \mapsto \mathbb R$. This
is analogous to our definition of the topology of local convergence
on $\mathcal P(\Omega, \mathcal F^{\tau})$.  The reader may easily
verify the following:

\begin{lem} The map $\mu \mapsto \overline{\mu}$ described above gives a one-to-one correspondence
between $\mathcal P(\Omega, \mathcal F^{\tau})$ and $\mathcal
P(\overline{\Omega}_G, \overline{\mathcal F}) \subset \mathcal
P(\overline{\Omega}, \overline{\mathcal F})$.  Moreover, the
topology of local convergence on $\mathcal P(\Omega, \mathcal
F^{\tau})$ (as defined in the previous section) is equivalent to the
topology of local convergence on $\mathcal P(\overline{\Omega},
\overline{\mathcal F})$, restricted to $\mathcal
P(\overline{\Omega}_G, \overline{\mathcal F})$. \end{lem}

We extend the definition of specific free energy to $\mathcal P(\overline{\Omega}_G,
\overline{\mathcal F})$ by writing $SFE(\overline{\mu}) = SFE(\mu)$ whenever $\mu \in \mathcal
P(\Omega, \mathcal F^{\tau})$.  Citing Lemma \ref{levelsetcompactness} and Lemma
\ref{levelsetmetric}, we have:

\begin{cor} \label{SFEislowersemicontinuous} $SFE$ is a lower semi-continuous function on $\mathcal
P(\overline{\Omega}_G, \overline{\mathcal F})$ with respect to the topology of local convergence;
moreover, the level sets $M_C =\{ \overline{\mu} | SFE(\overline{\mu}) \leq C \}$ are compact and
metrizable. \end{cor}

If $A \in \mathcal F^{\tau}$, denote by $\overline{A}$ the
corresponding subset of $\overline{\Omega}_G$. Then we can extend
the kernels $\gamma^{\Phi}_{\Lambda}(A, \phi)$ (defined in the first
chapter) to this context by writing
$\gamma^{\Phi}_{\Lambda}(\overline{A}, \nabla \phi) =
\gamma^{\Phi}_{\Lambda}(A, \phi)$.  (Since the latter term is
unchanged when a constant is added to $\phi$, the kernels are well
defined.)

We would like to argue that, in {\em some} sense, $\overline{\mu}$
is a {\em Gibbs measure} if and only if it is preserved by these
kernels (or, equivalently, if the corresponding measure $\mu$ is a
gradient Gibbs measure).   However, the following fact suggests that
this is impossible with the definition of Gibbs measure we presented
in the introduction:

\begin{prop} When $E = \mathbb R^m$ and $\Phi$ is a Gibbs potential (which
admits at least one Gibbs measure) there exists no Gibbs potential
$\overline{\Phi}$ such that $\mu \in \mathcal G^\Phi(\Omega,
\mathcal F^\tau)$ if and only if $\overline{\mu} \in \mathcal
G^{\overline{\Phi}}(\overline{\Omega}, \overline{\mathcal F})$.
\end{prop}

\begin{proof} If $\overline{\mu} \in \mathcal
G^{\overline{\Phi}}(\overline{\Omega}, \overline{\mathcal F})$, and
$\psi$ is sampled from $\overline{\mu}$, then by definition, the law
of $\psi(x)$ --- given the values of $\psi$ at the neighbors of $x$
--- is absolutely continuous with respect to the underlying measure
on $\overline E$, with Radon-Nikodym derivative given by the
Hamiltonian $H_{\{x\}}^{\overline{\Phi}}(\psi)$.  On the other hand,
if $\overline{\mu}$ is supported on gradient functions, then
$\psi(x)$ is completely determined (by the plaquette condition) from
the value of $\psi$ at the neighbors of $x$; thus the conditional
distribution of $\psi(x)$ is supported on a single point.  The only
way the above conditional measure can be absolutely continuous with
respect to an underlying measure $\overline{\lambda}$ on $\overline
E$ is if $\overline \lambda$ has a point mass at that point.  This
cannot happen if $\overline{\lambda}$ is Lebesgue measure. Moreover,
switching to another underlying measure $\overline{\lambda}$ does
not solve the problem.  Since the law of $\nabla \phi(x)$, when
$\phi$ is sampled from a gradient Gibbs measure, is absolutely
continuous with respect to Lebesgue measure, the measure
$\overline{\lambda}$ would have to point masses on a set of positive
Lebesgue measure (and in particular could not be $\sigma$-finite).
But any finite measure of the form $e^{-H_{\Lambda}^\Phi} \prod_{x
\in \Lambda} d\overline\lambda(\psi(x))$---in which the individual
random variables $\psi(x)$ are supported on point masses of
$\overline{\lambda}$---is supported on point masses of the product
$\prod_{x \in \Lambda} \overline\lambda$, and hence is supported on
a countable set. This is not the case for general $\Lambda$ when
$\psi$ is the gradient of a gradient Gibbs measure. \qed
\end{proof}

We will now expand our definition of Gibbs measure.  First, if
$(X,\mathcal X)$ is any probability space and $\mathcal B$ a
sub-$\sigma$-algebra of $\mathcal X$, then a probability kernel
$\pi$ from $\mathcal B$ to $\mathcal X$ is {\em proper} if $\pi(B|
\cdot) = 1_B$ for each $B \in \mathcal B$.  The Gibbs
re-randomization kernels $\gamma_\Lambda$ from $\mathcal T_\Lambda$
to $\mathcal F$ or from $\mathcal T_\Lambda^\tau$ to $\mathcal
F^\tau$, as defined in the introduction, are examples of proper
kernels.

Most of the theorems in \cite{G} are proved for a more general class
of families of proper Gibbs re-randomization kernels called {\it
specifications} (in the sense of sections 1.1 and 1.2 of \cite{G};
see also \cite{FS}) on $\overline{\Omega}_G$.  The following is
Definition 1.23 of \cite{G}:

{\flushleft{\bf Definition:} A {\it specification} on $\mathbb Z^d$
with state space $(\overline E, \overline {\mathcal E})$, is a family $\delta =
\{\delta_\Lambda\}_{\Lambda \subset \subset \mathbb Z^d}$ of proper
probability kernels $\delta_\Lambda$ from $\overline {\mathcal
T}_\Lambda$ to $\overline {\mathcal F}$ which satisfy the
consistency condition $\delta_\Delta \delta_\Lambda = \delta_\Delta$
when $\Lambda \subset \Delta$.  The random fields in the set

\begin{eqnarray}
\mathcal G(\delta) & := & \{ \mu \in \mathcal
P(\overline{\Omega},\overline{\mathcal F}):\mu(A|\overline{\mathcal
T}_\Lambda)
= \delta_\Lambda(A|\cdot)\mu \text{a.s.}\} \\
& & \text{ for all $A \in \mathcal F$ and $\Lambda \subset \subset
\mathbb Z^d$}
\end{eqnarray}
are called {\em Gibbs measures with respect to the specification
$\delta$}.}

Let $\Lambda'$ be the set $\Lambda \cap \left(\Lambda + \sum_{i=1}^d
e_i \right)$.  Note that if the $\nabla \phi$ is known at all $x
\not \in \Lambda$, and the value of $\phi$ is known at some $x \not
\in \Lambda$, then we can deduce the value of $\phi$ at $x \not \in
\Lambda'$; however, this information tells us nothing about the
value of $\phi$ at vertices in $\Lambda'$.

Now, define $\delta_\Lambda$ to be the kernel from $(\overline
\Omega, \overline {\mathcal T}_\Lambda)$ to $(\overline \Omega,
\overline {\mathcal F})$ corresponding to the kernel
$\gamma^\Phi_{\Lambda'}$ on $(\Omega, \mathcal F^\tau)$.  The reader
may verify that these kernels form a specification on $\mathbb Z^d$
with state space $(\overline E, \overline {\mathcal E}) := (E^d,
\mathcal E^d)$.

The following statement follows immediately from the definitions:

\begin{lem} $\mu$ is a Gibbs measure (ergodic measure, extremal Gibbs measure) on $(\Omega,
\mathcal F^{\tau})$ (with respect to the Gibbs specification
$\gamma^\Phi$) if and only if $\overline{\mu}$ is a Gibbs measure
(resp., ergodic measure, extremal Gibbs measure) on
$(\overline{\Omega}, \overline{\mathcal F})$ with respect to the
specification $\delta$.
\end{lem}

\section{Extremal and ergodic decompositions} \label{decompositionsection}

In this section, we cite several results from Chapter 7 and Chapter 14 of \cite{G} (e.g., existence
of extremal and ergodic decompositions), all of which apply to both ordinary measures and gradient
measures.  Throughout this section, we assume that a gradient potential $\Phi$ is given.  In each
case, the gradient analog of the statement follows from the cited non-gradient result by the
correspondence described in the previous section.

As in the first chapter, we say a measurable subset $A$ of $\Omega$
is a {\it tail event}, if $A \in \mathcal T = \cap_{\Lambda \subset
\subset \mathbb Z^d} \mathcal F_{\mathbb Z^d - \Lambda}$.  We say
$A$ is an {\it $\mathcal L$-invariant event} if it is preserved by
translations by members of $\mathcal L$; both $\mathcal T$ and the
set $\mathcal I_{\mathcal L}$ of $\mathcal L$-invariant events are
$\sigma$-algebras. (See Proposition 7.3, Corollary 7.4, and Remark
14.3 of \cite{G}.) We denote the sets of extremal and ergodic Gibbs
measures by $\ex\mathcal G(\Omega, \mathcal F)$ and $\ex\mathcal
G_{\mathcal L}(\Omega, \mathcal F)$, respectively; we sometimes
abbreviate $\mathcal G(\Omega, \mathcal F)$ and $\mathcal
G_{\mathcal L}(\Omega, \mathcal F)$ by $\mathcal G$ and $\mathcal
G_{\mathcal L}$, respectively.  Similarly, the set of extremal and
ergodic gradient Gibbs measures are written, respectively,
$\ex\mathcal G(\Omega, \mathcal F^{\tau})$ and $\ex\mathcal
G_{\mathcal L}(\Omega, \mathcal F^{\tau})$; we abbreviate $\mathcal
G(\Omega, \mathcal F^{\tau})$ and $\mathcal G_{\mathcal L}(\Omega,
\mathcal F^{\tau})$ by $\mathcal G^{\tau}$ and $\mathcal G_{\mathcal
L}^{\tau}$.

We say $\mu$ is {\it trivial} on a $\sigma$-algebra $\mathcal A$ if $\mu(A) \in \{0, 1 \}$ for all
$A \in \mathcal A$.  Now, we cite the following characterization of extremal and $\mathcal
L$-ergodic measures in terms of their behavior on tail and $\mathcal L$-invariant events,
respectively (Theorems 7.7 and 14.5 of \cite{G}):

\begin{lem} \label{zeroone} The following hold for all Gibbs measures on $(\Omega, \mathcal F)$
\begin{enumerate} \item A probability measure $\mu \in \mathcal P_{\mathcal L}(\Omega, \mathcal F)$
is extreme in $\mathcal P_{\mathcal L}(\Omega, \mathcal F)$ if and only if $\mu$ is trivial on
$\mathcal I_{\mathcal L}$. \item A Gibbs measure $\mu\in \mathcal G$ is extreme in $\mathcal G$ if
and only $\mu$ is trivial on $\mathcal T$.  Distinct extremal Gibbs measures $\mu_1$ and $\mu_2$
are mutually singular in that there exists an $A \in \mathcal T$ with $\mu_1(A)=0$ and
$\mu_2(A)=1$.
\item A Gibbs measure $\mu \in \mathcal G_{\mathcal L}$ is $\mathcal L$-ergodic
if and only if $\mu$ is trivial on $\mathcal I_{\mathcal L}$. Distinct $\mathcal L$-ergodic
measures are mutually singular in that there exists an $A \in \mathcal I_{\mathcal L}$ with
$\mu_1(A) = 0$ and $\mu_2(A)=1$. \end{enumerate} Analogous statements are true for gradient
measures $\mu \in \mathcal P_{\mathcal L}(\Omega, \mathcal F^{\tau})$. \end{lem}

There may exist extremal Gibbs measures that are not $\mathcal L$-ergodic and $\mathcal L$-ergodic
Gibbs measures that are not extremal.  However, as the following lemma shows, every {\it extremal}
$\mathcal L$-invariant measure is necessarily $\mathcal L$-ergodic.  (Theorem 14.15 of \cite{G}.)
\begin{lem} \label{ergodiciffextreme} A Gibbs measure $\mu \in \mathcal G_{\mathcal L}$ is an extreme point of the convex set
$\mathcal G_{\mathcal L}$ if and only if $\mu$ is $\mathcal
L$-ergodic.  That is, $$\ex \mathcal G_{\mathcal L} = \mathcal
G_{\mathcal L} \cap \ex \mathcal P_{\mathcal L}(\Omega, \mathcal
F).$$ Furthermore, $\mathcal G_{\mathcal L}$ is a face of $\mathcal
P_{\mathcal L}(\Omega, \mathcal F)$. That is, if $\mu, \nu \in
\mathcal P_{\mathcal L}(\Omega, \mathcal F)$ and $0<s<1$ are such
that $s \mu + (1-s)\nu \in \mathcal G_{\mathcal L}$, then $\mu, \nu
\in \mathcal G_{\mathcal L}$. Analogous statements are true for
gradient measures. \end{lem}

Given a single observation $\phi$ from an extremal Gibbs measure or
an $\mathcal L$-ergodic measure $\mu \in \mathcal P(\Omega, \mathcal
F)$, it is $\mu$-almost-surely possible to reconstruct $\mu$ from
$\phi$ a way that we will now describe.  Whenever
$\gamma_{\Lambda_n}(\cdot|\phi)$ has a limit in the topology of
local convergence as $n$ tends to $\infty$, we denote this limit by
$\pi^{\phi}$. Let $\{\Lambda_n\}$ be any increasing sequence of
cubes in $\mathbb Z^d$ such that $|\Lambda_n| \rightarrow \infty$.
We denote by $\pi_{\mathcal L}^{\phi}$ the ``shift-averaged''
measure given by $$\pi_{\mathcal L}^{\phi}(A) = \lim_{n \rightarrow
\infty} |\Lambda_n \cap \mathcal L|^{-1} \sum_{x \in \Lambda_n \cap
\mathcal L} 1_A(\theta_x \phi)$$ when this limit exists.  We can
extend the functions $\pi^\cdot$ and $\pi_{\mathcal L}^\cdot$ to
functions from $\Omega$ to $\mathcal G(\Omega, \mathcal F)$ and
$\mathcal P_{\mathcal L}(\Omega, \mathcal F)$, respectively, by
setting them equal to some arbitrary $\nu_0$ in (respectively)
$\mathcal G(\Omega, \mathcal F)$ and $\mathcal P_{\mathcal
L}(\Omega, \mathcal F)$ when these limits fail to exist. The
following lemma makes precise our ability to recover $\mu$ from a
single observation.  (The first half is \cite{G}, Theorem 7.12. The
second follows from \cite{G}, Theorem 14.10.)

\begin{lem} \label{extremalconvergence}The following are true: \begin{enumerate} \item If $\mu \in
\ex\mathcal G$, then $\mu (\{ \phi \in \Omega : \pi^{\phi} = \mu \}) = 1$.
\item If $\mu \in \ex \mathcal P_{\mathcal L}(\Omega, \mathcal F)$, then $\mu (\{ \phi \in \Omega :
\pi^{\phi}_{\mathcal L} = \mu \}) = 1$.  \end{enumerate}  Analogous statements are true of gradient
measures. \end{lem}

Next, we would like to say that each measure in $\mathcal G$ (respectively, $\mathcal G_{\mathcal
L}$) is a weighted average of extremal (respectively, ergodic) measures.  In order to precisely
define a ``weighted average'' of elements in $\ex\mathcal G$ and $\ex\mathcal G_{\mathcal L}$, we
need to define $\sigma$-algebras on these sets of measures. To do this, for each $A \in \mathcal
F$, consider the evaluation map $e_{A}: \mu \mapsto \mu(A)$. Denote by $e(\ex \mathcal G)$ the
smallest $\sigma$-algebra on $\ex\mathcal G$ with respect to which each $e_A$ is measurable. Define
$e(\ex \mathcal G_{\mathcal L})$ similarly.  The following decomposition theorem shows that
$\mathcal G$ and $\mathcal G_{\mathcal L}$ are isomorphic to the simplices of probability measures
on $\ex\mathcal G$ and $\ex\mathcal G_{\mathcal L}$, respectively.  (See \cite{G}, Theorem 7.26 and
Theorem 14.17.)

\begin{lem} \label{Gibbsdecomposition} For each $\mu \in \mathcal G$ there exists a unique weight
$w_{\mu} \in \mathcal P(\ex\mathcal G, e(\ex\mathcal G))$ such that for each $A \in \mathcal F$,
$$\mu(A) = \int_{\ex \mathcal G} \nu(A) w_{\mu}(d \nu).$$ The mapping $\nu \mapsto w_{\mu}$ is a
bijection between $\mathcal G$ and $\mathcal P(\ex \mathcal G, e(\ex \mathcal G))$. Furthermore,
$w_{\mu}$ has the same law as the image of $\mu$ under the mapping $\phi \mapsto \pi^{\phi}$. These
results remain true when $\mathcal G$ is replaced by $\mathcal G_{\mathcal L}$ and $\pi^{\phi}$ is
replaced by $\pi_{\mathcal L}^{\phi}$. Analogous decompositions exist for gradient measures.
\end{lem}

\begin{lem} \label{decompositionergodic} For each $\mu \in \mathcal P_{\mathcal L}(\Omega, \mathcal
F)$ there exists a unique weight $$w_{\mu} \in \mathcal P(\ex \mathcal P_{\mathcal L}(\Omega,
\mathcal F), e(\ex \mathcal P(\Omega, \mathcal F)))$$ such that for each $A \in \mathcal F$,
$$\mu(A) = \int_{\ex \mathcal P_{\mathcal L}(\Omega,
\mathcal F)} \nu(A) w_{\mu} (d\nu).$$ The mapping $\mu \mapsto
w_{\mu}$ is an bijection between $\mathcal P_{\mathcal L}(\Omega,
\mathcal F)$ and $\mathcal P(\ex \mathcal P_{\mathcal L}(\Omega,
\mathcal F), e(\text{ex} \mathcal P(\Omega, \mathcal F)))$.
Furthermore, $w_{\mu}$ gives the law for the image of $\mu$ under
the mapping $\phi \mapsto \pi^{\phi}_{\mathcal L}$. Analogous
decompositions exist for gradient measures. \end{lem}

In less formal terms, the lemmas state that sampling $\phi$ from
$\mu \in \mathcal G$ is equivalent to: \begin{enumerate} \item First
choosing an extremal measure $\mu_0$ from an extremal decomposition
measure.
\item Then choosing $\phi$ from $\mu_0$. \end{enumerate} Similarly,
sampling $\phi$ from $\mu \in \mathcal P_{\mathcal L}(\Omega,
\mathcal F)$ is equivalent to: \begin{enumerate} \item First
choosing an $\mathcal L$-ergodic measure $\mu_0$ from an ergodic
decomposition measure.
\item Then choosing $\phi$ from $\mu_0$.
\end{enumerate} Note that an $\mathcal L$-ergodic Gibbs measure $\mu \in \mathcal G_{\mathcal L}$
may or may not be extremal. In Chapter \ref{clusterswappingchapter} and \ref{discretegibbschapter},
we will be interested not only in classifying $\mathcal L$-ergodic gradient Gibbs measures but also
in determining how each $\mathcal L$-ergodic gradient Gibbs measure decomposes into extremal
components.

\section{Decompositions and $SFE$} \label{SFEdecompositionsection}

The fact that $SFE\[\int_{\ex\mathcal P_{\mathcal L}(\Omega,
\mathcal F)} \nu w_{\mu} (d\nu)\] = w_{\mu}(SFE)$ is essential to
our proof of the large deviations principle and the variational
principle.  A longer, more detailed proof of this result---which
follows Chapter 15 of \cite{G}---is given in Appendix
\ref{SFEdecompositionchapter}. An advantage of the longer proof is
that it also yields results of independent interest, including one
interpretation of $SFE(\mu)$ in terms of the ``conditional free
energy'' of one fundamental domain of $\mathcal L$---conditioned on
its ``lexicographic past''---and another interpretation involving
discrete derivatives.  We do not use these interpretations anywhere
else (hence, their relegation to the appendix). The proof described
here follows \cite{DS}.

\begin{lem} \label{SFEisaffine} The function $\mu \mapsto SFE(\mu)$ is affine. \end{lem}
\begin{proof} We follow the proof given in Exercise 4.4.41 of \cite{DS}.  Suppose $\rho = a\mu +
(1-a)\nu$ with $0<a<1$. Recall the definition
$$SFE(\rho) = \lim_{n \rightarrow \infty} |\Lambda_n|^{-1}\lim_{n \rightarrow \infty}
FE_{\Lambda_n}(\rho) = |\Lambda_n|^{-1} \mathcal H
(\rho_{\Lambda_n}, e^{-H^o_{\Lambda_n}}\lambda^{|\Lambda_n -1|}).$$

For this proof only, denote $\pi_n =
e^{-H^o_{\Lambda_n}}\lambda^{|\Lambda_n -1|}$.  Taking the limit,
the convexity of $SFE$ follows immediately from the convexity of
$\mathcal H(\cdot,\pi_n)$. To prove concavity, it is enough to show
that $$\frac{1}{|\Lambda_n|} \mathcal H(\rho, \pi_n) \geq
\frac{1}{|\Lambda_n|} \[a \mathcal H(\mu, \pi_n) + (1-a) \mathcal
H(\nu, \pi_n) \] + o(1).$$ If either of $\mu_{\Lambda_n}$ or
$\nu_{\Lambda_n}$ fails to be absolutely continuous with respect to
$\pi_n$, then both sides of the inequality are equal to infinity.
Otherwise, let $f_n$ and $g_n$ be the Radon-Nikodym derivatives of
$\mu_{\Lambda_n}$ and $\nu_{\Lambda_n}$, respectively, with respect
to $\pi_n$.  Then we can rewrite the inequality:
$$ - \pi_n \[(a f_n + (1-a) g_n) \log (a f_n+ (1-a) g_n) \] \geq $$ $$- \pi_n \[ a f_n \log f_n +
(1-a)g_n \log g_n \] + o(|\Lambda_n|).$$  In fact, a simple identity
(see Exercise 4.4.41 of \cite{DS}) states that if $f_n$ and $g_n$
assume any values in $[0,\infty)$ and $0 < a < 1$, then
$$-(a f_n + (1-a) g_n) \log (a f_n+ (1-a) g_n) \geq - \[ a f_n \log f_n + (1-a)g_n \log g_n \] +
 \frac{|f_n - g_n|}{e}.$$  Since $\pi_n(f_n) = \pi_n(g_n) = 1$, the integral of the ``error
term'' (the last term on the right in the above expression) is at
most $2$; in particular, it is $o(|\Lambda_n|)$. \qed \end{proof}

In fact, $SFE$ is also ``strongly affine'' in the following sense:

\begin{thm} \label{SFEisstronglyaffine} If $\mu$ can be written $$\mu = \int_{\ex\mathcal
P_{\mathcal L}(\Omega, \mathcal F)} \nu w_{\mu} (d\nu),$$ then
$$SFE(\mu) = \int_{\text{ex}\mathcal P_{\mathcal L}(\Omega, \mathcal
F^{\tau})} SFE(\nu) w_{\mu} (d\nu) = w_{\mu}(SFE).$$ \end{thm}

The above Theorem clearly follows from Lemma \ref{SFEisaffine},  Lemma
\ref{SFEislowersemicontinuous}, and the following lemma (applied to the gradient field
$\overline{\mu}$).

\begin{lem} If $F: \mathcal P_{\mathcal L}(\Omega, \mathcal F) \mapsto [0, \infty]$ is
lower-semicontinuous (with respect to a topology in which the level
sets $\{ \mu : F(\mu) \leq C \}$ are metrizable) and affine, then it
is strongly affine.  That is, $$F \[ \int_{\text{ex}\mathcal
P_{\mathcal L}(\Omega, \mathcal F)} \nu w_{\mu}(d\nu) \] =
\int_{\text{ex}\mathcal P_{\mathcal L}(\Omega, \mathcal F)}
F(\nu)w_{\mu}(d\nu).$$ \end{lem}

\begin{proof} The proof is identical to the proof of Lemma 5.2.24 of \cite{DS}. \qed \end{proof}

\begin{cor} \label{SFEisanexpectation} The function $\alpha:\Omega \mapsto \mathbb R$, defined by
$\alpha(\phi) = SFE(\pi_{\mathcal L}^{\phi})$ is bounded below, and satisfies $$SFE(\mu) =
\mu(\alpha) = \mu(SFE(\pi_{\mathcal L}^\cdot))$$ for all $\mu \in \mathcal P_{\mathcal L}(\Omega,
\mathcal F)$. \end{cor}

\begin{proof} By Lemma \ref{extremalconvergence}, the function $\mu \mapsto \mu(\alpha)$ agrees
with the function $\mu \mapsto SFE(\mu)$ when $\mu$ is $\mathcal L$-ergodic.  Since $\alpha$ is
$e(\mathcal P_{\mathcal L}(\omega, \mathcal F^{\tau}))$ measurable (it is a limit of the
$e(\mathcal P_{\mathcal L}(\omega, \mathcal F^{\tau}))$-measurable functions $\mu \mapsto
|\Lambda_n|^nFE_{\Lambda}(\mu_n)$), it follows from Lemma \ref{Gibbsdecomposition} that the
statement proved in Theorem \ref{SFEisstronglyaffine} for $\mu \mapsto SFE(\mu)$ applies to the
function $\mu \mapsto \mu(\alpha)$ as well.  \qed \end{proof}

\chapter{Surface tension and energy} \label{surfacetensionchapter} Define the {\it surface tension}
$\sigma:\mathbb R^{d \times m} \mapsto \mathbb R$ by writing
$$\sigma(u) = \inf_{\mu \in \mathcal P_{\mathcal L}(\Omega, \mathcal
F^{\tau}), S(\mu) = u} SFE(\mu).$$ We will give another equivalent
definition of surface tension (as a normalized limit of log
partition functions on tori) in Chapter
\ref{LDPempiricalmeasurechapter}. In this chapter, we will make
several elementary observations about the function $\sigma$ and the
set of slopes on which $\sigma$ is finite.  We also discuss some
basic facts about the existence of finite energy functions
satisfying boundary conditions. Further discussion of this and
related problems can be found in many texts on linear programming
and network flows; see , e.g., \cite{AMO}.  Unless otherwise noted,
we will assume throughout this chapter, for notational simplicity,
that $m=1$.  The extensions to higher dimensions $m > 1$ are in all
cases straightforward.

Throughout this chapter, we assume that an $\mathcal L$-invariant
perturbed simply attractive potential $\Phi$ is given.  Denote by
$U_{\Phi}$ the interior of the region on which $\sigma$ is finite.
Since $SFE$ is affine, it is clear that $\sigma$ is convex; in
particular, this implies that $U_{\Phi}$ is a convex open set and
that $\sigma$ is continuous on $U_{\Phi}$ and (provided $U_{\Phi}$
is non-empty) equal to $\infty$ outside of the closure
$\overline{U_{\Phi}}$ of $U_{\Phi}$. We refer to $U_{\Phi}$ as the
space of {\it allowable} slopes.  Some of the results in subsequent
chapters (including the second half of the variational principle)
will apply only to slopes in $U_{\Phi}$.  In the next section, we
make some preliminary constructions that will be necessary in the
continuous cases $E = \mathbb R$ and $E = \mathbb R^m$, and which we
will apply mainly to the case that $\Phi$ is a perturbed ISAP.

\section{Energy bounds for ISAPs when $E = \mathbb R$} \label{overlineV}
\subsection{Bounding free energy and surface tension in terms of $V$}

We will assume throughout this section that $m=1$ and $V: \mathbb R
\mapsto \mathbb R^+$ is a convex, symmetric (i.e., $V(x) = V(-x)$)
difference potential, and that $\Phi_V$ is the corresponding ISAP.
We begin with following question: what does the shape of $V$ tells
us about the shape of $\sigma$?

To begin to answer this, we extend the definition of $V$ to vectors
by writing $V(u) = \sum_{i=1}^d V(u)$.  Note that $V(u)$ gives the
specific energy of the deterministic singleton measure $\mu \in
\mathcal P(\Omega, \mathcal F^\tau)$ supported on the plane
$\phi_u(x) = (x,u)$ (defined up to additive constant). Thus, if we
replaced $SFE$ by $SE$ in the definition of $\sigma$, we would have
exactly $\sigma(u) = V(u)$ (by Jensen's inequality).

However, $SFE(\mu)$ for the $\mu$ defined above is infinite, since
the entropy of $\mu$ is $-\infty$.  Nonetheless, if $\mu$ is instead
as the law of $\phi = \phi_u + \phi'$, where the law $\mu'$ of
$\phi'$ has slope zero, specific entropy $\log \epsilon$, and
$|\phi_1(x)-\phi(y)| \leq \epsilon$ $\mu'$-a.s.\ (for some constant
$\epsilon$; see the next subsection for the existence of such a
measure $\mu'$), then $SFE(\mu) \leq \sup \{V(u'):|u' - u|_\infty
\leq \epsilon\} - \log \epsilon$, where $|u'-u|_\infty =
\sup_{i=1}^d |u_i - u'_i|$. Applying this to a particular choice of
$\mu'$ yields the following:

\begin{lem} \label{sigmalessthanV}
When $\Phi_V$ is an ISAP, $E = \mathbb R$, and $V$ grows at most
exponentially fast, there exist constants $C_1>0$ and $C_2$
(depending only on $V$) such that $\sigma(u) \leq C_1V(u)+C_2$.  If
$\Phi_V$ is any ISAP, then $\sigma(u) \leq \sup \{V(u'):|u' -
u|_\infty \leq \epsilon \} - \log \epsilon$ for all $\epsilon$.
\end{lem}

To get a bound in the other direction, note that $SFE(\mu) =
\sigma(u)$ for some $\mu$ with $S(\mu) = u$, and Lemma \ref{SEbound}
then implies $\mu \[V(\phi(y) - \phi(x))\] \leq C_1\sigma(u) + C_2$
for some constants $C_1$ and $C_2$ when $x$ and $y$ are adjacent
(with the constants depending on $V$, but not on $\mu$, $x$, or $y$)
which in turn implies (by Jensen's inequality) the following:

\begin{lem} \label{Vlessthansigma}
When $\Phi_V$ is an ISAP and $E = \mathbb R$, there exist constants
$C_1>0$ and $C_2$ (depending only on $V$) such that $V(u) \leq
C_1\sigma(u)+C_2$.
\end{lem}

Combining the lemmas implies that $\sigma$ and $V$ generate
identical Orlicz-Sobolev norms (discussed in detail in Chapter
\ref{orliczsobolevchapter}) when $V$ grows at most exponentially
fast. When $V$ increases super-exponentially, however, this need not
be the case. One of the most important cases in which the first part
of Lemma \ref{sigmalessthanV} fails is the hard constraint model in
which $V(\eta)=0$ if $\eta \in [-1,1]$ and $\infty$ otherwise. It is
easy to see that $\sigma(u)$ tends to $\infty$ as $u$ tends to the
boundary of $[-1,1]^d$, even though $V(u)$ is constant for $u \in
[-1,1]^d$.

We will now define a function $\overline{V}$, called the {\bf
wedge-normalization of $V$} such that $\overline{V}(u)$ and
$\sigma(u)$ do agree up to a constant factor.  In the hard
constraint model above, it would be natural to guess that the second
estimate in Lemma \ref{sigmalessthanV} is tight, so that on the
interval $[-1,1]$, $\overline{V}(\eta)$ should be approximately
(some constant times) the log of the distance from $\eta$ to the
nearest of the endpoints $-1$ and $1$. In fact, the particular
expression for $\overline{V}(\eta)$ we use is $V(\eta)$ minus four
times the log of the distance from $F(\eta)$ to either $0$ or $1$,
where $F(\eta)$ is the fraction of the mass of the measure
$e^{-V(\eta)}d\eta$ (where $d\eta$ is Lebesgue measure on $\mathbb
R$) that lies to the left of $\eta$. (The particular form is
convenient because it allows us express free energy with respect to
$\overline{V}$ in terms of relatively entropy with respect to
certain ``wedge measure,'' defined below.)

In addition to obtaining information about $\sigma$ (which will
later be useful in determining the topologies in which surface-shape
large deviations principles apply), we will show that if
$\phi:\Lambda \rightarrow \mathbb R$ is such that the
nearest-neighbor sum $\sum\overline{V}\phi(x)-\phi(y))$ is equal to
$C$, then there is measure describing a random small perturbation of
$\phi$ whose free energy (with respect to $V$) is at most a constant
times $C$.

To motivate the construction of these random perturbations, we
mention that they will be useful in Chapter \ref{LDPchapter} when we
engage in the mechanical process of proving lower bounds on the
Gibbs measure of the set of functions $\phi: D_n \rightarrow \mathbb
R$ that approximate a particular function $f$ on $D$; one technique
will involve constructing a $\phi_0$ that approximates $f$, is
piecewise linear on large pieces of $D_n$, and has a $\overline{V}$
energy that we can control. Then we can get a lower bound by
estimating the Gibbs measure of the set of $\phi$ that are
``pertubations'' of $\phi_0$ in the non-linear regions and allowed
to vary more freely on the linear regions (which are dealt with
separately using the various subadditivity limits of Chapter
\ref{LDPempiricalmeasurechapter}), see Figure \ref{mostlylinear}.

The reader who is only interested in the case that $V$
increases at most exponentially fast may skip the proofs of the
lemmas in the remainder of this section (since in this case, one may
take $\overline{V} = V$ and the results are still obviously true).
The reader who is only interested in the case $E = \mathbb Z$ need
not read this section at all.

\subsection{Defining box measures and $\overline{V}$}

Given $\Lambda$, a {\it box measure} $\mu$ on the space of functions
$\phi$ on $\Lambda$ is a uniform measure on a set $B = \{ \phi|
\phi_1 < \phi < \phi_2 \}$ where $\phi_1$ and $\phi_2$ are also real
functions on $\Lambda$. (We write $f < g$ if $f(x) < g(x)$ for each
$x\in\Lambda$.)  An upper bound on the free energy of $\mu$ is given
by $\sup_{\phi \in B} H_{\Phi}(\phi) - \sum_{\eta \in \Lambda} \log
(\phi_2(\eta) - \phi_1(\eta))$. In the remainder of this text, when
we need to prove that a reasonably low (or at least non-infinite)
free-energy Gibbs measure exists with certain properties, we will
sometimes construct them explicitly using box measures.  In this
section we will construct a convex function $\overline{V}$, based on
$V$, and use the shorthand $\Phi = \Phi_V$ and $\overline{\Phi} =
\Phi_{\overline{V}}$.  Both $\Phi$ and $\overline{\Phi}$ will be
ISAPs.  We derive an upper bound on
$FE_{\Lambda}^{\overline{\Phi}}(\mu)$ in terms of
$FE_{\Lambda}^{\Phi}(\mu)$ and show that, for constants $C_1$ and
$C_2$, there always exists a box measure centered at $\phi$ with
$\Phi$ free energy at most $C_1 H^{\overline{\Phi}}_{\Lambda}(\phi)
+ C_2$.  The first part of our construction involves showing that
the relative entropy of a measure $\nu$ on $[0,1]$ with respect to a
certain ``wedge-density'' measure is not much more than twice the
relative entropy of $\nu$ with respect to a uniform measure.

\begin{lem} \label{wedgebound} Let $\mu_1$ be the uniform distribution on $[0,1]$ and $\mu_2$ the
measure with density given by the wedge-shaped function $g(\eta) = 2-4|\eta - \frac{1}{2}|$.  If
$\mathcal H(\nu,\mu_1)$ is finite, then $\mathcal H(\nu,\mu_2)$ is finite also and furthermore, for
any $C_1 > 2$, there exists an $C_2$ for which $\mathcal H(\nu,\mu_2) \leq C_1 \mathcal
H(\nu,\mu_1) + C_2$ for any Lebesgue-measurable $\nu$ on $[0,1]$. \end{lem}

\begin{proof} Given $\beta$, by Lemma \ref{relativeentropylemma}, we can compute the density
function $f$ on $[0,1]$ that minimizes $\int(f(\eta) \log f(\eta)) d\eta - \beta \int f(\eta) \log
g(\eta)d\eta$: we use the fact that this expression is equal to $\mathcal H(\mu,\nu) + \log C$
where $\mu$ has density $f$ and $\nu$ has density $C g^{\beta}$ with $C^{-1}= \int g^{\beta}(\eta)
dx$ (which makes sense provided $\beta > -1$). The minimum occurs when $\mu = \nu$, i.e., when $f =
g^{\beta}$ (up to constant multiple).  In particular, this minimum is finite when $-1 < \beta < 0$.
We can rewrite this statement (using the fact that $(1-\beta) + \beta = 1$) as follows: whenever
$-1 < \beta < 0$, there exists an $\alpha$ for which
$$(1-\beta) \int f(\eta) \log f(\eta)d\eta + \beta\[\int f(\eta) \log f(\eta)d\eta - \int
f(\eta)\log g(\eta)d\eta \] \geq \alpha.$$ Dividing by $\beta$ (and
recalling that $\beta < 0$) we have:
$$ \int f(\eta) \log \frac{f(\eta)}{g(\eta)}d\eta \leq \frac{\alpha}{\beta} + \frac{\beta-1}{\beta}
\int f(\eta) \log f(\eta)d\eta $$
$$\mathcal H(\nu, \mu_2) \leq \frac{\alpha}{\beta} + \frac{\beta-1}{\beta}\mathcal H(\nu, \mu_1).$$
By taking $\beta$ close to $-1$, we can make $\frac{\beta-1}{\beta}$ arbitrarily close to $2$. \qed
\end{proof}

Now let $Z = \int_{-\infty}^{\infty} e^{-V(\eta)}d\eta$ and $F(\eta)
= Z^{-1} \int_{-\infty}^{\eta}e^{-V(\eta)}d\eta$.  If $X$ is the
random variable on $\mathbb R$ with distribution given by $F$, then
$F(X)$ is uniform on $[0,1]$.  Now, we define $\overline{V}(\eta) =
V(\eta) - \log g(F(\eta))$.  Writing $\overline{Z} =
\int_{-\infty}^{\infty}e^{-\overline{V}(\eta)}d\eta$, it is also
easy to check that $\overline{Z}^{-1}e^{-\overline{V}}$ is the
density function for the random variable $F^{-1}(G^{-1} F(X))$,
where $G(\eta) = \int_0^\eta g(\zeta)d\zeta$ is the distribution
function for $g$.

\begin{lem} \label{someoverlineVfacts}
The following are true for the $\overline{V}$ defined above: \begin{enumerate} \item
$\overline{V}(\eta) \geq V(\eta) - \log 2$. \item If $(a,b)$ is the
interval on which $V$ is finite (here $a,b \in \mathbb R \cup
\{-\infty, \infty \}$), then $$\lim^+_{\eta \mapsto a}
\overline{V}(\eta) = \lim^-_{\eta \mapsto b} \overline{V}(\eta) =
\infty.$$ \item $\overline{V}$ is convex.
\end{enumerate} \end{lem}

\begin{proof} Recall that $\overline{V}(\eta) - V(\eta) = -\log g(F(\eta))$. The first item follows
because $g \leq 2$, the second because $g(F(\eta))$ tends to zero as
$\eta$ tends to $a$ or $b$.

For the third item, since we are assuming that $V$ is convex, it is
enough to prove that $$A(\eta) = \overline{V}(\eta) - V(\eta) =
-\log g(F(\eta)) = - \log g \int_{-\infty}^\eta e^{-V(\zeta)}dy$$ is
also convex. From the symmetries of $V$ and $g$, we can see that
$A(\eta) = A(-\eta)$; $A$ is strictly increasing for $\eta < 0$ and
decreasing for $\eta > 0$. It is therefore sufficient to show that
$A$ is convex on the interval $(0, \infty)$.  On this interval,
using the definition of $g$ given above, we can write $$A(\eta) = -
\log 4(1 - F(\eta)) = - \log 4 - \log \int_\eta^{\infty}
e^{-V(\zeta)}dy.$$ Since this function is continuous, it is enough
to check that for any $0 < \epsilon < \eta$, we have $$2A(\eta) \leq
A(\eta + \epsilon) + A(\eta - \epsilon),$$ or equivalently (by
Fubini's theorem) we must show $$ \int_\eta^{\infty}
\int_\eta^{\infty} e^{-V(a)-V(b)}dadb -
\int_{\eta+\epsilon}^{\infty} \int_{\eta-\epsilon}^{\infty}
e^{-V(a)-V(b)}dadb = \geq 0.$$ Canceling the common region of
integration, we can rewrite the left hand side as
$$\int_\eta^{\infty} \int_\eta^{\eta+\epsilon} e^{-V(a)-V(b)}dadb - \int_{\eta-\epsilon}^{\eta}
\int_{\eta+\epsilon}^{\infty} e^{-V(a)-V(b)}dadb.$$ Relabeling
variables, this becomes:
$$\int_\eta^{\infty} \int_\eta^{\eta+\epsilon} e^{-V(a)-V(b)}dadb - \int_{\eta+\epsilon}^{\infty}
\int_{\eta-\epsilon}^\eta e^{-V(a)-V(b)}dadb = $$
$$\int_\eta^{\infty} \int_\eta^{\eta+\epsilon} e^{-V(a)-V(b)} -
e^{-V(a+\epsilon)- V(b-\epsilon)}dadb.$$ Since $a - b + \epsilon
\geq 0$ throughout the region of integration, we claim that
$V(a+\epsilon)+V(b-\epsilon)-V(a)-V(b) \geq 0$.  To see this, we can
write by convexity: $$V(a) \leq \frac{a-b+\epsilon}{a-b + 2
\epsilon} V(a+\epsilon) +
\frac{\epsilon}{a-b+2\epsilon}V(b-\epsilon)$$ $$V(b) \leq
\frac{a-b+\epsilon}{a-b + 2 \epsilon} V(b-\epsilon) +
\frac{\epsilon}{a-b+2\epsilon}V(a+\epsilon).$$ Summing these two
lines gives the claim. From this, we conclude that $$e^{-V(a)-V(b)}
- e^{-V(a+\epsilon) - V(b-\epsilon)} \geq 0$$ throughout the region
of integration, and the lemma follows. \qed
\end{proof}

Since we will use it again in Chapter \ref{clusterswappingchapter},
we state the simple fact we used at the end of the proof as a lemma.
\begin{lem} \label{convexidentity} If $a-b + \epsilon > 0$, and
$V:\mathbb R \mapsto \mathbb R$ is convex, then $V(a+\epsilon) +
V(b-\epsilon) \geq V(a) + V(b)$.
\end{lem}

\subsection{Bounding $FE^{\overline{\Phi}}$ in terms of $FE^{\Phi}$}
\begin{lem} \label{edgefebound} For each $C_1 > 2$, there exists a $C_2$ such that for any edge $e$
and measure $\mu$ on $\Omega$, we have $$FE^{\Phi}_e(\mu) - C_2 \leq
FE^{\overline{\Phi}}_e(\mu) \leq C_1 FE^{\Phi}_e(\mu) + C_2.$$
\end{lem}

\begin{proof}  The first half of the inequality follows immediately
from \ref{someoverlineVfacts}.  Next, since $F$ is a continuously
differentiable, increasing one-to-one function, one easily checks
that if $Y$ is any other real-valued random variable, then the
relative entropy $\mathcal H(Y,X)$ is equal to $\mathcal
H(F(Y),F(X))$.  (Here, when $X$ and $Y$ are real-valued random
variables, we write $\mathcal H(Y,X)$ to mean $\mathcal H(\mu_Y,
\mu_X)$, where measures $\mu_Y$ and $\mu_X$ are the laws of $Y$ and
$X$ on $\mathbb R$.)  Define $X$ and $Z$ as in the previous section
and recall by definition of free energy that if $Y= \phi(y) -
\phi(x)$, then $FE_e^{\Phi}(\mu) = FE^{V}(\mu_Y)$, where
$FE^{V}(\mu_Y) + \log Z = \mathcal H(Y,X)$; and
$FE^{\overline{\Phi}}_e (\mu) = FE^{\overline{V}}(\mu_Y)$, where
$FE^{\overline{V}}(\mu_Y) + \log \overline{Z} = \mathcal
H(Y,F^{-1}G^{-1}F(X))$.  We conclude from Lemma \ref{wedgebound}
that
$$FE_{\overline{V}}(\mu) + \log \overline{Z} = \mathcal H(Y,F^{-1}G^{-1}FX) = \mathcal H(F(Y),
G^{-1}F(X)) \leq$$ $$ C_1 \mathcal H(F(Y), F(X)) + C_2 = C_1\mathcal
H(Y,X) + C_2 = C_1 FE_{V}(\mu) + C_2 + \log Z.$$ \qed \end{proof}

\begin{lem} \label{vbarfebound} For any $C_1 > 2d+1$, there is a constant $C_2$ such that for any
finite connected $\Lambda \subset \mathbb Z^d$ and any measure $\mu$
on $\Omega$, we have $$FE^{\Phi}_{\Lambda}(\mu) - C_2|\Lambda| \leq
FE^{\overline{\Phi}}_{\Lambda}(\mu) \leq C_1
FE^{\Phi}_{\Lambda}(\mu) + C_2|\Lambda|.$$ \end{lem}

\begin{proof} The first half of the inequality follows immediately
from \ref{someoverlineVfacts}.  Next, we claim that that $\sum_{e
\subset \Lambda} FE^V_e(\mu) \leq 2d FE^V_{\Lambda}(\mu)$.  We can
prove this by invoking Lemma \ref{subadditive} if we assume that
$\alpha$ used to define $FE'$ in the setup to that lemma is equal to
zero; for simplicity, we will assume this to be the case (the reader
may check that without the assumption, the result remains true with
a different value of $C_2$).  Given Lemma \ref{subadditive}, it is
enough to show that there is a sequence $e_1, e_2, \ldots,
e^{|\Lambda|-1}$ of edges that forms a spanning tree of $\Lambda$
and satisfies
$$\sum_{i=1}^{|\Lambda|-1} FE^V_{e_i}(\mu) \geq \frac{1}{2d} \sum_{e
\subset \Lambda} FE^V_e(\mu).$$ The latter point is a simple graph
theoretic fact that is easily verified with a greedy algorithm:
choose $e_1$ to be an edge $e$ for which $FE^V_{e}(\mu)$ is maximal.
Choose each subsequent $e^i$ to be an edge $e$ for which
$FE^V_{e}(\mu)$ is maximal among those edges that do not, together
with $e_1, \ldots, e_{i-1}$, form cycles.  The number of edges which
form cycles at the $i$th step is at most the number of edges that
are incident to vertices of $e_1, \ldots, e_{i-1}$---this number is
bounded above by $2d(i-1)$.  Hence $FE^{V}_{e_i} \geq
FE^{V}_{e_{(2di+1)}}(\mu)$ where the latter represents the
$(2di+1)$th edge when the edges are listed in rank order of free
energy values. Since $\sum_{i=1}^{(|\Lambda|-1)/2d)}
FE^{V}_{e_{(2di+1)}}(\mu) \geq \frac{1}{2d} \sum_{e \subset \Lambda}
FE^V_e(\mu)$, this fact follows.

Using this fact in the last step, and taking $C_1'$ and $C_2'$ to be
the constants from Lemma \ref{edgefebound}, we now see that
$$FE^{\overline{\Phi}}_{\Lambda}(\mu) - FE^{\Phi}_{\Lambda}(\mu) =
\mu H^{\overline{\Phi}}_{\Lambda} - \mu H^{\Phi}_{\Lambda} = \sum_{e
\subset \Lambda} \mu(\overline{V}(e) - V(e)) = \sum_{e \subset
\Lambda} FE^{\overline{V}}_e(\mu) - FE^V_{e}(\mu) \leq $$ $$\sum_{e
\subset \Lambda} ((C_1'-1)FE^V_e(\mu) + C_2') \leq 2d(C_1'-1)
FE^V_{\Lambda}(\mu) + 2d |C_2'| |\Lambda|.$$ We conclude the proof
by adding $FE_{\Lambda}^{\Phi}(\mu)$ to both sides of the above
equation sequence and setting $C_1 = 2d(C_1'-1)+1$ and $C_2 =
2d|C_2'|$, recalling that $C_1'$ can be chosen arbitrarily close to
$2$. \qed \end{proof}

Similarly, we have:

\begin{lem} \label{overlineVSFEbound} There exist $C_1$ and $C_2$ such that
$$SFE^{\Phi}(\mu) - C_2 \leq SFE^{\overline{\Phi}}(\mu) \leq C_1 SFE^{\Phi}(\mu) + C_2$$ for all
$\mu$. \end{lem}

\begin{proof} This follows immediately from Lemma \ref{vbarfebound}; simply take limits, using the
definition of $SFE$. \qed \end{proof}

\begin{cor}
\label{overlineVsigmabound} There exist $C_1$ and $C_2$ such that
$$C_1 \sigma^{\Phi}(u) - C_2 \leq \sigma^{\overline{\Phi}}(u) \leq C_1 \sigma^{\Phi}(u) + C_2$$ for all
$u$.
\end{cor}

\subsection{Using $\overline{\Phi}$ to bound box measure free energies}
\begin{lem} \label{boxmeasureexistence} There exist constants $C_1$ and $C_2$ such that whenever
$\Lambda$ is a connected subset of $\mathbb Z^d$ and
$H^{\overline{\Phi}}_{\Lambda}(\phi) = \epsilon$, there exists a box
measure $\overline{\mu}$ on $\mathbb R^{|\Lambda|}$ centered at
$\phi$ whose gradient measure $\mu$ satisfies
$FE^{\Phi}_{\Lambda}(\mu) \leq C_1 \epsilon + C_2|\Lambda|$.
\end{lem}

\begin{proof} Fix $\eta>0$ and let $(-a_\eta,a_\eta)$ be the largest interval for which $\zeta \in
(-a_\eta,a_\eta)$ implies $V(\zeta) \leq V(\eta) + 1$.  Then we
claim that for some $C_1>0$, independent of $\eta$ and $V$, we have
$-\log|a_\eta - \eta| \leq C_1 A(\eta) =
C_1(\overline{V}(\eta)-V(\eta))$. For $\eta>0$, recall from that
$A(\eta) = - \log 4 - \log \int_\eta^{\infty} e^{-V(\zeta)}d\zeta$.
Now, by convexity and the fact that $V$ is positive, for any
$\zeta>a_\eta$, we have $V(\zeta) \geq \frac{1}{a_\eta -
\eta}(\zeta-a_\eta)$. This implies that
$$\int_{a_\eta}^{\infty}e^{-V(\zeta)}d\zeta \leq
\int_{0}^{\infty}e^{-\frac{\zeta}{a_\eta-\eta}}d\zeta =
(a_\eta-\eta).$$

Furthermore, since $V$ is positive, $\int_{\eta}^{a_\eta}
e^{-V(\zeta)} \leq a_\eta-\eta$.  It follows that $$e^{-A(\eta)} =
4\int_\eta^{\infty} e^{-V(\zeta)}d\zeta \leq 8(a_\eta - \eta).$$
Hence, $$A(\eta) \geq - \log 8 - \log(a_\eta -\eta).$$ Now, for each
$v \in \Lambda$, choose $\epsilon(v)$ to be the minimum of
$\frac{1}{16}e^{-A(\phi(v') - \phi(v))}$ for all $v'$ adjacent to
$v$.  By the construction of $a_x$, the definition of $\epsilon(v)$,
and the bound in the previous paragraph, for all $\psi$ in the box
$\{\psi|\phi(v)-\epsilon(v) < \psi < \phi(v) + \epsilon(v), v \in
\Lambda \}$, and any edge $e$, we have $H^{\Phi}_{e}(\phi) <
H^{\Phi}_e(\psi)$. Summing over all edges, we have
$H^{\Phi}_{\Lambda}(\psi) < H^{\Phi}_{\Lambda}(\phi) + 2d|\Lambda|$.
Now, let $\mu$ be the box measure that corresponds to fixing
$\psi(v_0)=\phi(v_0)$ for a reference vertex $v_0$ and choosing the
other $\psi(v)$ independently, uniformly in $(\phi(v)-\epsilon(v),
\phi(v) + \epsilon(v))$. Since the widths of the intervals are
$2\epsilon(v)$, this energy bound and Lemma
\ref{freeenergymonotonicity} together imply
$$FE^{\Phi}_{\Lambda}(\mu) \leq H_{\Lambda}^{\Phi}(\phi)+ 2d|\Lambda|
- \sum_{v \in \Lambda, v \not = v_0} \log (2\epsilon(v)).$$

Now, for each $v$, we have $-\log 2(\epsilon(v)) \leq
\sup_{|v-v'|=1, v' \in \Lambda} A(\phi(v')-\phi(v))$.  How does the
sum of these supremum values compare to the sum over all edges? We
know that $A(\eta) + \log 2 \geq 0$, for any $\eta$.  Thus,
$$FE^{\Phi}_{\Lambda}(\mu) - H^{\Phi}_{\Lambda}(\phi) -2d|\Lambda|
\leq$$ $$\sum_{v \in \Lambda, v \not = v_0} \sup_{|v-v'|=1, v' \in
\Lambda} A(\phi(v')-\phi(v))+ \log 2 \leq $$ $$\sum_{v \in \Lambda,
v \not = v_0} \sum_{|v-v'|=1, v' \in \Lambda}
A(\phi(v')-\phi(v))+\log 2 \leq $$ $$2\sum_{e=(v_1,v_2) \subset
\Lambda}[A(\phi(v_1)-\phi(v_2)] + 2|\Lambda|\log 2 = $$
$$2[H^{\overline{\Phi}}_{\Lambda}(\phi) - H^{\Phi}_{\Lambda}(\phi)]
+ 2|\Lambda|\log 2.$$ The lemma thus holds with $C_1 = 2$ and $C_2 =
2 \log 2 + 2d$. \qed \end{proof}

Applying similar analysis to the case $\phi = \phi_u$ and $\Lambda =
\mathbb Z^d$ yields the following analog of Lemma
\ref{sigmalessthanV}.

\begin{lem} \label{infiniteboxmeasureexistence} There exist constants $C_1$ and $C_2$
such that for all $u$ there exists a box measure $\overline{\mu}$
centered at $\phi_u$ whose gradient measure $\mu$ satisfies
$SFE^{\Phi}(\mu) \leq C_1 V(u) + C_2$.  In particular, $\sigma^\Phi
(u) \leq C_1 \overline{V}(u) + C_2$.
\end{lem}

Lemma \ref{Vlessthansigma} and \ref{overlineVsigmabound} imply the
other direction (that $\overline{V}(u) - C_1 \leq \sigma(u)$).
Hence, we have the following:

\begin{lem} \label{overlineVsigmaVbound} There exist $C_1$ and $C_2$ depending only on $V$ such that
$$\overline{V}(u) - C_2 \leq \sigma^{\Phi}(\mu) \leq
C_1 \overline{V}(u) + C_2.$$
\end{lem}

\section{Torus approximations} \label{torussection}
In this section, we allow $\Phi$ to be an SAP or a perturbed SAP and
allow either $E = \mathbb R$ or $E = \mathbb Z$. One way to compute
$U_{\Phi}$ and approximate $\sigma$ is using the finite torus $T$
given by $\mathbb Z^d$ modulo $\mathcal L$.  A function $g$ from the
directed edges of $T$ to $\mathbb R$ is called a {\it local gradient
function} if the sum of $g$ along any directed null-homotopic cycle
in $T$ is zero. If $g$ were extended periodically to $\mathbb Z^d$,
then it would be the discrete derivative of a height function
$\phi_g$ (defined up to additive constant); we define the {\it
slope} (or {\it homology class}) of $g$, written $S(g)$, to be the
asymptotic slope of the corresponding function $\phi_g$.  Observe
that $\psi_g(x) = \phi_g(x) - (S(g),x)$ is $\mathcal L$-periodic.

Write $H^{\Phi}(g) = \sum V(\phi_g(y) - \phi_g(x))$ where the sum
ranges over a set of edges in $\mathbb Z^d$ corresponding to the
edges $T$.  Write $\chi(u) = \inf_{S(g) = u} H^{\Phi}(g)$ and
similarly $\overline{\chi}(u) = \inf_{S(g) = u}
H^{\overline{\Phi}}(g)$. From Lemma \ref{boxmeasureexistence} we can
deduce the following:

\begin{lem} There exist constants $C_1$ and $C_2$ such that $\sigma(u) \leq C_1 \overline{\chi}(u)
+ C_2$. \end{lem}

In particular, $\sigma(u)$ is finite whenever $\overline{\chi}(u)$ is finite (which is true
whenever $\chi(u)$ is finite). It is easy to see that the converse holds.  Hence, we have the
following:

\begin{lem} \label{finiteenergyU} If $E = \mathbb R$ and $\Phi$ is an ISAP, then the set $U_{\Phi}$ is
the interior of the set of slopes $u$ for which $\chi(u)$ is finite.\end{lem}

We can easily generalize the above result the case $E = \mathbb Z$.
Let $\Phi$ be the nearest neighbor potential with edge potentials
given by $V_{x,y}$.  In this case, we let define
$\overline{V}_{x,y}: \mathbb R \mapsto \mathbb R$ to be the largest
convex function which agrees with $V_{x,y}$ on the integers. This is
a linear interpolation of $V$ in between integers at which $V_{x,y}$
is finite; it is equal to $\infty$ outside of the interval spanned
by the integers at which $V_{x,y}$ is finite.

Now, given any (real) $g$ on $G$ for which $H^{\overline{\Phi}}$ is
defined, we can extend this to a real-valued $\phi_g$.  Now, let
$\mu$ be the measure on $(\Omega, \mathcal F)$ that returns $\lfloor
\phi_g + \epsilon \rfloor$ where $\epsilon$ is chosen uniformly in
$[0,1]$.  It is easy to see that the restriction of this measure to
$(\Omega, \mathcal F^{\tau})$ is an $\mathcal L$-invariant measure
of slope $u$, and that $\mu(\overline{\Phi}) =
H^{\overline{\Phi}}/I$ where $I$ is the index of $\mathcal L$ in
$\mathbb Z^d$.  In the discrete setting, the specific relative
entropy of a singleton measure is zero and the specific entropy of
any non-singleton measure is strictly greater than zero.  We may
conclude that $SFE(\mu) < H^{\overline{\Phi}}/I$ and $\sigma(u) <
\chi(u)$ whenever $\chi(u)$ is finite.  Also, if $\sigma(u)$ is
finite, then there exists a slope-$u$ measure $\mu$ for which
$SFE(\mu)$ is finite; let $g$ be the local gradient function $T$ for
which $\phi_g(y) - \phi_g(x) = \mu(\phi(y) - \phi(x))$.  It is
easily seen that $H^{\overline{\Phi}}(g)$ must be finite in this
case.  Thus $\sigma(u)$ is infinite whenever $\chi(u)$ is infinite.
The reader may verify the generalization to all simply attractive
potentials.

\begin{lem} If $E = \mathbb R$, and $\Phi$ is an SAP, then the set $U_{\Phi}$ is the interior of the
set of slopes $u$ for which $\chi(u)$ is finite.  If $E = \mathbb
Z$, then $U_{\Phi}$ is the interior of the set of slopes $u$ for
which $\overline{\chi}(u)$ is finite. \end{lem}

When $E = \mathbb R$ and $\Phi$ is an ISAP, then $\chi(u)$ is finite
precisely when $\overline{\chi}(u)$ is finite; but if $\Phi$ is not
an ISAP, then we have not defined $\overline{\chi}(u)$.  However, if
$E = \mathbb Z$, then our definition of $\overline{\chi}(u)$ applies
for all SAPs.

The reader may also check the following: \begin{lem}
\label{convexhull} If $E = \mathbb Z$ and $\Phi$ is a Lipschitz
simply attractive potential, then set $U_{\Phi}$ is the interior of
the convex hull of the set $\mathcal S$ of integer slopes $u$ for
which $\chi(u)$ is finite. \end{lem}

The following is also trivial to verify: \begin{cor} If all of the
convex nearest neighbor potentials $V_{x,y}: \mathbb R \mapsto
\mathbb R$ are everywhere finite, then $\sigma$ is everywhere finite
and $U_{\Phi} = \mathbb R^d$. \end{cor}

Finally, we would like to define homology-class-restricted Gibbs
measures on a torus, using a method similar to the one presented in
\cite{FS}.  Let $x_0,x_1,\ldots, x_{j-1}$ be representatives of a
fundamental domain of $\mathbb Z^d$ modulo $\mathcal L$; these
vertices are in one-to-one correspondence with the elements of $T$.
Every local gradient function $g$ is determined by its homology
class and the values of the function $\phi_g$ (with additive
constant chosen so that $\phi_g(x_0)=0$) at $x_1, x_2, \ldots,
x_{j-1}$.   Now, we can define a measure on gradient $\mu^u$ on
local gradient functions $g$ of homology class $u$ by
$$Z_T^{-1} e^{-H^{\Phi}(g)} \prod_{i=1}^{j-1} d\phi_g(x_i),$$ where $Z_T$ is
the normalizing constant that makes the above a
probability measure, and the connection between $g$ and $\phi_g$ is
determined by the value $u$.

Now, choose $k$ so that $k\mathbb Z^d \subset \mathcal L$, and
define $\mu^n_u$ to be the measure produced as above on the torus
$T_n = [0,kn-1]^d$ that one gets by replacing $\mathcal L$ by
$nk\mathbb Z^d$.  (If $E = \mathbb Z^d$, then when defining
$\mu^n_u$, we replace $u$ with $\frac{1}{n} \lfloor un \rfloor$ in
order to ensure that the slope $u$ homology class of finite-energy
local gradient functions is non-empty.)

\begin{lem} \label{torusconvergence} If $u \in U_{\Phi}$, then some subsequence of measures
$\mu^n_u$ converges in a topology of local convergence to a $\mathcal L$-invariant gradient Gibbs
measure $\mu \in \mathcal P_{\mathcal L}(\Omega, \mathcal F^{\tau})$ with finite specific free
energy, which is less than or equal to $\liminf_{n \rightarrow \infty} - |T_n|^{-1} \log Z_{T_n}$.
\end{lem}

\begin{proof} The proof is very similar to the one in the proof of Lemma
\ref{levelsetcompactness}---namely, one observes that there is an
upper bound (independent of $n$) on the free energy of $\mu^n_u$
restricted to a set $\Lambda$, then uses Lemma \ref{compact} and
Lemma \ref{weaktauequivalent} to prove convergence for the
$\mu^n_u$'s restricted to a set $\Lambda \subset \mathbb Z^d$ (which
can be treated as a subset of any sufficiently large torus), and
then uses the standard diagonalization argument to get convergence
for all $\Lambda$.  \qed \end{proof}

Statements similar to the above are proved for Ginzburg-Landau
models in \cite{FS} and for domino tiling models in \cite{CKP} and
\cite{CEP}.

\section{Distance functions and interpolations}
In this section, we give another perspective on the set $U_{\Phi}$
and discuss the problem of interpolating functions defined on
subsets of $\mathbb Z^d$ to finite-energy functions defined on all
of $\mathbb Z^d$.  The basic ideas are simple and very standard (see
\cite{AMO} for more exposition and many additional references; look
for the keywords {\it cycle decomposition} and {\it Dijkstra's
algorithm}).

If $x$ and $y$ are neighboring vertices (as always, $x$ preceding $y$ lexicographically) and
$V_{x,y}$ is the pair potential connecting the two, then define $d(x,y) = \sup_{\eta \in E}
V_{x,y}(\eta) < \infty$ and $d(y,x) = -\inf_{\eta \in E} V_{x,y}(\eta) < \infty$.  By definition,
$d(x,y)$ is an upper bound on the value that $\phi(y)-\phi(x)$ can take for any $\phi$ with
$\Phi_{\{x,y\}}(\phi) < \infty$. For simplicity, will assume here that the interval on which
$V_{x,y}$ is finite is closed (i.e., $V_{x,y}$ is lower semicontinuous), so the upper bound is
always achievable (this assumption does not affect the Gibbs kernels, since changing a convex
$V_{x,y}$ to make it lower semi-continuous---by altering its values at the endpoints---only alters
the potential on a set of Lebesgue measure zero). If $x$ and $y$ are not adjacent vertices but are
connected by a path $P = (x=p_0, p_1, p_2, \ldots, p_k = y)$ of vertices, we define $d_P(x,y)$ to
be $\sum_{j=0}^{k-1} d(p_j, p_{j+1})$ and define $D(x,y)$ to be the minimal value of $d_P(x,y)$ as
$P$ ranges over all paths connecting $x$ and $y$. We assume that there are no negative
cycles---i.e., no paths $P$ from a vertex $x$ to itself with $d_P(i,x)<0$. (Otherwise, all height
functions on regions containing $P$ would have infinite energy.)

A classical theorem (see \cite{AMO}) is the following:

\begin{lem} \label{extension} If we fix the values of a function $\phi'$ (real if $E$ is $\mathbb R$,
discrete if $E = \mathbb Z$) at all vertices in $\Lambda \subset
\mathbb Z^d $, then there exists a finite-energy height function
$\phi$ defined on all of $\mathbb Z^d$ {\it extending $\phi'$}
(i.e., satisfying $\phi(x) = \phi'(x)$ for $x \in \Lambda$) if and
only if $D(x,y) \geq \phi'(y) - \phi'(x)$ for each $x,y \in
\Lambda$.
\end{lem}

\begin{proof} The proof is straightforward---simply observe that $\phi(y) = \inf_{x \in \Lambda}
\phi'(x) + D(x,y)$ is such an extension function; in fact, this is the maximal such function. \qed
\end{proof}

A similar argument gives the following:

\begin{lem} \label{torusexistence} There exists a finite energy local gradient function $g$ on a
torus $T$ with slope $u$ if and only if there exists no path $P =
\{p_0, p_1, \ldots, p_k = p_0 \}$ in $T$ (with no vertices repeated
except the first/last value) such that the lifting $\tilde{P} =
\{\tilde{p}_0, \tilde{p}_1, \ldots, \tilde{p}_k \}$ of $P$ to the
covering space $\mathbb Z^d$ satisfies $(u,\tilde{p}_k -
\tilde{p}_0) < D_{\tilde{P}}(\tilde{p}_0, \tilde{p}_k)$.  \end{lem}

\begin{proof} Note that if $E = \mathbb R$, then the gradient function $g$ has slope $u$ on $T$ if and
only the function $f(x) = g(x) - (u ,x)$ has slope zero (and hence
extends to a periodic function on $\mathbb Z^d$); also, $g$ has
finite energy with respect to $\Phi$ if and only if $f$ has finite
energy with respect to the nearest neighbor potential $\Phi'$
defined by $V'_{x,y} (\eta) = V_{x,y}
\[ (u,y - x) + \eta \]$.  Now, the paths $P$ described above correspond precisely to negative
cycles with respect to $d$ defined using the potential $\Phi'$; it follow from Lemma
\ref{extension} that they fail to exist if and only if there exists a finite energy, zero-slope
height function $f$ (and hence a finite energy local gradient function $g$ of slope $u$).

If $E = \mathbb Z$, then the integer-valued gradient function $g$
has slope $u$ on $T$ if and only the integer-valued function $f(x) =
g(x) - \lfloor (u ,x) \rfloor $ has slope zero (and hence extends to
a periodic function on $\mathbb Z^d$); also, $g$ has finite energy
with respect to $\Phi$ if and only if $f$ has finite energy with
respect to the nearest neighbor potential $\Phi'$ defined by
$V'_{x,y} (\eta) = V_{x,y} \[ \lfloor (u,y) \rfloor - \lfloor (u,x)
\rfloor + \eta \]$.  The argument proceeds as in the real case. \qed
\end{proof}

Let $P_1, P_2, \ldots, P_k$ be the finitely many non-self-intersecting cycles in $T$.  The above
lemma and Lemma \ref{finiteenergyU} imply the following:

\begin{lem} \label{halfspace} The set $U_{\Phi}$ is interior of the set of values $u$ for which
$(u,x_i-y_i) \leq D_{P_i}(x_i,y_i)$ for $1 \leq i \leq k$, where
here $x_i$ and $y_i$ are the starting and ending points of $P_i$.
In particular, if $U_{\Phi}$ is not the entire space $\mathbb R^d$,
then it is the intersection of finitely many half spaces.  If each
$V_{x,y}$ is equal to $\infty$ outside of a finite interval, then
$U_{\Phi}$ is the interior of a convex polyhedron.
\end{lem}

A finite energy function $\phi$ defined on all of $\mathbb Z^d$ is
said to be {\it upward taut} if for some $x \in \mathbb Z^d$, and
some $C$, every finite energy $\phi':\mathbb Z^d \mapsto E$ which
agrees with $\phi'$ on all but finitely many places satisfies
$\phi'(x) \leq C$; $\phi$ is {\it downward taut} if every such
$\phi'$ satisfies $\phi'(x) \geq C$.

\begin{lem} A finite energy $\phi: \mathbb Z^d \mapsto E$ is {\it upward taut} if and only if for some
$P_i$ and some $x$, the expression $k D_{P_i}(x_i,y_i) - \phi\[x + k(y_i - x_i) \]$ is bounded
below for all $k \geq 0$.  Similarly, $\phi$ is {\it downward taut} if and only if for some $P_i$,
the expression $k D_{P_i}(y_i,x_i) - \phi\[x + k(y_i - x_i) \]$ is above for all $k \geq 0$.
\end{lem}

\begin{proof} By Lemma \ref{extension}, $\phi$ is upward taut if and only if for some $x$, there
exists a sequence $y_j$ arbitrarily far away from $x$ for which $\phi(y) - \phi(x) \geq D(x,y_j) -
C$.  Let $Q_j$ be the corresponding paths connecting $x$ and $y_j$.  Taking a subsequential limit
of these paths gives a path $Q$ from $x$ to $\infty$ with $\phi(y) - \phi(x) \geq D(x,q_j) - C$ for
every $q_j$ in $Q$.  The path $Q$, projecting down to the torus $T$, must completely traverse at
least one of the paths $P_i$ infinitely often.  We refer to such a path as a {\it $C$-taut path}.
(A zero-taut path is called simply a {\it taut path}.) Let $x'$ be the first time in $Q$ a vertex
equal to the initial vertex of $P_i$, modulo $T$, occurs.  Let $z_1$ and $z_2$ be the first and
last points of the first time the path $P_i$ occurs. Then we can replace $Q$ with a new $C$-taut
path by moving the translating the path $P_i$ so that it begins at $x$, and translating the path
between $x'$ and $z_1$ so that it begins at $x + (y_i - x_i)$.  Repeating this process, we can find
arrange for $Q$ to begin with an arbitrarily long sequence of $P_i$'s.  Again, we find that the
path $P$ consisting of an infinite sequence of $P_i$'s is also a $C$-taut path.  A similar argument
applies to downward taut functions. \qed \end{proof}

\begin{lem} \label{ergodicnottaut} If $\Phi$ is simply attractive and $\mu$ is an $\mathcal
L$-ergodic gradient Gibbs measure with finite specific free energy and $S(\mu) = u \in U_{\Phi}$,
then $\mu$-almost all $\phi$ are not taut. \end{lem}

\begin{proof} If $\phi$ is taut, then for some $P_i$ and $C$, there the ergodic theorem implies
that there almost surely exist a positive fraction of vertices which
are the beginnings of infinite $C$-taut paths formed by
concatenating $P_i$'s end to end.   We claim that each of these
paths must in fact be $0$-taut; suppose otherwise.  Then along one
of the infinite sequence $x + \mathbb Z (y_i - x_i)$, there would be
have to be a place where, for some $\epsilon$, $x$ was the last
value in the sequence to be the beginning of an $\epsilon$-taut
path. But only one such event can occur in each path; such events
cannot occur with $\mu$ positive probability because $\mu$ is
invariant under translation by $k(y_i-x_i)$, for $k \in \mathbb Z$.

A similar argument shows that in fact every $x' \in x + \mathcal L$ must belong to an infinite taut
path of this form; it follows that the slope $u$ of $\mu$ satisfies $(u,x_i-y_i) =
D_{P_i}(x_i,y_i)$---and hence, it lies on the boundary of $U_{\Phi}$, by Lemma \ref{halfspace}.
\qed \end{proof}

Note that from the proof, we also have:

\begin{lem} \label{tautpaths} Let $E = \mathbb Z$, and $u \in \partial U_{\Phi}$ and let $\mu$ be a
$\mathcal L$-invariant measure of finite specific free energy and slope $u$.  Then for some $x$ and
some $P_i$, the path formed by concatenating infinitely many copies of $P_i$ starting at any $x+y$,
with $y \in \mathcal L$, is almost surely taut.  In fact, for each face $X$ of the polyhedron
$\partial U_{\Phi}$, there exists at least one $P_i$ such that the preceding statement is true for
$P_i$ whenever $u \in X$. \end{lem}

\begin{proof}  Let $P_i$ be a path which determines the face $X$ (as in Lemma \ref{halfspace}).
\qed \end{proof}

We can also conclude that if $\mu \in \mathcal P_{\mathcal
L}(\Omega, \mathcal F^{\tau})$ is any gradient measure supported on
finite energy functions $\phi$, then $\phi(x_2)-\phi(x_1)$ is $\mu$
almost surely constant whenever $x_1$ and $x_2$ lie in a path $P$ of
the type described above.  It is easy to see that such a measure
must have infinite specific free energy if $E = \mathbb R$; but it
will have finite specific free energy if $E = \mathbb Z$.  We can
lift the path $P$ defined on $T$ to an infinite path in $\mathbb
Z^d$ along which all height differences are ``frozen'' by $\mu$.
Whenever $\mu$ has slope $u$ which lies on the boundary of
$U_{\Phi}$, we say that $\mu$ is {\it taut} and that $u$ is a {\it
taut slope}.  It is not hard to check the following:

\begin{lem} \label{conttobound} Let $\Phi$ be a simply attractive potential. If $E = \mathbb Z$, then
$\sigma$ extends continuously to a real-valued function on the
closure of $U_{\Phi}$. If $E = \mathbb R$, then $\sigma(u)$ tends to
infinity as $u$ approaches the boundary of $U_{\Phi}$ and $\sigma(u)
= \infty$ on the boundary itself. \end{lem}

We can approximate $D(x,y)$ as follows: write $D_{\Phi}(x) = \sup_{u \in U_{\Phi}} (u,x)$.

\begin{lem} \label{distancenorm}  There exists a constant $C$ such that $|D_{\Phi}(y-x) - D(x,y)| <
C$ for all $x$ and $y$ in $\mathbb Z^d$. \end{lem}

\begin{proof} Using the above arguments of Lemma \ref{torusexistence}, for each $u$, we can
construct a finite-energy function which, restricted to any offset $b + \mathcal L$, is equal to a
plane of slope $u$.  It is easy to see that for any $b \in \mathcal L$, there exists some $a$ in a
fundamental domain of $\mathcal L$ for which $D(a,b + a) = D_{\Phi}(b)$. The result follows by
taking $C = 2 \sup |D(z_1,z_2)| + 2\sup |D_{\Phi}(z_2-z_1)|$ where $z_1$ and $z_2$ are members of
that fundamental domain. \qed \end{proof}

\section{Gradient phase existence}
\begin{lem} \label{exposedexistence} If $u$ is either an {\it exposed} point or a convex
combination of exposed points of $\sigma$, then there exists a measure $\mu \in \mathcal
P_{\mathcal L}(\Omega, \mathcal F^{\tau})$ with $SFE(\mu) = \sigma(u)$. If $u$ is an exposed point
of $\sigma$, then $\mu$ can be taken to be ergodic. \end{lem}

\begin{proof} For $w \in \mathbb R^d$, let $\Psi^w$ be the nearest neighbor potential defined by
$$\Psi^w_{\{x,y \}}(\phi) = [\phi(x)-\phi(y)](w,x-y).$$ Define $\Phi^w = \Phi + \Psi^w$. It is not
hard to verify (using the discrete fundamental theorem of calculus) the equality of the Gibbs
kernels:
    $$\gamma^{\Phi^w}_{\Lambda} = \gamma^{\Phi}_{\Lambda}$$ as well as the identity $SFE^{\Phi^w}(\mu) =
SFE^\Phi(\mu) + (w, S(u))$, which in turn implies $\sigma_{\Phi^w}(u) = \sigma_{\Phi}(u) + (w, u)$.
Now, if $\sigma_{\Phi}$ has an exposed point at $u$, then we can replace $\Phi$ with a $\Phi^w$ for
which $\sigma_{\Phi^w}$ has a unique minimum at $u$.  By Lemma \ref{minimizerexists}, there exists
a measure $\mu$ with $SFE(\mu) = \sigma_{\Phi^w}(u)$. Clearly, such a measure has finite slope, so
(by definition of $\sigma_{\Phi^w}$) it must have slope $u$ and specific free energy equal to
$\sigma(u)$.  Furthermore, Jensen's inequality, convexity of $SFE$, and the ergodic decomposition
theorem imply that the ergodic components of $\mu$, with probability one, must have slope $u$. Now,
suppose $u$ is not an exposed point but $u$ is a convex combination of exposed points $u = \sum a_i
u_i$.  Let $\mu_u$ denote an ergodic measure of slope $u$ with $SFE(\mu_u) = \sigma(\mu_u)$.  Then
Lemma \ref{SFEisanexpectation} implies that $\mu = \sum a_i \mu_{u_i}$ is a Gibbs measure of slope
$u$ with $SFE(\mu) = \sigma(u)$.  \qed \end{proof}

We say $V:\mathbb R \mapsto \mathbb R$ is {\it super-linear} if for
all $c>0$, there exists a $b>0$ such that $V(\eta) > c |\eta|$
whenever $|\eta| \geq b$.  Say $\Phi$ is a {\it super-linear simply
attractive potential} if $\Phi$ is simply attractive and all of the
nearest neighbor potentials of $\Phi$ grow super-linearly.  The
following is a simple consequence of Lemma \ref{SEbound} and Lemma
\ref{levelsetcompactness}.

\begin{lem} If $\Phi$ is a super-linear simply attractive potential, then $\mu_i \rightarrow \mu$
(in the topology of local convergence) and $SFE^{\Phi}(\mu_i) \leq C$ for all $i$ together imply
$S(\mu_i) \rightarrow S(\mu)$ and $SFE^{\Phi}(\mu) \leq C$. \end{lem}

\section{A word on Lipschitz potentials}
Assume that $E = \mathbb Z$.  Many naturally arising discrete random
surfaces (height functions for domino tilings, square ice, etc.)
have the property that the nearest neighbor potential functions are
equal to infinity outside of a finite interval.  A {\it Lipschitz
potential} is a potential $\Phi$ such that for each adjacent pair of
vertices $x$ and $y$, the potential $\Phi_{\{x, y \}}(\phi)$ is
equal to $\infty$ when $\phi(x)-\phi(y)$ lies outside of a bounded
interval of $\mathbb Z$.  The reader may easily check the following:
\begin{lem} A gradient, nearest neighbor, $\mathcal L$-invariant
potential $\Phi$ is {\it Lipschitz} if and only if it is a perturbed
LSAP. \end{lem}

By Lemma \ref{halfspace}, $U_{\Phi}$ is the interior of a convex
polyhedron and by Lemma \ref{conttobound}, $\sigma$ is a bounded,
continuous function on the closure of $U_{\Phi}$. Discrete Lipschitz
potentials (i.e., Lipschitz potentials in the case $E = \mathbb Z$)
are convenient to work with for many reasons.  First of all, an
$\mathcal L$-periodic Lipschitz potential $\Phi$ can be completely
described with a finite amount of data. And although it is generally
not known how to compute $\sigma$ exactly, one can use simulations
(and the alternate definition of $\sigma$ given in Chapter
\ref{LDPempiricalmeasurechapter}), to approximate $\sigma$ to
arbitrary precision. Also, for discrete LSAPs, Propp and Wilson give
an algorithm (called ``coupling from the past'') for perfect
sampling from $\gamma_{\Lambda}^{\Phi}(*,\phi)$, where $\Lambda$ is
a finite subset of $\mathbb Z^d$ and $\phi$ is an arbitrary
``boundary function'' (which need only be defined on the boundary of
$\Lambda$) \cite{PW}. Finally, given any discrete SAP $\Phi$ and
$C>0$, we can approximate $\Phi$ with a Lipschitz potential,
$\Phi^C$, defined by
$$\Phi^C_{\Lambda}(\phi) =
\begin{cases} \Phi_{\Lambda}(\phi) & \Phi_{\Lambda}(\Phi) \leq C \\
\infty & \text{ otherwise.} \end{cases} $$  It is not hard to show
that as $C$ gets large, the variational distance between the
measures $\gamma_{\Lambda}^{\Phi}(*,\phi)$ and
$\gamma_{\Lambda}^{\Phi^C}(*,\phi)$ decays exponentially.

\section{Examples: layered surfaces, non-intersecting lattice paths} Fix $d$, and let $\Phi$ be a simply
attractive potential.  We can use a Gibbs measure of $\Phi$ to
choose a single random surface $\phi:\mathbb Z^d \mapsto \mathbb Z$.
Now, we would like a way choose a sequence of ``layered surfaces''
$\phi_i$, defined for $i \in \mathbb Z$, in such a way that

\begin{enumerate}
\item For each $x \in \mathbb Z^d$ and $i \in \mathbb Z$, $\phi_i(x) < \phi_i(x+1)$.
\item Given $\phi_{i-1}$, $\phi_{i-2}$, and the values $\phi_i(x)$ for $x$ outside of a finite set
$\Lambda$, the conditional distribution on the values $\phi_i(x)$
for $x \in \Lambda$ can be described as follows: it is the measure
$\gamma^{\Phi}_{\Lambda}(*,\phi_i)$ conditioned on $\phi_{i-1}(x) <
\phi_i(x) < \phi_{i+1}(x)$ for all $x \in \Delta$. \end{enumerate}
The natural way to do this is to replace $\Phi$ by a
$d+1$-dimensional $\Phi'$.  Write $\phi'(x,i) = \phi_i(x)$, where $x
\in \mathbb Z^d$ and $i \in \mathbb Z$.  Then we write
$\Phi'_{\{(x,i),(y,i)\}} = \Phi_{x,y}$ and
    $$\Phi'_{\{(x,i),(x,i+1)\}}(\phi') = \begin{cases} \beta (\phi'(x,i+1) - \phi'(x,i)) & \phi'(x,i+1) >
    \phi'(x,i)\\
    \infty & \text{ otherwise.}\end{cases}$$  Intuitively, we would like to set $\beta = 0$.  In this case,
    however, $\Phi'$ would not be simply attractive.  But since the kernels $\gamma_{\Phi'}$ are
    independent of $\beta$, we can get around this by fixing $\beta$ to be an arbitrary positive
    constant.

We can think of a Gibbs measure in $\mathcal G^{\Phi'}$ as a way of
choosing a sequence of layered surfaces.  This construction also
makes sense if we replace $E= \mathbb Z$ by $E = \mathbb R$.  It is
now easy to check that $U_{\Phi'} = U_{\Phi} \times (1, \infty)$
when $E = \mathbb Z$ and $U_{\Phi'} = U_{\Phi} \times (0, \infty)$
when $E = \mathbb R$.  If $u = (u_1, u_2)$, with $u_1 \in \mathbb
R^d$ and $u_2 \in \mathbb R$, then we can think of a random function
$\phi'$ of slope $u$ as a sequence of layered slope $u_1$ functions
$\phi_i$, spaced apart with density $1/u_2$.  It is a consequence of
the results in Chapter \ref{clusterswappingchapter} that, at least
when $E = \mathbb R$, the minimal-$SFE$ ergodic slope-$u_1$, density
$1/u_2$ layered surface Gibbs measure is unique. In the special case
$d=1$ and $E = \mathbb Z$, we can think of layered surfaces as
non-intersecting paths in a lattice.  (The latter appear frequently
in the theory of random matrices.)

\chapter{Analytical results for Sobolev spaces} \label{orliczsobolevchapter} In order to prove our
large deviations principles on surface shapes in later chapters, we will need to construct a
topology in which certain sets of bounded-average-energy surfaces are compact.  To this end, we
will ultimately convert some of our questions about discretized surfaces into analogous questions
about continuous functions.  The topologies on the continuous function spaces will be generated by
the {\it Orlicz-Sobolev norms} defined in this chapter, which are generalizations of $L^p$ norms
that arise when the function $|\cdot|^p$ is replaced by a more general positive, convex, symmetric
function $A$; it will generally be desirable to choose $A$ in such a way that this topology is as
strong as possible, since this will lead to the strongest large deviations principles and
concentration inequalities. Thus, for a given potential $\Phi$, we would like to determine the
strongest topologies in which the necessary compactness results hold.

In Section \ref{orliczsobolevsection}, we will cite several known
results ({\it compact imbedding theorems} and various bounds) about
Orlicz-Sobolev spaces on a bounded domain $D \subset \mathbb R^d$.
We cite these results without proof from \cite{A}, \cite{C2}, and
\cite{MP}.  While much of this theory is classical, some of the
results we employ have been proved during the last few
years---including the strongest versions of the imbedding theorems
and some of the bounds we need. A summary of these recent results
(many by Cianchi) can be found in \cite{C2}.

The results in this chapter will assist us primarily when $\Phi$ is
an ISAP (or a perturbed ISAP). When $\Phi$ is an LSAP (or perturbed
LSAP) the variational and large deviations principles can be derived
without the theorems of this chapter.  If $V: \mathbb R \mapsto
[0,\infty]$ is a non-constant, even, convex function, as always, we
will denote by $\Phi_V$ the ISAP whose nearest neighbor gradient
potentials are given by $V$. In Section \ref{simplexsection}, we
will describe how the Hamiltonians $H^{\Phi_V}(\phi)$ and
$H^{\Phi_{\overline{V}}}(\phi)$ can be approximated via energy
integrals of continuous interpolations of $\phi$; we also define the
``good approximations'' of bounded, sufficiently regular domains $D$
(by subsets of $\frac{1}{n} \mathbb Z^d$) that we will use in the
large deviations principle of Chapter \ref{LDPchapter}.  We extend
the compactness results of Section \ref{orliczsobolevsection} to
families of functions defined on approximating subsets of $D$ (and
not necessarily all of $D$).

In Section \ref{interpolationsection}, we derive a technical result
we need for both the proof of the variational principle in Chapter
\ref{LDPempiricalmeasurechapter} and the proof of the large
deviations principle in Chapter \ref{LDPchapter}.  Finally, in
Section \ref{mlapproximationsection}, we show how to approximate
certain functions $f:D \mapsto \mathbb R$ by functions that are
``mostly linear'' and whose ``energy'' is not much larger than the
energy of $f$; this result will be useful in proving the lower bound
on probabilities in the large deviations principle of Chapter
\ref{LDPchapter}.

\section{Orlicz-Sobolev spaces} \label{orliczsobolevsection}
\subsection{Orlicz-Sobolev space definitions}
We define $$W^{j,p}(D) = \{ f \in L^p(D): D^{\alpha}f \in L^p(D)
\text { for } 0 \leq |\alpha| \leq j \},$$ where $\alpha =
(\alpha_1, \ldots, \alpha_d)$ is a multi-index, with $0 \leq
\alpha_i \in \mathbb Z^d$, and $D^{\alpha}$ is the distributional
derivative. (Here, we use the definition $|\alpha| = \sum_{i=1}^d
\alpha_i$.) Define a norm by $$||f||_{j,p} =\[ \sum_{0 \leq |\alpha|
\leq j} ||D^{\alpha}f||^p_p \]^{1/p}$$ for $1 \leq p < \infty$ and
$||f||_{j,\infty} = \sup_{0 \leq |\alpha| \leq j}
||D^{\alpha}f||_{\infty}$.  These are {\it Sobolev spaces}.

We define a {\it Young function} to be a convex, even function $A:
\mathbb R \mapsto \mathbb R^+\cup \{\infty\}$ for which $A(0)=0$,
$A$ is finite on some open interval $(-a,a)$, and $A$ is not
identically zero (which implies that $A(\eta)$ grows at least
linearly in $\eta$ for $\eta$ large). Observe that for an isotropic
simply attractive potential $\Phi_V$, we may and will assume (adding
a constant to $V$ if necessary to make $V(0)=0$) that $V$ is a Young
function.  When $A$ is a Young function, the {\it Orlicz} space
$L^A(D)$ is the space of functions $f: D\mapsto \mathbb R$ for which
the norm
$$||f||_{A,D} = \inf \{ k| \int_D A(\frac{f(\eta)}{k})d\eta \leq 1
\}$$ is finite. We write $||f||_{A,D}$ simply as $||f||_A$---and
$L^A(D)$ as $L^A$---when the choice of $D$ is clear from context.
The {\it Orlicz-Sobolev} spaces are the spaces $$W^{j,A}(D) = \{ f
\in L^A(D): D^{\alpha}f \in L^A(D) \text { for } 0 \leq |\alpha|
\leq j \}.$$ Clearly, Sobolev spaces are Orlicz-Sobolev spaces
defined using the Young functions $A(\eta) = |\eta|^p$ for some $p$.
From here on, we will assume $j=1$ and deal only with the spaces
$W^{1,A}$.

If $v = (v_1, \ldots, v_d) \in \mathbb Z^d$, then we write $A(v) =
\sum_{i=1}^d A(v_i)$.  We write $|v|$ for the Euclidean norm of $v$.
Since $\sup_{i=1}^d |v_i| \leq |v| \leq d \inf_{i=1}^d |v_i|$ and
$A$ is convex (with $A(0)=0$), it is clear that $A(\frac{v}{d}) \leq
A(|v|) \leq A(d v)$.  Now, there are two natural ways to define
$||\nabla f||_{A,D}$; one is $$||\nabla f||_{A,D} = \inf \{ k|
\int_D A\[\frac{|\nabla f(\eta)|}{k}\] d\eta \leq 1 \}.$$ The second
is the same but without the absolute value sign on $\nabla f$, i.e.,
$$||\nabla f||_{A,D} = \inf \{ k| \int_D A\[\frac{ \nabla f(\eta)
}{k}\] d\eta \leq 1 \}.$$  The above discussion implies that the two
norms thus defined are equivalent---they differ from one another by
at most a factor of $d$.  Unless otherwise specified, we will always
use the second definition.  Some of the theorems we cite (compact
imbeddings, etc.) are proved in papers which use the first
definition; however, all these cited results obviously remain true
when the norm is replaced by an equivalent one.

We now cite the following (Section 2.3 of \cite{C2}):

\begin{lem} The spaces $W^{1,A}$, equipped with the norm $||f||_{W^{1,A}} = ||f||_A + ||\nabla f
||_A$, are Banach spaces. \end{lem}  We define $L^{1,A}$ to be the
set of weakly differentiable functions $f:D \mapsto \mathbb R$ with
$|\nabla f| \in L^{A}(D)$.  Note that $W^{1,A}(D) = L^A(D) \cap
L^{1,A}(D)$.  We will later see that $W^{1,A}(D) = L^{1,A}(D)$ on
all domains $D$ of interest to us.

\subsection{Regular domains: the domain class $\mathbb G(\frac{d-1}{d})$} Denote by $f_D$ the mean
$|D|^{-1} \int_{D} f(\eta) d\eta$, where $|D|$ is the Lebesgue volume of $D$. (Unless otherwise
specified, all integrals over $D$ are with respect to Lebesgue measure.) The following kinds of
assertions are starting points for the Orlicz-Sobolev theory: \begin{enumerate}
\item For some constant $C$, $||f-f_D||_{A^*} \leq C||\nabla f||_A$ (for appropriately chosen Young
functions $A$ and $A^*$). \item For some constant $C$, $||f||_{A^*}
\leq C(||f||_{A} + ||\nabla f||_A)$ (where $A$ is a Young function
and $A^*$ is some appropriately chosen Young function that grows
more rapidly than $A$).  In other words, the {\it imbedding} of
$W^{1,A}(D)$ into $L^{A^*}(D)$ exists and is continuous.
\end{enumerate}

Both kinds of results depend on regularity properties of the domain
$D$. Even without specifying $A$ and $A^*$, we can imagine what
might go wrong if we did not have regularity conditions.  Suppose
$D$ consists of the interior of $D_b \cup D_1 \cup D_2$ where $D_1$
and $D_2$ are very large closed cubes (not necessarily of equal
volume) and $D_b$ is a long but very skinny rectangular tube (a
``bottleneck'') connecting them.  Then take $f$ to be a function
that is equal to $0$ on $D_1$ and $M$ on $D_2$, and varies linearly
between between $0$ and $M$ within $D_b$.  By making $M$
sufficiently large and making $D_b$ sufficiently long and skinny, we
can make $||f-f_D||_{A^*}$ arbitrarily large even while making
$||\nabla f||_A$ arbitrarily small.

We can imagine that if $D$ had infinitely many (increasingly skinny)
bottlenecks of this form, it might be possible---by choosing $f$
equal to zero on one side of some bottleneck and some value $M$ on
the other side---to produce functions $f$ on $D$ for which
$\frac{||f-f_D||_{A^*}}{||\nabla f||_A}$ is arbitrarily large,
contradicting the first assertion. Similar problems arise with the
second assertion.  Both assertions require a bound on the volume
separation of the bottleneck (i.e., the amount of mass on the
smaller side of the bottleneck) in terms of its width.

To be precise, define $P(E;D)$ to be the {\it perimeter} of $E$
relative to $D$, i.e., the total variation over $D$ of the gradient
of the characteristic function of $E$.  (When $E$ has a smooth
boundary, then $P(E,D)$ is simply the $(n-1)$-dimensional Hausdorff
measure of $\partial E \cap D$. See Section $6.1.1$ of \cite{M} for
more details.) Let $\mathbb G(z)$ be the set of bounded domains $\{
D \subset \mathbb R^n \}$ for which there exists a constant $C$ such
that
$$[\min \{|E|, |D-E| \}]^{z} \leq CP(E;D)$$ for all Lebesgue
measurable subsets $E$ of $D$, where $|\cdot|$ denotes Lebesgue
measure. In particular, when $z = \frac{d-1}{d}$, this is a very
natural restriction, which we write
$$\min \{|E|, |D-E| \} \leq CP(E;D)^{\frac{d}{d-1}}.$$ This is natural because we can interpret
$CP(E;D)^{\frac{d}{d-1}}$ is a constant times the volume contained
in a sphere with surface area $P(E,D)$. The infimum of the $C$ for
which this inequality holds is called the {\it isoperimetric
constant} of $D$.  The space $\mathbb G(\frac{d}{d-1})$ includes
many families of domains for which Sobolev-type results were proved
classically, including bounded sets satisfying the {\it cone
property} (see Corollary $3.2.1/3$ of \cite{M} and Section 2.4 of
\cite{C2}), the {\it strong local Lipschitz property}, and the {\it
uniform $C^m$-regularity property} for $m\geq 1$ (see $4.3$ to $4.7$
of \cite{A}).

There is a range of weaker Orlicz-Sobolev theorems that apply when
weaker regularity conditions are placed on $D$ (see, e.g., Remark
3.12 of \cite{C2}); it may be possible to use these more general
results to prove weaker large deviations principles for random
surfaces on mesh approximations of these more general domains.  (See
Chapter \ref{openproblemschapter}.)  However, we will limit our
attention to the domains $\mathbb G(\frac{d-1}{d})$.  Random
surfaces on what we will call ``good approximations'' of these
domains are especially convenient because they have the same large
deviations properties as random surfaces on the standard
approximations of the unit cube $[0,1]^d$.

\subsection{Comparing Young functions}
When $A$ and $B$ are Young functions, we say that $B$ {\it dominates} $A$ near infinity if there
exist positive constants $c_1$ and $c_2$ such that $A(\eta) \leq B(c_1 \eta)$ for $\eta > c_2$. If
for every $c_1 > 0$, there is a $c_2$ for which this holds, we say $A$ increases {\it essentially
more slowly} than $B$. We call $A$ and $B$ {\it equivalent near infinity} if each dominates the
other near infinity.  The following fact is a motivation for this definition.  (See Remark $3.3$ of
\cite{C2}.)

\begin{lem} If $D$ has finite volume and the Young functions $A$ and $B$ are equivalent near
infinity, then the Luxemburg norms $||f||_A$ and $||f||_B$ are equivalent norms. \end{lem}

\subsection{Sobolev conjugates} \label{sobconjsection}
Following Section $3.2$ of \cite{C2}, given $n\geq 2$ and a Young function $A$, we define an
increasing function $H: [0,\infty) \mapsto [0,\infty)$ by $$H(r) = \[\int_{0}^r
\[\frac{t}{A(t)}\]^{\frac{1}{d-1}} dt\]^{\frac{d-1}{d}},$$ and $A_d: [0,\infty) \mapsto [0,\infty)$
by $$A_d = A \circ H^{-1},$$ where $H^{-1}$ is the left-continuous
inverse of $H$.  We extend $A_d$ to $\mathbb R$ by $A_d(\eta) =
A_d(-\eta)$. Note: in order for $H$ to be finite and well-defined,
we have to assume that for $c>0$, we have $$\int_0^{c} \[
\frac{t}{A(t)} \]^{\frac{1}{d-1}} < \infty.$$ If this is the case,
we say that $A$ {\it has a conjugate} $A_d$.  If this is not the
case, we can replace $A$ by any equivalent Young function for which
the integral {\it does} converge and define $H$ and $A_d$ using that
Young function instead; we call this $A_d$ an {\it equivalency
conjugate} for $A$. (Such an equivalency conjugate always exists by
Remark 3.3 of \cite{C2}.) Also, note that if $C = \int_c^{\infty}
(\frac{t}{A(t)})^{\frac{1}{d-1}} = \infty$, then $A_d(\eta)$ is
everywhere finite. Otherwise, it is infinite for $|\eta| > C$.  In
particular (as the reader may easily check), $A_d(\eta)$ is
everywhere finite if $A(\eta) \leq \eta^d$ and infinite outside an
interval if $A(\eta) \geq \eta^{d+\epsilon}$ and $\epsilon >0$.

We cite a concrete example from Example 3.17 of \cite{C2}:  If
$A(\eta)$ is equivalent near infinity to $\eta^p(\log(\eta))^q$,
then $A_d(\eta)$ is equivalent near infinity to

$$\begin{cases} \eta^{dp/(d-p)}(\log \eta)^{dq/(d-p)} & \text{ if } 1 \leq p \leq d. \\
\exp \[(\eta^{d/(d-1-q)}\] & \text { if } p=d, q<d-1. \\
\exp \[\exp \[\eta^{d'}\]\] & \text { if } p=d, q=d-1. \\
\end{cases}$$
If either $p>n$, or $p=n$ and $q>n-1$, then $A_d$ is equal to
$\infty$ outside of a finite interval.  We use $A^*$ to denote a
{\em sub-conjugate} of $A$, i.e., any Young function that increases
essentially more slowly than $A_d$.

\subsection{Imbeddings}
We cite the following (Theorem $3.9$, Remark $3.10$, and Theorem
$3.13$ of \cite{C2}): \begin{thm} \label{imbedding} Let $A$ be any
Young function with conjugate $A_d$ and $D \in \mathbb
G(\frac{d-1}{d})$.  Then the following are true: \begin{enumerate}
\item There exists a constant $K$, depending only on $A$, the volume of $D$, and the isoperimetric constant
$C$ of $D$ such that for any $f \in W^{1,A}$, $$||f-f_D||_{A_d} \leq
K ||\nabla f||_{A}.$$ \item The imbedding $$W^{1,A}(D) \mapsto
L^{A_d}(D)$$ is well-defined and continuous.  Given $A$, if $A$ has
a conjugate (not merely an equivalency conjugate) $A_d$, then $K$
depends only on $C$ (and the same $K$ holds if $C$ is replaced by
any $C_0 < C$).
\item In each of the two previous items, the Young space $L^{A_d}$ is the smallest Orlicz space for
which the result is true. \item If $A^*$ is any Young function
increasing essentially more slowly near infinity than $A_d$, then
the imbedding $$W^{1,A}(D) \mapsto L^{A^*}(D)$$ is compact.
\end{enumerate} \end{thm}

(The fact that the $K$ in the first statement that holds for $C$
also holds for any $C_0<C$ is not stated explicitly in \cite{C2},
but it appears to be understood in the context.  In any case, it is
not hard to see that if a counterexample to the statement
exists---in the form of a function $f$ on a domain $D$---for some
$C$ and $K$, then a counterexample also exists for $K$ and any
larger $C$. The counterexample can be obtained by removing a
zero-volume subset of $D$ to produce a new set $D'$ with the
appropriate larger isoperimetric constant---and then replacing $f$
with its restriction to the $D'$.)

The following corollary is useful, as it implies that whenever we can prove $f \in L^{1,A}(D)$, we
can apply Theorem \ref{imbedding} to produce a bound on the $L^{A_d}$ norm.

\begin{cor} \label{imbeddingcor} If $D \in \mathbb G(\frac{d-1}{d})$, then the spaces $L^{1,A}(D)$ and
$W^{1,A}$ are equivalent. \end{cor}

\begin{proof} Suppose that $||\nabla f||_A$ is finite but $||f||_A$ infinite. Then we can take a
truncation $f_M(\eta)$ equal to $f(\eta)$ when $|f(\eta)| \leq M$, $M$ for $f(\eta) > M$, and $-M$
for $f(\eta) < -M$.  Write $g_M = f_M + c_M$ where the constants $c_M$ are chosen in such a way
that $g_M$ has mean zero.  As $M$ tends to $\infty$, $||\nabla g_M||_A$ tends to $||\nabla f||_A$
and $||g_M||_A$ tends to infinity.  To see the latter fact, observe that if $f$ has finite mean
$f_D$, then constants $c_M$ converge to $f_D$ and $||g_M||_A$ converges to $||f||_A$.  If $f$ does
not have finite mean, then $||g_M||_1$ tends to infinity; hence $||g_M||_A$ tends to infinity for
any Young function $A$.  These facts imply that the ratio $||g_M||_A/||\nabla g_M||_A$ grows
arbitrarily large, contradicting Theorem \ref{imbedding}. \qed \end{proof}

Now we can make a statement closer to the form we will actually need it to be in for our proofs:

\begin{cor} \label{energyimbedding} Let $A$ be any Young function which has a conjugate $A_d$ and
$D \in \mathbb G(\frac{d-1}{d})$, and let $A^*$ be any function that
$A_d$ dominates near infinity, and $C$ a positive constant.  Then
there exists a constant $K$, depending only on $A$, $A^*$, and $C$,
such that whenever the isoperimetric constant of $D$ is less than
$C$ and $f \in W^{1,A}$,
$$||f-f_D||_{A^*} \leq K \int A(\nabla f(\eta)) d\eta.$$ Moreover, if $A^*$ increases essentially more
slowly than $A_d$, then the mapping from $L^{1,A}(D)$ to $L^{A^*}$
that sends $f$ to $f - f_D$ is compact. \end{cor}

\begin{proof}  This follows from Theorem \ref{imbedding} and a couple of simple observations.  The
first is that the convexity of $A$, the fact that $A(0)=0$, and Jensen's inequality imply that
$$\int A(\nabla f(\eta)) d\eta \leq ||f||_A.$$ The second is that the image of the unit ball in
$L^{1,A}$ is a subset of the image of the unit ball of $W^{1,A}$; since the latter image is
precompact, the former is also.  \qed \end{proof}

\section{Connection to the discrete settings} \label{simplexsection}
\subsection{Simplex interpolations and discrete norms} Given $w= (w_1, w_2, \ldots, w_d) \in \mathbb R^d$,
write $\lfloor w \rfloor = (\lfloor w_1 \rfloor, \lfloor w_2
\rfloor, \ldots, \lfloor w_d \rfloor)$, where for any real $\eta$,
$\lfloor \eta \rfloor$ is the integer part of $\eta$. Also, let
$s(w) \in S^d$ be the permutation (uniquely defined for almost all
$w$) that gives the rank ordering of the components of $w - \lfloor
w \rfloor$. For each vertex $v \in \mathbb Z^d$ and $s \in S^d$, we
denote by $C(v,s)$ the closure of the simplex of vertices $w$ with
$\lfloor w \rfloor = v$ and $s(w) = s$.

We say that a domain in $\mathbb R^d$ is a {\it simplex domain} if
it is the interior of a finite union of simplices of the form
$C(v,s)$. We say that a subset $\Delta \subset \mathbb Z^d$ is a
{\it simplex boundary set} if it is the union of the corner sets of
the simplices in a simplex domain (denote this simplex domain by
$\hat \Delta$) and if every adjacent pair of vertices $\Delta$ forms
an edge of at least one of these simplices. In a sense we describe
precisely in the next subsection, we will often fix a domain $D
\subset \mathbb R^d$ and choose a sequence of simplex boundary sets
$D_n$ so that the domains $\tilde{D}_n = \frac{1}{n} \hat{D}_n$ are
increasingly close ``approximations'' of $D$. We will assume that
the volume of the $\tilde{D}_n$, written $|\tilde{D}_n|$, is bounded
between positive constants $c_1$ and $c_2$, independently of $n$. To
summarize:
\begin{enumerate}
\item $D$ is a domain in $\mathbb R^d$.
\item $D_n \subset \mathbb R^d$ is an ``approximation'' (to be precisely defined later) of
$nD$.
\item $\hat D_n$ is a simplex domain in $\mathbb R^d$ derived from $D_n$; it is an ``approximation''
of $nD$.
\item The normalized domain $\tilde D_n = \frac{1}{n} \hat{D}_n$ is an ``approximation'' of
$D$.
\end{enumerate}

Given a function $\phi: D_n \mapsto E$, denote by
$\hat{\phi}:\hat{D}_n \mapsto \mathbb R$ the unique function that
extends to the closure of $\hat{D}_n$ in such a way that it agrees
with $\phi$ on $D_n$ and is linear on the closure of each simplex of
$\hat{D}_n$; define a rescaled function $\tilde{\phi}:\tilde{D}_n
\mapsto \mathbb R$ by $\tilde{\phi}(\eta) = \frac{1}{n}
\hat{\phi}(n\eta)$. Again, to summarize:\begin{enumerate}
\item Begin with $\phi: D_n \rightarrow E$.
\item Then $\hat{\phi}_n \rightarrow E$ is linear interpolation of $\phi$ to
$\hat{D}_n$.
\item $\tilde{\phi}:\tilde{D}_n \mapsto \mathbb R$ by $\tilde{\phi}(\eta) = \frac{1}{n}
\hat{\phi}(n\eta)$ is a rescaling of $\hat{\phi}$. \end{enumerate}

Write $$\mathcal E_{n,A}(\phi) = |D_n|^{-1} \sum_{y \in D_n}
A\[\frac{\phi(y)}{n}\] $$ and similarly $$\mathcal E_{n,A}(\nabla
\phi) = |D_n|^{-1}H^o_{\Lambda}(\phi) = |D_n|^{-1}\sum_{\{x,y\}
\subset D_n, |x-y|=1} A(\phi(y) - \phi(x)).$$

Here $H$ is the Hamiltonian of the ISAP $\Phi_A$.  Note the
normalization by the size of $D_n$ built into these definitions.
Roughly speaking, the first gives the {\it average} value of
$A(\phi(\cdot)/n)$ on $D_n$; the second gives the ``energy per
site'' of $\phi$.  The following simple lemma gives a connection
between continuous and discrete energy integrals.  For any
continuous $f: D \rightarrow \mathbb R$, write $\mathcal E_{D}
(\nabla \phi) = \int_{\tilde{D}_n} A (\nabla
\tilde{\phi}(\eta))d\eta$ and $\mathcal E_{D} (\phi) =
\int_{\tilde{D}_n} A(\tilde{\phi})$.

\begin{lem} \label{discretenormbound} Assume (as above) that the volume of $\tilde D_n$ is bounded
between two positive constants $c_1$ and $c_2$.  Then there exist
positive constants $C_1$ and $C_2$ (independent of $n$) such that
for any $\phi: D_n \mapsto \mathbb R$, we have $$C_1 \mathcal
E_{n,A}(\nabla \phi) \leq \mathcal E_{D} (\nabla \phi) \leq C_2
\mathcal E_{n,A}(\nabla \phi)$$ for all $\phi$.   There also exist
positive constants $C_1$ and $C_2$ (independent of $n$) for which
$$E_{n,A} (C_1 \phi) \leq \mathcal E_{D} (\phi) \leq \mathcal E_{n,A} (C_2 \phi).$$ \end{lem}

\begin{proof} For the first statement, let $C(v,s)$ be a simplex of $\hat{D_n}$ and let $v^0, v^1,
v^2,\ldots, v^d$ be the edges in the path starting at $v^0=v$ and
stepping at the $i$th step on unit in the $s(i)$th, so that $v^d =
v+(1,1, \ldots, 1)$.  The vertices of this path are the vertices of
$C(v,s)$. Observe that $$\int_{\frac{1}{n}C(v,s)} A(\nabla
\tilde{\phi}(\eta))d\eta = \frac{1}{n^d} \int_{C(v,s)} A(|\nabla
\hat{\phi}|) = \frac{1}{n^dd!} \sum_{i=1}^d A(\phi(v_i) -
\phi(v_{i-1})).$$  Since $\tilde{D_n}$ has volume bounded between
two positive constants, we may deduce that the ratio of
$\frac{1}{n^dd!}$ and $\frac{1}{|D_n|}$ is also bounded between two
positive constants.  The result now follows from the fact that every
edge of $D_n$ is an edge of at least one simplex and at most $d!$
simplices of $\hat{D}_n$---and each edge of a simplex of $\hat{D}_n$
is also an edge of the graph $D_n$.

For the second statement, by arguments similar to those above, it is
enough to prove the result for the case that $\hat{D_n}$ consists of
a single simplex $C(v,s)$.  That is, it is enough to show that when
$D_n$ is the set of vertices of a single simplex, then for some
$C_1$ and $C_2$ $$ \sum_{x \in D_n} A(C_1 \phi(x)) \leq
\int_{C(v,s)}A(\hat{\phi}(\eta))d\eta \leq  \sum_{x \in D_n} A(C_2
\phi(x)),$$ and this is a simple exercise.  \qed \end{proof}

Given the sets $D_n$ as above, we can now define discrete norms
$$||\phi||_{A, n} = \inf \{k : \mathcal E_{A, n} (\phi/k) \leq 1
\}$$ and $$||\nabla \phi||_{A, n} = \inf \{k : \mathcal E_{A, n}
(\nabla \phi/k) \leq 1 \}.$$  The following is an immediate
consequence of Lemma \ref{discretenormbound}.

\begin{cor} \label{discretenormcor} Let $D_n$ be a sequence of simplex boundary sets and assume
that the volume of $\tilde{D}_n$ lies between two positive
constants. Then there exist positive constants $C_1$ and
$C_2$---independent of $n$---for which the following hold for any
$\phi:D_n \mapsto \mathbb R$:

$$C_1 ||\phi||_{A, D_n} \leq ||\tilde{\phi}||_{A,\tilde{D}_n} \leq C_2 ||\phi||_{A, D_n}$$

$$C_1 ||\phi||^1_{A, D_n} \leq ||\nabla \tilde{\phi}||_{A,\tilde{D}_n} \leq C_2 ||\phi||^1_{A,
D_n}.$$ \end{cor}

In particular, this result implies that Lemma \ref{imbedding} and Corollary \ref{imbeddingcor}
remain true for spaces of functions defined on $\Lambda_n$ if we replace the continuous norms with
discrete ones (and the relevant constants in these bounds are independent of $n$).

\subsection{Good approximations of domains $D \in \mathbb G(\frac{d-1}{d})$} \label{goodapproximationsection}
Let $\partial D$ denote the boundary of $D$.  Given a sequence of
subsets $D_n$ of $\mathbb Z^d$, we say that the sets
$\frac{1}{n}D_n$ are {\it good approximations} of $D$ if:
\begin{enumerate}
\item $D_n \subset Z^d \cap nD$ for all $n$.
\item $\lim_{n \rightarrow \infty} \sup \{|\frac{1}{n}x-y| : x \in [\mathbb Z^d \cap nD]\backslash D_n, y \in \partial
D\} = 0$.
\item For each $x \in D_n$ and $y \in \partial D$, $|x-ny|_{\infty} \geq 1$, where $|\cdot|_{\infty}$ is the
supremum norm of a vector.
\item Each $D_n$ is simplex boundary set.
\item For some $C$, each $\tilde{D}_n$ has isoperimetric constant less than or
equal to $C$. \end{enumerate}

The first item states that $\frac{1}{n}D_n$ approximate $D$ from within.  The second one states
that the approximation gets progressively better as $n$ gets large; it implies, that every compact
subset of $D$ is contained in $\tilde{D}_n$ for all sufficiently large $n$. The third is a
technical condition that requires that we exclude from $D_n$ points that are {\it too} close to the
boundary of $nD$. This condition ensures, for example, that every edge in $D_n$ will be completely
contained in $nD$. The fourth condition makes it possible to use Section \ref{simplexsection} to
approximate energies of functions on $D_n$ with continuously defined energies of continuous
functions on $\tilde{D}_n$.  The fifth condition prevents us from using a sequence of
approximations in which the bottlenecks become increasingly severe with $n$; in particular, it
implies that each $D_n$ is connected.

\begin{figure} \begin{center} \leavevmode \epsfbox[0 0 220 260]{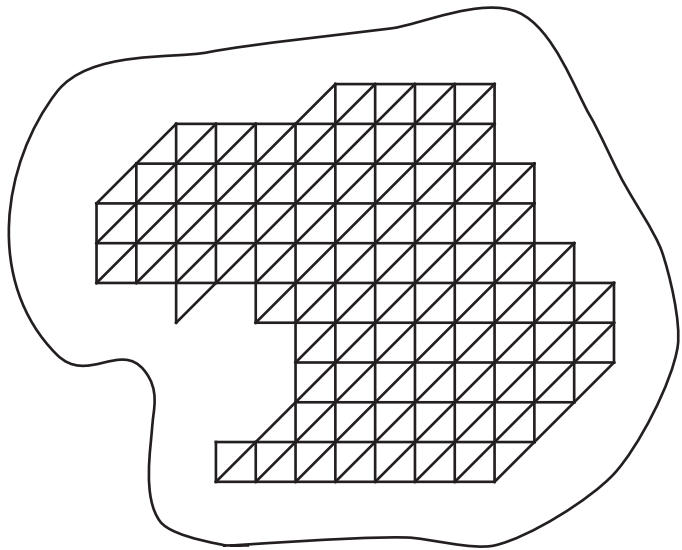} \end{center} \caption{A
possible simplex domain approximation for a domain $D$} \label{simplex} \end{figure}

Finally, for the purposes of our large deviations principles, we
will need a topology---similar to the $L^{A^*}$ topology, for an
appropriate Young function $A$---on a space which includes the
functions from $\tilde{D}_n$ to $\mathbb R$ for all $n$.  Suppose
that $\frac{1}{n}D_n$ are good approximations of $D$. Let
$$L^{A^*}_0(D)= L^{A^*}(D)\cup_{n=1}^{\infty} L^{B}(\tilde{D}_n).$$ For
notational convenience, we write $\tilde{D}_{\infty} = D$.  Assume
$A$ has Sobolev conjugate $A_d$ and that $A^*$ increases essentially
more slowly than $A_d$ near infinity.  For any $\phi_1$ and $\phi_2$
in $L^{A^*}_0(D)$ (defined on $\tilde{D}_i$ and $\tilde{D}_j$
respectively), write $\delta(\phi_1, \phi_2) =
||\phi_1-\phi_2||_{A^*, \tilde{D}_i \cap \tilde{D}_j} + |\tilde{D}_i
- \tilde{D}_j|$, where $|\tilde{D}_i - \tilde{D}_j|$ denotes the
Lebesgue measure of the symmetric difference of these two sets.   We
treat $L^{A^*}_0(D)$ as being endowed with the topology generated by
this metric.  If $\phi$ is defined on $\tilde{D}_n$, write $\mathcal
E(\nabla \phi) = \int_{\tilde{D}_n} A(\nabla (\phi_i))$.

\begin{lem} \label{goodcompact} The set $X_C$ of zero-mean functions $\phi \in L^{A^*}_0(D)$ with
$\mathcal E(\nabla \phi) \leq C$ is compact in $L^{A^*}_0(D)$.
\end{lem}

\begin{proof} Since $L^{A^*}_0(D)$ is a metric space, it is sufficient to prove that $X_C$ is
sequentially compact.  Let $\phi_i$ be any sequence of functions in
$L^{A^*}_0(D)$; if an infinite number of the $\phi_i$ are defined on
the same $\tilde{D}_n$, then the existence of a convergent
subsequence follows from Lemma \ref{imbedding}.  Otherwise, we may
assume that $\phi_i$ are defined on $\tilde{D}_n$ where $n$ is
increasing in $i$.  By a diagonalization argument, we can choose a
subsequence $m(i)$ of the integers, indexed by $i$, whose
restrictions to $\tilde{D}_n$ (which are defined when $i$ is large
enough so that $\phi_i$ is defined on a superset of $\tilde{D}_n$)
converge, for each $n$, to a limit in $L^{A^*}(\tilde{D}_n)$ (again,
by Lemma \ref{imbedding}). Thus, there is a unique function $\phi$
such that $\phi_{m(i)}$ converges to $\phi$ on every compact subset
of $D$. We need only to check that $\phi \in L^{A^*}(D)$ and that
the $\phi_{m(i)}$ also converge to $\phi$ with respect to the metric
$\delta$.

Since there is a uniform upper bound on the isoperimetric constants
of the $\tilde{D}_n$, by Theorem \ref{imbedding}, there is a uniform
upper bound $k$ on $||\phi_i||_{A^*, \tilde{D}_i}$, and hence on
$||\phi||_{A^*, \tilde{D}_n}$ as well.  If we had $\int_D
A^*(\frac{\phi(\eta)}{k})d\eta > 1$, then for some $i$, we would
have $\int_{\tilde{D}_i} A^*(\frac{\phi(\eta)}{k})d\eta > 1$, a
contradiction. Hence $||\phi||_{A^*, D} \leq k$.

Next, we know that for ever $\epsilon$ and $n$ there exists an $N$
such that for all $m \geq N$, we have $||\phi_m - \phi||_{A^*,
\tilde{D}_n} < \epsilon$.  Let $\psi^{\epsilon}_n = 1_{D \backslash
\tilde{D}_n} \phi_m - \phi$; it will be clear that $\delta(\phi_n,
\phi) \rightarrow 0$ if we can show that
$||\psi^{\epsilon}_n||_{A^*, D}$ tends to zero for every fixed
$\epsilon$.  This follows easily from the fact that
$||\psi^{\epsilon}_n||_{A_d,D}$ is uniformly bounded (independently
of $n$), that the volume of the sets on which $\psi^{\epsilon}_n$
are supported is tend to zero in $n$, and that $A^*$ increases
essentially more slowly near infinity than $A_d$. \qed
\end{proof}

The following corollary is the only compactness result we will
actually need in our proof of the large deviations principle (in the
space $L^{A^*}_0(D)$) for random surfaces defined using isotropic
convex nearest neighbor potentials.

\begin{cor} \label{discretegoodcompact} Let $Y_C$ be the set of zero-mean functions of the form
$\tilde{\phi}$, where $\phi:D_n \mapsto \mathbb R$ satisfies
$$H^o_{D_n}(\phi) \leq C|D_n|.$$  Here $H$ is the Hamiltonian for the
isoperimetric convex nearest neighbor potential determined by $A$.
If $A$ has a conjugate $A_d$ and $A^*$ increases essentially more
slowly than $A_d$, then $Y_C$ is precompact in $L^{A^*}_0(D)$.  The
same is true if $A$ does not have a conjugate but $A_d$ is a
conjugate of a Young function equivalent to $A$.\end{cor}

\begin{proof} When $A$ has conjugate $A_d$, the statement follows immediately from Lemma
\ref{discretenormbound} and Lemma \ref{goodcompact}.  When $A$ does
not have a conjugate, we can replace $A$ with an $A_0$ that does
have a conjugate and for which $|A(\eta)-A_0(\eta)| \leq b$ for
$\eta \in \mathbb R$.  (See Remark 3.3 of \cite{C2}.) Thus,
$|D_n|^{-1} H^o_{D_n}(\phi)$ defined using $A$ differs only by a
constant amount from the analogous expression defined using $A_0$.
Since $D$ has finite volume, the $Y_C$ defined using $A$ is a subset
of some $Y_{C_0}$, defined using $A_0$ in place of $A$, and the
corollary follows. \qed \end{proof}

\section{Low energy interpolations from $L^1$ and $L^{1,A}$ bounds} \label{interpolationsection}
\subsection{$L^A$ bounds from $L^1$ and $L^{1,A}$ bounds}
Let $C \in \mathbb R$ be an arbitrary constant and let $S_{A,C}(D)$
be the set of weakly differentiable functions $f$ on $D$ with
$||\nabla f||_A \leq C$.  The following theorem states that the
imbedding of $L^1(D) \cap S_{A,C}(D)$ (with the $L^1$ norm) into
$L^A \cap S_{A,C}$ (with the $L^A$ norm) is continuous. \begin{lem}
\label{fancybound} Let $A$ be a Young function and $D \in \mathbb
G(\frac{d-1}{d})$ and $C>0$. For every $\epsilon$, there exists a
$\delta = \delta(C,\epsilon,A)>0$ such that the following two
conditions \begin{enumerate}
\item $||\nabla f||_A \leq C$ \item $||f||_1 \leq \delta$ \end{enumerate} together imply $||f||_A
\leq \epsilon$.  The lemma remains true if we replace the first statement with the modified
statement: $\int_D A (\nabla f(\eta))d\eta \leq C$. \end{lem}

\begin{proof} First of all, if $||\nabla f||_A \leq C$, then $g = f/C$ satisfies $\int_D (\nabla
f(\eta))d\eta \leq 1$---i.e., $g$ satisfies the modified statement with $C=1$.  It will therefore
be enough to show that the modified first statement and the second statement together imply the
conclusion.  Statements $2.21$ (and subsequent discussion) and $2.25$ of \cite{C2} indicate that
for some constant $K$ (depending only on $D$ and $A$), we have $||\nabla f||_A \geq ||K \nabla
f^*||_A$ whenever $f^*$ is a positive, radially decreasing, spherically symmetric ``rearrangement''
of $f$ defined on a sphere $D'$ with the same volume as $D$.  Though we do not define
rearrangements precisely here, it suffices to note that such rearrangements always exist for
measurable $f$, and whenever $f^*$ is a rearrangement of $f$, we have $||f||_A = ||f^*||_A$ and
$||f||_1 = ||f^*||_1$ and $\int_{D} A(|\nabla f(\eta)|)d\eta= \int_{D'} A(|\nabla
f^*(\eta)|)d\eta$.  Thus, if $f$ violates the theorem conclusion for a particular choice of $C$,
$\delta$, and $\epsilon$, then $f^*$ violates the statement for the same values.

It is therefore sufficient to prove the result for positive, radially decreasing spherically
symmetric functions $f(\eta) = g(|\eta|)$---defined on a sphere $D'$ of radius $R$) with $g:(0, R]
\mapsto [0, \infty)$ decreasing. Thus, the theorem is implied by (the $K=1$ case of) the following
assertion: For every $\epsilon>0$ and $K>0$ and $C>0$ there exists a $\delta>0$ such that
\begin{enumerate}
\item $\int_0^R \[A(|g'(r)|)d\Gamma_d r^{d-1}\] dr \leq C$ \item $\int_0^R \[|g(r)|d\Gamma_d r^{d-1}\]dr
\leq \delta$ \end{enumerate} together imply $\int_0^R \[A\[\frac{g(r)}{\epsilon}\]d\Gamma_d
r^{d-1}\]dr \leq K$.

Here $\Gamma_d$ is the volume of the unit sphere of radius $1$ in $d$ dimensions (and $d\Gamma_d$
the $n-1$-dimensional volume of its boundary). Replacing $C$, $\delta$, and $K$ with their values
divided by $d\Gamma_d$, we eliminate that constant: it is now enough for us to prove that for any
$\epsilon>0$, $C>0$, and $K>0$ there exists $\delta>0$ such that \begin{enumerate} \item $\int_0^R
A(|g'(r)|)r^{d-1} dr \leq C$ \item $\int_0^R |g(r)| r^{d-1}dr \leq \delta$ \end{enumerate} together
imply $\int_0^R A\[\frac{g(r)}{\epsilon}\] r^{d-1}dr \leq K$.

We say the theorem ``holds for $(R,C, K,\epsilon)$'', if for those
fixed values, the above implication is true for some sufficiently
small $\delta$.  Let $\alpha>0$ and $\beta>0$ be arbitrary
constants.  We now show that to prove that the theorem holds for all
positive $(R,C,K,\epsilon)$, it is sufficient to prove that the
theorem holds for all positive $(R,C,K,\epsilon)$ for which $R \leq
\alpha$; furthermore, it is enough to prove the implication for the
positive, decreasing functions $g$ on $(0,R]$ for which $|g(R)| \leq
\beta$. To just these two reductions, first suppose $R>\alpha$ but
$0<R'<\alpha$. Since $g$ is decreasing, we have $|g(r)| \leq
|g(R')|$ for $r > R'$ and $\delta \geq \int_0^R |g(r)| r^{d-1}dr
\geq (R')^d g(R')$. We conclude that for any fixed value of $R'$, we
can assume (when $\delta$ is small enough) that $g(R') \leq \beta$.
Moreover,
    $$\int_{R'}^R A\[\frac{g(r)}{\epsilon}\] r^{d-1}dr
    \leq \int_{R'}^R A\[\frac{\delta}{(R')^d}\].$$
For fixed $R'$, $R$, and $A$, the latter term clearly tends to zero as $\delta$ tends to zero.
Thus, it is enough for us to bound $\int_{0}^{R'} A(\frac{g(r)}{\epsilon}) r^{d-1}dr$; and we can
do this by solving the modified problem with $R = R' < \alpha$.

Thus, we may assume $R \leq \frac{\epsilon}{2}$ throughout the
remainder of the proof. Now, putting $t = R-r$, and noting
$\frac{t}{\epsilon} \leq \frac{1}{2}$ (since $R \leq
\frac{\epsilon}{2}$), we can write
$$\frac{g(r)}{\epsilon}=\frac{g(R)}{\epsilon}+\frac{1}{\epsilon}\int_{r}^{R}
|g'(u)|du =
\[1-\frac{t}{\epsilon}\] \frac{g(R)}{\epsilon\[1 - \frac{t}{\epsilon}\]} + \frac{t}{\epsilon}
\int_r^R |g'(u)|d(u/t).$$ Applying Jensen's inequality repeatedly
(see explanation below), we conclude that
$$A\[\frac{g(r)}{\epsilon}\] \leq
\[1-\frac{t}{\epsilon}\]
A\[\frac{g(R)}{\epsilon\[1-\frac{t}{\epsilon}\]}\] +
\frac{t}{\epsilon} A\[\int_r^{R} |g'(u)|d(u/t)\] \leq $$
$$A\[\frac{g(R)}{\epsilon}\]+ \frac{1}{\epsilon}\int_{r}^{R}
A(|g'(u)|) du.$$ The second inequality uses the fact that $A$ is
convex and $A(0)=0$ along with the bound $\frac{t}{\epsilon}<1/2$
for the first term and a second application of Jensen's inequality
to the probability measure $d(u/t)$ on $(r, R)$ for the second term.
We conclude that
$$A\[\frac{g(r)}{\epsilon}\] \leq A\[\frac{g(R)}{\epsilon}\] + \frac{1}{\epsilon}\int_r^R
A(|g'(u)|) du.$$ Now we integrate both sides of the inequality,
using Fubini's theorem on the second term of the right hand side to
change order of integration.  We then get the result $$\int_0^R
A\[\frac{g(r)}{\epsilon}\] r^{d-1} dr \leq \int_0^R r^{d-1}
A\[\frac{g(R)}{\epsilon}\]dr + \frac{1}{\epsilon}\int_0^R \int_r^R
r^{d-1} A\[|g'(u)|\]dudr = $$ $$A\[\frac{g(R)}{\epsilon}\]
\frac{R^d}{d} + \frac{1}{\epsilon}\int_0^{R}\int_0^u r^{d-1}
A\[|g'(u)|\] drdu = $$
$$A\[\frac{g(R)}{\epsilon}\] \frac{R^d}{d} + \frac{1}{\epsilon}\int_0^{R}\frac{u^d}{d} A(|g'(u)|)
du.$$ Since $u^d < Ru^{d-1}$ we have the bound (from the modified
first assumption in the lemma statement) $$\int_0^R
A\[\frac{g(r)}{\epsilon}\] r^{d-1}dr \leq A\[\frac{g(R)}{\epsilon}\]
\frac{R^d}{d} + \frac{CR}{\epsilon}.$$ For any fixed $A$,
$\epsilon$, and $C$, we can assume $\alpha$ and $\beta$ are small
enough so that the latter expression is less than one (since $R \leq
\alpha$ and $g(R)\leq \beta$). \qed \end{proof}

\subsection{Interpolations from $d-1$-dimensional $L^1$ and $L^{1,A}$ bounds}

Let $G \subset \mathbb R^{d-1}$ be the $d-1$-dimensional unit cube.
Suppose we are given a linear function $f_u:\mathbb R^d \mapsto
\mathbb R$ defined by $f_u(\eta) = (u,\eta)$ (where $u \in \mathbb
R^d$), and a function $f$ defined on $G$.

For any $\epsilon$, we can construct an {\it interpolation between
$f_u$ and $f$} on $G \times [0,\epsilon]$ as follows.  Then write
$$g(\eta) = \frac{\eta_n}{\epsilon} f_u(\eta_1,\eta_2,\ldots, \eta_{n-1},\epsilon) +
\frac{\epsilon-\eta_n}{\epsilon}f(\eta_1,\eta_2, \ldots,
\eta_{n-1}).$$ The following lemma implies that, under appropriate
conditions, we can make $\epsilon$ small and still maintain an upper
bound on the ``average energy'' of the interpolation function $g$
over $G \times [0,\epsilon]$.

\begin{lem}  \label{sliceinterpolation} Fix a constant $C > 0$.  For every $\epsilon>0$, there
exists a $\delta>0$ so that the following three statements \begin{enumerate} \item $A(u) \leq C$
\item $\int_G A(\nabla f(\eta))d\eta \leq C$
\item $||h(\eta)||_1 \leq \delta$, where $h(\eta) =
f_u(\eta_1,\ldots,\eta_{n-1},0)-f(\eta_1,\ldots,\eta_{n-1})$ \end{enumerate} together imply that
$\frac{1}{\epsilon} \int_{[G \times \epsilon]}A(g(\eta))d\eta \leq 4C$. \end{lem}

\begin{proof} If $g_{\eta_1},\ldots, g_{\eta_d}$ are the partial derivatives of $g$, we can write
$$\frac{1}{\epsilon} \int_{G \times [0,\epsilon]} A(\nabla g(\eta))d\eta = \sum_{i=1}^d
\frac{1}{\epsilon} \int_{G \times [0,\epsilon]} A(\nabla
g_{\eta_i}(\eta))d\eta.$$ Since the first $d-1$ derivatives are
weighted averages of derivatives of $f$ and $f_u$, the sum of the
first $d-1$ terms is bounded by $\int_G A(\nabla f(\eta))d\eta +
\sum_{i=1}^{d-1} A(u_i) \leq C+C = 2C$.  The $d$th term is given by
$$\int_G
A\[\frac{f_u(\eta_1,\ldots,\eta_{n-1},\epsilon)-f(\eta_1,\ldots,\eta_{n-1})}{\epsilon}\]d\eta
= \int_G A\[\frac{h(\eta)}{\epsilon}+u_n\]d\eta.$$ Write $B(\eta) =
\frac{A(\eta)}{2C}$.  For any $\delta_0 >0$, Lemma \ref{fancybound}
implies that for $\delta$ small enough we have $||f||_B \leq
\delta_0$.  Note that $||u_n||_{B,G} \leq 1/2$ (where $u = (u_1,
\ldots, u_n)$, and we treat $u_n$ as a constant function).  If
$\delta_0$ is small enough so that $||\frac{h}{\epsilon}||_{B,G}
\leq \frac{1}{2}$, then $||\frac{h}{\epsilon} + u_n||_B \leq 1$ and
hence $\int_G A\[\frac{h(\eta)}{\epsilon}+u_n\]d\eta \leq 2C$. \qed
\end{proof}

Now, suppose that $G_1,\ldots, G_{2d}$ are the $2d$ faces of the
$d$-dimensional unit cube $D \subset \mathbb R^d$ and fix some
$\epsilon < 1/2$.  Let $D_{\epsilon}$ be the cube of side length
$1-2\epsilon$, positioned in $\mathbb R^d$ so as to be concentric
with $D$.  Given $f:D \mapsto \mathbb R$, we can interpolate between
$f$ and $f_u$ as follows.

For $1\leq \leq 2d$, let $G_i'$ be the face obtained by shifting $G_i$ by $\epsilon$ units in a
perpendicular direction, so that it borders $D'$.  Let $g_i$ be the linear interpolation between
between $f$ (on $G_i$) and $f_u$ (on $G_i'$), similar to the linear interpolation on $G \times
[0,\epsilon]$ described above; we can extend $g_i$ to all of $D$ by writing $g_i(\eta) = f_u(\eta)$
whenever $\eta$ and $G_i$ are on opposite sides of $G_i'$.

Now, we define an {\it inner interpolation function} $$g_I(\eta) = \begin{cases}
     \inf_{1 \leq i \leq 2d}\{g_i(\eta)\}, & f(\eta) \leq \inf_{1 \leq i \leq 2d} g_i(\eta)\\
     \sup_{1 \leq i \leq 2d}\{g_i(\eta)\}, & f(\eta) \geq \sup_{1 \leq i \leq 2d} g_i(\eta)\\
     f(\eta), & \text{otherwise.} \\
\end{cases}$$ Since $f(\eta)=g_i(\eta)$, for some $i$, for all $\eta$ on the boundary of $D$, we
have $g(\eta) = f(\eta)$ on the boundary of $D$.  Moreover, since all of the $g_i$'s are equal to
$f_u$ inside $D'$, we have $g(\eta) = f_u(\eta)$ for $\eta \in D'$.  Similarly, we define an {\it
outer interpolation function} by letting $g_i'$ be the interpolation between $f_u$ (on $G_i$) and
$f$ (on $G_i'$ and the portion of $D$ that lies on the opposite side of $G_i'$ from $G_i$). Then
$$g_O(\eta) = \begin{cases}
     \inf_{1 \leq i \leq 2d}\{g_i(\eta)\}, & f_u(\eta) \leq \inf_{1 \leq i \leq 2d} g_i(\eta)\\
     \sup_{1 \leq i \leq 2d}\{g_i(\eta)\}, & f_u(\eta) \geq \sup_{1 \leq i \leq 2d} g_i(\eta)\\
     f_u(\eta), & \text{otherwise.} \\
\end{cases}$$

This $g_O(\eta)$ is equal to $f_u(\eta)$ on the boundary of $D$ and $f(\eta)$ inside of $D'$.  The
following lemma gives energy bounds on $g_O(\eta)$ and $g_I(\eta)$.

\begin{lem} \label{preboxinterpolation} Let $G$ be the outer surface of a unit cube $D \subset \mathbb
R^d$. Fix a constant $C > 0$.  For every $\epsilon>0$, there exists a $\delta>0$ for which the
following four statements \begin{enumerate} \item $A(u) \leq C$ \item $\int_{D \backslash D'}
A(\nabla f(\eta))d\eta \leq C |D \backslash D'|$
\item $\int_{\partial D} A(\nabla(f(\eta))d\eta$ is well defined as a $d-1$-dimensional integral and is less than or
equal to $2d C$
\item $||f - f_u||_{1, \partial D} \delta$ \end{enumerate} together imply $$|D \backslash
D'|^{-1} \int_{D \backslash D'}A(g_I(\eta))d\eta \leq (4C)(2d+1) = 8dC + 4C$$ and the analogous
statement for $g_O$. \end{lem}

\begin{proof} The main observation is that $\nabla g_I(\eta) \leq \max \{ f(\eta), \nabla
g_1(\eta), \nabla g_2(\eta), \ldots, \nabla g_{2d}(\eta) \}$ for almost all $\eta$ (and the same is
true for $g_O$). Arguments similar to Lemma produce a bound of $4C$ on $|D \backslash D'|^{-1}
\int_{D \backslash D'}A(g_i(\eta))d\eta$ for each $i$, and the result follows. \qed \end{proof}

We can actually derive a similar version of the above lemma in which
the bound on $\int_{\partial D} A(\nabla(f(\eta))d\eta$ is omitted.
Define $g_I^{\gamma, \epsilon}$ as follows.  First, for some $0 <
\gamma < 1$, we define $f_{\gamma} D \mapsto \mathbb R$ by
$f_{\gamma}(\eta) = \frac{1}{\gamma}f(\gamma \eta)$. (Assume here
that $D$ is centered at the origin; so $f_{\gamma}$ is essentially
``zooming in'' on the portion of $f$ defined on the cube
$D_{\gamma}$ of side length $\gamma$ that is concentric with $D$.)
The construct the interpolations $g_I$ and $g_O$ (as described
above) using $f_{\gamma}$ instead of $f$.  Then we ``zoom back out''
by writing
$$g_I^{\gamma, \epsilon}(\eta) = \begin{cases} \gamma g_I\[\frac {\eta}{\gamma}\] & \eta \in D_\gamma\\
f(\eta) & \text{otherwise.} \\ \end{cases}$$

Now we have:

\begin{lem} \label{boxinterpolation} Let $D$ be the unit cube.  There exists a constant $C_0>0$ for
which the following is true. Fix a constant $C > 0$. For all sufficiently small $\epsilon>0$, there
exists a $\delta>0$ for which the following three statements \begin{enumerate}
\item $A(u) \leq C$
\item $\int_D |D|^{-1} A(\nabla f(\eta))d\eta \leq C$
\item $||f - f_u||_{1,D} \leq \delta$ \end{enumerate} together imply that for some
$\gamma$ with $1-\epsilon < \gamma < 1$, $$\int_{D \backslash
D_{1-2\epsilon}}A(g_I^{\gamma, \epsilon}(\eta))d\eta \leq C_0 C|D
\backslash D_{1-2\epsilon}|.$$ In particular, the average value of
$A(\nabla g_I^{\gamma, \epsilon})$, over the region where it fails
to be equal to $f_u$, is bounded above by $C_0 C$. A similar
statement is true for $g_O^{\gamma, \epsilon}$, where the roles of
$f_u$ and $f$ are reversed. Also, the conclusion of the lemma
remains true if we replace $D$ by $\frac{1}{n}D$, provided we
$\delta$ with $\frac{\delta}{n^{-d-1}}$. \end{lem}

\begin{proof} Since $\int_D A(\nabla f(\eta))d\eta \leq C$, then we must have $\int_{\partial
D_{\gamma}} A(\nabla(f(\eta))d\eta \leq 2C/\epsilon$ for a set of $\gamma$ values in $[1-\epsilon,
1]$ with measure at least $\epsilon/2$.  Since $||f-f_u||_1 \leq \delta$, then we must have
$$\int_{\partial D_{\gamma}} ||f-f_u||_1 \leq 2\delta/\epsilon$$ for at least one of the $\gamma$
values for which $\int_{\partial D_{\gamma}} A(\nabla(f(\eta))d\eta
\leq 2C/\epsilon$.  For $\epsilon$ sufficiently small, the volume of
$D_{\gamma} \backslash D_{\gamma (1-\epsilon)}$ is bounded above by
$3d \epsilon$. The result then follows by applying Lemma
\ref{preboxinterpolation} to $f_{\gamma}$ (as defined above). The
argument for $g_O^{\gamma, \epsilon}$ is similar.  The analogous
statement for $\frac{1}{n}D$ follows immediately from the rescaling
$f_n:\frac{1}{n}D \mapsto \mathbb R$ given by $f_n(\eta) =
\frac{1}{n}f(n\eta)$. \qed \end{proof}

\subsection{Discrete interpolation lemma}
The proof of the following discrete analog of Lemma
\ref{boxinterpolation} is virtually identical to the proof of its
continuous counterpart.  Let $\Delta_n = [0,n]^d - v$, where $v \in
\mathbb Z^d$ is the vector with all components equal to $\lfloor n/2
\rfloor$.  If $a \not \in \mathbb Z$, write $\Delta_a =
\Delta_{\lfloor a \rfloor}$.

\begin{lem} \label{insidebound} Fix a constant $C>0$.  There exists a constant $C_0$ such that the
following is true: For all sufficiently small $\epsilon > 0$, there
exists a $\delta > 0$ such that for all sufficiently large $n$ and
any $\phi:\Delta_n \mapsto \mathbb R$, the following three
statements \begin{enumerate} \item $A(u) \leq C$
\item $H^o_{\Delta_n}(\phi) \leq C |\Delta_n|$ (where $\Phi$ has nearest neighbor gradient potentials
determined by $A$)
\item $\sum_{x \in \Delta_n} |\phi(x) - (x,u)| \leq \delta |\Delta_n|$ \end{enumerate} together imply that there
exists an {\it inside interpolation function} $\psi_I$ for which \begin{enumerate}
    \item $\psi_I (x) = \phi(x)$ for all $x \in \partial \Delta_n$.
    \item $\psi_I(x) = (u,x)$ for all $x$ contained in
$\Delta_{(1-\epsilon)n}$.
    \item $H^{\overline{\Phi}}_{\Delta_n \backslash \Lambda_{(1-\epsilon)n}}(\psi_I) \leq
        C_0 C\epsilon n^d$. \end{enumerate}The same criteria imply the
existence of an {\it outside interpolation function} $\psi_O :
\Delta_n \mapsto \mathbb R$ for which
\begin{enumerate}
    \item $\psi_O(x) = (u,x)$ for all $x \in \partial\Delta_n$.
    \item $\psi_O(x) = \phi(x)$ for all $x$ contained in
$\Delta_{(1-\epsilon)n}$.
    \item $H^{\overline{\Phi}}_{\Delta_n \backslash \Delta_{(1-\epsilon)n}}(\psi_O) \leq
        C_0 C \epsilon n^d$. \end{enumerate}\end{lem}

\section{Approximation by ``mostly linear'' functions} \label{mlapproximationsection}
In this section, we will see that every $f \in L^{1,A}(D)$ can be approximated by a function which
``mostly agrees'' with one of the linear approximators $F_n$ (defined in the following lemma) and
whose total energy outside of the areas on which it agrees with $F_n$ is small.  This construction
will be useful to us in Chapter \ref{LDPchapter}.

\begin{lem} \label{piecewiselinear} Fix $f \in L^{1,1}(D)$.  Let $F_n(\eta)$ be the piecewise
linear (linear on simplices) function (not continuous in general)
whose mean and mean gradient are equal to those of $f$ on each
simplex of $\frac{1}{n} \mathbb Z^d$. Then
$$\lim_{n \rightarrow \infty} ||\nabla f - \nabla F_n||_1 = 0.$$ \end{lem}

\begin{proof} This merely says $\nabla f$ can be approximated by step functions. \qed \end{proof}

\begin{lem} For any $f \in L^{1,1}(D)$, $||f - F_n||_1 = o(\frac{1}{n})$. \end{lem}

\begin{proof} This follows from applying Lemma \ref{piecewiselinear} and Theorem \ref{imbedding} to
$f-F_n$. \qed \end{proof}

\begin{lem} \label{mostlylinearapproximator} Suppose $f \in L^{1,A}(D)$ and fix $\epsilon > 0$.
Then for all $n$ sufficiently large, there exists a function
$F_{\epsilon} \in L^{1,A}(D)$ that is equal to $F_n$ on a closed
subset $D'$ of $D$, where the volume of $D \backslash D'$ is less
than $\epsilon$ and $\int_{D-D'} A(\nabla F_{\epsilon}) \leq
\epsilon$. \end{lem}

\begin{proof} First, for any $\epsilon_0 > 0$, choose $C \geq 2\frac{\int_D A(\nabla
f(\eta))d\eta}{\epsilon_0}$. For any $\delta > 0$, the fraction of
boxes $b$ (of side length $1/n$) that satisfy $||f-F_n||_{1,b} \leq
\frac{\delta}{n^{d+1}}$ tends to zero in $n$.  It follows that for
sufficiently large $n$ and any $\delta$, the three conditions of
Lemma \ref{boxinterpolation} will apply to at least a $1-\epsilon_0$
fraction of the boxes of the mesh $\frac{1}{n} \mathbb Z^d$. We can
now define our interpolation function $F_{\epsilon}$ to be equal to
$f$ on the boxes for which the conditions of Lemma
\ref{boxinterpolation} do not apply and equal to the interpolation
described by Lemma \ref{boxinterpolation} on the boxes on which they
do apply.  It is clear from the lemma that if we take a sufficiently
small $\epsilon_0 < \epsilon$, we can also arrange to have
$\int_{D-D'} A(\nabla F_{\epsilon}) \leq \epsilon$. \qed \end{proof}

In Chapter \ref{LDPchapter}, given a function $f$ defined on $D$,
and good approximation $D_n$ of $f$, we will sometimes use a very
naive approach (described in the following lemma) for defining
functions $\phi_n: D_n \mapsto \mathbb R$ which (appropriately
rescaled) approximate $f$.

\begin{lem} \label{discretefunctionapproximator} Fix: $f: D \mapsto \mathbb R$ and suppose that
$\int_{D} A(\nabla f(\eta))d\eta$ is finite. Let the sequence $D_n$
be a good approximation of $D$. Then the sequence of functions
$\phi_n: D_n \mapsto \mathbb R$ (approximating $f$) defined by
$\phi_n(y) = n^{d+1} \int_{\{\eta : \lfloor n\eta \rfloor = y\}}
f(\eta)d\eta$ satisfies $|D_n|^{-1}H^o_{D_n}(\phi_n) \leq C_0
\int_{D} A(\nabla f(\eta))d\eta$ for some $C_0$ independent of $n$
and $f$. \end{lem}

\begin{proof} Write $s_n(\eta)$ for the cube $\{\eta : \lfloor n\eta \rfloor = y \}$. Because $D_n$
is a good approximation, all of the cubes $s_n(\eta)$, for $\eta \in
D_n$, are subsets of $D$. Note further that $$\phi_n(\eta) -
\phi_n(\eta+e_i) = n^{d+1} \int_{s_n(\eta)} \int_{t=0}^{\frac{1}{n}}
\frac{\partial}{\partial \eta_i} f(\eta+te_i))dt d\eta.$$ By
Jensen's inequality (and the fact that $n^{d+1} \int_{s_n(\eta)}
\int_0^{1/n} dt d\eta = 1$), we have
$$A(\phi_n(\eta) - \phi_n(\eta+e_i)) \leq n^{d+1} \int_{s_n(\eta)} \int_{t=0}^{\frac{1}{n}}
A\[\frac{\partial}{\partial \eta_i}f(\eta)\]d\eta. $$ The result
follows by summing these integrals over all edges of $D_n$. \qed
\end{proof}

\chapter{Limit equalities and the variational principle} \label{LDPempiricalmeasurechapter}

The first goal of this chapter is to derive some equivalent
definitions of specific free energy (making use of the notion of
``empirical measure'' of a configuration $\phi: \Lambda_n \mapsto
E$) and an equivalent definition of the surface tension.  In Section
\ref{limitequalitysection}, we prove the equivalence of these
definitions in the simplest setting: when $\Phi$ is an SAP and $(E,
\mathcal E)$ is $\mathbb R$, endowed with the Borel
$\sigma$-algebra.  In Section \ref{limitequalityother}, we describe
some (relatively minor) technical modifications to the arguments of
Section \ref{limitequalitysection} that enable us to extend these
results to all perturbed SAPs, as well as to higher dimensional real
and discrete models ($E = \mathbb R^m$ or $E = \mathbb Z^m$).

The second goal of the chapter is to prove the second half of the
variational principle---namely, that {\it every} $\mathcal
L$-invariant gradient phase $\mu$ of slope $u \in U_{\Phi}$ has
minimal specific free energy among $\mathcal L$-invariant measures
of slope $u$ and hence satisfies $SFE(\mu) = \sigma(u)$.  In Section
\ref{variationalprinciplesection} we prove this when $\Phi = \Phi_V$
is an ISAP and $(E, \mathcal E)$ is $\mathbb R$, endowed with the
Borel $\sigma$-algebra. In Section \ref{variationalother}, we
describe the (again, relatively minor) technical modifications to
the arguments of Section \ref{variationalprinciplesection} that
enable us to extend them, for $E = \mathbb R^m$ or $E = \mathbb
Z^m$, to all perturbed ISAP's and perturbed LSAP's.   (It is not
known whether the second half of the variational principle---like
the first half---holds for all perturbed SAPs.)

\section{Limit equalities: $PBL(\mu)$, $FBL(\mu)$, and $SFE(\mu)$} \label{limitequalitysection}
Throughout this section, we will assume that $E = \mathbb R$ and
$\Phi$ is an SAP.  Denote by $\mathcal A$ the topology of local
convergence on $\mathcal P(\Omega, \mathcal F^{\tau})$ and by
$\mathcal B$ the basis for that topology consisting of finite
intersections of sets of the form $\{ \mu : \mu(F) < \epsilon \}$
where $0<\epsilon <1$ and $F: \Omega \rightarrow \mathbb R$ is
$\mathcal F_{\Lambda}^{\tau}$-measurable for some $\Lambda \subset
\subset \mathbb Z^d$, and $F$ is bounded between $0$ and $1$.

For notational convenience, throughout the following two chapters,
we will also augment the space $\Omega$ to be the space of functions
from $\mathbb Z^d$ to $\mathbb R \cup \{\infty\}$ instead of merely
from $\mathbb Z^d$ to $\mathbb R$. Whenever we deal with a
real-valued function $\phi$ that we have defined only on a subset
$\Lambda$ of $\mathbb Z^d$, we tacitly assume that it extends to a
function in $\Omega$ for which $\phi(x) = \infty$ for $x \not \in
\Lambda$. (So in this context, the statement $\phi(x) = \infty$ is a
way of saying that $\phi(x)$ has not been defined.)

Given $u \in \mathbb R^d$, let $\phi_u: \mathbb Z^d \mapsto \mathbb
R$ be a function with an $\mathcal L$-periodic gradient and
$\epsilon$ a positive constant for which \begin{enumerate}
\item $\phi_u$ has asymptotic slope given by $u$.
\item There exists a finite constant $C$ such that for any edge $(x,y)$ in $\mathbb Z^d$ and any
$\phi$ bounded between $\phi_u-2\epsilon$ and $\phi_u + 2\epsilon$, we have: $V_{x,y}(\phi(y)-
\phi(x)) \leq C$. \end{enumerate}

Arguments similar to those in Lemma \ref{finiteenergyU} and Lemma
\ref{torusexistence} imply that---whenever $\Phi$ is an SAP and $u
\in U_{\Phi}$---there exists at least one such pair
$(\phi_u,\epsilon)$.  As in previous chapters, we take $\Lambda_n =
[0,kn - 1]^d$ where $k$ is chosen so that $k \mathbb Z^d \subset
\mathcal L$.  Let $C^u_n$ be the subset of maps $\phi:\Lambda_n
\mapsto \mathbb R$ for which $|[\phi(x) - \phi(x_0)] - [\phi_u(x) -
\phi_u(x_0)]| \leq \epsilon$ for all $x \in \partial \Lambda_n
\backslash \{x_0\}$. (Here, $x_0$ is the origin, so it lies on one
corner of the boundary of $\Lambda_n = [0,kn-1]^d$.) The functions
in $C^u_n$ are those functions which (up to additive constant)
closely approximate the periodic function $\phi_u$ on the boundary
of $\Lambda_n$.

Given a $B \in \mathcal B$, let $B_n$ be the set of functions $\phi:
\Lambda_n \mapsto \mathbb R$ whose empirical measures lie in $B$. To
say this precisely, let $\theta_x$ denote translation by $x \in
\mathcal L$.  Define $L_n(\phi)$, a measure on $(\Omega, \mathcal
F)$, by $$L_n(\phi) = |\Lambda_n \cap \mathcal L|^{-1} \sum_{x \in
\Lambda_n \cap \mathcal L} \delta_{\theta_x \phi}.$$ Then write $B_n
= \{ \phi | L_n(\phi) \in B \}$.  In the following definitions, we
will assume that $u \in U$ and that $\phi_u$ and $\epsilon$ are
fixed.  After we prove Theorem \ref{limitequality}, it will be clear
that the definition of $PBL$ is independent of the particular choice
of $\phi_u$ and $\epsilon$.  We use the initials $PBL$ to mean
``pinned boundary limit'' and write
$$PBL^u_B(\mu) = \limsup_{n \rightarrow \infty} - |\Lambda_n|^{-1} \log \[ \int 1_{C^u_n \cap B_n}
e^{-H^o_{\Lambda}(\phi)} \prod_{x \in \Lambda_n \backslash \{x_0\}
}d\phi(x) \].$$ In the integrals in the above limiting sequence, we
assume $\phi(x_0)$ is set to $0$ and $\phi(x) = \infty$ for $x \not
\in \Lambda_n$.  We also write
$$PBL^u(\mu) = \sup_{B \ni \mu, B \in \mathcal B} PBL^u_B(\mu)$$
$$PBL(\mu) = PBL^{S(\mu)}(\mu).$$ If $u \not \in U$, then write $PBL^u(\mu) = \infty$; in
particular, if $S(\mu) \not \in U$, then $PBL(\mu) = \infty$.  We now define the ``free boundary
limit'' as follows:
$$FBL_B(\mu) = \liminf_{n \rightarrow \infty} -|\Lambda_n|^{-1} \log \[ \int 1_{B_n}
e^{-H^o_{\Lambda}(\phi)} \prod_{x \in \Lambda_n \backslash \{x_0\} }d\phi(x) \]$$
$$FBL(\mu) = \sup_{B \ni \mu, B \in \mathcal B} FBL_B(\mu)$$

\begin{thm} \label{limitequality} If $\Phi$ is an SAP, $E = \mathbb R$, and $\mu \in \mathcal
P_{\mathcal L}(\Omega, \mathcal F^{\tau})$, then $$SFE(\mu) =
FBL(\mu) = PBL(\mu).$$ \end{thm}

It is obvious that $FBL(\mu) \leq PBL^u(\mu)$ for any $u$; in particular $FBL(\mu) \leq PBL(\mu)$.
In the following three subsections, we will prove Theorem \ref{limitequality} in three steps:
\begin{enumerate} \item $PBL(\mu) \leq SFE(\mu)$ when $\mu$ is ergodic
\item $PBL(\mu) \leq SFE(\mu)$ for any $\mu$
\item $SFE(\mu) \leq FBL(\mu)$
\end{enumerate}

\subsection{$PBL(\mu) \leq SFE(\mu)$ for $\mu$ ergodic} \label{ergodicPBLSFEsection}

First, as a simple consequence of the ergodic theorem, we show that if the gradient of $\phi$ is
chosen from an $\mathcal L$-ergodic gradient measure $\mu$ with finite slope, then $\phi$ closely
approximates a plane, in an $L^1$ sense, with high probability.  Precisely, we show the following:

\begin{lem} \label{ergodicL1} Let $\mu \in \mathcal P_{\mathcal L}(\Omega, \mathcal F^{\tau})$ be
$\mathcal L$-ergodic with finite slope $u \in \mathbb R^d$.  Then
for any fixed $\epsilon$, we have
$$\lim_{n \rightarrow \infty} \mu \{ |\Lambda_n|^{-1} \sum_{x \in \Lambda_n} |\phi(x) - (x,u) -
\phi_{\Lambda_n}| \geq \epsilon n^{d+1} \} = 0,$$ where
$\phi_{\Lambda_n}$ is the mean value of $\phi(x) - (x,u)$ on
$\Lambda_n$. \end{lem}

\begin{proof} First, we may assume without loss of generality that when $\phi$ is chosen (defined
up to additive constant) from $\mu$, we have $u=0$ and $\mu (\phi(y)
- \phi(x)) = 0$ for each $x,y \in \mathbb Z^d$. If this is not the
case, we can let $f$ be a function (determined from $\mu$ up to
additive constant) for which $f(y) - f(x) = \mu (\phi(y) - \phi(x))$
for each $x,y \in \mathbb Z^d$, and then replace $\mu$ by $\mu - f$.
Thus, we need only show in this case that for $\epsilon>0$, we have
$$\lim_{n \rightarrow \infty} \mu \{ |\Lambda_n|^{-1} \sum_{x \in
\Lambda_n} |\phi(x) - \phi_{\Lambda_n}| \geq \epsilon n^{d+1} \} =
0,$$ where $\phi_{\Lambda_n}$ is the mean value of $\phi$ on
$\Lambda_n$.

When $\phi$ is chosen from $\mu$, the classical ergodic theorem
(see, e.g., Chapter 14 of \cite{G}) states that for any $\gamma > 0$
$$\lim_{n \rightarrow \infty} \mu \{ |\frac{\sum_{x \in \Lambda_n, |x-y|=1}\nabla
\phi(x)}{|\Lambda_n|}| \geq \gamma \} = 0.$$

A trivial corollary of the ergodic theorem states that if $h:D
\rightarrow \mathbb R^d$ is any continuous, bounded function on the
unit cube $D$, then
$$\lim_{n \rightarrow \infty} \mu \{ |\frac{\sum_{x \in \Lambda_n, |x-y|=1}(\nabla
\phi(x),h(x/n))}{|\Lambda_n|}| \geq \gamma \} = 0.$$ Now, write
$\Delta_n = (nk) e_1 + \Lambda_n$, where $e_1$ is a basis vector of
$\mathbb Z^d$.  So $\Delta_n$ and $\Lambda_n$ are adjacent blocks. A
discrete integration by parts gives $$\sum_{y \in \Delta_n} \phi(y)
- \sum_{x \in \Lambda_n} \phi(x) = (nk)^{d+1} \sum_{x \in \Delta_n
\cup \Lambda_n} [\phi(x) - \phi(x - e_1)] h(x/n),$$ where $h(v) =
v_1 $ for $v = (v_1, \ldots, v_d) \in D$ and $h(v) = 2 - v_1 $ for
$v \in D + e_1$. It follows that for any $\gamma$, the probability
that the mean value of $\phi$ on $\Lambda_n$ differs from the mean
value on $\Delta_n$ by more than $\gamma n^{d+1}$ tends to zero in
$n$.

Similarly, if we take a large integer $c$ and form a $(ckn)^d$ cube by joining $c^d$ translated
$\Lambda_n$ blocks, then for any fixed $\gamma$, the probability that the mean value on any one of
these blocks differs from the value on any other by $\gamma n^{d+1}$ will tend to zero in $n$.

Let $C$ be supremum of $\mu ( |\nabla \phi(x)| )$ for $x \in \mathbb
Z^d$.  The probability that the mean value of $|\nabla (\phi)|$ on
one of these blocks is greater than $2C$ also tends to zero in $n$.
It follows from the first part of Theorem \ref{imbedding} and
Corollary \ref{discretenormcor} that for some $C_0$ (independent of
$\delta$) we have, for any $\delta_0$, $$\lim_{n \rightarrow \infty}
\mu \{ |\Lambda_{n}|^{-1} \sum_{x \in \Lambda_{n}} |\phi(x) -
\phi_{\Lambda_{n}}| \geq C_0 n^{d+1} \} = 0 .$$

Now, taking $N = cn$, the above claims imply that (when $N$ is restricted to multiples of $c$)
$$\lim_{N \rightarrow \infty} \mu \{ |\Lambda_{N}|^{-1} \sum_{x \in \Lambda_{N}} |\phi(x) -
\phi_{\Lambda_{N}}| \geq C_0 c^d n^{d+1} = (C_0/c) N^{d+1} \} = 0.$$

It is not hard to see (by applying the same result and restricting
to a slightly smaller box) that the above remains true if $N$ is not
restricted to multiples of $c$.  Since the above is true for any
$c$, we have
$$\lim_{N \rightarrow \infty} \mu \{ |\Lambda_{N}|^{-1} \sum_{x \in \Lambda_{N}} |\phi(x) -
\phi_{\Lambda_{N}}| \geq \epsilon N^{d+1} \} = 0,$$ for all
$\epsilon > 0$.  \qed \end{proof}

A simple corollary is the following:

\begin{lem} \label{ergodicL1cor} Let $\mu \in \mathcal P_{\mathcal L}(\Omega, \mathcal F^{\tau})$
be $\mathcal L$-ergodic with finite slope $u \in \mathbb R^d$ and
finite specific free energy.  Then for any fixed $\epsilon$,
$$\lim_{n \rightarrow \infty} \mu \{ |\Lambda_n|^{-1} \sum_{x \in
\Lambda_n} |\phi(x) - (x,u) - \phi(0)| \geq \epsilon n^{d+1} \} =
0.$$  \end{lem}

\begin{proof} To deduce this from Lemma \ref{ergodicL1}, we need only show that the probability
that $\phi(x_0)-(x_0,u)$ differs from $|\Lambda_n|^{-1} \sum_{x \in
\Lambda_n} \phi(x) - (x,u)$ on $\Lambda_n$ by more than $\epsilon n$
tends to zero in $n$. By shift-invariance, we can show equivalently
that if we choose uniformly an $x$ in $\Lambda_{2n}$, then the
probability that $f(x) =  \phi(x)-(x,u)$ differs from $g(x) =
|\Lambda_n|^{-1} \sum_{y \in \theta_x \Lambda_n} \phi(y)-(y,u)$ by
$\epsilon n$ tends to zero in $n$.  However, an immediate
consequence of Lemma \ref{ergodicL1} is that the probability that
either $f(x)$ or $g(x)$ differs from $|\Lambda_{2n}|^{-1} \sum_{y
\in \Lambda_{2n}} \phi(y) - (y-u)$ by $\epsilon n$ tends to zero in
$n$. \qed \end{proof}

Now, we continue with the proof that $PBL(\mu) \leq SFE(\mu)$ when
$\mu$ is $\mathcal L$-ergodic and $S(\mu) = u \in U$.  Recall that
$SFE(\mu) = \lim_{n \rightarrow \infty} |\Lambda_n|^{-1}
FE_{\Lambda_n}(\mu)$ and $$PBL(\mu) = \sup_{B \ni \mu, B \in
\mathcal B} \limsup_{n \rightarrow \infty} |\Lambda_n|^{-1}
FE_{\Lambda_n}(\nu^B_n),$$ where $\nu^B_n$ is the Gibbs measure
$$Z^{-1} 1_{A_n \cap B_n} e^{-H^o_{\Lambda_n}} \prod_{x \in
\Lambda_n \backslash \{x_0\}} d \phi(x),$$ and $Z$ is chosen to make
the above a probability measure.  Let $\mu_n$ be the restriction of
$\mu$ to $\mathcal F^{\tau}_{\Lambda_n}$. It is now enough for us to
show that for each $B \ni \mu, B \in \mathcal B$, we have
$FE_{\Lambda_n}(\nu^B_n) \leq FE_{\Lambda_n}(\mu_n) +
o(|\Lambda_n|)$. Since $\nu^B_n$ has minimal free energy among
measures on $(\Omega, \mathcal F^{\tau}_{\Lambda_n})$ that are
supported on $A_n \cap B_n$, it will be sufficient to generate a
measure $\mu'_n$ on $(\Omega, \mathcal F^{\tau}_{\Lambda_n})$---also
supported on $A_n \cap B_n$ for sufficiently large $n$---which
satisfies $FE_{\Lambda_n}(\mu'_n) \leq FE_{\Lambda_n}(\mu_n) +
o(|\Lambda_n|)$.

We take $B$ to be the set of measures $\pi$ on $(\Omega, \mathcal
F^{\tau})$ for which $|\pi(F_i) - \mu(F_i)| \leq \epsilon$, for some
finite sequence of cylinder functions $F_i: \Omega \rightarrow
[0,1]$ and some constant $\epsilon$. Write $B'$ to be the set of
measures $\pi$ for which $|\pi(F_i) - \mu(F_i)| \leq \epsilon/2$ for
each $i$.  Now, $\mu_n$ is not necessarily supported on $A_n \cap
B_n$.  However, the ergodic theorem implies that the probability
that a sample from $\mu_n$ lies in $B'_n \subset B_n$ tends to one
as $n$ tends to $\infty$.

To sample from $\mu'_n$, we will first sample $\phi$ from $\mu_n$ (conditioned on $\phi$ lying in
$B'_n$), and then use a ``random truncation'' to alter $\phi$ in a way that forces it to lie in
$A_n \cap B_n$. In a separate step, we will check that $FE_{\Lambda_n}(\mu'_n) \leq
FE_{\Lambda_n}(\mu_n) + o(|\Lambda_n|)$ by showing that these random truncations change the free
energy by at most $o(|\Lambda_n|)$.

Let $D$ be the unit cube $[0,1]^d$ and let $f:D \rightarrow \mathbb
R$ be the piecewise-linear ``pyramid-shaped'' function for which
$f(z) = 0$ for $z \in \partial D$, $f(z_0)=1$ where $z_0$ is the
center point of $D$, and $f$ is linear on each line segment
connecting $z_0$ to a point in $\partial D$.  Now, given $\epsilon >
0$, we define a pyramid-shaped, $\epsilon$-sloped function on
$\Lambda_n$ by $p_n(x) = \phi_u(x) + n\delta f(\frac{x}{kn})$. If
$\epsilon$ is as described in the definition of $\phi_u$, then $p_n$
will have the property that $V_{x,y}(\phi(y) - \phi(x)) \leq C$
whenever $x,y \in \Lambda_n, |x-y|=1$ and $|(\phi(y) - \phi(x)) -
(p_n(y) - p_n(x))| < \epsilon$.

Now, we first define an ''upper truncated'' measure $\mu_n''$ as follows: to sample from $\mu''_n$,
first choose $\phi$ from $\mu$, taking the additive constant so that $\phi(x_0) = 0$.  Let $A^+$ be
the set of vertices $x \in \Lambda_n$ for which $\phi(x) > p_n(x)$ and let $A^-$ the set of $x \in
\Lambda_n$ for which $\phi(x) < - p_n(x)$. Then for each $x \in A^+$, we re-sample $\phi(x)$ from
$$Z^{-1} 1_{ \{p_n(x) \leq \phi(x) \leq p_n(x) + \epsilon \}} \exp\[ \sum_{y \in \Lambda_n
\backslash A^+_n, |x-y|=1} V_{x,y}(\phi(y) - \phi(x))\] d\phi(x),$$
where $Z$ is the appropriate normalizing constant.

\begin{figure} \begin{center} \leavevmode \epsfbox[0 0 372 225]{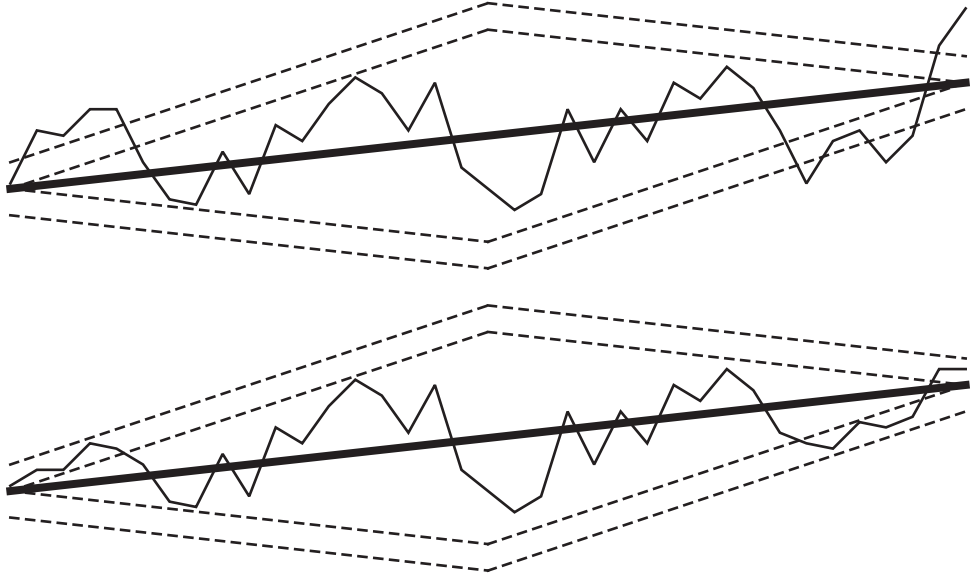} \end{center}
\caption[A random truncation] {Shown here are $\phi_u$ (dark lines); $\phi_u \pm \epsilon nf(nx)$
and $\phi_u(x) \pm (\epsilon nf(nx) + \epsilon)$ (dotted lines); and $\phi$ before and after a
random truncation (where $d=1$).} \label{truncate} \end{figure}

Now, we claim that $FE_{\Lambda_n}(\mu_n'') \leq
FE_{\Lambda_n}(\mu_n) + o(|\Lambda_n|)$.  To see this, note that by
Lemma \ref{conditionalentropy}, we can compute each relative entropy
in stages: to sample from $\mu_n$ or $\mu_n''$, we can first choose
$A^+$, then choose the values of $\phi$ conditioned on $A^+$.  The
relative entropy is the relative entropy of the first random process
(the choice of $A^+$) plus the expected relative entropy of the
second; similarly, the relative entropy of the choice of $\phi$
(conditioned on $A^+$) can be written as the relative entropy of the
choice of $\phi$ for $x \not \in A^+$ plus the expected relative
entropy of the choice of $\phi(x)$ for $x \in A^+$.  Since the
latter step is the only one which is different in the two process,
the difference in the two relative entropies depends only on the
expected difference in the relative entropy of the last two choices.
For a given choice of $A^+$ and $\phi$ (defined for $x \not \in
A^+$), write $$Z^{-1} \pi_0 = 1_{p_n(x) \leq  \phi(x), x \in A^+}
e^{H^o_{\Lambda_n}(\phi)} \prod_{x \in A^+} d \phi(x),$$ where $Z$
makes the above a probability a measure; let $\pi$ )(resp., $\pi''$)
be the conditional distribution of $\mu_n$ (resp. $\mu_n''$)
conditioned on $A^+$ and $\phi$ (for $x \not \in A^+$).  Since, by
Corollary \ref{ergodicL1cor}, the expected size of $|A^+|$ is
$o(|\Lambda_n|)$, it is enough for us to show that $\mathcal
H(\pi'', \pi_0) \leq \mathcal H(\pi, \pi_0) + c|A^+|$ for some
constant $c$. Since $\mathcal H(\pi, \pi_0)$ is positive, it is
enough to show $\mathcal H(\pi'', \pi_0) \leq c|A^+|$.  The reader
may check this last fact, again, by writing the relative entropy as
a sum of a sequence of expected conditional relative entropies, one
for each $x \in A^+$.  The key observation is that when one chooses
$\phi(x)$---and $\phi(y)$ have already been chosen for some nonempty
subset $S_x$ of the neighbors $y$ of $x$, and each $\phi(y) \leq
\phi_u(y) + \epsilon n f(y/n)$, then the measure $e^{-\sum_{y \in
S_x} V((\phi(y) - \phi(x)) s(y))} 1_{\phi(x) \leq \phi_u(x)+\epsilon
n f(y/n)} d\phi(x)$ (where $s(y)$ is $1$ if $x$ precedes $y$
lexicographically, and $-1$ otherwise) obtains a bounded fraction of
its mass in the interval $[\phi_u(x)+\epsilon n f(y/n),
\phi_u(x)+\epsilon n f(y/n)+1]$, independently of the precise values
of $\phi(y)$.

To sample $\phi$ from $\mu'_n$, first sample $\phi$ from $\mu_n''$, then apply a ``lower
truncation'' (analogous to the upper truncation described above, using $A^-$ instead of $A^+$) and
condition on the output lying in $B_n$. Since, by the ergodic theorem, the probability of the
latter event tends to one, we have $FE_{\Lambda_n}(\mu_n') \leq FE_{\Lambda_n}(\mu_n) +
o(|\Lambda_n|)$, as desired.

\subsection{$PBL(\mu) \leq SFE(\mu)$ for general $\mu$}

Since $SFE$ is affine, it is enough to show that $PBL(\mu)$ is ``strongly convex''---i.e., that
when $\mu = \int \nu w_{\mu}(\nu) d \nu$, we have $PBL(\mu) \leq \int PBL(\nu) w_{\mu}(\nu) d \nu$.

First we will show that $PBL$ is convex.  Suppose that $\mu_1$ and
$\mu_2$ have slopes $u_1$ and $u_2$ respectively, $0<a<1$, $\mu = a
\mu_1 + (1-a)\mu_2)$, and $u = au_1 + (1-a)u_2$.  It is enough to
show that $PBL^u_B{\mu}$ is less than or equal to $PBL^{u_1}(\mu_1)
+ PBL^{u_2}(\mu_2)$ for each $B \ni \mu, B \in \mathcal B$.  We may
assume $B$ is a finite intersection of sets of the form $\{\pi :
|\pi(F_j) - \mu(F_j)| \leq \epsilon \}$, where each $F_j : \Omega
\rightarrow [0,1]$ is a cylinder function.  For $i \in \{1,2\}$,
take $B^i$ to be the intersection of the sets $\{ \pi: |\pi(F_j) -
\mu_i(F_j)| \leq \epsilon/2 \}$.  By the definition of $PBL$, it is
now enough to show that $PBL^u_B{\mu}$ is less than or equal to $a
PBL^{u_1}_{B^1}(\mu_1) + (1-a) PBL^{u_2}_{B^2}(\mu_2)$. And to prove
this, it is enough to us to show that for any fixed $n$, letting $M$
get large, we have
$$- |\Lambda_M|^{-1} \log \[ \int 1_{A_M \cap B_M} e^{-H^o_{\Lambda}(\phi)} \prod_{x \in \Lambda_n
\backslash \{x_0\} }d\phi(x) \] \leq $$ $$- a |\Lambda_n|^{-1} \log \[ \int 1_{A_n \cap B^1_n}
e^{-H^o_{\Lambda}(\phi)} \prod_{x \in \Lambda_n \backslash \{x_0\} }d\phi(x) \] - $$
$$(1-a)|\Lambda_n|^{-1} \log \[ \int 1_{A_n \cap B^2_n} e^{-H^o_{\Lambda}(\phi)} \prod_{x \in
\Lambda_n \backslash \{x_0\} }d\phi(x) \] + o(1).$$

\begin{figure} \begin{center} \leavevmode \epsfbox[0 0 225 160]{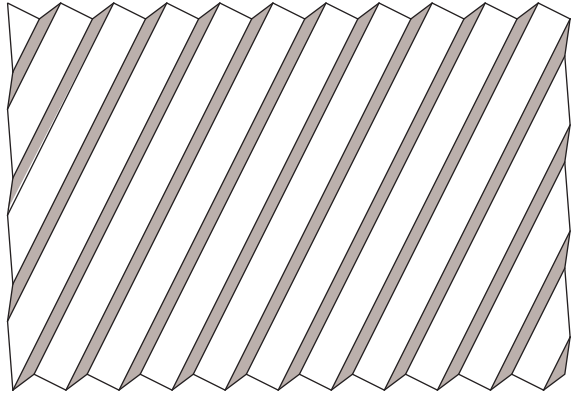} \end{center}
\caption[A ``washboard'' shaped surface] {A ``washboard'' shaped surface $\psi$ when $d=2$.  The
shaded regions have slope $u_1$ and the unshaded regions have slope $u_2$.} \label{washboard}
\end{figure}

Next, roughly speaking, we would like to combine $\phi_{u_1}$ and $\phi_{u_2}$ to form a
``washboard'' function of slope $u$ whose gradient agrees with that of $\phi_{u_1}$ on an $a$
fraction of the points and that of $\phi_{u_2}$ on a $(1-a)$ fraction of points.  (See Figure
\ref{washboard}.)

Write $p_1 = (1-a)$ and $p_2 = a$.  Fix some large integer $N$ and
consider the layered sequence of surfaces $\phi^j_{u_i}(x) =
\phi_{u_i}(x) + p_i Nj$ (for $j \in \mathbb Z$); let
$\tilde{\phi}_{u_i}(x,\eta)$ give the index $j$ of the lowest layer
which lies beneath the point $(x, \eta)$, i.e., the smallest $j$ for
which $\phi^j_{u_i}(x) \leq \eta$.

Now, write $\psi(x) = \inf \{ \eta | \tilde{\phi}_{u_1}(x,\eta) +
\tilde{\phi}_{u_2}(x,\eta) \leq 0$.  Now, fix $n$ and tile $\mathbb
Z^d$ with a sequence $\Delta_j$ of cubes of side-length $n$.  It is
not hard to see, for fixed $n$, that $\psi$ has asymptotic slope
given by $u$ and that the fraction of cubes on which the gradient of
$\psi$ fails to agree with that of either $\phi_1$ or $\phi_2$
throughout the cube is $o(\frac{1}{N})$.  Moreover, if $\epsilon$ is
the minimum of the $\epsilon$ values in the definitions of
$\phi_{u_1}$ and $\phi_{u_2}$, then the energy at any edge of $\phi$
will be less than or equal to $C$ whenever $\phi$ is sampled from
the box measure of radius $\epsilon$ centered at $\psi$.

Define $\nu_n$ to be the measure $1_{A_n \cap
B_n}e^{H^o_{\Lambda_n}} \prod_{x \in \Lambda_n \backslash \{x_0 \} }
d\phi(x)$; define $\nu^i_n$ analogously using $B^i_n$ instead of
$B$.  We can restate our goal as follows: we need to prove that as
$N$ and $M$ get large, $$|\Lambda_M|^{-1} FE (\nu_M) \leq a
|\Lambda_n|^{-1} FE(\nu^1_n) + (1-a) |\Lambda_n|^{-1} FE(\nu^2_n) +
o(1).$$  Since $\nu_M$ has minimal free energy among measures
supported on $A_M \cap B_M$, it is enough to generate some measure
$\nu'_M$ supported on $A_M \cap B_M$ for which the analogous
expression holds.

Now, we define a measure $\nu'_M$ on $(\Omega, \mathcal
F^{\tau}_{\Lambda})$ as follows. To sample from $\nu'_M$, first
choose $\phi$ from the radius-$\epsilon$ box measure centered at
$\psi$ on $\Lambda$. Then, for all cubes $\Delta_j$ on which the
gradient of $\psi$ is identically equal to that of either
$\phi_{u_1}$ or $\phi_{u_2}$ (say $\phi_{u_i}$), we re-sample $\phi$
from $\nu^i_n$ (when $\Delta$ is translated to coincide with
$\Lambda_n$ for the purposes of the definition). Letting $N$ and $M$
get large in such a way that $M/N \rightarrow \infty$, we see that
$$|\Lambda_M|^{-1} FE (\nu_M) \leq a |\Lambda_n|^{-1} FE(\nu^1_n) +
(1-a) |\Lambda_n|^{-1} FE(\nu^2_n) + o(1).$$  Moreover, if we modify
$\nu'_M$ by adding a truncation (as in the previous section) this
modification will also change the normalized free energy by at most
$o(1)$, and the result follows.

It remains to check that $PBL(\mu)$ is ``strongly convex''---i.e.,
that when $\mu = \int \nu w_{\mu}(\nu) d \nu$, we have $$PBL(\mu)
\leq \int PBL(\nu) w_{\mu}(\nu) d \nu.$$  However, it is obvious
from the definition of $PBL$ that it is lower semi-continuous in the
topology of local convergence.  So it is enough to observe that we
can approximate $\mu$ in this topology by a sequence of slope-$u$
weighted averages of the form $\mu^k = \sum_{i=1}^k a_{i,k}
\mu_{i,k}$ where $0 \leq a_i \leq 1$, the $\mu_i$ are ergodic, the
$\mu^k$ converge to $\mu$ in the topology of local convergence, and
$\limsup SFE(\mu^k) \leq SFE(\mu)$.

If one samples a sequence $\mu_i$ of ergodic components
independently from $w_{\mu}$, then the law of large numbers implies
that $\nu_n = \frac{1}{n} \sum_{i=1}^n \mu_i$ converges to $\mu$ in
the topology of weak local convergence and that $SFE(\nu_n)
\rightarrow SFE(\mu)$ almost surely.  The desired approximation is
now easily obtained by altering the coefficients of the $\nu_n$
slightly (so that the slope is exactly $u$ instead of approximately
$\nu$).

\subsection{Empirical measure argument: $SFE(\mu) \leq FBL(\mu)$}

Recall the definition $$FBL_B(\mu) = \liminf_{n \rightarrow \infty}
-|\Lambda_n|^{-1} \int 1_{B_n} e^{-H^o_{\Lambda}(\phi) \prod_{x \in
\Lambda_n \backslash \{x_0\} }}d\phi(x).$$

We can normalize the measure $1_{B_n} e^{-H^o_{\Lambda}(\phi)
\prod_{x \in \Lambda_n \backslash \{x_0\} }}d\phi(x)$ to produce a
probability measure on the set of functions from $\Lambda_n
\backslash \{x_0\}$ to $\mathbb R$.  We can extend this to a
probability measure $\nu^n_B$ on $\Omega$ in which the events
$\phi(x_0) = 0$ and $\phi(x = \infty), x \not \in \Lambda_n$ have
probability one.  Let $\Lambda_n^m$ be the subset of $\Lambda_n$
containing vertices which are at least $m$ units in distance from
the boundary of $\Lambda_n$.  Now, let $m(n)$ be some function of
$n$ for which $m(n)$ tends monotonically to $\infty$ in $n$ but
$m(n) = o(n)$.  Let $\mu^n_B = |\Lambda_n^{m(n)} \cap k \mathbb
Z^d|^{-1} \sum_{y \in \Lambda_n^{m(n)} \cap k \mathbb Z^d} \theta_y
\nu^n_B$.

\begin{lem} \label{fblsfelemma} As $n$ tends to infinity along an increasing sequence for which
$$-|\Lambda_n|^{-1} \int 1_{B_n} e^{-H^o_{\Lambda}(\phi) \prod_{x \in \Lambda_n \backslash \{x_0\}
}} d\phi(x) \rightarrow FBL_B(\mu),$$ at least one subsequential
limit $\mu_B$ of the measures $\mu^n_B$ exists; any such limit
satisfies $$SFE(\mu_B) \leq FBL_B(\mu).$$ \end{lem} \begin{proof}

Fix any integer $r>1$ and suppose that $\Delta_1, \ldots \Delta_j$
are disjoint translations of $\Lambda_r$ contained in $\Lambda_n$.
Then Lemma \ref{subadditive} implies that $$FE_{\Lambda_n} (\nu^n_B)
\geq \sum_{i=1}^j FE_{\Delta_i} + c\[|\Lambda_n \backslash
\cup_{i=1}^j \Delta_i| + j\],$$ where $c$ is a fixed
constant---namely, the minimum possible value of $FE_e(\mu)$ where
$e$ is a single edge.  Now, if $n$ is large enough so that $m(n) >
r$, then $\mu_n$, restricted to $\Lambda_r$, is a weighted average
of the restrictions of $\nu_n$ to translations $\Delta_1, \ldots,
\Delta_{(n-m(n))^d}$ of $\Lambda_r$.  We can divide $\{\Delta_i\}$
into $r^d = |\Lambda_r \cap k \mathbb Z^d|$ sub-collections
according to the values of the components of their lexicographically
minimal corners modulo $rk$; each such sub-collection consists of
disjoint copies of $\Lambda_r$.  Using convexity of free energy for
the first step, we may conclude: $$|\Lambda_r \cap k \mathbb
Z^d|^{-1} FE_{\Lambda_r} (\mu^n_B) \leq r^{-d} $$ $$
\sum_{i=1}^{(n-m(n)^d)} (n-m(n))^{-d} FE_{\Delta_i} \nu^n_B \leq $$
$$(1 - o(n))|\Lambda_n \cap k\mathbb Z^d|^{-1} FE_{\Lambda_n}(\nu^n_B)
- o(n) + c(r^{-d})$$ This implies that $$\limsup |\Lambda_r|^{-1}
FE_{\Lambda_r}(\mu_n) \leq \limsup |\Lambda_n|^{-1}
FE_{\Lambda_n}(\nu^n_B) = FBL_B(\mu).$$  By Lemma \ref{compact},
this also implies that $\mu_n$ restricted to $\Lambda_r$ has a
subsequential limit.  By a diagonalization argument, we can can
choose a subsequence for which such a limit exists for every $r$;
and in this case, and this subsequence converges to a measure
$\mu_B$ on $(\Omega, \mathcal F^{\tau})$, which is easily seen to be
$\mathcal L$-invariant and satisfy $FE_{\Lambda_r}(\mu_B) \leq
FBL_B(\mu) - cr^{-d}$. Since this is true for any $r$, we have
$SFE(\mu_B) \leq FBL_B(\mu)$. \qed
\end{proof}

Now, by Lemma \ref{weaktauequivalent}, the set $M_{FBL(\mu)} = \{
\nu: SFE(\nu) \leq FBL(\mu) \}$ is metrizable in the topology of
weak local convergence.  Hence, we can choose a sequence of sets
$B^m \in \mathcal B$ in such a way that each is contained in the
ball of radius $1/m$ about $\mu$ with respect to this metric.  It is
then clear that $\mu_{B^m} \rightarrow \mu$; hence, by Lemma
\ref{fblsfelemma} and Theorem \ref{levelsetcompactness},
$$SFE(\mu) \leq \liminf_{m \rightarrow \infty} SFE(\mu_{B^m}) \leq \limsup_{m \rightarrow \infty}
FBL_{B^m}(\mu) = FBL(\mu).$$

\subsection{Alternate definition of $\sigma$}
Theorem \ref{limitequality} also implies the following alternative
definition of $\sigma$:

\begin{cor} \label{newsigmadef} If $\Phi$ is an SAP and $E = \mathbb R$ or $E = \mathbb Z$, then
$$\sigma(u) = \liminf_{n \rightarrow \infty} - |T_n|^{-1} \log
Z_{T_n} \geq SFE(\mu),$$ where $Z_{T_n}$ are as defined in Section
\ref{torussection} and $\mu$ is the subsequential limit of the
measures $\mu^n_u$ described in Lemma \ref{torusconvergence}.
\end{cor}

\begin{proof}   Suppose $E = \mathbb R$.  By tiling $T_n$ with large cubes and considering measures
obtained by taking a box measure centered at $\phi_u$ outside of the
cubes and then sampling the interiors according to a Gibbs measure,
it is easy to see that for any $\mu$ with slope $u$,
$$\liminf_{n \rightarrow \infty} - |T_n|^{-1} \log Z_{T_n} \leq
PBL(\mu) = SFE(\mu).$$ Hence,
$$\sigma(u) \geq \liminf_{n \rightarrow \infty} - |T_n|^{-1} \log
Z_{T_n}.$$

For the other direction, it is enough to construct a measure $\mu$
with $$SFE(\mu) \leq \liminf_{n \rightarrow \infty} - |T_n|^{-1}
\log Z_{T_n},$$ and we can choose a subsequential limit $\mu$ of the
torus measures which has this property.  We will discuss the case $E
= \mathbb Z$ in the next section.  \qed \end{proof}

\section{Limit equalities in other settings} \label{limitequalityother}
\subsection{Discrete systems}

In this section, we describe the modifications to the proof of Theorem \ref{limitequality}
necessary for the following discrete analog of the theorem.

\begin{thm} \label{discretelimitequality} If $\Phi$ is an SAP, $E = \mathbb Z$, and $\mu \in \mathcal
P_{\mathcal L}(\Omega, \mathcal F^{\tau})$, then $$SFE(\mu) =
FBL(\mu) = PBL(\mu).$$ \end{thm}

If $E = \mathbb Z$, then the only relevant values of $V_{x,y}:
\mathbb R \mapsto \mathbb R$ are the values assumed at integers.
Thus, we lose no generality in making the convenient assumption that
each $V_{x,y}$ is linearly interpolated between integer values ---
that is, $V_{x,y}$ is linear on $(j,j+1)$ for each $j \in \mathbb
Z$, lower semi-continuous, and convex.  (Thus, $V_{x,y}$ is
continuous and finite on some closed interval $[i,j]$ with $i,j \in
\mathbb Z \cup \{-\infty, \infty\}$, and infinite elsewhere.)

If $\phi: \mathbb Z^d \rightarrow R$ is any continuous
configuration, then we can define a ``randomly rounded'' discrete
configuration by $\phi^{\epsilon} = \lfloor \phi + \epsilon \rfloor$
where $\epsilon$ is chosen uniformly from $[0,1)$.  A key
observation is that the expected value of
$V_{x,y}(\phi^{\epsilon}(y) - \phi^{\epsilon}(x))$ is equal to
$V_{x,y}(\phi(y) - \phi(x))$.

Now, we take $\phi_u$ to be any real-valued function of slope $u$, with $\mathcal L$-periodic
gradient, which has finite $\Phi$ energy.  Now, we can define $FBL(\mu)$ and $PBL(\mu)$ precisely
as in the continuous case except that when defining $PBL(\mu)$, we fix the boundary conditions by
randomly rounding $\phi_u$.  So $PBL^B_{\Lambda_n}(\mu) = FE_{\Lambda_n}(\mu_n)$ where $\mu_n$ is
define as follows: to sample from $\mu_n$, first choose $\phi$ on the boundary via a random
rounding $\phi^{\epsilon}_u$ of $\phi_u$, where $\epsilon$ is chosen uniformly in $[0,1)$; then
choose $\phi$ on the inside according to the Gibbs measure conditioned on the empirical measure
lying in $B$ (i.e., on $\phi$ lying in $B_n$, as defined above).  We define $FBL$ exactly as
before.

Except for this change in setup, all of the arguments and definitions are essentially identical to
the continuous cases.  As before, it is obvious that $FBL(\mu) \leq PBL(\mu)$, and the argument
that $SFE(\mu) \leq FBL(\mu)$ is the same as before.  In the proof that $PBL(\mu) \leq SFE(\mu)$,
however, there is a slightly difference in that we no longer need to define box measures, since
singleton measures themselves have finite free energy, and the ``random truncation'' can be
replaced by a non-random one.  For the alternate definition of $\sigma$, since we cannot
necessarily cause the slope of a configuration on the torus $T_n$ to be exactly $u$, we have to
round $u$ to an integer vector multiple of $1/n$ (as described in Lemma \ref{torusconvergence}).
However, the remainder of the proof is the same as the continuous case.

\subsection{Higher dimensions and PSAPs}
Finally, in this section, we give the most general form of Theorem \ref{limitequality}:

\begin{thm} \label{generallimitequality} If $\Phi$ is a perturbed SAP, $E = \mathbb Z^m$ or $E = \mathbb
R^m$, and $\mu \in \mathcal P_{\mathcal L}(\Omega, \mathcal
F^{\tau})$, then $$SFE(\mu) = FBL(\mu) = PBL(\mu).$$ \end{thm}

First, note that $\Phi$ is an SAP (which we can write $\Phi =
\sum_{i=1}^m \Phi_i$, where each $\Phi_i$ is a one-dimensional SAP),
then the definitions of $PBL$ and $FBL$ are the same as those given
at the beginning of Section \ref{limitequalitysection} except that
$\phi_u = \prod_{i=1}^m \phi_u^i$, $A_n = \prod_{i=1}^m A^i_n$, and
$B_n = \prod_{i=1}^m B^i_n$ where the $A_i^n$ and $B^i_n$ are
subsets of functions from $\Lambda_n$ to $\mathbb R$, defined as
they would be if we were in the one-dimensional setting using
$\Phi_i$ and $u_i$ (the $d$-dimensional slope determined by the
$i$th row of $u$) instead of $\Phi$. If $\Phi + \Psi$ is a perturbed
SAP, where $\Phi$ strictly dominates $\Psi$, then we will also
define $A_n$ and $B_n$ using $\Phi$.

\begin{proof}  In both the $m>1$ and perturbed settings, the arguments for $SFE(\mu) \leq FBL(\mu)$
is exactly the same as before; in each setting, as in the simplest one-dimensional case, it is
obvious from the definitions that $FBL(\mu) \leq PBL(\mu)$.

In the un-perturbed case when $m>1$, the proof that $PBL(\mu) \leq
SFE(\mu)$ for $\mu$ ergodic is the same as in the one-dimensional
case except that we apply a separate randomized truncation in each
of the $m$ coordinate directions. The generalization of this result
to non-ergodic $\mu$ is exactly the same as the proof given for $E =
\mathbb R$.

In the perturbed case where the potential is defined to be $\Phi +
\Psi$ (where $\Phi$ is simply attractive and strictly dominates
$\Psi$, and $\Psi$ has range $k < \infty$), instead of defining
$A_n$ to be the set for which the boundary values $\phi(x) -
\phi(x_0)$ (for $x \in \partial \Lambda_n$) are within $\epsilon$ of
$\phi_u(x) - \phi_u(x_0)$, we take $A_n$ to be the set in which
$|(\phi(x) - \phi(x_0)) - (\phi_u(x) - \phi_u(x_0))| < \epsilon$ for
all $x \in \Lambda_n$ which are within $k$ units of distance from
$\partial \Lambda_n$.  Then we use the one-dimensional argument
(just as before) to produce a random truncation in which the
expected energy added each time we decide $\phi(x)$ for another $x
\in A^+$ is finite; then we observe that since the expected combined
$\Phi$ energy that occurs in edges within $k$ units of points in
$A^+$ is $o(|\Lambda_n|)$, the expected combined $\Psi$ energy is
$o(|\Lambda_n|)$ as well. \qed \end{proof}

\section{Variational principle}\label{variationalprinciplesection}
We have already proved in Lemma \ref{minimizersaregibbs} that if $\mu \in \mathcal P_{\mathcal
L}(\Omega, \mathcal F^{\tau})$ has finite slope $u$ and $SFE(\mu) = \sigma(u)$, then $\mu$ must be
a Gibbs measure.  This, together with the converse, is called the {\it variational principle}.
\begin{thm} \label{variationalprinciple} Let $\Phi$ be a perturbed ISAP (when $E = \mathbb R^m$ or
$\mathbb Z^m$) or a perturbed LSAP (when $E = \mathbb Z^m$).  Then
if $\mu$ is an ergodic gradient measure of finite slope $u$, then
$\mu$ is a Gibbs measure if and only if $SFE(\mu) = \sigma(u)$. If
$\mu$ is a not-necessarily-ergodic gradient Gibbs measure with
ergodic decomposition given by
    $$\mu = \int_{\text{ex}\mathcal P_{\mathcal L}(\Omega, \mathcal F^\tau)} w_{\mu}(\nu)d\nu,$$
    then (letting $S(\nu)$ denote the slope of $\nu$) $\mu$ is a Gibbs measure if and only if
    $$SFE(\mu) = w_{\mu}(\sigma(S*))$$ (using the abbreviation $w_{\mu}(\sigma(S*)) := \int_{\text{ex}\mathcal P_{\mathcal L}(\Omega, \mathcal F^\tau)}
    w_{\mu}\sigma(S(\nu))d\nu$).  \end{thm}
We will first observe here that the second statement follows from
the first; we will then prove the first statement in case that $E =
\mathbb R$ and $\Phi = \Phi_V$ is an ISAP, delaying the more general
discussion until Section \ref{variationalother}.

\begin{proof} The second statement in the lemma follows from the first, using Lemma
\ref{SFEisanexpectation} and the fact that slope is
$e(\text{ex}\mathcal P_{\mathcal L}(\Omega, \mathcal F^\tau))$
measurable. For the first statement, by Lemma
\ref{minimizersaregibbs}, we need only show that if $\mu$ is an
ergodic Gibbs measure of slope $u$, then $SFE(\mu) \leq
\sigma(\mu)$. By Lemma \ref{limitequality} (noting that in the
definition of $PBL$ and any $\mu$ of slope $u$, the limit defining
$PBL^{u}_B(\mu)$ only gets smaller if we replace $A_n \cap B_n$ with
$A_n$---i.e., if we take $B = \Omega$), it is enough to show
$$SFE(\mu) \leq \liminf_{n \rightarrow \infty} -|\Lambda_n|^{-1}\log \int_{A_n}
e^{-H^o_{\Lambda_n}(\phi)} \prod_{x \in \Lambda_n \backslash
\{x_0\}} d \phi(x).$$

By Lemma \ref{levelsetcompactness}, if we can construct a sequence
$\mu_n$ of slope $u$ measures with $$\lim_{n \rightarrow \infty}
SFE(\mu_n) = \liminf_{n \rightarrow \infty} -|\Lambda_n|^{-1}\log
\int_{A_n} e^{-H^o_{\Lambda_n}(\phi)} \prod_{x \in \Lambda_n
\backslash \{x_0\}} d \phi(x),$$ and $\mu_n$ converging weakly to
$\mu$, then it will follow that $SFE(\mu) \leq \sigma(u)$.

We define measures $\mu_n^{\epsilon}$ as follows: to sample from
such a measure, first tile $\mathbb Z^d$ with cubes of size
$\Lambda_n$; on each cube, we will independently choose a
configuration belonging to $A_n$ in the following way. First, sample
$\phi$ from $\mu$ and consider its restriction to $\Lambda_n$; we
condition on the event that $\overline{H}^o_{\Lambda_n}(\phi) \leq
C|\Lambda_n|$ where $C$ is twice the specific $\overline{V}$ energy
of $\mu$ (a value which is finite by Lemma \ref{overlineVSFEbound}
and Lemma \ref{SEbound}). The ergodic theorem implies the
probability of this event tends to one as $n$ tends to $\infty$.
Then, we let $\phi_O$ be the the outside interpolation function
whose existence is guaranteed by Lemma \ref{insidebound} and we
re-sample the values in $\Lambda_n \backslash
\Lambda_{(1-\epsilon)n}$ according to the box measure (whose
existence is guaranteed by Lemma \ref{boxmeasureexistence}) centered
at $\phi_O$. (Throughout this argument, we tacitly assume that
$(1-\epsilon)n$ is rounded down to the nearest integer.)

Now, we would like to show that for any $\delta$, if $\epsilon$ is
small enough, and $n$ large enough, then the specific free energy
will be less $\delta + SFE(\mu)$.  To see this, we will prove the
result for measures $\nu_n^{\epsilon}$ whose specific free energies
are clearly higher than those of the $\mu_n^{\epsilon}$.  To sample
from $\nu_n^{\epsilon}$, first sample form $\mu_n^{\epsilon}$, and
then on each box, let $\phi_I$ be the interpolation between
$\Lambda_{(1-\epsilon)n}$ and $\Lambda_{(1-\epsilon)^2n}$ guaranteed
by Lemma \ref{insidebound}; then re-sample using the values in
$\Lambda_{(1-\epsilon)n} \backslash \Lambda_{(1 - \epsilon)^2n}$
from the box measure whose existence is guaranteed by Lemma
\ref{boxmeasureexistence}.  This measure is equal to the Gibbs
measure on $A_{(1-\epsilon)n}$ in the blocks forming a
$(1-\epsilon)^2$ fraction of $\mathbb Z^d$, and outside, equal to a
box measure with specific free energy bounded by a constant (again,
by Lemma \ref{boxmeasureexistence}).  Thus, for a sufficiently small
$\epsilon$, we have that for sufficiently large $n$,
$SFE(\nu^{\epsilon}_n) \leq SFE(\mu) + \delta$.  It is now trivial
to check that $\mu$ is a limit point of the measures
$\mu^{\epsilon}_n$. \qed \end{proof}

\section{Variational principle in other settings} \label{variationalother}
\subsection{Discrete models}

To prove Theorem \ref{variationalprinciple} in the setting $E =
\mathbb Z$, we take $\overline{V} = V$ (as before, assuming $V$ is
linear between integers and lower semi-continuous) and construct
continuous interpolations exactly as in the $E = \mathbb R$ setting.
Our definition of $\mu_n^{\epsilon}$ and $\nu_n^{\epsilon}$ is
essentially the same as in the previous case.  We tile by
$\Lambda_n$ blocks just as before, and choose the interpolations
just as before; the only difference is that instead of using a box
measure, we use a random rounding: i.e., we add a variable uniformly
distributed in $[0,1]$ to the entire choice of $\phi$ and then round
down (in the set $R$ consisting of points outside of the
$(1-\epsilon)^2$ boxes), after which we sample the remainder of
$\phi$ inside each of the boxes according to the appropriate Gibbs
measure (with boundary conditions given by the values of $\phi$ in
$R$).

\subsection{Lipschitz simply attractive potentials}

In the Lipschitz models---when $\Phi$ is not necessarily an ISAP,
and we may not have a good definition of $\overline{V}$---we define
$\mu_n^{\epsilon}$ and $\nu_n^{\epsilon}$ slightly differently: in
this case, we use a truncation of the kind used in the proof that
$PBL(\mu) \leq SFE(\mu)$ in Section \ref{ergodicPBLSFEsection}.
Because $\mu$ is Lipschitz, Lemma \ref{ergodicL1} actually implies
that as $n$ tends to infinity, the probability that $\phi$ sampled
from $\mu$ differs from the plane of slope $S(\mu)$ (with
appropriate additive constant) by more than $\epsilon n$ in the
supremum norm tends to zero in $n$.  Thus, as $n$ gets large, the
probability tends to one that truncations of the type in used in
Lemma \ref{ergodicPBLSFEsection} will not affect the value of $\phi$
anywhere inside the box $\Lambda_{(1-\epsilon)n}$. These truncations
play the role of the interpolations used in the non-Lipschitz
setting.  It is not hard to define an ``inside truncation'' in a
similar fashion, to play the role of the inside interpolations.

\subsection{Higher dimensions and perturbed ISAPs and LSAPs}
The change to higher dimensions (without perturbations) is essentially trivial; we simply define
the interpolations separately on each coordinate as before.  The change to perturbed potentials
$\Phi + \Psi$ is also essentially trivial; we simply observe that if the expected $\Phi$ energy in
$\Lambda_n \backslash \Lambda_{(1-\epsilon)^2n}$ is $o(|\Lambda_n|)$, the expected $\Psi$ energy is
$o(|\Lambda_n|)$ as well.

\chapter{LDP for empirical measure profiles} \label{LDPchapter}
\section{Empirical measure profiles and statement of LDP}
Let $D$ be a bounded domain, and the sequence $D_n$ a good
approximation of $D$.  The aim of this chapter is to prove a large
deviations principle (of speed $n^d$) for the behavior of a random
gradient configuration $\phi_n$ on $D_n$. Instead of merely
considering the gradient empirical measures---as in the previous
section---we will consider {\it profiles}, which also contain
information about how the occurrence of various gradient events are
distributed throughout $D$. We will also investigate behavior of the
normalized functions $\frac{1}{n}\phi_n(nx)$ (interpolated to
functions on $D$) and examine their large deviations behavior with
respect to topologies induced by $L^p$ and Orlicz $L^A$ metrics. We
can define an {\it empirical profile measure} $R_{\phi_n, n} \in
\mathcal P(D \times \Omega)$ by $$R_{\phi_n, n} = \int_{D}
\delta_{x, \mathcal \theta_{\lfloor nx \rfloor} \phi_n}dx. $$
Informally, to sample a point $(x, a)$ we choose $x$ uniformly from
$D$, and then take $a = \theta_{\lfloor nx \rfloor} \phi_n$.  (As in
the previous chapter, it is convenient to write $\phi_n(x) = \infty$
for $x \not \in D_n$.)  Also, using $\phi_n$, we will define a
function $\tilde \phi_n$ by interpolating the function
$\frac{1}{n}\phi_n(nx)$ to a continuous, piecewise linear (on
simplices) function on the simplex domain corresponding to $D_n$;
each such $\tilde \phi_n$ is a member of the space
$L^{\overline{V}^*}_0(D)$ constructed in Section
\ref{goodapproximationsection} for Corollary
\ref{discretegoodcompact}, where we take $\overline{V}^*$ to be any
function increasing essentially more slowly that $\overline{V}_d$,
where $\overline{V}$ is defined as in Section \ref{overlineV}, and
$\overline{V}_d$ is the Sobolev conjugate as defined in Section
\ref{sobconjsection}.

Write $\mu_n = Z_n^{-1} e^{-H^o_{D_n}(\phi)} \prod_{x \in D_n
\backslash \{x_0 \}} d\phi(x)$ where the $Z_n$ are normalizing
constants chosen to make $\mu_n$ probability measures.  Let $\rho_n$
be the measure on $X = \mathcal P(D \times \Omega) \times
L^{\overline{V}^*}_0(D)$ induced by $\mu_n$ and the map $\phi_n
\mapsto (R_{\phi_n,n}, \tilde \phi_n)$. We say a measure $\mu \in
\mathcal P(D \times \Omega)$ is {\it $\mathcal L$-invariant} if $\mu
(\cdot, \Omega)$ is a Lebesgue measure on $D$ and for any $D'
\subset D$ of positive Lebesgue measure, $\mu (D', \cdot)$ is an
$\mathcal L$-invariant measure on $\Omega$. Given any subset $D'$ of
$D$ with positive Lebesgue measure, we can write $S(\mu(D',\cdot))$
for the slope of the measure $\mu(D', \cdot)/\mu(D' \times \Omega)$
(we have normalized to make this a probability measure) times
$\mu(D' \times \Omega)$; the map $D' \mapsto S(\mu(D',\cdot))$ is a
signed/vector-valued measure on $D$.   Let $\mathcal X$ be the
topology on $\mathcal P(D \times \Omega) \times
L^{\overline{V}^*}(D)$ which is the product of (on the first
coordinate) the smallest topology in which $\mu\rightarrow \mu(D'
\times f)$ is measurable for all rectangular subsets $D'$ of $D$ and
bounded cylinder functions $f$ and (on the second coordinate) the
$L^{\overline{V}^*}_0(D)$ topology.  We say an ISAP $\Phi = \Phi_V$
is {\it super-linear} if $\lim_{\eta \rightarrow \infty}
V(\eta)/\eta = \infty$; although many of our bounds still hold for
ISAPs which are not super-linear, the full LDP, which we prove
below, is false as stated when $V$ is not super-linear.

\begin{thm} \label{ldpprofiletheorem} If $\Phi$ is a super-linear ISAP.  The measures $\rho_n$
satisfy a large deviations principle with speed $n^d$ and rate function
    $$I(\mu,f) = \begin{cases}
        SFE (\mu(D, \cdot)) - P(\Phi)  &
        \text{$\mu$ is $\mathcal L$-invariant and $S(\mu(x,\cdot)) = \nabla f(x)$}\\
        & \text{\em as a distribution}\\
        \infty & \text{otherwise} \\
        \end{cases}$$ on the space $\mathcal P(D \times \Omega) \times L^{\overline{V}^*}(D)$ in the
         topology $\mathcal X$ described above.
\end{thm}

In the case that $\Phi$ is an SAP and $E = \mathbb R$, the
uniqueness of the minimizer of the rate function $I$ described above
is an immediate consequence of the uniqueness of the gradient Gibbs
measure of a given slope (Theorem \ref{minimalphaseuniqueness}) and
the strict convexity of $\sigma$ (Theorem \ref{sigmaconvex}---which
in particular implies that $\sigma$ has a unique minimum).  This
will also imply uniqueness of the rate function in the presence of
boundary conditions (see Section \ref{boundarycondsect}).

\section{Proof of LDP}
Recall that, in general, a sequence of measures $\rho_n$ on a
topological space $(X, \mathcal X)$ is said to {\it satisfy a large
deviations principle with rate function $I$ and speed $n^d$} if $I:
X \rightarrow [0, \infty]$ is lower-semicontinuous and for all sets
$B \in \mathcal X$, $$-\inf_{x \in B^o} I(x) \leq \liminf_{n
\rightarrow \infty} n^{-d} \log \rho_n(B) \leq \limsup_{n\rightarrow
\infty} n^{-d} \log \rho_n(B) \leq -\inf_{x \in \overline{B}}
I(x).$$ Let $\mathcal A$ be any basis of $\mathcal X$. It is not
hard to check (and is proved in \cite{DZ}, Theorem 4.1.11 and Lemma
1.2.18) that the large deviations principle is a consequence of the
following three statements:
\begin{enumerate}
\item \bf{Lower bound on probabilities:} $$\inf_{S \in \mathcal A; x \in S} \limsup
 \frac{1}{n^d} \text{log} \rho_n(S) \geq -I(x).$$
\item {\bf Upper bound on probabilities:} $$\inf_{S \in \mathcal A; x \in S} \liminf
 \frac{1}{n^d} \text{log} \rho_n(S) \leq -I(x).$$
\item {\bf Exponential tightness:} For every $\alpha < \infty$, there exists a compact set $K_{\alpha} \subset X$
for which $\limsup_{n \rightarrow \infty} n^{-d} \log \rho_n (X \backslash K_{\alpha}) < - \alpha$.
\end{enumerate}  The above conditions also imply that all the level sets $M_{\alpha} = \{x : I(x)
\leq \alpha \}$ are compact (by Lemma 1.2.18 of \cite{DZ}).  When
the first and second statement hold, we say that the $\rho_n$
satisfy a {\it weak large deviations principle with rate function
$I$}. When the third statement holds, we say that the $\rho_n$ are
{\it exponentially tight}. In the following subsections, we will
prove Theorem \ref{ldpprofiletheorem} by checking each of these
statements in turn.  We will do this first for the simplest
case---when $\Phi$ is an ISAP and $E = \mathbb R$---and address
generalizations in Section \ref{LDPothersettings}.

\subsection{Lower bounds on probabilities} \label{lowerboundonprob}
The lower bound follows almost immediately from Lemma \ref{mostlylinearapproximator} and Lemma
\ref{discretefunctionapproximator} and Theorem \ref{limitequality} ---in particular, the fact that
$SFE(\mu) = PBL(\mu)$.

Suppose we are given $I(f,\mu) < \infty$.  This implies that $f \in L^{1,A}(D)$, and hence, by
Lemma \ref{mostlylinearapproximator}, for all $n$ sufficiently large, there exists a function
$F_{\epsilon} \in L^{1,A}(D)$ that is equal to a piecewise linear approximator $F_n$ (as defined
prior to Lemma \ref{mostlylinearapproximator}) on a closed subset $D'$ of $D$, where the volume of
$D \backslash D'$ is less than $\epsilon$ and $\int_{D-D'} A(\nabla F_{\epsilon}) \leq \epsilon$.

Now, a basis set $S \subset A$ centered at $f, \mu$ can be written
as the set of pairs $(g,\nu)$ for which $\delta(f,g) < \gamma$ (for
some $\gamma > 0$, where $\delta$ is the distance function described
in Section \ref{goodapproximationsection}---roughly, $\delta(f,g) =
||f-g||_{\overline{V}^*}$) and $|\mu(1_{D_i} H_i) - \nu(1_{D_i}
H_i)| < \gamma$ for each of a finite set of cylinder functions $H_i$
(each bounded between $0$ and $1$) on $\Omega$ and rectangular
subsets $D_i$ of $D$.

Assume without loss of generality (rescaling if necessary) that the volume of $D$ is less than one.

Now, given any $F_{\epsilon}$, and a large $n$, we define a measure
$\mu_{\epsilon, n}$ on configurations on $D_n$ as follows: to sample
$\phi$ from $\mu_{\epsilon,n}$, first compute the approximation
$F_{\epsilon}^n:D_n \rightarrow \mathbb R$ of $F_{\epsilon}$
guaranteed by Lemma \ref{discretefunctionapproximator} (with $A =
\overline{V}$) and fix $\phi(x_0) = F_{\epsilon}^n(x_0)$ for some
reference vertex $v_0$.  Then sample $\phi(x)$ from the box measure
centered at $\mu_{\epsilon,n}$ (the one whose existence is given by
Lemma \ref{vbarfebound}) for all values of $x \in D_n \backslash
\{x_0\}$ for which $x$ does not lie on the interior of one of the
linear regions of $F_{\epsilon}$. The free energy of this process
is, by Lemma \ref{vbarfebound} and Lemma
\ref{discretefunctionapproximator} and the assumed bound on
$\int_{D-D'} A(\nabla F_{\epsilon})$, at most a constant times
$\epsilon n^d$.

\begin{figure} \begin{center} \leavevmode \epsfbox[0 0 225 160]{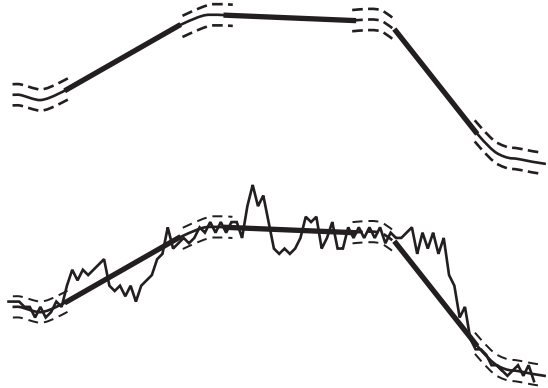} \end{center}
\caption[Random surface approximating a ``mostly linear function''] {Above, a mostly linear
function $F_{\epsilon}$, with dotted lines among the non-linear parts describing the bounds of a
corresponding box measure), in the trivial case $d=1$. Below, an instance of $\phi$ sampled from
$\mu_{n,\epsilon}$: in non-linear parts, $\phi$ is chosen from a box measure; in the linear parts,
$\phi$ is chosen with empirical measure constraints.} \label{mostlylinear} \end{figure}

Then, inside each large box $\Delta_n$ of $D_n$---assume it has size $\Lambda_k$---approximating a
box $\Delta$ on which $F_{\epsilon}$ is linear, we sample $\phi$ from (appropriately translated)
the Gibbs measure on $A_n \cap G_n$ where $G$ is the set of measures $\nu$ for which $|\nu(H_i)-
 \frac{\mu(H_i 1_{D_i})}{|D_i|}| < \gamma/2$ for each $i$.  See Figure \ref{mostlylinear}.

If $\epsilon < \gamma/2$ and $n$ is sufficiently large, then it is
not hard to see that $\mu_{\epsilon, n}$ is supported inside the set
$S$.  Since the free energy from the box measure choice in the
non-linear sections of $\mu_{\epsilon, n}$ is $O(\epsilon n^d)$, and
the free energy from the remaining choice is at most
$SFE(\mu(\overline{\Delta} \times \cdot)/\mu(\overline{\Delta}
\times \Omega)$ times $n^d |\overline{\Delta}|$, where
$\overline{\Delta}$ is the union of the square regions of $\mathbb
R^d$ on which $F_{\epsilon}$ is linear.  If follows from Theorem
\ref{limitequality}---since $SFE(\mu) = PBL(\mu)$---that $\limsup_{n
\rightarrow \infty} \frac{1}{n^d} FE(\nu_{\epsilon,n}) \leq I(f,\mu)
+ o(\epsilon)$, which implies the desired lower bound on
probabilities.

\subsection{Upper bounds on probabilities when $I(\mu,f) < \infty$}
The upper bound follows almost immediately from Lemma \ref{limitequality}---in particular, the fact
that $SFE(\mu) = FBL(\mu)$.  Since $I(\mu,f) < \infty$, $\mu$ is an $\mathcal L$-invariant gradient
measure. Now, we may partition $D$ into disjoint cubes $K_1, \ldots, K_k$ of equal size that cover
at least a $1-\epsilon$ fraction of the Lebesgue volume of $D$, and let $K = \cup_{i=1}^k K_i$. In
particular, for any basis set $S$, for a fine enough partition, we will have
$\mu(K,\cdot)/\mu(K,\Omega)$; the $\liminf$ of the average specific free energy within $K$---as we
choose $S$ small enough so that $\nu(K,\cdot)$ lies in a sufficiently small neighborhood of $\mu(K,
\cdot)$---is at least $SFE(\mu(K,\cdot))$.  As $\epsilon$ tends to zero, $\mu(K,\cdot)$ converges
to $\mu(D,\cdot)$ in the topology of local convergence.  The average specific free energy outside
of $K$ is at least $\alpha$ (as defined in Corollary \ref{alphabound}).  This gives the desired
bound.

\subsection{Upper bounds on probabilities when $I(\mu, f)=\infty$}  \label{upperinfty}
If $\mu$ is $\mathcal L$-invariant and $S(\mu(x,\cdot)) = \nabla f(x)$ as a distribution and
nonetheless $I(\mu)=\infty$ (because $SFE(\mu(D, \cdot) = \infty)$), then the argument is the same
as in the previous section. We have now to show the upper bound on probabilities when either
$S(\mu(x,\cdot)) \not = \nabla f(x)$ as a distribution or $\mu$ is not $\mathcal L$-invariant.

First, if $\mu$ is not $\mathcal L$-invariant, then there is an event $H$ and a rectangle $D'$ for
which $\mu(D' \times H) > \mu(D \times \theta_x H) + \delta$ for some $\delta > 0$.  We choose a
neighborhood $S$ to be such that for each $(g,\nu)$ in this neighborhood, $\nu(D' \times H) > \nu(D
\times \theta_x H) + \delta/2$.  Now, it is not hard to see that the probability of belonging to
this neighborhood becomes zero when $n$ is large enough.

Second, suppose that there exists a rectangle $D'$ for which $S(\mu(x, D'))$ is not equal to the
mean value of $\nabla f(x)$ on $D'$.  Suppose these two values disagree in the $i$th component
direction---and suppose that $D'$ has length $L$ in the $i$th direction and normal cross-sectional
area given by $\alpha$.  Then, integrating in the $i$th direction, this implies that the difference
$\delta$ of the mean value of $f$ on opposite sides of $D'$, divided by $L$, is different from
$\mu(D', \phi(e_i) - \phi(0))$---different by, say, $\gamma > 0$.  Now, let $$H_C(\phi) =
\begin{cases}
\phi(e_i) - \phi(0) & |\phi(e_i) - \phi(0)| \leq C \\
-C & \phi(e_i) - \phi(0) < -C \\
C & \phi(e_i) - \phi(0) > C \\
\end{cases} $$ By letting $C$ get large and choosing a neighborhood that ensures that the average
value of $H_C$ on $D'$ tends to the $S(\mu(x, D'))$; by continuity
of the average cross-sectional area, in $f$, if we take $S$ also to
include a sufficiently tight restriction on $g$, then we can force
the sum of $\phi(x+e_i) - \phi(x)$ over an arbitrarily small
fraction of the points $x$ to grow in $n$ like $\gamma n^d$.  Since
$V$ is super-linear, taking this fraction arbitrarily small implies,
as desired, that $$\inf_{S \in \mathcal A; x \in S} \liminf
\frac{1}{n^d} \text{log} \rho_n(S) \leq -\infty.$$

\subsection{Exponential tightness} We define the set $K_C \subset \mathcal P(D \times \Omega) \times L^{\overline{V}^*}_0(D)$ to
be the set of of profile/interpolation pairs corresponding to
functions from $D_n$ to $\mathbb R$ with average $\overline{V}$
energy per edge equal to or less than $C$.  Exponential tightness
will follow once we show that \begin{enumerate}
\item For each $C$, the set $K_C$ is pre-compact in $(X, \mathcal X)$.
\item For any $\alpha > 0$, if we choose $C$ large enough, we will have
$$\liminf_{n \rightarrow \infty} n^{-d} \log \mu_n(K_C^c) \leq - \alpha.$$ \end{enumerate}

For the first statement, it is enough to show that the projection of
$K_C$ onto each of its two components---in $\mathcal P(D \times
\Omega)$ and in $L^{\overline{V}^*}_0(D)$---is pre-compact in the
corresponding topology.  The first is a simple exercise; the second
follows from Lemma \ref{discretegoodcompact}.

For the second statement, by Lemma \ref{vbarfebound}, it is enough
to prove an analogous statement using $V$ instead of $\overline{V}$.
Suppose it were the case that for some $\alpha$, $$\liminf \log
\mu_n (1_{n^{-d}|H^o(\phi)| \geq C}) \geq -\alpha$$ for every $C$.
Then this would imply that if we replaced $\Phi$ with
$(1-\epsilon)\Phi$, for any small $\epsilon>0$, then we would have
$$\liminf n^{-d} \log \int e^{-H^o(\phi)} \prod_{x \in D_n \backslash
\{x_0\}} = \infty,$$ i.e., the log partition function growing
super-exponentially in $n$. This is a contradiction to Corollary
\ref{alphabound}.

\section{LDP in other settings} \label{LDPothersettings}
In this section, we briefly describe some variants Theorem
\ref{ldpprofiletheorem} (to $E = \mathbb Z$ or $m >1$ or $\Phi$ not
an ISAP) and describe the modifications to the proof required for
these settings.  However, we will not reproduce in detail the proof
of Theorem \ref{ldpprofiletheorem} in each of these settings, as
this would consume a good deal of space and the modifications are
all straightforward.

\subsection{Discrete models and Lipschitz models}
All of the arguments in this chapter carry through to higher
dimensional spin spaces, perturbed systems, and discrete models
(using similar rounding arguments to those described in the previous
chapter) with little or no modification.  However, some additional
care is required in the case that $E = \mathbb Z$ but $\Phi$ is
Lipschitz, so that $V$ fails to be everywhere finite.  In this case,
$\sigma$ does not approach infinity near the boundary of the space
of allowable slopes, so it is possible that the rate function
$I(f,\mu)$ will be non-infinite even if $\nabla f_{\mu}(x)$ lies on
this boundary for $x$ in a subset of $D$ with positive measure.  If
this is the case, we say that $f$ is a {\it taut} height function.
Since the variational principle and the uniqueness of the gradient
phase of slope $u$ (as shown in Chapter
\ref{clusterswappingchapter}) may not apply when $u$ is one of these
boundary slopes, we cannot expect the minimizer of the rate function
to be unique in this case. However, the large deviations principle
does go through. A simple analytical argument shows that we can
always approximate a taut $(f, \mu)$ by a sequence of pairs
$(f_i,\mu_i)$ for which the $f_i$ are not taut; this enables us to
deduce the necessary lower bound, and the upper bounds and
exponential tightness arguments are the same as before.

\subsection{LDP with boundary conditions} \label{boundarycondsect}

In the continuous setting, we sometimes wish to limit our attention
to functions $f$ that extend to the closure of $D$ and satisfy $f =
f_0$ on $\partial D$, where $f_0 :
\partial D \rightarrow \mathbb R^m$ is a continuous {\it boundary condition}.
Of course, elements of $L^{\overline{V}^*}(D)$ are not continuous
for general $\overline{V}^*$, and in particular need not be
continuous at the boundary of $D$. But we will say that a function
$f$ on $D$ {\em has $f_0$ as its boundary} if $f$ is a limit in
$L^{\overline{V}^*}(D)$ of functions in $L^{\overline{V}^*}(D)$,
each of which agrees with $f_0$ outside of a compact subset of
$\partial D$.

We would like to impose similar conditions on the discrete models.
However, since none of the elements of the boundaries of the $D_n$
actually lies on the boundary of $\partial D$, we cannot simply
require that $\phi(x) = f_0(x)$ for $x \in
\partial D_n$.  In fact, there are many ways to specify discrete
boundary conditions; we will choose the one that is most convenient
for us.

Assume that $f_0$ extends continuously to a function in
$L^{1,A}(D)$; then we define as in the previous section the
functions $f^n_0$, approximating $f_0$ on $D_n$, as in Lemma
\ref{discretefunctionapproximator}, and box measures
$\nu_{n,\epsilon}$ centered at these functions.   Now, take $\mu_n$
to be the sequence of measures $\mu_n = Z_n^{-1}
e^{-H^o_{D_n}(\phi)} \prod_{x \in D_n \backslash \{x_0 \}}
1_{\chi_n} d\phi(x)$ where the $Z_n$ are normalizing constants
chosen to make $\mu_n$ probability measures and each $\chi_n$ is the
set of $\phi$ for which $|\phi(x) - f^n_0(x)| \leq \epsilon$ for
each $x \in \partial D_n$.  These induce measures $\rho_n$, defined
as above.

\begin{thm} \label{ldpprofiletheoremwithboundary} If $\Phi$ is a super-linear ISAP.  The measures
$\rho_n$ satisfy a large deviations principle with speed $n^d$ and
rate function (up to an additive constant) given by
    $$I(\mu,f) = \begin{cases}
        SFE (\mu(D, \cdot)) - P(\Phi) &
        \text{$\mu$ is $\mathcal L$-invariant and $S(\mu(x,\cdot)) = \nabla f(x)$}\\
        & \text{ as a distribution and $f$ has $f_0$ as}\\
        & \text{ its boundary} \\
        \infty & \text{otherwise} \\
        \end{cases}$$ on the space $\mathcal P(D \times \Omega) \times L^{\overline{V}^*}(D)$ in the
         topology $\mathcal X$ described above.
\end{thm}

\begin{proof} We assume $m=1$ (the extension to $m>1$ being straightforward).
Exponential tightness and the upper bounds on probabilities in the
case that $f$ has $f_0$ as its boundary are the same as in the
no-boundary case. If $f$ agrees with $f_0$ outside of a compact
subset of $f$, then lower bound argument is also exactly the same as
before (noting that the approximation of $f$ defined in Lemma
\ref{mostlylinearapproximator} agrees with $f$---and hence
$f_0$---on the boundary of $D_n$ for all large enough $n$); since,
by definition, any $f$ that has $f_0$ as its boundary is a limit in
$L^{\overline{V}^*}(D)$ of functions with this property, this gives
the lower bound in general.

It remains only to check the upper bound on probabilities in the
case that $f$ does not have $f_0$ as its boundary (and hence $I(\mu,
f) = \infty$).  We know by compactness that the probability that the
probability that $\phi$ fails to have a subsequential limit in
$L^{\overline{V}^*}(D)$ tends to zero super-exponentially.  However,
in light of the discrete boundary conditions, it is not hard to see
that any subsequential limit in $L^{\overline{V}^*}(D)$ of $\phi$
chosen from $\rho_n$ must have $f_0$ as its boundary. \qed
\end{proof}

The generalization to LSAPs with $E = \mathbb Z$ is straightforward
when $f_0$ is not taut, since in this case, any $f$ extending $f_0$
can be approximated by functions which are not taut.  However, if
$f_0$ is taut---for example, if it is a plane of slope $u$---then
the empirical measure large deviations principle need not hold.  If
there are distinct ergodic Gibbs measures $\mu_1$ and $\mu_2$ of
slope $u$ that have different specific free entropies, and $\phi_1$
and $\phi_2$ are samples from $\mu_1$ and $\mu_2$, then the large
deviations behavior of the sequence of Gibbs measures with boundary
conditions given by $\phi_1$ outside of $\Lambda_n$ will be
different from the one with boundary conditions given by $\phi_2$
outside of $\Lambda_n$.



\subsection{Gravity and other external fields} \label{externalfields}
Let $h:X \rightarrow \mathbb R$ be any continuous function on $X$;
we would now like to replace the measure $\rho_n$ with $\pi_n =
e^{-h}\rho_n$ (times the normalizing constant that makes $\pi_n$ is
a probability measure). These new measures $\pi_n$ clearly satisfy
upper and lower bounds on probabilities described above when $I(x)$
is replaced by $I_0(x)= I(x) + h(x)$. When $m=1$, typical example of
a continuous function $h$ might be (in the presence of boundary
conditions) $h(f,\mu) = \int_D f(\eta) d\eta$--- this corresponds to
weighting a configuration with additional energy proportional to the
``gravitational potential'' energy of the surface (causing typical
surfaces to sag lower in the interior of $D$). Another example of
such an $h$ might be $h(f,\mu) = \mu(H)$ for some function $H:\Omega
\rightarrow \{0,1\}$; this corresponds to weighting by the number of
times a particular local configuration appears.

We would like to argue that for some constant $C$, the $\pi_n$
satisfy a large deviations principle with rate function $I(x) + h(x)
+ C$.  But this follows immediately from Varadahn's Integral Lemma
(Theorem 4.3.1 of \cite{DZ}) provided that a certain tightness
conditions holds.  Namely, we require
$$\lim_{M \rightarrow \infty} \limsup n^d \log E\[ e^{n^d h(x)} 1_{h(x) \geq M } \] = - \infty.$$ This
holds for the gravitational potential described above (when $\sigma$
is super-linear) and many other kinds of external fields.

\chapter{Cluster swapping} \label{clusterswappingchapter} {\it Cluster swapping} is a simple
geometric operation that we will use to prove strict convexity of
the surface tension function $\sigma$ and to classify $\mathcal
L$-ergodic gradient Gibbs measures and their extremal decompositions
(which may be non-trivial if $E=\mathbb Z$) whenever $\Phi$ is an
$\mathcal L$-invariant simply attractive potential. Throughout this
chapter and Chapter \ref{discretegibbschapter}, we will assume that
$m=1$, so that either $E = \mathbb Z$ or $E = \mathbb R$, and that
$\Phi$ is an $\mathcal L$-invariant SAP.

The prerequisites for this chapter are Chapters \ref{SFEchapter},
\ref{decompositionchapter}, and \ref{surfacetensionchapter}. The
only results from Chapters \ref{orliczsobolevchapter},
\ref{LDPempiricalmeasurechapter}, and \ref{LDPchapter} that we will
even mention in either this chapter or Chapter
\ref{discretegibbschapter} are Corollary \ref{newsigmadef} (which
gives the alternate definition of $\sigma$ using limits of log
partition functions on tori) and the second half of the variational
principle, Theorem \ref{variationalprinciple}; and the latter we
mention only in the following remark.

Recall that the variational principle has two parts: the first,
which we will use frequently, is Theorem \ref{minimizersaregibbs},
which states that whenever $SFE(\mu) = \sigma(S(\mu))$, the measure
$\mu$ is a gradient Gibbs measure.  This result holds for all simply
attractive potentials. Recall also that a {\it minimal $\mathcal
L$-ergodic gradient phase} of slope $u$ is defined to be a slope-$u$
$\mathcal L$-ergodic gradient Gibbs measure $\mu$ on $(\Omega,
\mathcal F^{\tau})$ with minimal specific free energy; i.e.,
$SFE(\mu) = \sigma(u)$. In this chapter, we will classify the
minimal $\mathcal L$-ergodic gradient phases of slope $u$ for any $u
\in U_{\Phi}$.

The second half of the variational principle, Theorem \ref{variationalprinciple}, states that {\it
every} $\mathcal L$-ergodic gradient phase $\mu$ of slope $u \in U_{\Phi}$ is in fact a {\it
minimal} $\mathcal L$-ergodic gradient phase.  We have proved the second half only for perturbed
isotropic and (discrete) Lipschitz simply attractive potentials---it is not known whether Theorem
\ref{variationalprinciple} can be extended to all simply attractive potentials.  In the cases where
Theorem \ref{variationalprinciple} is true, our classification of the minimal $\mathcal L$-ergodic
gradient phases of slope $u$ may be interpreted as a classification of {\it all} $\mathcal
L$-ergodic gradient phases of slope $u$.

\section{Introduction to cluster swapping} \label{clusterintrosection}

\subsection{Review: Fortuin-Kasteleyn; Swendsen-Wang;\\ Edwards-Sokal updates}

Before describing cluster swapping, we review some facts about the
related Swendsen-Wang algorithm, introduced in 1987 \cite{SW} and
generalized to the form we present below by Edwards and Sokal in
1989 \cite{ES}. For this subsection only, we let $(E, \mathcal E)$
be a finite set endowed with counting measure $\lambda$ (for
example, $E$ could be $\{-1,1\}$, as in the Ising model), $\Lambda$
any finite graph with a subset $\partial \Lambda$ of its vertices
designated as boundary vertices, and $\Phi = \{\Phi_\Lambda \}$ any
Gibbs potential on functions $\xi:\Lambda \rightarrow E$.  The
crucial idea is the introduction of an independent random auxiliary
function called the {\em residual energy}.

Let $(\xi, r)$ be a random pair in which $\xi:\Lambda \rightarrow E$
is sampled from the Gibbs measure $e^{-H^{\Phi}_\Lambda} \prod_{x
\in \Lambda \backslash
\partial \Lambda} d \xi(x)$ (times a normalizing constant), with
boundary values $\xi = \xi_0$ fixed on $\partial \Lambda$, and $r$
is an independent function from the subsets of $\Lambda$ to
$[0,\infty)$ where the values $r(\Delta)$ are all independent
exponentials with parameter $1$, i.e., distributed according to the
measure $e^{-x}dx$ on $[0, \infty)$.  We refer to the quantity
$r(\Delta)$ as the {\em residual energy} in $\Delta$,
$\Phi_\Delta(\xi)$ as the {\em potential energy} in $\Delta$, and
$t(\Delta) := r(\Delta) + \Phi_\Delta(\xi)$ as the {\em total
energy} in $\Delta$.  Note that the probability density of the pair
$(\xi, r)$ with respect to the natural underlying measure (i.e.,
$\lambda^{|\Lambda \backslash
\partial \Lambda|}$ times the product---over all $\Delta \subset
\Lambda$---of Lebesgue measure on $[0,\infty)$) is proportional to
$e^{-|t|}$ where $|t| = \sum_{\Delta \subset \Lambda} t(\Delta)$.

A general version of the random {\em Swendsen-Wang update} \cite{SW}
(as described by Edwards and Sokal \cite{ES}) to the pair $(\xi, r)$
is the following: first re-sample all of the residual energies
$r(\Delta)$ for $\Delta \subset \Lambda$ from the marginal law of
$r$ (i.e., $e^{-|r|}\prod_{\Delta \subset \Lambda} dr(\Delta)$,
where each $dr(\Delta)$ is Lebesgue measure on $[0,\infty)$). Then
re-sample the pair $(\xi, r)$ {\em conditioned} on the total energy
function $t$.

If the latter step happens to be computationally easy (as it turns
out to be for Ising and Potts models and spin glasses, see below),
then one can often efficiently generate (approximately) random
samples from the Gibbs measure by beginning with a deterministic
pair $(\xi, r)$ and applying the Swendsen-Wang update repeatedly
\cite{ES, SW}. This method of Monte Carlo sampling, called the {\em
Swendsen-Wang algorithm}, is the subject of a large literature. (See
also \cite{H} for an exact sampling analog of Swendsen-Wang.)
Although we will not discuss sampling problems in this paper, we
will use related random updates to generate couplings and to prove
other results.

\begin{flushleft} {\bf Remark:} The Edwards-Sokal formulation is cosmetically different
from ours.  They first add additive constants to the functions
$\Phi_\Delta$ if necessary so that they are all non-negative and
replace what we call ``total energy'' $t(\Delta)$ with the quantity
$t'(\Delta) = e^{-t(\Delta)}$.  They then study the joint law of
$(\xi, t')$ instead of $(\xi, t)$ or $(\xi, r)$.  We will use the
fact that conditioned on $\xi$, the law of $t(\Delta)$ is
$\Phi_\Delta(\xi)$ plus an independent parameter one exponential.
Edwards and Sokal use the fact that conditioned on $\xi$, the law of
$t'(\Delta)$ is uniform on $[0, e^{-t'(\Delta)}]$. Edwards and Sokal
also do not interpret $t'$ as an ``energy.'' For our purposes, it is
more natural to deal with $r$ or $t$ and interpret them as an
energies, given the role they will later play in variational
principles, $SFE$-preserving updates on infinite systems, etc.
\end{flushleft}

From here on, we will specialize to the case that $\Phi$ is a
nearest-neighbor pair potential.  In this case, the values of $r$ on
sets $\Delta$ that are not endpoint pairs of edges of $\Lambda$ are
irrelevant to the way $\xi$ is updated in the Swendsen-Wang
algorithm, since the total energy in such a set is $r(\Delta)$
independently of $\xi$.  From here on, we will ignore these
$r(\Delta)$, and think of $r$ as a function on the edges of
$\Lambda$ only.

The simplest and most well-studied example---and the one that will
turn out to be most relevant to random surfaces---is the Ising model
with coupling constants $K_e \in \mathbb R$ on the edges $e$ of
$\Lambda$, i.e, $E = \{-1,1\}$ and $\Phi_{(x,y)}(\xi) =
K_e\xi(x)\xi(y)$ for each edge $e=(x,y)$. Then the potential energy
$K_e\xi(x)\xi(y)$ of $e$ takes on only the two values $\pm K_e$, and
the total energy in $e$ is given by $t(e) = K_e\xi(x)\xi(y) + r(e)$.
The expression $r(e) = t(e) - K_e\xi(x)\xi(y)$ is always
non-negative.  If we fix $t(e)$ and $t(e) < |K_e|$, then this can
only be the case if $K_e\xi(x)\xi(y)$ is negative
--- i.e., the edge energy is the lower of its two possible energies
(in which case the edge is said to be {\em satisfied} by $\xi$). If
$t(e) \geq |K(e)|$ then this will be the case for both possible
values of the potential energy. Let $\mathcal S$ be the set of edges
of $\Lambda$ at which $t(e) \geq |K(e)|$.  An edge is said to be
{\em open} if it lies in the complement of $\mathcal S$.

The reader may verify the following:

\begin{lem} \label{Isingclusterexpansion}
Let $\xi$ be a random function on $\Lambda$ with law given by the
Ising model with coupling constants $\{K_e\}$, with boundary
conditions $\xi_0$ on $\partial \Lambda$.  Let $r$ be an independent
product of parameter one exponentials.  If we condition on the total
energy function $t$ (which determines $\mathcal S$), then (any
regular version of) the conditional law of $\xi$ satisfies the
following almost surely:
\begin{enumerate}
\item All open edges are a.s.\ satisfied by $\xi$. This fact and the boundary conditions uniquely determine $\xi$ on
each open cluster that contains a vertex of $\partial \Lambda$.
\item In each open cluster that does not contain a vertex of $\partial \Lambda$, there are exactly two
ways, differing by a sign change, of defining $\xi$ on that cluster
so that all the edges in the cluster are satisfied.
\item The law of $\xi$ conditioned on $t$ is given by tossing an
independent fair coin to determine the sign of each open cluster of
$\Lambda$ that does not contain a boundary vertex.
\end{enumerate}
\end{lem}

In other words, conditioned on $t$, there are $2^\text{(number of
open clusters)}$ ways to choose a pair $(\xi, r)$ with total energy
$t$, and each of them is equally likely. (Note that isolated vertices---
i.e., vertices all of whose edges are in $\mathcal S$---are also clusters,
so the coin toss applies to these sites as well.) In particular, all the
information needed for determining the law of $\xi$ conditioned on
$t$ is contained in $\mathcal S$. Note that $\mathcal S$ is the set
of edges $e$ that are either unsatisfied or are satisfied and have
$r(e) \geq 2|K_e|$. Thus, conditioned on $\xi$, the law of $\mathcal
S$ is given by a Bernoulli percolation on $\Lambda$, with an edge
belonging to $\mathcal S$ with probability $1$ if $e$ is unsatisfied
and with probability $\int_{2|K_e|}^\infty e^{-x}dx = e^{-2|K_e|}$
if $e$ is satisfied. The Swendsen-Wang update to $\xi$ described in
\cite{SW} is performed by first sampling the set of open edges
(using Bernoulli percolation on the satisfied clusters of $\xi$ with
the probabilities given above) and then tossing a fair coin to
decide the sign of $\xi$ on each open cluster that does not contain
a boundary vertex. (The law of the set of open clusters of $(\xi,
r)$---called {\em Fortuin-Kasteleyn clusters}---is also simple and
was described by Fortuin and Kasteleyn in 1972 \cite{FK}. A similar
analysis applies to Potts models.)

\subsection{Random surfaces and Ising models}
We now return to the main setting of this chapter: $\Phi$ is an SAP,
$E = \mathbb Z$ or $\mathbb R$, and $\Lambda \subset \subset \mathbb
Z^d$.   Given a pair $\phi_1, \phi_2 \in \Omega$ of admissible
functions (as defined in Section \ref{gibbsdefinition}), define a
non-decreasingly-ordered-pair valued function $$\xi(x) = (\min
\{\phi_1(x), \phi_2(x) \}, \max \{\phi_1(x), \phi_2(x) \})$$ and a
$\{-1,0,1\}$-valued function $\zeta(x) = 1_{\phi_1(x) > \phi_2(x)} -
1_{\phi_1(x) < \phi_2(x)}$.  We will also sometimes interpret
$\xi(x)$ as representing the unordered set $\{\phi_1(x), \phi_2(x)
\}$ since it contains no information about which of the two values
came from which function.  Note that $\xi_1(x)$ refers to the
smaller of $\phi_1(x)$ and $\phi_2(x)$ and $\xi_2(x)$ to the larger.
In light of the following trivial result, we may think of the map
$(\phi_1, \phi_2) \rightarrow (\xi, \zeta)$ as a measure-preserving
change of coordinates.

\begin{lem} \label{xizetameasurepreserving}
The map $E^2 \rightarrow E^2 \times \{-1,0,1\}$ that sends $(\eta_1,
\eta_2)$ to $$((\min \{\eta_1, \eta_2 \}, \max \{\eta_1, \eta_2 \}),
1_{\eta_1 > \eta_2} - 1_{\eta_1 < \eta_2})$$ is injective.  The
$\lambda \times \lambda$ measure of any measurable subset of $E
\times E$ is equal to the $\lambda \times \lambda$-times-counting
measure of its image under this map.
\end{lem}

We now make the connection between this setup and the previous
section:

\begin{lem} \label{xizetalaw}
Let $\phi_1, \phi_2$ be an independent pair of random functions
sampled from $\gamma_\Lambda(\cdot, \phi_1^0)$ and
$\gamma_\Lambda(\cdot, \phi_2^0)$ (where $\phi_1^0, \phi_2^0 \in
\Omega$ are admissible functions that determine the boundary
conditions outside of $\Lambda$). If $(\xi, \zeta)$ are constructed
from $(\phi_1, \phi_2)$ as above, then conditioned on $\xi$, it is
almost surely the case that (any regular version of) the conditional
law of $\zeta$ in $\Lambda \backslash \{x \in \Lambda: \phi_1(x) =
\phi_2(x) \}$ is given by a ferromagnetic Ising model with coupling
constants $K_{(x,y)}\leq 0$ given by the potential edge energy
$\Phi_e(\phi_1)+\Phi_e(\phi_2)$ when one of the $\phi_i$ is greater
than or equal to the other at both endpoints minus the potential edge energy when this is not the
case.  Explicitly,
\begin{eqnarray} K_{(x,y)} &=& \left[ V_{x,y} (\xi_1(y) - \xi_1(x)) + V_{x,y}( \xi_2(y) -
\xi_2(x))\right]\\ & & - \left[V_{x,y} (\xi_1(y) - \xi_2(x)) + V_{x,y}( \xi_2(y) -
\xi_1(x))\right].\\
\end{eqnarray}

\end{lem}

\begin{proof}
If $E = \mathbb R$, then $\{x \in \Lambda: \phi_1(x) = \phi_2(x) \}$
is almost surely empty. For each possible value of $\xi$ restricted
to $\Lambda$, the map $(\phi_1, \phi_2) \rightarrow \xi$ has
$2^{|\Lambda|}$ inverses, and the map is Lebesgue measure preserving
(in fact, affine) in a neighborhood of each of them. We conclude
that the conditional law of each possible inverse $(\phi_1, \phi_2)$
of $\xi$, given $\xi$, is proportional to its Gibbs weight
$e^{-H_\Lambda(\phi_1) - H_\Lambda(\phi_2)}$, and this proves the
lemma.  A similar argument applies when $E = \mathbb Z$.  In this
case, there are $2^{|\Lambda \backslash \{x \in \Lambda: \phi_1(x) =
\phi_2(x) \}|}$ possible inverses, and the probability of each of
them is proportional to its Gibbs weight.  The fact that the value
of $K_e$ described above is always non-positive (and hence the model
is ferromagnetic) follows immediately from the convexity of the
$V_{x,y}$ and Lemma \ref{convexidentity}. \qed
\end{proof}

\begin{flushleft}{\bf Remark:} Lemma \ref{xizetalaw} also applies
when the $V_{x,y}$ are non-convex, but the $K_e$ may be positive in
this case, i.e., the corresponding Ising model may not be
ferromagnetic. It also applies if we replace $E = \mathbb Z$ or $E=
\mathbb R$ with another choice of $E$ (e.g., a finite set) and
define $\xi$ and $\zeta$ in terms of an arbitrary ordering on $E$
(which may not come with a canonical ordering, as $\mathbb Z$ and
$\mathbb R$ do).
\end{flushleft}

\subsection{Defining cluster swaps}

Now we will formally define some cluster swapping maps.  Write
$$\overline{\Omega} = \Omega \times \Omega \times \Sigma,$$ where
$\Sigma$ is the set of functions from $\mathbb E^d$ to $[0, \infty)$
and $\mathbb E^d$ is the set of edges of the lattice $\mathbb Z^d$.
Let $\overline{\mathcal F}$ be the product $\sigma$-algebra on
$\overline{\Omega}$ and let $\overline{\mathcal F}^{\tau}$ be the
$\sigma$-algebra generated by $\mathcal F^{\tau} \times \mathcal
F^{\tau} \times$ times the product $\sigma$-algebra on $\Sigma$. Let
$\pi$ be the measure on $\Sigma$ describing an independent product
of parameter one exponentials (i.e., the Gibbs measure on $\Sigma$
in which $H_\Lambda(r) = \sum_{x \in \Lambda}r(x)$ for each $\Lambda
\subset \subset \mathbb E^d$ and each $r \in \Sigma$).

Given a triple $(\phi_1, \phi_2, r)$, we can define $(\xi, \zeta,
t)$, where $t$ is the total energy per edge, as defined in the
previous section.  First we observe that $(\phi_1, \phi_2, r)
\rightarrow (\xi, \zeta, t)$ is a measure-preserving coordinate
change:

\begin{lem} \label{xizetatmeasurepreserving}
The map $(\phi_1, \phi_2, r) \rightarrow (\xi, \zeta, t)$ on any
finite graph $\Lambda$ is an injective, measure preserving
coordinate change.  That is, the measure of a measurable set in the
natural underlying measure on admissible $(\phi_1, \phi_2, r)$
configurations (i.e., $\lambda^{|\Lambda|} \times
\lambda^{|\Lambda|}$ times Lebesgue measure on the product of
$[0,\infty)$ over the edges of $\Lambda$) is equal to the measure of
its image in the natural underlying measure on $(\xi, \zeta, t)$
configurations (i.e., $(\lambda^2)^{|\Lambda|}$ times counting
measure times Lebesgue measure on the product of $\mathbb R$ over
the edges of $\Lambda$).
\end{lem}
\begin{proof}

The map $(\phi_1, \phi_2, r) \rightarrow (\xi, \zeta, r)$ is
injective and measure preserving by Lemma
\ref{xizetameasurepreserving}.  The map $(\xi, \zeta, r) \rightarrow
(\xi, \zeta, t)$ is injective and measure preserving because $t-r$
is a continuous function of $(\xi, \zeta)$.
\end{proof}

Combining Lemmas \ref{Isingclusterexpansion}, \ref{xizetalaw}, and
\ref{xizetatmeasurepreserving} gives the following:

\begin{lem} \label{clustersignsindependent}
Let $(\phi_1, \phi_2, r) \in \overline{\Omega}$ be a random triplet
with law $\gamma_\Lambda(\cdot, \phi_1^0) \otimes
\gamma_\Lambda(\cdot, \phi_2^0) \otimes \pi$ (so $\phi_1^0, \phi_2^0
\in \Omega$ determine the boundary conditions outside of $\Lambda$),
and define $(\xi,\zeta, t)$ from $(\phi_1, \phi_2,r)$ as above. Then
conditioned on $\xi$ and the total energy $t$ (which determine the
$K_e$ and $\mathcal S$), the conditional law of $\zeta$ is as
follows: throughout each component of the complement of $\mathcal S$
containing a vertex outside of $\Lambda$, $\xi$ is equal to its
value at that vertex. On each component of the complement of
$\mathcal S$ that is strictly contained in $\Lambda$, $\eta$ is a.s.\
constant, and the law of the values on these components is given by
an independent fair coin toss on each component.
\end{lem}

We now define the maps that we will call {\em cluster swaps}. Let
$R_\Lambda^x:\overline{\Omega} \rightarrow \overline{\Omega}$ be the
map such that $R_\Lambda^x(\phi_1, \phi_2, r) = (\phi_1', \phi_2',
r')$ where --- if $(\xi, \zeta,t)$ and $(\xi', \zeta',t')$ are
defined from the triplets as above --- we have $\xi' = \xi$, $t' =
t$, and $\zeta' = \zeta$ unless the vertices in the open cluster
containing $x$ (as defined by $(\phi_1, \phi_2, r)$) are all
contained within $\Lambda$, in which case $\zeta' = -\zeta$ on that
cluster and $\zeta' = \zeta$ everywhere else.  We write $R^x =
R_{\mathbb Z^d}^x$.

Informally, $R^x_\Lambda$ swaps values of $\phi_1$ and $\phi_2$ on
the open cluster containing $x$ (provided that cluster is contained
in $\Lambda$) and then adjusts the residual energy in such a way
that the total energy on each edge is unchanged.  Clearly,
$R^x_\Lambda$ is an involution.  When we use $(\xi, \zeta, t)$
coordinates, it is obvious that $R^x_\Lambda$ preserves both the
underlying measure and the Hamiltonian $|t|$.  We conclude the
following:

\begin{lem} \label{Rxfixesmu}
If $\phi^0_1$ and $\phi^0_2$ are admissible and $\overline{\mu} =
\gamma_\Lambda(\cdot|\phi^0_1) \otimes
\gamma_\Lambda(\cdot|\phi^0_2) \otimes \pi$, where $\Lambda \subset
\subset \mathbb Z^d$ and $x \in \Lambda$, then
$\overline{\mu}R^x_\Lambda = \overline{\mu}$.
\end{lem}

We refer to $\mathcal S$ as the set of {\em closed} or {\em
swappable} edges, meaning that there is enough total energy on the
edge to make it possible to swap the values of $\phi_1$ and $\phi_2$
at one endpoint of the edge and not the other while (after adjusting
$r$) preserving the total energy on that edge.  Edges in $\mathbb
E^d \backslash \mathcal S$ are called {\em open} or {\em
unswappable}. Observe in particular that whenever $\phi_1$ and
$\phi_2$ agree on one of the endpoints of an edge $e$, we have $e
\in \mathcal S$. Thus, each of the points on which $\zeta= 0$ is its
own cluster of $\mathbb E^d \backslash \mathcal S$.

\subsection{Perfect matching example and uniqueness proof overview}

One simple setting for cluster swapping is domino tiling or
perfect-matching-on-$\mathbb Z^2$ model described in Section
\ref{dominosubsection}.  We interpret cluster swaps in this setting
and, as a preview of later sections, sketch the arguments that show
the uniqueness of the gradient Gibbs measure of slope $u \in
U_\Phi$.

Recall that there was a one-to-one correspondence between perfect
matchings of $\mathbb Z^2$ and finite energy height functions on the
faces of $\mathbb Z^2$ (defined up to additive constant), with
respect to the appropriate potential. (Although the rest of our
exposition assumes the height functions are defined on vertices, it
will be simpler to visualize the correspondence in this section if
we adopt the dual perspective and consider the functions to be
defined on faces.)

\begin{center}\begin{figure} \label{dominoswap} \begin{center} {
\resizebox{16
cm}{!}{\includegraphics[angle=-90]{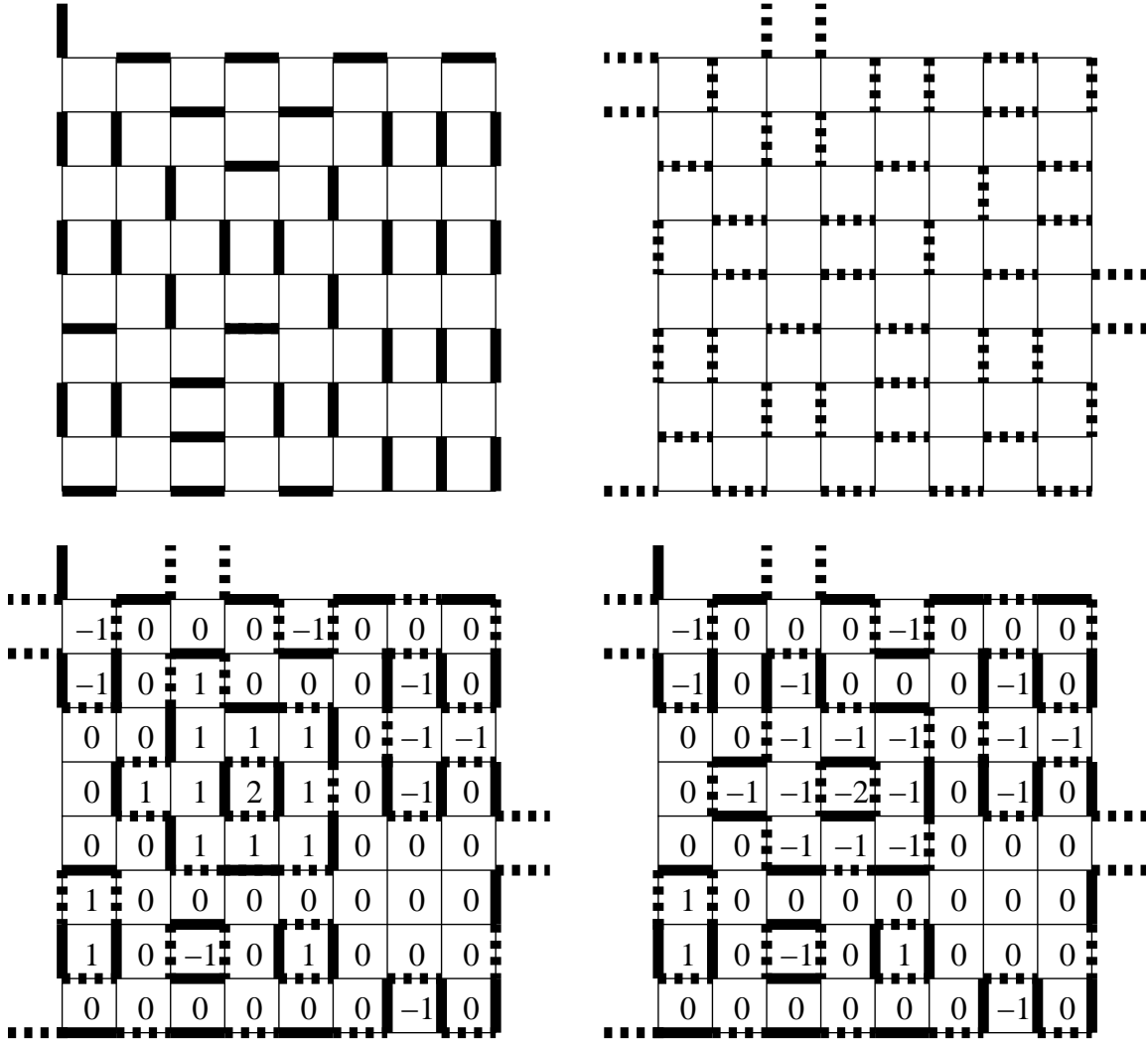}}}
\caption{Above: the edges of $\mathcal T_1$ and $\mathcal T_2$
intersecting the grids shown. Bottom left: the edges of $\mathcal
T$, the symmetric difference of $\mathcal T_1$ and $\mathcal T_2$,
together with the height difference $\phi' = \phi_1 - \phi_2$.  The
closed-edge set $\mathcal S$ consists of edges (dual to those shown)
whose endpoints are squares at least one of which has height zero.
The open clusters are the islands on which $\phi'$ is positive or
negative.  If $\Lambda$ is the $8 \times 8$ collection of square
faces shown, then there are four open clusters strictly contained in
$\Lambda$. Bottom right: the height difference $\psi' = \psi_1' -
\psi_2'$ and the corresponding tilings, where $(\psi_1, \psi_2, s) =
R_x(\phi_1, \phi_2, r)$ and $x$ is the square with $\phi'(x)=2$.
(The values of $r$ and $s$ are irrelevant in this model.)}
\end{center}
\end{figure}
\end{center}

If $\phi_1$ and $\phi_2$ are two such height functions,
corresponding to perfect matchings $\mathcal T_1$ and $\mathcal
T_2$, then $\phi' := \phi_2-\phi_1$ is a function on the square
faces in the $\mathbb Z^2$ lattice with the following properties
(see Figure \ref{dominoswap}):
\begin{enumerate}
\item If $x$ and $y$ are adjacent squares and the edge between them
lies in both or neither of $\mathcal T_1$ and $\mathcal T_2$, then
$\phi'(x) = \phi'(y)$.
\item If $x$ and $y$ are adjacent squares and the edge between them
lies in exactly one of $\mathcal T_1$ and $\mathcal T_2$, then
$|\phi'(x) - \phi'(y)|=1$.
\end{enumerate}

Let $\mathcal T$ be the set of edges that belong to exactly one of
$\mathcal T_1$ and $\mathcal T_2$.  Since every vertex is incident
to exactly zero or two edges in $\mathcal T$, $\mathcal T$ is a
disjoint union of finite-length cycles and infinite paths, which
partition the squares of $\mathbb Z^2$ into regions on which $\phi'$
is constant. The value of $\phi'$ changes by $\pm 1$ as one crosses
one of these cycles.  Recall also that $\mathcal L$ is the set of
elements of $\mathbb Z^2$ such that translation by these elements
preserves the standard ``chessboard'' coloring of the squares of
$\mathbb Z^2$. The following simple proposition is illustrated in
Figure \ref{dominoswap}.

\begin{prop} In the domino tiling setting, $\mathcal S$ consists of the
set of all edges such that $\phi' = 0$ on at least one endpoint of
that edge. The open clusters are the connected components of $\{x :
\phi'(x) \not = 0 \}$ and the boundary of each open
non-boundary-intersecting cluster is a cycle of edges in $\mathcal
T$.  The cluster swap $R_x$ reverses the sign of $\phi'$ on the
component of $\{x : \phi'(x) \not  = 0 \}$ containing zero and
leaves $\phi'$ unchanged elsewhere.
\end{prop}

A cluster swap in this context---as described in the previous
section---amounts to swapping the edge sets of $\mathcal T_1$ and
$\mathcal T_2$ that lie in the interior of one of the cycles.  When
one swaps the edge sets of $\mathcal T_1$ and $\mathcal T_2$ within
a region, this does not alter $\mathcal T$, but it {\em reverses}
whether the value of $\phi'$ changes by $1$ or $-1$ as one crosses
each cycle in that region. We refer to the swapping of the edges
within a single cycle of $\mathcal T$ (i.e., swapping which of the
two alternating sets of edges in the cycle belongs to $\mathcal T_1$
and which belongs to $\mathcal T_2$) as {\em reversing the
orientation} of the cycle.


Since in this model, each $V_{x,y}(\eta)$ is equal to $0$ or
$\infty$ for all $\eta$, the values of $r$ are in fact irrelevant.
That is, whenever $R^x(\phi_1, \phi_2, r) = (\psi_1, \psi_2, s)$, we
have $r=s$, and the values of $\psi_1$ and $\psi_2$ do not depend on
the value of $r$.

Now, Lemma \ref{Rxfixesmu} implies that conditioned on a finite
cycle of $\mathcal T$---separating a height zero region outside from
a height $\pm 1$ region inside---and on all the edges of $\mathcal
T_1$ and $\mathcal T_2$ outside of that cycle, each of the two
orientations of the cycle is equally probable.  Applying the same
argument using $R^x(\phi_1, \phi_2+c, r)$, where $c \in \mathbb Z$,
we see that the same is true for all cycles, and from this it is not
hard to prove the following:

\begin{lem}
Conditioned on the set $\mathcal T$, the orientations of the finite
cycles of $\mathcal T$ have the law of independent fair coin tosses.
\end{lem}

Of course, this is essentially obvious even without cluster
swapping, but the cluster swapping argument will be useful in more
general settings. Using this lemma, the reader may be able to
mentally prove the following (a more general version of which we
prove later):

\begin{lem} \label{dominoequality} If $\mu_1$ and $\mu_2$ are distinct $\mathcal L$-invariant Gibbs
measures on perfect matchings and $\mu_1 \otimes \mu_2$-almost
surely the symmetric difference $\mathcal T$ of a pair $(\mathcal
T_1,\mathcal T_2)$ contains no infinite paths, then $\mu_1 = \mu_2$.
\end{lem}

The next few sections will use the variational principle and a
variety of cluster swapping arguments to prove a version of the
following which holds for SAPs in any dimension:

\begin{lem} \label{atmostoneinfinitepath}
If $\mu_1$ and $\mu_2$ are ergodic gradient Gibbs measures of the
same slope $u \in U_\Phi$, then $\mu_1 \otimes \mu_2$ almost surely,
the symmetric difference $\mathcal T$ of $\mathcal T_1$ and
$\mathcal T_2$ contains {\em at most one} infinite path.
\end{lem}

Chapter \ref{discretegibbschapter} will then---in a much more
general but strictly two-dimensional context---rule out the case of
one infinite path:

\begin{lem} \label{noinfinitepath}
In the setting of the previous lemma, $\mathcal T$ almost surely
does not contain a single infinite path.
\end{lem}

The lemmas above will together imply the uniqueness of the ergodic
Gibbs measure of slope $u \in U_\Phi$.  We now return to the more
general setting in which $\Phi$ is any SAP.

\section{Monotonicity and log concavity}
\subsection{Stochastic domination via cluster swapping}
The following ``monotonicity'' property is very well known for many
systems with convex difference potentials; it is used, for example,
in \cite{FS} and \cite{CEP} for Ginzburg-Landau and domino tiling
models, respectively. Cluster swapping is one convenient way of
proving this fact.

Recall, first that if $\mu$ and $\nu$ are probability measures on an
arbitrary measure space $(X, \mathcal X)$ and $\leq$ is a partial
ordering $X$, then we say that $\mu \prec \nu$ or {\it $\nu$
stochastically dominates $\mu$} if there exists a measure $\rho$ on
$X \times X$ (with the product $\sigma$-algebra) such that on the
set of pairs $(a,b) \in X \times X$, it is $\rho$ a.s. the
case that $a \leq b$, and the first and second marginals of $\rho$
are respectively $\mu$ and $\nu$. When $\phi$ and $\psi$ are real or
integer valued functions with common domains, we use the partial
ordering $\phi \leq \psi$ to mean that $\phi(x) \leq \psi(x)$ for
all $x$ in the domain.

\begin{lem} \label{stochasticdomination} Suppose that $\phi^0_1,
\phi^0_2 \in \Omega$ are admissible and $\phi^0_1 \leq \phi^0_2$.
Then for any $\Lambda \subset \subset \mathbb Z^d$, we have
$\gamma_{\Lambda}(\cdot|\phi^0_1) \prec
\gamma_{\Lambda}(\cdot|\phi^0_2)$. \end{lem}

\begin{proof} The measure $\gamma_{\Lambda}(\cdot|\phi^0_1)
\otimes \gamma_{\Lambda}(\cdot|\phi^0_2) \otimes \pi$ on triplets
$(\phi_1,\phi_2,r)$ induces a corresponding measure on triplets
$(\xi,\zeta,t)$.  Clearly, $\zeta \geq 0$ at all vertices outside of
$\Lambda$.  Since $\zeta$ is a.s.\ constant on each component of the
complement of $\mathcal S$ and either $0$ or $1$ at all vertices
outside of $\Lambda$, this implies that any open cluster that is not
strictly contained in $\Lambda$ has $\zeta \geq 0$.

By Lemma \ref{clustersignsindependent}, the sign of $\zeta$ on each
open cluster on which $\zeta \not = 0$ may be determined by an
independent coin toss---in other words, on each open cluster the
coin toss decides whether $\phi_1 = \xi_1$ and $\phi_2 = \xi_2$ or
$\phi_1 = \xi_2$ and $\phi_2 = \xi_1$.

Suppose that instead, for each cluster we toss an independent fair
coin and, depending on the outcome, either take $\phi_1 = \phi_2
=\xi_1$ or $\phi_1 = \phi_2 = \xi_2$.  Clearly, this change does not
affect the marginal distributions of $\phi_1$ and $\phi_2$. However,
it does guarantee that we will have $\phi_1 = \phi_2$ at all
vertices that are not on open clusters containing vertices outside
of $\Lambda$.  Since $\phi_1 \leq \phi_2$ on such clusters, $\phi_1
\leq \phi_2$ throughout $\Lambda$. \qed
\end{proof}

The following is immediate:

\begin{cor} \label{boundarycontinuity} If $\phi^0_1(x) \leq \phi^0_2(x) \leq \phi^0_1(x) + c$ for all
vertices $x \in \mathbb Z^d \backslash \Lambda$ which are adjacent
to a vertex in $\Lambda$, then
$$\gamma_{\Lambda}(\cdot|\phi^0_1) \prec \gamma_{\Lambda}(\cdot|\phi^0_2) \prec
\gamma_{\Lambda}(\cdot|\phi^0_1+c).$$ \end{cor}

We will also use the following as a technical lemma.

\begin{cor} \label{boundarycontinuity2} Let $\phi_c$ be the function which is equal to an admissible
function $\phi_0$ everywhere except at one vertex $x \in \mathbb Z^d
\backslash \Lambda$ where it is equal to $c$; when $c$ is in the
interval for which $\phi_c$ is admissible, let $F(c)$ be the
$\gamma_{\Lambda}(\cdot|\phi_c)$-expected value of $\phi(y)$, where
$y \in \Lambda$. Then $F(c)$ is monotone increasing and $F(c_2) -
F(c_1) \leq c_2 - c_1$ for all $c_1, c_2 \in \mathbb E$. In
particular, if $c$ is chosen from a distribution $\nu$ on $E$
(supported on $c$ for which $\phi_c$ is admissible), then the
variance of $F(c)$ is less than or equal to the variance of $c$.
\end{cor}

\begin{proof} The first two claims follow immediately from Lemma \ref{stochasticdomination}.  All
that remains to prove is that if $c$ is chosen from $\nu$ and $F$ is
monotone Lipschitz (i.e., $F(c_2) - F(c_1) \leq c_2 - c_1$) then the
$\nu$-variance of the random variable $F(c)$ is less than or equal
to that of $c$. Equivalently, if $a_1$ and $a_2$ are sampled
independently from $\nu$, we would like to show that the variance of
$F(a_1) - F(a_2)$ is less than or equal to that of $a_2 - a_1$.
Since both variables have mean zero, and $(F(a_1) - F(a_2))^2 \leq
(a_1 - a_2)^2$ for all $a_1,a_2$, the result follows.  \qed
\end{proof}

\subsection{Log concavity via cluster swapping} \label{highdimensionswap}

A probability distribution on $\mathbb E$ is {\em log concave} if
the log of its Radon-Nikodym derivative $f$ with respect to $\mathbb
E$ is a concave function.  (In particular, $f$ is continuous on the
interval on which it is finite). On $\mathbb Z$, this is equivalent
to the statement that $f$ is continuous and $2\log f(a) \geq \log
f(a+1)+\log f(a-1)$ for all $a \in \mathbb Z$ (where we write $\log
0 = - \infty$), or equivalently $f(a)^2 \geq f(a+1)f(a-1)$ for all
$a \in \mathbb Z$.  Log concavity also implies that $f(a)^2 \geq
f(a+c)f(a-c)$ for all $c \in \mathbb Z$.

If $E =\mathbb R$ and $f$ is assumed to be continuous, then the log
concavity of $f$ is equivalent to the statement that $f(a)^2 \geq
f(a-c)f(a+c)$ for all $c \in \mathbb R$.  (If $f$ is continuous and
fails to be convex, then it is easy to see that there is some
arithmetic sequence along which it fails to be convex, and the
discrete characterization above applies to that sequence.)

\begin{lem}\label{logconcavefiniteset} Suppose $\Lambda \subset \subset \mathbb Z^d$, $x_0 \in
\Lambda$, and $\phi_0 \in \Omega$ is admissible.  If $\phi$ is a
random function chosen from $\gamma_{\Lambda}(\cdot|\phi_0)$, then the
law of $\phi(x_0)$ is log concave.  In fact, the same result holds
if $\Lambda$ is a finite subset of the vertices of any connected
graph and $\Phi$ a convex nearest neighbor potential on that graph,
with admissible boundary conditions fixed outside of $\Lambda$.
\end{lem}

\begin{proof}
When $E = \mathbb R$, this follows from a variant of the
Brunn-Minkowski inequality known as the Pr\'{e}kopa-Leindler
inequality; see Theorem 4.2 of \cite{Gar} for details.  We present a
simple argument that uses cluster swaps in the case $E = \mathbb Z$.

Let $\phi_1^0 = \phi_0$ and $\phi_2^0 = \phi_0 + c$ for some $c\in
\mathbb Z$ with $c>0$, and sample $(\phi_1, \phi_2, r)$ from
$\gamma_\Lambda(\cdot|\phi_1^0) \otimes
\gamma_\Lambda(\cdot|\phi_2^0) \otimes \Sigma$.  Then conditioned on
$\xi(x_0)=(a,a+c)$ (for some $a \in \mathbb Z$ such that this occurs
with positive probability) and $t$, the probability that
$\phi_2(x_0) \geq \phi_1(x_0)$ is $1/2$ if the open cluster
containing $x_0$ is contained in $\Lambda$ and $1$ otherwise, by
Lemma \ref{clustersignsindependent}. Thus, conditioned only on
$\xi(x_0)=(a,a+c)$ the probability that $\phi_2(x_0) \geq \phi_1(x_0)$
is between $1/2$ and $1$, which implies that $f(a)^2 \geq
f(a+c)f(a-c)$.

A similar argument---using regular conditional
probabilities---yields an alternate proof in the case $E = \mathbb
R$. The extension to general graphs is trivial.  \qed
\end{proof}

\subsection{Log concavity for extremal (non-gradient) Gibbs measures}
The log concavity arguments of Section \ref{highdimensionswap},
combined with Lemma \ref{extremalconvergence}, imply the following:

\begin{lem} \label{logconcave} If $\Phi$ is a simply attractive measure, and $\mu \in
\text{ex}\mathcal G(\Omega, \mathcal F)$, and $x \in \mathbb Z^d$
then the density of the height distribution of $\phi(x)$, for $\phi$
chosen from $\mu$, is log concave.  In particular, the random
variable $\phi(x)$ has finite mean, variance, and moments of all
orders. \end{lem}

\begin{proof} By Lemma \ref{extremalconvergence}, for $\mu$ almost all $\phi$, the measures
$\gamma_{\Delta_n}(\cdot|\phi)$ converge to $\mu$ in the topology of
local convergence, where the $\Delta_n$ are cubes of side length
$(2n+1)$, centered at the origin. Let $\nu_{\phi,n,x}$ be the law
for $\psi(x)$, when $\psi$ is chosen from
$\gamma_{\Delta_n}(\cdot|\phi)$; let $\nu_{x}$ be the law of
$\psi(x)$ when $\psi$ is chosen from $\mu$.  Then the preceding
statement implies that for $\mu$ almost all $\phi$, the measures
$\nu_{\phi,n,x}$ converge to $\nu_x$ in the $\tau$-topology (i.e.,
the smallest topology in which $\nu \mapsto \nu(A)$ is open for each
Borel set $A \subset E$). The reader may check that the
$\tau$-topology limit of a sequence of log concave distributions on
$\mathbb R$ or $\mathbb Z$ is necessarily log concave, and the
result then follows. The fact that a distribution $\nu$ is log
concave implies that the log probabilities must decrease at least
linearly; thus, the tails of $\nu$ decrease exponentially, and
moments of all orders exist. \qed
\end{proof}

A similar argument yields a characterization of smooth gradient Gibbs measures:

\begin{lem} \label{smoothchar} If $\Phi$ is a simply attractive measure, and $\mu \in
\text{ex}\mathcal G(\Omega, \mathcal F^{\tau})$, and $x \in \mathbb
Z^d$, then $\mu$ is a smooth gradient Gibbs measure (i.e., a
restriction to $\mathcal F^{\tau}$ of a Gibbs measure on $(\Omega,
\mathcal F)$) if and only if, for $\mu$ almost every $\phi$, the
measures $\nu_{\phi,n,x}$ converge to a non-zero limit in the
$\tau$-topology. \end{lem}

\begin{proof} The proof of Lemma \ref{logconcave} implies that $\nu_{\phi,n,x}$ almost surely has a
limit whenever $\mu$ is smooth.  (This also follows from Lemma
\ref{extremalconvergence}.)  For the converse, let $M_n$ be the
median of the measure $\nu_{\phi,n,x}$, and note that the value $M =
\lim_{n \rightarrow \infty} M_n$ converges for $\mu$ almost every
$\phi$ and is tail measurable. So, we define $\overline{\mu}$ as
follows: to sample $\phi$ from $\overline{\mu}$, first sample $\phi$
from $\mu$ (determined only up to additive constant) and then choose
the constant in such a way that $M = 0$ (when $E = \mathbb R$) or $M
\in [0,1)$ (when $E = \mathbb Z$). \qed \end{proof}

\subsection{Existence of minimal gradient Gibbs measures of a given
slope} Log concavity also yields a simple proof of the existence of
an ergodic $\mathcal L$-invariant gradient Gibbs measure of a given
slope $u \in U_{\Phi}$ with $\sigma(u) < \infty$.  Lemma
\ref{exposedexistence} implies the existence of an $\mathcal
L$-ergodic gradient Gibbs measure slope $u$ provided that $u$ does
not lie in an unbounded subset of $\mathbb R^d$ on which $\sigma$ is
linear.  We will now strengthen that result to all $u \in U_\Phi$.
(Note: this is a preliminary result that we will use to prove that
$\sigma$ is strictly convex.  If we could prove strict convexity of
$\sigma$ without this result, then this result would follow from
Lemma \ref{exposedexistence}.)  As always in this chapter, we assume
that $m=1$ and $\Phi$ is simply attractive.

\begin{lem} \label{ergodicexistence} There exists an ergodic gradient phase $\mu_u \in \mathcal
G^{\tau}_{\mathcal L}$ on $(\Omega, \mathcal F^{\tau})$ of slope $u$ for every $u \in U_{\Phi}$.
\end{lem}

\begin{proof} First, fix $k$ so that $k\mathbb Z^d \subset \mathcal L$.  By Lemma
\ref{torusconvergence}, some subsequence of the measures $\mu^n_u$
(defined in Section \ref{torussection}) on functions on $(nk)^d$
tori converges in the topology of local convergence to an $\mathcal
L$-invariant gradient Gibbs measure $\mu \in \mathcal P_{\mathcal
L}(\Omega, \mathcal F^{\tau})$.

Now, using the notation of Section \ref{torussection}, Lemma
\ref{logconcavefiniteset} implies that the probability density of
$\phi_g(y) - \phi_g(x)$ induced by $\mu_n$ is log concave for any
$x,y \in \mathbb Z^d$.  Also, if $x-y \in \mathcal L$, then this
density has expectation equal to the inner product
$(\frac{1}{n}\lfloor u \rfloor, y-x)$ when $E = \mathbb Z$ and
$(\frac{u}{n}, y-x)$ if $E = \mathbb R$.

For any $C_1>0$ and $C_2$, let $S_{C_1,C_2, C_3}$ be the set of log
concave probability densities on $E$ which are laws for random
variables $Y$ for which the expectation of $|Y|$ is less than or
equal to $C_1$ and the expectation of $Y$ is contained in
$[C_2,C_3]$. The reader may easily verify (by deriving a uniform
exponential bound on the decay of the tail probability densities)
that the sets $S_{C_1,C_2,C_3}$ are compact in the $\tau$-topology.
In particular, this implies that for all $x$ and $y$ in $\mathcal
L$, the $\mu$ probability density of the random variable
$\phi(y)-\phi(x)$ is log concave and $\mu (\phi(y) - \phi(x)) =
(u,y-x)$.  (One can see this by choosing $C_2 = (u,y-x-\epsilon)$
and $C_3 = (u,y-x+\epsilon)$ for arbitrarily small $\epsilon$.)

By Lemma \ref{decompositionergodic}, we can write $\mu = \int_{\ex\mathcal G^{\tau}_{\mathcal L}}
\nu w_{\mu} (d\nu)$ for some probability measure $w_{\mu}$ on the space of ergodic gradient Gibbs
measures.  By Theorem \ref{SFEisstronglyaffine}, $w_{\mu}$ is supported on the space of ergodic
gradient Gibbs measures with finite specific free energy.  By Lemma \ref{SEbound}, $w_{\mu}$ is
also supported on the space of gradient Gibbs measures with finite slope, and $S(\mu) =
\int_{\ex\mathcal G^{\tau}_{\mathcal L}} S(\nu) w_{\mu} (d\nu)$.  We claim that the random variable
$S(\nu)$---where $\nu$ is sampled from $w_{\mu}$---is equal to $u$ with probability one, and hence
$w_{\mu}$-almost every $\nu$ is an ergodic gradient Gibbs measure of slope $u$ and finite specific
free energy (as is needed for the lemma).

To this end, we first observe that there is a uniform bound
(independent of $n$) on the expected difference $\phi(y) - \phi(x)$
for any neighboring $x$ and $y$.  Choose a subsequence of the $n$
values along which $\lim - |T_n|^{-1} \log Z_{T_n}$ converges to the
value $ \liminf_{n \rightarrow \infty} - |T_n|^{-1} \log Z_{T_n}.$
The uniform bound on specific free energy implies a uniform bound
$C$ on the $\mu^u_n$ expected values of $|\phi(y) - \phi(x)|$ (as in
Lemma \ref{SEbound}) which in turn implies a uniform bound on the
$\mu^u_n$ variance of $\phi(y) - \phi(x)$ for any $n$ and any pair
of neighboring points $x$ and $y$ (as in Lemma \ref{SEbound}).

Now, we use a martingale/monotonicity argument identical to the one
in \cite{CEP} to show that if $x$ and $y$ are $j$ units apart in
$T_n$, then the variance of $\phi(y) - \phi(x)$ is bounded above by
$jC$.  First construct a path $x=a_0, a_1,a_2,\ldots, a_j = y$. Add
a constant to $\phi$ so that $\phi(a_0)=0$ and write $\delta_i =
\phi(a_i) - \phi(a_{i-1})$.  Let $\mathcal I_n$ be the smallest
$\sigma$ algebra in which $\delta_i$ is measurable for $i \leq n$.

Write $b_n = \mathbb E_{\mathcal I_n} (\phi(y) - \phi(x)) - \mathbb
E (\phi(y) - \phi(x))$, where $\mathbb E_{\mathcal A}$ represents
conditional expectation with respect to a $\sigma$ algebra $\mathcal
A$. Clearly, the sequence $b_n$ is a martingale. Writing $b_n =
b_{n-1} + \left(b_n - b_{n-1} \right)$, we have
$$\mathbb E b_n^2 = \mathbb E \left(b_{n-1}^2 + 2 b_{n-1} \left( b_n
- b_{n-1} \right) + (b_n - b_{n-1})^2 \right) = \mathbb E b_{n-1}^2
+ \mathbb E(b_n - b_{n-1})^2.$$  Inducting on $n$ gives the
following standard fact about martingales: $\mathbb E b_n^2 =
\sum_{i=1}^n \mathbb E(b_i - b_{i-1})^2$.  We will now derive bounds
on the individual terms $\mathbb E(b_i - b_{i-1})^2$

If $\delta_1, \ldots, \delta_{i-1}$ are fixed, then we may view
$b_i$ as a function of $\delta_i$.  By Corollary
\ref{boundarycontinuity2} (and its obvious analog on the torus),
this function is Lipschitz with Lipschitz constant one, and the
expected variance of $b_i$ (conditioned on $\delta_1, \ldots,
\delta_{i-1}$) is less or equal to that of $\delta_i$, i.e.,

$$\mathbb E b_i^2 - \mathbb E \left( \mathbb E_{I_{i-1}}b_i \right)^2 \leq \mathbb E \delta_i^2
- \mathbb E\left( \mathbb E_{\mathcal I_{i-1}} \delta_i\right)^2.$$

The left hand side is equal to $\mathbb E b_i^2 - b_{i-1}^2 =
\mathbb E (b_i- b_{i-1})^2$, and the right hand side is less than or
equal to $\mathbb E \delta_i^2$.  Since these $\delta_i$ are finite,
periodic functions of the edges in $\mathbb Z^d$ they correspond to,
summing over $i$ yields that that the variance of $\phi(y) -
\phi(x)$ is indeed bounded above by $jC$, where $C = \sup \{\mathbb
E(\phi(x_1) - \phi(x_2))^2 : |x_1- x_2|=1;x_1,x_2 \in \mathbb
Z^d\}$.

Together with log concavity and the compactness of the set of log
concave probability densities with given bounds on their variances
and the expectations, this implies that (for each basis vector $e_i$
of $\mathbb Z^d$), the $\mu$ probability that $\phi(je_i) -
\phi(e_i)$ differs from its expected value by more than $\epsilon j$
decays exponentially in $j$.  If, with some $w_{\mu}$-positive
probability $\delta$, the $i$th component of the slope of an ergodic
component of $\mu$ is greater than $u_i + 2\epsilon$, then the
one-dimensional ergodic theorem implies that $\liminf_{j \rightarrow
\infty} \mu (\{ \phi(je_i) - \phi(e_i) \geq (u_i + \epsilon) j \})
\geq \delta$, contradicting this exponential decay. \qed\end{proof}

Lemma \ref{torusconvergence} and Corollary \ref{newsigmadef} now imply the following:

\begin{lem} \label{ergodicisminimal} The ergodic measure $\mu_u$ of slope $u$ constructed above is
a minimal gradient phase---i.e., it satisfies $SFE(\mu_u) = \sigma(u)$. \end{lem}

\section{Measures on triplets}  \label{tripletsection}
\subsection{Definitions}
The proofs of the main theorems in this chapter all rely on
``infinite cluster swapping'' maps, which we have yet to define.  In
this section, we define and make some simple observations about
measures defined on the space $\overline{\Omega} = \Omega \times
\Omega \times \Sigma$ of infinite triplets.  Let $\sigma$-algebra
$\overline{\mathcal F}$ be the product $\sigma$-algebra on
$\overline{\Omega}$ and let $\overline{\mathcal F}^{\tau}$ be the
$\sigma$-algebra generated by $\mathcal F^{\tau} \times \mathcal
F^{\tau}$ times the product topology on $\Sigma$.

We think of $\Phi$ as extending to this space by writing:
$$\Phi_{\Lambda} (\phi_1, \phi_2, r) = \Phi_{\Lambda}(\phi_1) +
\Phi_{\Lambda}(\phi_2) + \sum_e r(e),$$ where the latter sum is over
all edges $e$ which contain at least one vertex of $\Lambda$.  We
similarly extend $H_{\Lambda}$ and the probability kernels to
triplets.  We defined these kernels in Section \ref{gibbsdefinition}
to be $$\gamma^{\Phi}_{\Lambda}(A, \phi) = Z_{\Lambda}(\phi)^{-1}
\int \prod_{x \in \Lambda} d \phi(x) \exp[-H_{\Lambda}(\phi)]
1_A(\phi).$$ We now write $$\gamma^{\Phi}_{\Lambda}(A, (\phi_1,
\phi_2, r)) = Z_{\Lambda}(\phi_1)^{-1} Z_{\Lambda}(\phi_2)^{-1} $$
$$\int \prod_{x \in \Lambda} d \phi_1(x) \prod_{x \in \Lambda} d
\phi_2(x) \prod_{e} d r(e) \exp[-H_{\Lambda}(\phi_1, \phi_2, r)]
1_A(\phi_1, \phi_2, r),$$ where again, the products over $e$ are
taken over edges with at least one vertex of $\Lambda$.  (Note that we
do not need the term $Z_\Lambda(r)^{-1}$, since this value is
identically one regardless of the size of $\Lambda$ and the value of $r$
on edges not intersecting $\Lambda$.)  A {\it
Gibbs measure} on $(\overline{\Omega}, \overline{\mathcal F})$ is a
measure on $(\overline{\Omega}, \overline{\mathcal F})$ which is
preserved by these kernels. A {\it gradient Gibbs measures} is
defined accordingly, replacing $\overline{\mathcal F}$ with
$\overline{\mathcal F}^{\tau}$.  Note that our definition implies
that in any Gibbs measure or gradient Gibbs measure on
$(\overline{\Omega}, \overline{\mathcal F})$, the random variables
$r(e)$ are independent of $\phi_1$, $\phi_2$, and independently
identically distributed according to a parameter one exponential
distribution on $[0,\infty)$.  We will denote the latter measure on
$\Sigma$ by $\pi$ throughout this chapter.

We say that a (gradient) Gibbs measure on triplets is {\it $\mathcal L$-invariant} if it is
invariant under the shifts $\theta_v, v \in \mathcal L$ which move the three components $\psi_1,
\psi_2, r$ in tandem; it is {\it $\mathcal L$-ergodic} if it is extremal in the set of $\mathcal
L$-invariant measures and {\it extremal} if it is extremal in the set of Gibbs measures (gradient
Gibbs measures) on $(\overline{\Omega}, \overline{\mathcal F})$ (resp., $(\overline{\Omega},
\overline{\mathcal F}^{\tau})$).

We can also define Gibbs measures, free energy, and specific free
energy as we did in Chapter \ref{SFEchapter}, replacing
$H_{\Lambda}$ by $\overline{H}_{\Lambda}:\overline{\Omega} \mapsto
\mathbb R$ defined by $$\overline{H}_{\Lambda}(\phi_1, \phi_2, r) =
H_{\Lambda}(\phi_1) + H_{\Lambda}(\phi_2) + \sum_{e} r(e).$$
In particular, if $\mu \in \mathcal P_\mathcal L(\overline \Omega, \overline {\mathcal F}^\tau)$,
then we write

$$SFE(\mu) = \lim_{n \rightarrow \infty} |\Lambda_n|^{-1} \mathcal H \left( \mu_{\Lambda_n},
e^{-\overline H^o_{\Lambda_n}}\lambda^{|\Lambda_n -1|}\otimes \lambda^{|\Lambda_n -1|} \otimes[0, \infty)^{|\Sigma_n|}
 \right),$$ where the expressions $\lambda^{|\Lambda_n -1|}$ are interpreted
the same way as in Section \ref{SFEdefsection}, and $\overline H^o_\Lambda$ is defined
analogously to $H^o_\Lambda$ (i.e., it is the sum of the energy contributions from edges strictly contained
in $\Lambda$), and $\Sigma_n$ is the set of edges with both endpoints in $\Lambda_n$.

We can also define the slope $S(\mu) = (u,v)$ to be the two slopes of the
marginal distributions of $\mu$.  We write $S_a(\mu) = \frac{u +
v}{2}$ for the {\it average slope} of $\mu$.

\subsection{Extremal decompositions of Gibbs measures on triplets}

The following simple fact will be frequently useful: \begin{lem}
\label{independentextremalcomponents} If a gradient measure $\mu$ on
triplets is extremal, then its three marginals are also extremal.
Also, any independent product of the form $\mu_1 \otimes \mu_2
\otimes \rho$ where $\mu_1$, $\mu_2$, and $\rho$ are extremal, is
itself extremal.
\end{lem}

(We remark that not every extremal measure on triplets is an independent of
its three components; for example, an extremal measure could have extremal marginals
but have the first two components coupled in such a way that they are almost
surely equal.)

\begin{proof} If $\mu$ is extremal, then it is clear that its
marginals must be extremal, since otherwise there would be a tail
event (an event involving only one of the three components) with
non-trivial $\mu$ probability.

To prove that the product of extremal measures is extremal suppose
otherwise, i.e., that that there exists a tail-measurable set $A
\subset \overline{\mathcal F}$ with $0 < \mu_1 \otimes \mu_2 \otimes
\rho(A) < 1$.  Then the conditional probability of $A$ given the
first component $\phi_1$ is a tail measurable function of $\Omega$
which is {\it not} $\mu_1$-almost surely constant, a contradiction.
\qed \end{proof}

Note that the analogous result for ergodic gradient measures is false.  To see that independent
products of ergodic measures can fail to be ergodic, consider the following example: let $d = 1$
and $\phi_0(i) = i \pmod 2$, and let $\mu$ be the probability measure on $\Omega$ that puts half
its mass on $\phi_0$ and half on $1 - \phi_0$. Clearly, $\mu$ is ergodic, as is its restriction to
$\overline{\mathcal F}^{\tau}$. However, $\mu \otimes \mu$ is {\it not} ergodic on $\Omega \times
\Omega$, since $\{ (\phi_0, \phi_0), (\phi_1, \phi_1) \}$ and $\{(\phi_0, \phi_1), (\phi_1, \phi_0)
\}$ are both shift-invariant events with probability $\frac{1}{2}$.

\subsection{First half of variational principle for triplets}

We will need the obvious analog of Theorem \ref{minimizersaregibbs}; its proof is identical to that
of Theorem \ref{minimizersaregibbs}.

\begin{lem} \label{tripletminimizersaregibbs}  If $\mu \in \mathcal P_{\mathcal
L}(\overline{\Omega}, \overline{\mathcal F}^{\tau})$ has minimal specific free energy among
$\mathcal L$-invariant measures with slope $(u,v)$, then $\mu$ is a gradient Gibbs measure.
\end{lem}

\begin{lem} \label{SFEminimum} The minimal specific free energy among measures of slope $(u,v)$ is
equal to $\sigma(u) + \sigma(v)$ and is obtained by an independent
product $\mu_u \otimes \mu_v \otimes \nu$ (as defined in Lemma
\ref{ergodicexistence}, where $\nu$ is an independent product of
parameter one exponentials). \end{lem}

\begin{proof} It is obvious that $SFE(\mu) \geq SFE(\mu_1) + SFE(\mu_2) + SFE(\nu)$ where $\mu_1$,
$\mu_2$, and $\nu$ are the marginals (use Lemma \ref{productentropy} and take limits).  Since the
latter term is at most zero, the statement now follows from Lemma \ref{ergodicisminimal}. \qed
\end{proof}

\section{Height offset variables} \label{HOVsection} Consider a measure $\mu \in
\mathcal P_{\mathcal L}(\Omega, \mathcal F^{\tau})$ with finite
specific free energy.  A function $h:\Omega \mapsto \mathbb R \cup
\{ \infty \}$ is called a {\it height offset variable} for $\mu$ if
the following are true: \begin{enumerate}
\item $h(\phi+c) = h(\phi)+c$ for all $\phi \in \Omega$ and $c \in
E$.
\item $h$ is tail-measurable on $\Omega$.
\item $h$ is $\mu$-almost surely finite.
\item If $v \in \mathcal L$, then, $\mu$-almost surely, $h(\phi) = h(\theta_v \phi) + (u,v)$, where
$u$ is the slope of the ergodic component of $\mu$ from which $\phi$
was chosen, i.e, $u = S(\pi^\phi_\mathcal L)$.  (Recall Lemma
\ref{decompositionergodic}.)
\end{enumerate} Although a sampling from a gradient Gibbs measure is
defined only up to additive constant, height offset variables, when
they exist, provide a canonical way of choosing that additive
constant.

Our main motivating example is when $h$ is the limit of the average
value of $\phi$ on increasingly large cubes centered at the
origin---and $h$ could be defined to be infinity if no such limit
exists; we will show in Section \ref{spectrumsection} that if $\mu$
is a smooth minimal phase, then the $h$ thus defined does satisfy
the above criteria.

\begin{lem} \label{offsetsmooth} If $\mu$ is a gradient Gibbs measure and $h$ is a height offset
variable for $\mu$, then $\mu$ is smooth---i.e., $\mu$ is the
restriction to $\mathcal F^{\tau}$ of a Gibbs measure $\mu'$ on
$(\Omega, \mathcal F)$ \end{lem} \begin{proof} We define $\mu'$ as
follows: to sample from $\mu'$, first choose $\phi$ (defined up to
additive constant) from $\mu$, and then output $\phi - h(\phi)$ (if
$E =\mathbb R$) or $\phi - \lfloor h(\phi) \rfloor$ (if $E = \mathbb
Z$). Since any function of $\phi - h(\phi)$ can be written as an
$\mathcal F^{\tau}$-measurable function of $\Omega$, this $\mu'$ is
well-defined.  Since $h$ is tail-measurable---and hence its value is
almost surely unchanged by the transitions kernels
$\gamma_{\Lambda}$---it is now straightforward to check that $\mu'$
is a Gibbs measure on $(\Omega, \mathcal F)$. \qed
\end{proof}

If $E = \mathbb R$, then $\mu'$-almost surely, $h(\phi) = 0$.  If $E
= \mathbb Z$, then $\mu'$-almost surely $h(\phi) \in [0,1)$.  When
$h$ is given, we refer to the measure $h(\mu)$, a measure on
$[0,1)$, as the {\it height offset spectrum} of $\mu$. The $\mathcal
L$-ergodicity of $\mu$ now implies the following:

\begin{lem} \label{offsetspectrum} If $E = \mathbb Z$ and $h$ is a height offset variable for an
$\mathcal L$-ergodic gradient measure $\mu$ of slope $u$, then the
height offset spectrum $h(\mu)$ (the law of $h(\phi)$ if $\phi$ is
chosen from $\mu$) is a measure on $[0,1)$ which is ergodic with
respect to the maps $x \mapsto x+(u,y) \pmod{1}$ for $y \in \mathcal
L$. In particular, if one of the components of $u$ is irrational,
then $h(\mu)$ is uniformly distributed.  Also, if $\mu$ is
extremal---so that $h(\mu)$ is a point measure---then we must have
$u \in \tilde{\mathcal L}$, where $\tilde{\mathcal L}$ is the dual
lattice of $\mathcal L$. \end{lem}

We will discuss the existence of height offset variables and their spectra in more detail in
Section \ref{spectrumsection}.

Next, we will need some analogous definitions for $\mathcal
L$-invariant, finite specific free energy measures $\mu$ on
$(\overline{\Omega}, \overline{\mathcal F}^{\tau})$.  In this
context, a function $h:\Omega \mapsto \mathbb E \cup \{ \infty \}$
is called a {\it height difference variable} for $\mu$ if the
following are true: \begin{enumerate}
\item $h(\phi_1+c_1, \phi_2+c_2, r) = h(\phi_1,\phi_2)+c_2-c_1$ for all $(\phi_1, \phi_2,r) \in \overline{\Omega}$
and $c_1,c_2 \in E$.
\item $h$ is tail-measurable on $\overline{\Omega}$.
\item $h$ is $\mu$-almost surely finite.
\item If $v \in \mathcal L$, then, $\mu$-almost surely, $h(\phi_1,\phi_2,r) = h(\theta_v \phi_1,
\theta_v\phi_2, \theta_v r)$. \end{enumerate}

(We will primarily use this definition primarily for measures $\mu$ almost all of whose
ergodic components have the same slope; hence the last requirement in the definition
does have need a term depending on slope like we have in the definition of a height offset variable.)
Again, our motivating example is that $h$ is the limit of the
average difference between $\phi_2 - \phi_1$ on large cubes centered
at the origin---if such a limit exists $\mu$ almost surely and
satisfies the above criteria.  Now, let $\overline{\mathcal
F}^{\tau}_0$ be the smallest $\sigma$-algebra in which, for any $x,y
\in \mathbb Z^d$ and $e \in \mathbb E^d$, the functions $r(e)$,
$\phi_1(x) - \phi_1(y)$, $\phi_2(x) - \phi_2(y)$, and $\phi_1(x) -
\phi_2(y)$ are all measurable functions on the set
$\overline{\Omega}$; $\overline{\mathcal F}^{\tau}_0$ differs from
$\overline{\mathcal F}^{\tau}$ in that differences {\it between}
$\phi_1$ and $\phi_2$ are $\overline{\mathcal
F}^{\tau}_0$-measurable.  Note the proper inclusions
$\overline{\mathcal F}^{\tau} \subset \overline{\mathcal F}^{\tau}_0
\subset \overline{\mathcal F}$.  We define Gibbs measures on
$(\overline{\Omega}, \overline{\mathcal F}^{\tau}_0)$ analogously to
Gibbs measures on $(\overline{\Omega}, \overline{\mathcal
F}^{\tau})$.

\begin{lem} \label{differenceoffset} Let $\mu$ be an $\mathcal L$-invariant gradient Gibbs measure
on $(\overline{\Omega}, \overline{\mathcal F}^{\tau})$ and $h$ a
height difference variable for $\mu$.  Then $\mu$ is the restriction
to $\overline{\mathcal F}^{\tau}$ of $\mathcal L$-invariant Gibbs
measure $\overline{\mu}$ on $(\overline{\Omega}, \overline{\mathcal
F}^{\tau}_0)$. Moreover, $\mu$ is $\mathcal L$-ergodic if and only
if $\overline{\mu}$ is $\mathcal L$-ergodic. \end{lem}

\begin{proof} We define $\overline{\mu}$ as follows: to sample from $\overline{\mu}$, first choose
$(\phi_1, \phi_2, r)$ (defined up to additive constant for each of $\phi_1$ and $\phi_2$) from
$\mu$, and then output $(\phi_1 + h(\phi_1), \phi_2, r), \phi_2, r)$ (defined up to a single
additive constant $c$ for both $\phi_1$ and $\phi_2$).  Since $h$ is tail-measurable---and hence
its value is almost surely unchanged by the transitions kernels $\gamma_{\Lambda}$---it is now
straightforward to check that $\overline{\mu}$ is a Gibbs measure on $(\overline{\Omega},
\overline{\mathcal F}^{\tau}_0)$.  Since $h$ is $\mathcal L$-invariant and $\overline{\mathcal
F}^{\tau}$-measurable, any $\mathcal L$-invariant, $\overline{\mathcal F}^{\tau}_0$-measurable
function of $(\phi_1 + h(\phi_1), \phi_2, r)$ can be written as an $\mathcal L$-invariant
$\overline{\mathcal F}^{\tau}$-measurable function of $(\phi_1, \phi_2, r)$---and vice versa.  It
follows that $\mu$ is ergodic if and only if $\overline{\mu}$ is ergodic. \qed \end{proof}

\begin{lem} \label{roughoffset} If $\mu$ is a gradient Gibbs measure on $(\overline{\Omega},
\overline{\mathcal F}^{\tau})$ and $h$ is a height difference variable for $\mu$, then both of the
first two marginals $\mu_1$ and $\mu_2$ of $\mu$ are smooth. \end{lem}

\begin{proof} Define a measure $\overline{\mu}$ on $(\overline{\Omega}, \overline{\mathcal F})$ as
follows: to sample from $\overline{\mu}$, first choose $(\phi_1,
\phi_2, r)$ (defined up to additive constants for $\phi_1$ and
$\phi_2$) from $\mu$.  Then pick the additive constant for $\phi_2$
in such a way that $\phi_2(0)=0$ and the additive constant for
$\phi_1$ in such a way that $h(\phi_1,\phi_2,r) = 0$.  Although
$\overline{\mu}$ is not a Gibbs measure on $(\overline{\Omega},
\overline{\mathcal F})$, its first marginal {\it is} a Gibbs measure
on $(\Omega, \mathcal F)$. This follows from the fact that for any
$\Lambda \subset \subset \mathbb Z^d$, we have $\mu = \mu
\gamma_{\Lambda}$ (where $\gamma_{\Lambda}$ is interpreted as a
transition kernel on the first coordinate of $\overline{\Omega}$
only), and that applying such a transition kernel to the first
coordinate of $(\phi_1,\phi_2,r)$ almost surely does not change
either $\phi_2$ or the $h(\phi_1,\phi_2,r)$.  A similar argument
holds for the second marginal of $\mu$.  \qed \end{proof}

\section{Infinite cluster classifications} Given $(\phi_1, \phi_2, r)$ chosen from a Gibbs measure on
$(\overline{\Omega}, \overline{\mathcal F})$, what can we say about
the infinite clusters of $\mathbb E^d \backslash \mathcal S$? How
many such clusters are there? How do the clusters change if one adds
a constant to $\phi_1$ or $\phi_2$?  In this section, we will
explore these and similar questions for Gibbs measures and gradient
Gibbs measures on triplets.

\subsection{More definitions} \label{moredefinitionssection}
For any $T \subset \mathbb Z^d$, we write the following:
\begin{enumerate} \item $T$ is {\it sparse} if $\limsup_{n
\rightarrow \infty} \frac{|\Lambda_n \cap T|}{|\Lambda_n|} = 0$.
(Throughout this chapter, we assume that the cubes $\Lambda_n$ are
centered at the origin.)  \item If $\lim_{n \rightarrow \infty}
\frac{|\Lambda_n \cap T|}{|\Lambda_n|} = \alpha$ for some $\alpha >
0$, then we say that $T$ is {\it $\alpha$-dense} or {\it has density
$\alpha$}.
\item An {\it island} of $T$ is a finite component of the complement of $T$ in $\mathbb Z^d$. \item
$\tilde{T}$ is the union of $T$ and all of the islands of $T$.
\end{enumerate} We will also apply the first two definitions to
subsets of the edges of $\mathbb Z^d$. If $T$ is chosen from an
$\mathcal L$-invariant measure on the space of subsets of $\mathbb
Z^d$, then the ergodic theorem implies that, with probability one,
it will almost surely either be empty or have some positive density
$\alpha$. Given $(\phi_1, \phi_2, r) \in \overline{\Omega}$, we
define the following variables (whose dependence on $(\phi_1,
\phi_2, r)$ will always be understood):

\begin{enumerate}
\item  $\mathcal S_c$ is the set of all edges for which $(\phi_1 + c, \phi_2, r)$ is swappable.
(Note: $\mathcal S_c$ is not $\overline{\mathcal
F}^{\tau}$-measurable, since it depends on the arbitrary constants
used to define $\phi_1$ and $\phi_2$.) \item $T^+_c$ is the union of
all vertices $v$ in infinite clusters of $\mathbb E^d \backslash
\mathcal S_c$ for which $\phi_2(v) > \phi_1(v) + c$ throughout the
cluster. Similarly, $T^-_c$ contains the vertices $v$ in infinite
clusters of $\mathbb E^d \backslash\mathcal S_c$ for which
$\phi_2(v) < \phi_1(v) + c$ throughout the cluster. (Note: if
$\phi_1(v)+c = \phi_2(v)$, then every edge incident to $v$ is
contained in $\mathcal S_c$. Also, observe that, given $\phi_1$ and
$\phi_2$, $T^+_c$ is decreasing in $c$ and $T^-_c$ is increasing in
$c$.)
\item $B^+  = \inf \{c: \text{ $T^+_c$ is empty} \}$ We say $B^+ =
\infty$ if no such $c$ exists, i.e., if $T^+_c$ fails to be empty for any finite $c$.  $B^-$ is
defined similarly: $B^- = \sup \{c: \text{ $T^+_c$ is empty} \}$.  In this case, $B^- = -\infty$ if
no such $c$ exists. \end{enumerate}

For a concrete example, the reader may check that in the perfect
matching example described earlier, an edge is in $\mathcal S_c$ if
and only if has an endpoint $x$ that satisfies $\phi_2(x) -
\phi_1(x) = c$. Also, $T^+_c$ is the set of points in or surrounded
by infinite clusters on which $\phi_2-\phi_1
>c$ and $T^-_c$ the set of points in or surrounded by infinite
clusters on which $\phi_2 - \phi_1 < c$.  Moreover, $B^-$ and $B^+$
are simply the largest and smallest values of $c$ for which the
level set $(\phi_2 - \phi_1)^{-1} (c)$ has an infinite cluster.  In
this setting, $B^- = B^+$ if and only if the union of the
corresponding perfect matchings contains no infinite paths.

The following is now a clear consequence of the above definitions:

\begin{lem} \label{crossingintervalfacts} The set $\{c : T^+_c = \emptyset \}$ is equal to the
interval $[B^+, \infty)$ if $E = \mathbb Z$ and either $(B^+,
\infty)$ or $[B^+,\infty)$ if $E = \mathbb R$.  Similarly, $\{c :
T^-_c = \emptyset \}$ is equal to $(-\infty, B^-]$ if $E = \mathbb
Z$ and either  $(-\infty, B^-)$ or $(-\infty, B^-]$ if $E = \mathbb
R$. Note that if $c \in (B^+, B^-)$, then {\it neither} $T^-_c$ nor
$T^+_c$ is empty.

Also, although $B^-$ and $B^+$ are not $\overline {\mathcal F}^\tau$-measurable, the following events {\it are} tail events
in $\overline{\mathcal F}^{\tau}$: \begin{enumerate} \item  $\{(\phi_1, \phi_2, r): B^+ = \infty
\}$ or $\{(\phi_1, \phi_2, r): B^+ = -\infty \}$ (or similar sets produced by replacing $B^+$ with
$B^-$)
\item $\{(\phi_1, \phi_2, r): \text{$B^+$ and $B^-$ are both finite and } B^+ - B^- \in \mathcal
B \}$ (where $\mathcal B$ is a Borel subset of $E$) \end{enumerate} \end{lem}

If an $\mathcal L$-ergodic gradient Gibbs measure $\mu$ on $(\overline{\Omega}, \overline{\mathcal
F}^{\tau})$ admits a height difference variable $h$ (as defined in Section \ref{HOVsection}), then
the following follows from the definitions:

\begin{lem} \label{bplus}  The values $B^- - h$, and $B^+ - h$ are both tail-measurable,
$\overline{\mathcal F}^{\tau}$-measurable, $\mathcal L$-invariant functions of $\overline{\Omega}$;
if $\mu$ is $\mathcal L$-ergodic, then they are both $\mu$ almost surely constant. \end{lem}

\subsection{Coupling extremal smooth Gibbs measures} \label{crossingcouplingsection}
Our first application of the above definitions to the comparison of
Gibbs measures is the following lemma: \begin{lem}
\label{clustercoupling} Suppose that $\mu_1$ and $\mu_2$ are
extremal (non-gradient) Gibbs measures on $\Omega$. Then there exist
values $B^+_0$ and $B^-_0$ such that $\mu_1 \otimes \mu_2 \otimes
\pi$-almost surely, we have $B^+ = B^+_0$ and $B^- = B^-_0$.
Moreover, $\mu_1 + B^-_0 \prec \mu_2 \prec \mu_1 + B^+_0$. In
particular, $B^-_0 \leq B^+_0$ with equality if and only if, up to
additive constant, $\mu_1 = \mu_2$ (i.e., the restrictions of
$\mu_1$ and $\mu_2$ to $\mathcal F^{\tau}$ are equivalent).
\end{lem} \begin{proof} The first statement simply follows from the
fact that $B^+$ and $B^-$ are tail measurable functions. To prove
the stochastic domination, we construct a coupling explicitly.
Suppose that $T^+_c$ is empty; that is, there is no infinite cluster
in $\mathbb E^d \backslash\mathcal S_c$ on which $\phi_2 > \phi_1 +
c$. Fix $k$. Then given a large box $\Lambda_n$, we let $\Lambda_n'$
be the set of vertices in $\Lambda_n$ that are {\it not} connected
to $\partial \Lambda_n$ by paths in $\mathbb E^d \backslash\mathcal
S_c$. These vertices are ``isolated'' from the boundary $\partial
\Lambda_n$.

Now, consider the swapping map on triplets that swaps at all vertices inside of $\Lambda_n'$ and
fixes all vertices outside of $\Lambda_n'$.  As in the proof of Lemma \ref{stochasticdomination},
we can use this measure-preserving involution to define a coupling $\nu_n \in \mathcal P (\Omega
\times \Omega, \mathcal F \times \mathcal F)$ of $\mu_1$ and $\mu_2$: to sample from the $\nu_n$,
first choose $(\phi_1, \phi_2, r)$ from $\mu_1 \otimes \mu_2 \otimes \pi$. Then modify $\phi_1$ and
$\phi_2$ in a way that replaces {\it both} of the values $\phi_1+c$ and $\phi_2$ inside of
$\Lambda_n'$ with either the values of $\phi_1+c$ (with probability $1/2$) or the values of
$\phi_2$ (with probability $1/2$), and output the modified pair $(\phi_1,\phi_2)$.

Now, fix a smaller box $\Lambda_k$, centered inside of $\Lambda_n$;
as $n$ tends to $\infty$, the probability that a pair produced
coupling satisfies $\phi_2 > \phi_1+c$ at some point in $\Lambda_k$
is bounded above by the $\mu_1 \otimes \mu_2 \otimes \pi$
probability that there exists a path in $\mathbb E^d
\backslash\mathcal S_c$ from $\Lambda_k$ to $\partial
\Lambda_n$---along which $\phi_2
> \phi_1 + c$.  This probability tends to zero as $n \rightarrow
\infty$. Now we claim that for $\mu_1 \otimes \mu_2 \otimes
\pi$-almost all triplets, there is a subsequential limit (in the
topology of local convergence) of these couplings $\nu_n$ which is a
coupling $\nu$ of $\mu_1$ and $\mu_2$ in which, $\nu$ almost surely,
$\phi_2 \geq \phi_1+c$; hence $\mu_1 + c \prec \mu_2$. (That the
marginals have limits follows from Lemma \ref{extremalconvergence},
and this implies the necessary tightness to ensure existence of a
subsequential limit.) The first half of stochastic domination
statement in the lemma now follows by taking $c=B^-_0$ (if $E
=\mathbb Z$) or by taking limits of $\nu$ defined by taking $c=c_i$
where the $c_i$ converge to $B^-_0$ from below (if $E = \mathbb R$).
\qed \end{proof}

Note that if both $T_c^+$ and $T_c^-$ are non-empty for all values
of $c \in \mathbb R$ (as might occur, for example, when samples from
$\mu_1$ and $\mu_2$ approximate planes of different slopes), then
$B^-_0 = - \infty$ and $B^+_0 = \infty$ and the above lemma gives us
no information.  The lemma also implies that in the setting of Lemma
\ref{clustercoupling}, we cannot have either $B_0^- = \infty$ or
$B_0^+ = - \infty$ (as would occur if either $T^-_c$ or $T^c_+$ were
empty for all values of $c \in \mathbb R$); the former would imply
$\mu_1 + \infty \prec \mu_2$ (or, precisely, $\mu_1 + c \prec \mu_2$
for all $c \in \mathbb R$), which is impossible.  The latter gives a
similar contradiction.  The following two lemmas are simple
consequences of Lemma \ref{clustercoupling}:

\begin{lem} \label{ergodicgradientbounds} If $\mu_1$ and $\mu_2$ are distinct gradient phases, then
$\mu_1 \otimes \mu_2 \otimes \pi$-almost surely, $B^- < B^+$. \end{lem} \begin{proof} By Lemma
\ref{clustercoupling}, it is enough to observe that $w_{\mu_1}$ and $w_{\mu_2}$ are mutually
singular, which follows from Lemma \ref{extremalconvergence}. \qed \end{proof}

\begin{lem} \label{ergodicsamemeasurebounds} If $\mu$ is a non-extremal gradient phase, then with
$\mu \otimes \mu \otimes \pi$ positive probability, $B^- < B^+$. \end{lem} \begin{proof} By Lemma
\ref{clustercoupling}, it is enough to observe that when $(\nu_1, \nu_2)$ are chosen from $w_{\mu}
\times w_{\mu}$, there is a positive probability that $\nu_1 \not = \nu_2$. \qed \end{proof}

(Since in the perfect matching model, $B^- < B^+$ if and only if the
union of the two matchings contains an infinite path, we may view
Lemma \ref{ergodicgradientbounds} as a generalization of Lemma
\ref{dominoequality}.)

\subsection{Uniqueness of infinite clusters}

\begin{lem} \label{onecluster}  Let $\mu$ be an $\mathcal L$-ergodic gradient Gibbs measure on
$(\overline{\Omega}, \overline{\mathcal F}^{\tau})$ with slope
$(u,v)$, where $u,v\in U_{\Phi}$; suppose that $h$ is a height
difference variable for $\mu$. Then there exists no $c \in \mathbb
R$ for which, with $\mu$-positive probability, either $T^+_{c-h}$ or
$T^-_{c-h}$ consists of more than one infinite component. \end{lem}

\begin{proof}  The number of infinite clusters in $T^+_{c-h}$ is a $\mathcal L$-invariant, tail
measurable, and $\overline{\mathcal F}^{\tau}$-measurable event.  As
such, it is almost surely constant.  Suppose that the number of
infinite clusters is almost surely $k$ for some $1 < k < \infty$,
and that the triplet $(\phi_1, \phi_2, r)$ is sampled from $\mu$.

Let $P$ be a path connecting points in two distinct infinite
clusters of $T^+_{c-h}$.  Observe that the set $T^+_{c-h}$ only
increases in size if we increase the value of $\phi_2$ at any finite
set of points. For each edge $e=(x,y)$ in $P$, there exists some
value $a_e$ such that if $\phi_2$ is modified so that $\phi_2(x)>a$
and $\phi_2(y) > a$, and $r$ and $\phi_1$ are left unchanged, then
we cannot have $(x,y) \in \mathcal S_{c-h}$.  By Lemma
\ref{ergodicnottaut}, it is thus almost surely possible to alter
$\phi_2$ in a finite number of places---keeping the energy
finite---in a way that connects two of the clusters of $T^+_{c-h}$.
It follows that for some $n$, there are {\it not}
$\gamma_{\Lambda_n}(\cdot|\phi_1, \phi_2, r)$-almost surely exactly
$k$ distinct clusters of $T^+_{c-h}$.  Since this is $\mu$ almost
surely the case for $(\phi_1, \phi_2, r)$ and {\it some} $n$, there
must exist an $n$ for which this is the case with positive $\mu$
probability.  But then it cannot be true that there are $\mu
\gamma_{\Lambda_n}$-almost surely $k$ infinite clusters of
$T^+_{c-h}$, so this is a contradiction.

Now, it remains only to rule out the case of infinitely many
clusters.  The following argument is due to Burton and Keane (see
\cite{BK} or \cite{Gr}).  By similar arguments to the above, we see
that for some $\epsilon$ and $\Lambda_n$, there is a $\mu$ finite
probability that applying the transition kernel $\Lambda_n$ to a
configuration has an $\epsilon$ probability of joining three or more
infinite cycles together. By the same token, there is a finite
probability that applying the transition kernel $\Lambda_n$ breaks a
single infinite cluster into three or more infinite clusters. Tile
all of $\mathbb Z^d$ with boxes of size $\Lambda_n$.  A given box is
called a {\it trifurcation box} if removing a connected cluster of
vertices inside of the box causes a single infinite cluster to break
into three or more pieces.

Now, let $Y$ be any finite set with $|Y| \geq 3$.  A {\it 3-partition} of $Y$ is a partition $\{
P_1, P_2, P_3 \}$ of $Y$ with exactly three non-empty sets $P_1$, $P_2$, and $P_3$.  Two
3-partitions $\{P_1, P_2, P_3 \}$ and $\{ Q_1, Q_2, Q_3 \}$ are {\it compatible} if there is an
ordering of their elements such that $P_1 \supset Q_2 \cup Q_3$ (or, equivalently, such $Q_1
\supset P_2 \cup P_3$).  We cite the following fact from Burton and Keane (Lemma 8.5 of \cite{G}):
If $\mathcal P$ is a family of distinct $3$-partitions of $Y$ such that each pair of elements in
$\mathcal P$ is compatible, then $|\mathcal P| \leq |Y| - 2$.

Now, for a large value of $k$, let $\Lambda_{kn}$ be a box of side length $kn$; observe that for
each trifurcation box of an infinite cluster $C$, we can choose a partition of the $\partial
\Lambda_n \cap C$ into three sets, each of which is the intersection of $\partial \Lambda_n$ with
one of the three components of the infinite cluster that is broken apart.  In fact, as Burton and
Keane observe (again, see Lemma 8.5 of \cite{G} or \cite{BK}), the set of partitions corresponding
to the trifurcation points of the intersection with any particular cluster forms a compatible
family of partitions of the intersection of $\partial \Lambda_n$ with that infinite cluster.  This
implies that the total number of trifurcation points in $\Lambda_n$ is less than $|\partial
\Lambda_n|$. Since the expected number of such points must grow linearly in $|\Lambda_n|$, this is
a contradiction.  \qed \end{proof}

\section{Gradient phase uniqueness and $\sigma$ strict convexity}
\subsection{Statement of main uniqueness and convexity results} In order to state the
main results of this section, we will need the following definition.
If $E = \mathbb Z$, then we say that a pair of ergodic gradient
phases $\mu_1$ and $\mu_2$ on $(\Omega, \mathcal F^{\tau})$are {\it
quasiequivalent} to one another if the following are true:
    \begin{enumerate}
    \item Each $\mu_i$ is a smooth phase; i.e., it is a restriction to $\mathcal F^{\tau}$ of a Gibbs measure $\mu_i'$ on
    $(\Omega, \mathcal F)$.
    \item $\mu'_1 \otimes \mu'_2 \otimes \pi$-almost surely, we have $B^+ - B^- \in
\{0,1\}$.
    \item $S(\mu_1) = S(\mu_2)$ (a simple consequence of the previous item when each $\mu_i$ has a well-defined
slope).
    \end{enumerate}
Note that by Lemma \ref{clustercoupling}, the second item implies that a pair of extremely
components $(\nu_1, \nu_2)$ sampled from $w_{\mu_1} \otimes w_{\mu_2}$ almost surely satisfies (if
additive constants are chosen correctly) $\nu_1 \prec \nu_2 \prec \nu_1 + 1$.

The following theorem is central to this section:

\begin{thm} \label{smallerSFEaverage} Suppose that $\mu$ is a measure on $(\overline{\Omega},
\overline{\mathcal F}^{\tau})$ whose first two marginals $\mu_1$ and
$\mu_2$ are minimal $\mathcal L$-ergodic gradient phases on
$(\Omega, \mathcal F^{\tau})$ with slopes in $U_{\Phi}$.  Suppose
further that with $\mu$ positive probability, we have $B^+ > B_-$
(when $E = \mathbb R$) or $B^+ > B^- + 1$ (when $E = \mathbb Z$).
Then for some appropriately defined ``infinite cluster swapping
map'' $R$, we have $SFE(R(\mu)) \leq SFE(\mu)$ and $S_a(\mu) =
S_a(R(\mu))$; moreover, $R(\mu)$ is an $\mathcal L$-invariant
gradient measure on $(\overline{\Omega}, \overline{\mathcal
F}^{\tau})$ which is {\it not} a gradient Gibbs measure. \end{thm}

From this theorem, we can immediately deduce surface tension strict
convexity, and ergodic gradient Gibbs phase uniqueness as
corollaries:

\begin{thm} \label{sigmaconvex} The surface tension $\sigma$ is strictly convex in $U_{\Phi}$.
\end{thm}

\begin{proof} Pick distinct slopes $u_1$ and $u_2$ in $U_{\Phi}$. By Lemma \ref{ergodicexistence},
there exist $\mathcal L$-ergodic Gibbs measures $\mu_1$ and $\mu_2$
of slopes $u_1$ and $u_1$ and $SFE(\mu_1) = \sigma(u_1)$,
$SFE(\mu_2) = \sigma(u_2)$.  Write $\mu = \mu_1 \otimes \mu_2
\otimes \pi$.  From Lemma \ref{clustercoupling}, it cannot be the
case that, $\mu$-almost surely, $B_+ \leq B_- +1$ (since this would
imply that the conclusion of Lemma \ref{clustercoupling} applies to
measures of different slopes).  Thus, $\mu$ satisfies the
requirements of Theorem \ref{smallerSFEaverage}. Let $(v_1,v_2)$ be
the slope of $R(\mu_1 \otimes \mu_2 \otimes \pi)$. Since $R(\mu)$ is
not a gradient Gibbs measure, by Lemma
\ref{tripletminimizersaregibbs} and Lemma \ref{SFEminimum}, the
convexity of $\sigma$, and the fact that $S_a(\mu) = S_a(R(\mu))$,
$$\sigma(u_1) + \sigma(u_2) \geq SFE(R (\mu_1 \otimes \mu_2 \otimes \pi)) > \sigma(v_1) +
\sigma(v_2) \geq 2\sigma(\frac{v_1 + v_2}{2}) = 2\sigma(\frac{u_1 +
u_2}{2}).$$

Since $\sigma$ is already known to be convex, it easily follows from
this that $\sigma$ is strictly convex on $U_{\Phi}$. \qed
\end{proof} Although $\sigma$ is also convex on the boundary of $U$,
it is not necessarily strictly convex there.  The surface tension
function for domino tilings described in \cite{CKP}, for example, is
constant on the boundary of $U$.  Using Theorem
\ref{smallerSFEaverage} we can deduce another key result.  In light
of Lemma \ref{onecluster}, we may view the $E=\mathbb Z$ part of
this statement as a generalization of Lemma
\ref{atmostoneinfinitepath}.

\begin{thm} \label{minimalphaseuniqueness} If $E = \mathbb R$, then for every $u \in U_{\Phi}$, there
exists a unique minimal gradient phase $\mu_u$ of slope $u$. If $E =
\mathbb Z$, and $\mu_1$ and $\mu_2$ are distinct minimal gradient
phases of slope $u$, then $\mu_1$ and $\mu_2$ are quasi-equivalent.
\end{thm}

\begin{proof} By Lemma \ref{ergodicexistence}, there exists at least one minimal gradient phase of
slope $u$.  Now, suppose that $\mu_1$ and $\mu_2$ are distinct
minimal gradient phases of slope $u$ which are not quasi-equivalent;
then $\mu = \mu_1 \otimes \mu_2 \otimes \pi$ satisfies the
requirements of Theorem \ref{smallerSFEaverage}.  Let $(v_1,v_2)$ be
the slope of $R(\mu_1 \otimes \mu_2 \otimes \pi)$. Since $R(\mu_1
\otimes \mu_2 \otimes \pi)$ is not a gradient Gibbs measure, by
Lemma \ref{tripletminimizersaregibbs} and Lemma \ref{SFEminimum},
the convexity of $\sigma$, and the fact that $S_a(\mu) =
S_a(R(\mu))$, we have the following contradiction: $$2\sigma(u) \geq
SFE(R (\mu_1 \otimes \mu_2 \otimes \pi)) > \sigma(v_1) + \sigma(v_2)
\geq 2\sigma(\frac{v_1 + v_2}{2}) = 2\sigma(u).$$ \qed \end{proof}

The remainder of this section is devoted to the proof of Theorem \ref{smallerSFEaverage} and the
definition of the cluster swapping map $R$ required by the theorem.  We will obtain $R(\mu)$ from
$\mu$ by either a ``single infinite swap'' or an ``infinitely repeated cluster swap'' as $\mu$
respectively does or does not admit a height difference variable $h$.

\subsection{Single infinite cluster swap} \label{singsection} Suppose that $h$ is a height difference variable for $\mu$;
write $B^-_0 = B^- - h$ and $B^+_0 = B^+ - h$.  We may assume that
either $B^- < B^+$ with positive probability and $E = \mathbb R$ or
$B^- < B^+ -1$ with positive probability and $E = \mathbb Z$.  In
particular, there exists a $c \in E$ such that with positive $\mu$
probability, $B^-_0 < c < B^+_0$, i.e., $B^- < c+h < B^+$.  Here,
both $T_{c+h}^+$ and $T_{c+h}^-$ are nonempty.

For any $c$, we define $R_c(\mu)$ to be the measure obtained as follows: to sample from $R_c(\mu)$,
first sample $(\phi_1, \phi_2, r)$ from $\mu$. Then add $c+h$ to $\phi_1$, do a cluster swap that
swaps everything outside of the new $T_0^-$ to get a new triple, and then subtract $c+h$ from first
coordinate of the new triple. To say this in precise terms, we define a map $R_c: \overline{\Omega}
\mapsto \overline{\Omega}$ (defined $\mu$-almost surely) by $R_c(\phi_1, \phi_2, r) = (\psi_1,
\psi_2, s)$ where
$$\psi_1(x) = \begin{cases} \phi_1(x) & x \in T^-_{c+h} \\
\phi_2(x)-c-h & \text{otherwise} \\
\end{cases} $$

$$\psi_2(x) = \begin{cases} \phi_2(x) & x \in T^-_{c+h} \\
\phi_1(x) + c + h & \text{otherwise} \\
\end{cases} $$ and for each $e \in \mathbb E^d$, define $s(e)$ in such a way that $(\phi_1, \phi_2,
r)$ and $(\psi_1, \psi_2, s)$ have the same total energy at the edge $e$.

Now, if $E = \mathbb Z$, then we write $R(\mu) = R_c(\mu)$, where
$c$ is some integer with the property that $T_{c+h}^+$ and
$T_{c+h}^-$ are both non-empty with positive probability (i.e.,
$B^-_0 < c < B^+0$). If $E = \mathbb R$, then (in order to simplify
a free energy computation) we will instead write $R(\mu) =
\frac{1}{c_1-c_2} \int_{c_1}^{c_2} R_c(\mu) dc$, where $B^-_0 < c_1
< c_2 < B^+_0$.

Clearly, $R(\mu)$ is $\mathcal L$-invariant.  Now, we need to argue
that $SFE(R(\mu)) \leq SFE(\mu)$, i.e., that $\lim_{n \rightarrow
\infty} |\Lambda_n|^{-1} FE_{\Lambda_n}(R(\mu)) \leq
|\Lambda_n|^{-1} FE_{\Lambda_n}(\mu)$. Recall the definition in this
context: $$FE_{\Lambda_n}(\mu) = \mathcal H_{\overline{\mathcal
F}^{\tau}_{\Lambda_n}} \[
\mu_{\lambda},e^{-H^o_{\Lambda}(\phi_1,\phi_2,r)} \prod_{x \in
\Lambda_n \backslash \{x_0\}} d[\phi_1(x) -
\phi_1(x_0)]d[\phi_2(x)-\phi_2(x_0)] \prod_e dr(e)\].$$  Now, if we
could show that $R(\mu)_{\Lambda_n}$ is the image of
$\mu_{\Lambda_n}$ under an injective map $R$ (from the space of
configurations on $\Lambda_n$ to itself) which preserves the above
measure, then the equivalence of $FE_{\Lambda_n}(\mu)$ and
$FE_{\Lambda_n}(R(\mu))$ would be obvious.  However, such a map
cannot quite be well-defined: to determine the gradient values of
the reflection map $R_c(\phi_1, \phi_2, r)$---restricted to a set
$\Lambda_n$---it is not quite enough to know $\phi_1(x) -
\phi_1(x_0)$, $\phi_2(x) - \phi_2(x_0)$ and $r(e)$ for all vertices
and edges in $\Lambda_n$.  It is also necessary to know the
difference between the additive constants of $\phi_1$ and $\phi_2$
(i.e., to know the value $\phi_2(x_0) - \phi_1(x_0) - c - h$ at some
reference vertex $x_0 \in \Lambda_n$) and to know which of the
values $x \in \partial \Lambda_n$ are members of $T^-_c$.  But if we
expand our definition of ``configuration on $\Lambda_n$'' to include
this additional information, then we can make the map well-defined.

To this end, we write $F(x) = 0$ if $x \in T^-_{c+h}$, and $F(x) =
1$ otherwise.  Given $\Lambda_n$, we will let $F_{\partial
\Lambda_n}$ be the restriction of $F$ to the boundary of
$\Lambda_n$.  Now, let $\mu'_{\Lambda_n}$ be the law of the
five-tuple $(\phi_1 - \phi_1(x_0), \nabla \phi_2 - \phi_2(x_0), r,
F_{\partial \Lambda_n}, c_0 = \phi_2(x_0) - \phi_1(x_0) - c - h)$,
where $x_0$ is a reference vertex defined on the boundary of
$\Lambda_n$, the $\phi_i$'s are defined on $\Lambda_n \backslash
\{x_0\}$, $r$ is defined on all edges within $\Lambda_n$, and $c_0
\in \mathbb R$.

Now, we can think of $R$ as a measure preserving map on the space of
five-tuples in the following way.  Given $(\phi_1 - \phi_1(x_0),
\phi_2 - \phi_2(x_0), r, F_{\partial \Lambda_n}, c_0)$, we can
compute a new five-tuple $R(\phi_1 - \phi_1(x_0), \phi_2 -
\phi_2(x_0), r, F_{\partial \Lambda_n}, c_0) = (\psi_1, \psi_2, s,
F_{\partial \Lambda_n}, c_0)$ as follows.  Define the cluster $C$ to
be the set of all vertices $v$ for which there is a $P_v$ from $v$
to a vertex $x_v \in \partial \Lambda_n$ such that $F(x_v) = 0$ and
no edge of $P_v$ is swappable in $(\phi_1, \phi_2, r)$ (where the
relative additive constants of $\phi_1$ and $\phi_2$ are chosen in
such a way that $c+h = 0$---i.e., $\phi_2(x_0) - \phi_1(x_0) =
c_0$).  The set $C$ is then simply $\Lambda_n \cap T^-_c$. Now
determine $(\psi_1,\psi_2,s)$ by fixing the values of $(\phi_1,
\phi_2, r)$ inside of $C$, swapping the values outside of $C$, and
adjusting $s$ so that the map is energy preserving on each edge;
leave the values of $F_{\partial \Lambda_n}$ and $c_0$ unchanged.
(Note: the value $F_{\partial \Lambda_n}$ is always defined based on
the infinite clusters of the original triplet $(\phi_1, \phi_2, r)$;
the value of $F_{\partial \Lambda_n}$ after swapping should not be
interpreted as referring to infinite clusters in an infinite
post-swapping configuration.  We include the same $F_{\partial
\Lambda_n}$ in both the pre-swap and post-swap five-tuples because
doing so allows us to make the swapping map invertible.)

Now, define $FE_{\Lambda_n}(\mu')$ to be the relative entropy of
$\mu'$ with respect to $\nu_1 \otimes \nu_2 \otimes \nu_3$ where
$$\nu_1 = e^{-H^o_{\Lambda_n}}(\phi_1, \phi_2, r) \prod_{x \in
\Lambda_n \backslash \{x_0 \}} d[\phi_1(x)- \phi_1(x_0)]
d[\phi_2(x)-\phi_2(x_0)] \prod_{e \in \Lambda_n} dr(e),$$ where each
$d[\phi_i(x) - \phi_i(x_0)]$ and each $dr(e)$ is Lebesgue measure
(and we write $e \in \Lambda_n$ when both endpoints of $e$ are in
$\Lambda_n$); $\nu_2 = dc_0$ is Lebesgue measure; and $\nu_3 =
dF(x)$ is counting measure.  We can now use arguments similar to
those given in Section \ref{clusterintrosection} to see that the
swapping map on five-tuples described above is invertible and
measure preserving with respect to the measure $\nu_1 \otimes \nu_2
\otimes \nu_3$.  (It is enough to observe that each $R_c$ is
invertible, well-defined, and $\nu_1 \otimes \nu_3$-measure
preserving on each of the regions $X_{F_0,C_0}$ on which $F = F_0$
and $C = C_0$.) It follows that $FE_{\Lambda_n}(\mu') =
FE_{\Lambda_n}(R(\mu'))$.

By Lemma \ref{conditionalentropy}, $FE_{\Lambda_n}(\mu')$ is equal to $FE_{\Lambda_n}(\mu)$ plus
the the $\mu$ expectation of the relative entropy of $c_0$ (with respect to Lebesgue measure) given
$(\phi_1,\phi_2,\mu)$, plus the expectation of the relative entropy of $F$ (with respect to
counting measure) given the four-tuple $(\phi_1, \phi_2, r, c_0)$.

Now, let $\mu^{\phi_1, \phi_2, r}$ be the regular conditional probability for $c_0$ given $\phi_1 -
\phi_1(x_0)$, $\phi_2 - \phi_2(x_0)$, and $r$.  Similarly, let $\mu^{\phi_1, \phi_2, r, c_0}$ be
the regular conditional probability for $F$ given the four-tuple $(\phi_1, \phi_2, r, c_0)$.  Now,
we can phrase Lemma \ref{conditionalentropy} as follows:
$$FE_{\Lambda_n}(\mu') = FE_{\Lambda_n}(\mu) + \mu \mathcal H(\mu^{\phi_1, \phi_2, r},\nu_2) + \mu'
\mathcal H(\mu^{\phi_1, \phi_2, r, c_0},\nu_3).$$

Note that the conditional distribution of $F$ with respect to
counting measure is bounded between between $-|\partial
\Lambda_n|\log(3)$ and $0$.  Also, since $c$ is chosen uniformly in
an interval of length $(c_2-c_1)$ independently of $\phi_1, \phi_2,
r$, it is clear that the second term on the righthand side is at
most $-\log(c_2-c_1)$ (or zero if $E = \mathbb Z$).  Thus,
$$FE_{\Lambda_n}(R(\mu')) = FE_{\Lambda_n}(\mu') \leq FE_{\Lambda_n}(\mu) + o(|\Lambda_n|).$$

Next, using similar notation to the above for $R(\mu)$ and $R(\mu')$
instead of $\mu$ and $\mu'$,
$$FE_{\Lambda_n}(R(\mu')) = FE_{\Lambda_n}(R(\mu)) + R(\mu) \mathcal H(R(\mu)^{\phi_1, \phi_2,
r}|\nu_2) + R(\mu') \mathcal H(R(\mu)^{\phi_1, \phi_2, r,
c_0}|\nu_3).$$

This time we need a lower bound on the second term; informally, we
must show that we do not expect the distribution of $c_0$, given
$\phi_1$ and $\phi_2$ and $r$, to be too ``spread out.''  The
easiest way to do this is to slightly alter our definition of $c_0$.
Recall that we define $c_0$ by $c_0 = \phi_2(x_0) - \phi_1(x_0) - c
- h$.  We then determined the clusters on which swapping occurred by
looking at edges at which $(\phi_1-\phi_1(x_0), \phi_2 - \phi_2(x_0)
+ c+0, r)$ is swappable.  For given values of the triplet $(\phi_1 -
\phi(x_0), \phi_2 - \phi_2(x_0),r)$, let $b_1$ and $b_2$ be the
lower and upper bounds on the set of choices of $c_0$ for which
$(\phi_1-\phi_1(x_0), \phi_2 - \phi_2(x_0) + c_0, r)$ has {\it any}
swappable edges on $\Lambda_n$; clearly, if $c_0$ lies outside of
the interval $[b_1, b_2]$, $F$ will be constant on $\Lambda_n$ and
the swapping map will either fix all of $\Lambda_n$ or swap all of
$\Lambda_n$. Thus, if $\phi_2(x_0) - \phi_1(x_0) - h \not \in
[b_1+c_1, b_2+c_2]$, then there will be no swappable edges
regardless of how $c$ is chosen in $(c_1,c_2)$.  Our new way of
choosing $c_0$ will be as follows:
first let $$B = \begin{cases} b_1+c_1 & \phi_2(x_0) - \phi_1(x_0) - h \leq b_1+c_1 \\
b_2+c_2 & \phi_2(x_0) - \phi_1(x_0) - h \geq b_2+c_2 \\
\phi_2(x_0) - \phi_1(x_0) - h & \text{otherwise} \\
\end{cases}$$ then, as before, choose $c$ uniformly in $[c_1, c_2]$ and write $c_0 = B-c$ when $E =
\mathbb R$, and simply $c_0 = B$ when $E = \mathbb Z$. Observe that
the above definitions and arguments above remain valid with this new
definition of $c_0$. Using the fact that $\mu$ has finite specific
free energy, it is not hard to show that the expected value of $b_2
- b_1$ is $O(|\Lambda_n|)$. Now, since the expected length of the
interval on which $R(\mu)^{\phi_1, \phi_2, r}$ is supported is
$O(\Lambda_n)$, and the minimal relative entropy with respect to an
interval of length $k$ is $-\log k$, Jensen's inequality implies
that $R(\mu) \mathcal H(R(\mu)^{\phi_1, \phi_2, r}|\nu_2) \geq -O
(\log |\Lambda_n))$.

Thus, we have $$FE_{\Lambda_n}(R(\mu)) \leq FE_{\Lambda_n}(R(\mu'))
+ o(|\Lambda_n|) \leq FE_{\Lambda_n}(\mu) + o(|\Lambda_n|),$$ and
$$SFE(\mu) = \lim_{n \rightarrow \infty} |\Lambda_n|^{-1}
FE_{\Lambda_n}(\mu) \leq \lim_{n \rightarrow \infty}
|\Lambda_n|^{-1} FE_{\Lambda_n}(R(\mu)) = SFE(R(\mu)).$$

Finally, note that both $T^-_{c+h}$ and $T^+_{c+h}$ become infinite clusters after the swap; thus,
it follows from Lemma \ref{onecluster}, that $R(\mu)$ is not a gradient Gibbs measure. Finally,
since the averages $\frac{[\phi_1(y) - \phi_1(x)] + [\phi_2(y) - \phi_2(x)]}{2}$ are always left
unchanged by swapping maps, it is clear that $S_a(\mu) = S_a(R(\mu))$; thus, Theorem
\ref{smallerSFEaverage} holds in this case.

\subsection{Infinitely repeated infinite cluster swaps} \label{infsection}
In this section we deal with the case that there exists no height difference variable $h$ for
$\mu$. Note that if either $B^+$ or $B^-$ were finite with positive probability, then this $B^+$ or
$B^-$ would itself be a height difference variable for the measure $\mu_0$ equal to $\mu$
conditioned on this event, and we could apply the reflection of the previous section to the measure
$\mu_0$.

We may thus assume that $B^+$ and $B^-$ are both almost surely not
finite.  Lemma \ref{ergodicgradientbounds} implies that we cannot
have either $B^+ = -\infty$ or $B^- = \infty$ with positive
probability.  Thus, $B^+ = \infty$ and $B^- = -\infty$ with positive
probability, in which case there exist infinite clusters $T^-_c$ for
all values $c \in \mathbb R$. Recall that the sets $T^-_c$ are
increasing in $c$. For any $x \in \mathbb Z^d$, let $F(x)$ be the
smallest value $c$ for which $x \in T^-_c$.

Now, given a set $\Lambda_n$, we can define a ``single cluster swapping'' operation $R_c$ that is
measure preserving on the four-tuple $(\phi_1, \phi_2, r, F_{\partial \Lambda_n})$ the same way we
did in the previous section (except this time setting $h=0$).  We write $R_c(\phi_1, \phi_2, r,
F_{\partial \Lambda_n}) = (\psi_1, \psi_2, s, F_{\partial \Lambda_n})$, where $F_{\partial
\Lambda_n}$ is left unchanged and $(\psi_1, \psi_2, s) = R_c(\phi_1, \phi_2,r)$ (as defined in the
previous section).  This map is measure preserving on $\nu_1 \otimes \nu_3$ (as defined in the
previous section).

Now, pick some positive value $M \in E$; we will consider maps of
the form $R_{kM}$ for $k \in \mathbb Z$.  Clearly, if $kM > \sup_{x
\in \partial \Lambda_n} F(x)$, then the map $R_{kM}$ leaves
$(\phi_1, \phi_2, r, F_{\partial \Lambda_n})$ unchanged.  Similarly,
if $kM < \inf_{x \in \partial \Lambda_n} F(x)$, then $R_{kM}$ (up to
additive constants) simply permutes $\phi_1 - \phi_1(x_0)$ and
$\phi_2 - \phi_2(x_0)$ and makes no change to $r$. Now, write $R =
\prod_{i=k_1}^{k_2} R_{kM}$ where $k_1$ is any {\it odd} integer for
which $k_1M < \inf_{x \in \partial \Lambda_n} F(x)$ and $k_2$ is any
integer for which $k_2 M > \sup_{x \in \partial \Lambda_n} F(x)$.
Note that up to additive constant, the map $R$ is independent of the
particular $k_1$ and $k_2$ we choose. (Because the maps $R_{kM}$
permutes $\phi_1$ and $\phi_2$ when $k < k_1$, this would not be
true if we did not fix the parity of $k_1$.)

Taking limits of the $R$ thus defined on increasingly large boxes, we can extend $R$ to a function
from $\Omega$ to $\Omega$ (see also the explicit description of $R$ below).  We now define the
measure $R(\mu)$ on $(\overline{\Omega}, \overline{\mathcal F}^{\tau})$ as follows: to sample from
$R(\mu)$, first choose $\phi_1$ and $\phi_2$ from $\mu$; pick the additive constants of $\phi_1$ so
that $\phi_1(x_0) = 0$, and choose the additive constant of $\phi_2$ so that $\phi_2(x_0)$ is
uniformly distributed in $[0, M)$.  The proof that $SFE(R(\mu)) \leq SFE(\mu)$ is now essentially
the same as the proof given in single infinite cluster swap case. We first observe that
$FE_{\Lambda_n}(\mu') = FE_{\Lambda_n}(R(\mu'))$ where $\mu'$ and $R(\mu')$ are measures on
five-tuples defined in the previous section; the only differences are first, that we may now assume
$c \in [0,M)$, since the map only depends on the value of $c$ modulo $M$, and second, that $F$ is
defined differently. However, it is still easy show that the growth of $|\mu' \mathcal
H(\mu^{\phi_1, \phi_2, r, c_0},\nu_3)|$ is $o(|\Lambda_n|)$ by using the fact that $SFE(\mu)$ is
finite to show that the discrete derivative of $F$ at every point in $\partial \Lambda_n$ has
finite expectation.

It now remains only to show that $R(\mu)$ is not a gradient Gibbs
measure.  We begin by giving a more explicit expression for the map
$R$,  We know that if $F(x) \leq kM$, then $R_{kM}$ fixes the pair
$(\phi_1(x), \phi_2(x))$ if $E = \mathbb Z$; if $F(x) < kM$, then
$R_{kM}$ fixes the pair $(\phi_1(x), \phi_2(x))$ when $E = \mathbb
R$. (When $E = \mathbb R$, the event that $F(x)$ is exactly equal to
$kM$ will always have measure zero, so we ignore this case.)  On the
other hand, if $F(x)
> kM$, then $R_{kM}$ sends the pair $(\phi_1(x), \phi_2(x))$ to the pair $(\phi_2(x) - kM,
\phi_1(x) + kM)$. And then $R_{(k-1)M}$ sends that pair to $(\phi_1(x) + M, \phi_2(x) - M)$. After
successively applying $R_{(k-2)M}, \ldots, R_{M}$, we thus end up with the pair $(\phi_1(x) +
\frac{k}{2}M, \phi_2(x) - \frac{k}{2}M)$ if $k$ is even and $(\phi_2(x) - \frac{k+1}{2}M, \phi_1(x)
+ \frac{k+1}{2}M)$ if $k$ is odd. Write $F_M(x) = \sup_{k : F(x) > kM}$. Then we have, for $\psi_1$
and $\psi_2$, $R(\phi_1, \phi_2, r) = (\psi_1, \psi_2, s)$ where

$$\psi_1(x) = \begin{cases} \phi_1(x) + \frac{F_M(x)}{2}M & \text{$F_M(x)$ is even} \\
\phi_2(x) - \frac{F_M(x)+1}{2}M & \text{$F_M(x)$ is odd} \\ \end{cases}$$

$$\psi_2(x) = \begin{cases}  \phi_2(x) - \frac{F_M(x)}{2}M) & \text{$F_M(x)$ is even} \\
\phi_1(x) + \frac{F_M(x)+1}{2}M & \text{$F_M(x)$ is odd} \\ \end{cases}$$

We can take this as the formal definition of the repeated swapping map $R$ on all of $\Omega$. As
always, $s(e)$ is determined on every edge by the requirement that $R$ preserve the total energy on
each edge. This gives us an explicit definition of $R: \overline{\Omega} \mapsto
\overline{\Omega}$. We now need the following:

\begin{lem} When $x$ and $y$ are neighboring vertices of $\mathbb Z^d$, the following are true:
\begin{enumerate}
\item If $F_M(y) = F_M(x)+1$ and $F_M(x)$ is even, then $(\psi_1+M, \psi_2,s)$ is swappable at
$(x,y)$.
\item If $F_M(y) = F_M(x)+1$ and $F_M(x)$ is odd, then $(\psi_1, \psi_2,s)$ is swappable at
$(x,y)$.
\item If $F_M(y) \geq F_M(x) + 2$, then both $(\psi_1, \psi_2,s)$ and $(\psi_1+M, \psi_2,s)$ are
swappable at $(x,y)$. \end{enumerate} \end{lem}

\begin{proof} This can be proved directly from the formal definition of $R$ given above; however,
it is most intuitively understood by tweaking the above to give yet another explicit formulation of
$R$. Consider first a triplet $(\phi_1, \phi_2, r)$ with additive constants fixed.  Now, it is
clear (e.g., from above description), that the average is unchanged by swapping, i.e., $a(x)
\frac{\phi_1 + \phi_2}{2} =\frac{\psi_1 + \psi_2}{2}$.  Thus, $\psi_1$ and $\psi_2$ are determined
by the difference: $\delta(x) = \psi_2- \psi_1$. We can write $\psi_2(x) = a(x) + \delta(x)/2$ and
$\psi_1(x) = a(x) - \delta(x)/2$, and thus the energy at contained in $\psi_1$ and $\psi_2$ at an
edge $e = (x,y)$ is $W_e = V(a(y)-a(x) + [\delta(y)-\delta(x)]/2) + V(a(y) - a(x) + [\delta(y) -
\delta(x)]/2)$. Even if $V$ is not a symmetric function, this expression is symmetric in $\delta(y)
- \delta(x)$. Denote by $\chi(e)$ the maximum value of $[\delta(y) - \delta(x)]$ for which the
above expression is less than or equal to the combined energy contained in the triplet $(\phi_1,
\phi_2,r)$ at the edge $e$. Whatever swaps we perform, $|\delta(y) - \delta(x)|$ may not exceed
$\chi(e)$.

Now, denote by $\gamma$ the function defined analogously to $\delta$
but using $\phi_1$ and $\phi_2$ instead of $\psi_1$ and $\psi_2$.
Now, it is easy to check, that the set $\mathcal S_c$ of places
where $\phi_1+c$ and $\phi_2$ are swappable is precisely the set of
points at which $|2c - \gamma(x) -\gamma(y)| \leq \chi(e)$.  If a
swapping $R_c$ swaps the values of $\phi_1 + c$ and $\phi_2$ at $y$
and fixes the values at $x$ then this has the affect of replacing
$\gamma(y)$ with $2c-\gamma(y)$ and leaving $\gamma(x)$ unchanged.
In other words, the act of ``swapping $\phi_1(y)+c$ and $\phi_2(y)$
becomes then the act of reflecting $\gamma(y)$ across the horizontal
axis of height $c$.  Note that whenever $\gamma(y)$ and $\gamma(x)$
lie on opposite sides of $c$ (or one of the values is equal to $c$)
then $e \in \mathcal S_c$. If $\gamma(x)$ and $\gamma(y)$ are on the
same side of $c$, then $e \in \mathcal S_c$ if and only if a string
of length $\chi(e)$ can stretch from $\gamma(x)$ to $c$ and back to
$\gamma(y)$: i.e., $|\gamma(x) - c| + |\gamma(y) - c| \leq \chi(e)$.

Now, we can extend $\gamma$ in a unique way to a continuous function
on each closed edge $e=(x,y)$ so that $\gamma$ is linear on each of
the two half segments of $e$, $\gamma$ achieves its maximum at the
midpoint $m$ of $e$, and $|\gamma(m) - \gamma(x)| + |\gamma(m) -
\gamma(y)| = \chi(e)$. Now, if either $\gamma(x) \leq c$ or
$\gamma(y) \leq c$, then we have $e \in \mathcal S_c$ if and only if
$\gamma$ assumes the value $c$ at some point along the edge $e$.

If $F_M(x) < F_M(y)$, then, by definition, $x$ and $y$ are not in
the same components of $\mathbb E^d \backslash\mathcal S_c$ where $c
\in \{ F_M(x)M + M, F_{M}(x)M + 2M, \ldots, F_{M}(y)M \}$; in
particular, this implies that $(\phi_1 + c, \phi_2, r)$ is swappable
at $e$ for $c = F_M(x)M+M$ and $c = F_{M}(y)M$.  Since $\gamma(x)
\leq F_M(x)$, this implies that $\gamma$ assumes all of the values
$F_M(x)M + M, F_{M}(x)M + 2M, \ldots, F_{M}(y)M$ at some point along
the edge $e$. Now, the map $R_c$ can be extended to the continuous
version of $\gamma$ as follows: first extend $T_c^-$ to continuous
points by letting it contain not only the points in $\mathbb Z^d$
defined to be in $T_c^-$ before but also those points $z$ on the
interior of an edge for which $\gamma(z) < c$ and there is a path
from $z$ to a point in $T_c^-$ along which $\gamma < c$ (or
equivalently, all points $z$ starting from which there exists an
infinite-length, non-self-intersecting path along which $\gamma<c$).
As before, we write $F(x)$ for the smallest value $c$ for which $x
\in T^-_c$, and $F_M(x) = \sup_{k : F(x) > kM}$. Then define
$R_c(\gamma) = \delta$ where $\delta(z) = \gamma(z)$ if $z \in
T_c^-$ and $\delta(z) = 2c-\gamma(z)$ otherwise.  We can similarly
extend the definition of $R$ to the interiors of the edges by
writing $R(\gamma) = \delta$ where
$$\delta(x) = \begin{cases} \gamma(x)  - F_M(x)M & \text{$F_M(x)$ is even} \\
(F_M(x)+1)M - \gamma(x) & \text{$F_M(x)$ is odd.} \\ \end{cases}$$

It is clear that the total variation of $\delta$ is equal to that of $\gamma$ along each edge $e$.
If $F_M(x) < F_M(y)$, and $z_{F_M(x)+1},\ldots, z_{F_M(y)}$ are the points along the edge from $x$
to $y$ at which $\gamma$ first assumes the values $F_M(x)M + M, F_M(x) + 2M, \ldots, F_M(y)M$, then
it is not hard to see that $\delta(z_i)$ will be $M$ if $i$ is even and $0$ if $i$ is odd (since
$F_M(z) = F(z) = \gamma(z)$ at these points).  From this, the lemma follows immediately.  \qed
\end{proof}

Now, if we define $\mathcal S_0$ and $\mathcal S_M$ using the triple
$(\psi_1, \psi_2, s)$, what are the infinite clusters of $\mathbb
E^d \backslash \mathcal S_0$? If $k$ is odd, then the above result
implies that each of the edges separating an element of $T^-{kM}$
(defined using the original $\phi_1, \phi_2, r$) from an element of
its complement is in $\mathcal S_0$; thus, for no odd $k$ does an
infinite cluster contain both a member of $T^-_{kM}$ and a member of
its complement. Since the $T^-_{kM}$ are nested sets, it follows
that the infinite cluster of the former must be contained in
$T^-_{kM} \backslash T^-_{(k+2)M}$ for some odd $k$. Now, if there
is one such cluster, there must almost surely be infinitely many,
since otherwise the minimum value of $k$ for which such a cluster
occurs in $(\psi_1, \psi_2, r)$ would be a height difference
variable for $\mu$ (and in this subsection, we are assuming that no
height difference variable exists).  If there are infinitely many
infinite clusters with positive probability, then Lemma
\ref{onecluster} implies that $R(\mu)$ is not a gradient Gibbs
measure.  A similar argument holds for the infinite clusters of
$\mathbb E^d \backslash \mathcal S_M$.  However, if there are no
infinite clusters in the complement of either $\mathcal S_0$ or
$\mathcal S_M$, then Lemma \ref{ergodicgradientbounds} implies that
$R(\mu)$ is not a gradient Gibbs measure. Finally, as in the
previous section, since the averages $\frac{[\phi_1(y) - \phi_1(x)]
+ [\phi_2(y) - \phi_2(x)]}{2}$ are always left unchanged by swapping
maps, we have $S_a(\mu) = S_a(R(\mu))$ and the statement of Theorem
\ref{smallerSFEaverage} follows.

\section{Height offset spectra} \label{spectrumsection}

From Theorem \ref{minimalphaseuniqueness}, we know that if $E=
\mathbb R$, then for $u \in U_{\Phi}$ and simply attractive $\Phi$,
there exists a unique minimal gradient phase $\mu_u$ on $(\Omega,
\mathcal F^{\tau})$.  In the case $E = \mathbb Z$, the theorem
implies that if there exists a rough measure $\mu_u$ of slope $u$,
it is also unique.  In fact, Theorem \ref{smallerSFEaverage} also
implies the following:

\begin{lem} If $u \in U_{\Phi}$ and either $E = \mathbb R$ or some minimal gradient phase $\mu_u$ of
slope $u$ is rough, then there is a unique minimal gradient phase
$\mu_u$ of slope $u$ and it is extremal. \end{lem} \begin{proof}
Theorem \ref{smallerSFEaverage} already gives uniqueness of $\mu_u$.
If $\mu_u$ fails to be extremal, then Lemma \ref{clustercoupling}
implies that with $\mu = \mu_u \otimes \mu_u \otimes \pi$ positive
probability, we have the strict inequality $B^+ > B^-$. If $E =
\mathbb R$ or if $E = \mathbb Z$ and $B^+ > B^- + 1$ with
$\mu$-positive probability, then Theorem \ref{smallerSFEaverage}
implies a contradiction (through the same argument as in the proof
of Theorem \ref{minimalphaseuniqueness}).  Suppose on the other hand
that $E = \mathbb Z$ and $B^+ - B^- \in \{0,1\}$ almost surely.
(Recall from Lemma \ref{clustercoupling} $B^+ \geq B^-$ almost
surely.) Then $B^+$ is a height difference variable for $\mu$ and
hence $\mu_u$ is smooth (by Lemma \ref{roughoffset}). \qed
\end{proof}

This section is devoted to the exceptional case that $E = \mathbb
Z$, $u \in U$, and every minimal gradient phase of slope $u$ is
smooth. In this case, $\mu = \mu_u$ may not be extremal, and we will
determine its extremal components.  Suppose that $\mu'$ is an
extremal component chosen from $w_{\mu}$.  Then since $\mu'$ is
($w_{\mu}$-almost surely) smooth, we can view $\mu'$ as a measure on
$\Omega$ (and we may choose the additive constant arbitrarily).
Since the additive constant is an integer, the average expected
value of $\mu'$ over any $\Lambda \subset \subset \mathbb Z^d$ ---
taken modulo $1$ --- is independent of the additive constant.

One way to extend $\mu$ to a Gibbs measure on $(\Omega, \mathcal F)$ is as follows; to sample from
$\mu$, first sample an extremal component $\mu'$ from $w_{\mu}$; then treat $\mu'$ as a measure on
$(\Omega, \mathcal F)$, adding an appropriate integer constant to cause the $\mu'$ expected height
of $\phi(0)$ to lie in $[0, 1)$.  It is then clear that $\mu(\phi(0)) \in [0,1)$; moreover, by the
definition of slope, $\mu(x) \in [(u,x), (u,x) + 1)$ for each $x \in \mathcal L$.

\begin{lem} \label{offsetdom}  If $\mu$ is minimal gradient phase with slope $u \in U_{\Phi}$,
extended as above to a measure $\mu'$ on $(\Omega, \mathcal F)$,
then for $w_{\mu} \otimes w_{\mu}$-almost all pairs of extremal
Gibbs measures $(\mu'_1, \mu'_2)$, we have either $\mu'_1 \prec
\mu'_2 \prec \mu'_1+1$ or $\mu'_2 \prec \mu'_1 \prec \mu'_2 + 1$.
\end{lem}

\begin{proof}  From Theorem \ref{smallerSFEaverage}, we have that $\mu \otimes \mu_u \otimes \pi$
almost surely $B^+ - B^- \in \{0,1\}$.  From Lemma
\ref{clustercoupling} (and Lemma \ref{extremalconvergence}) we have
that $w_{\mu} \otimes w_{\mu}$-almost surely, $\mu'_1 \prec \mu'_2+c
\prec \mu'_1 + 1$ for some value of $c \in \mathbb Z$.  But since
the expected value of $\phi(0)$ is in $[0,1)$ for $w_{\mu}$-almost
all measures, we may assume that either $c=0$ or $c=1$.  In the
former case, $\mu'_1 \prec \mu'_2 \prec \mu'_1 + 1$.  In the latter
case, $\mu_2' \prec \mu_1' \prec \mu_2' + 1$. \qed \end{proof}

\begin{thm} \label{avgHOV} Let $\mu$ be a smooth minimal phase of slope $u \in U_{\Phi}$. The
following is a height offset variable, as defined in Section
\ref{HOVsection}: $h(\mu) = \liminf |\Lambda_n|^{-1} \sum_{x \in
\Lambda_n} \phi(x)$. \end{thm}  \begin{proof}  Lemma \ref{offsetdom}
implies that for each $x \in \mathbb Z^d$, the $w_{\mu}$
distribution of the random variable $\mu'(\phi_0)$ is supported in
an interval of length one. In particular, for a point $x \in
\mathcal L$, since $\mu(\phi(x)) \in [(u,x), (u,x)+1)$, we may
conclude that $w_{\mu}$ almost surely (for $\mu'$ chosen from
$w_\mu$), $|\mu'(\phi(x)) - (u,x)| \leq 2$.

Now, taking $\Lambda_n$ to be the box of side length $2n+1$ centered
at the origin, it follows that $w_{\mu}$ almost surely, $$h(u') = -2
\leq \liminf |\Lambda_n|^{-1} \sum_{x \in \Lambda_n} \mu'(\phi(x))
\leq 2.$$  If we write $h(\phi) = h(\pi^{\phi})$ (as in Lemma
\ref{extremalconvergence}), then it is not hard to see that $h$ is a
height offset variable (as defined in Section \ref{HOVsection}).

We now claim that this $h$ is in fact equivalent to the $h$ given in
the statement of the theorem. To see this, first write $\phi_h(x) =
\phi(x) - h(\phi) - (u,x)$; note that $\phi_h$ is an $\mathcal
F^{\tau}$-measurable function.  Now, by Lemma \ref{logconcave},
$w_{\mu}$-almost surely, the $\mu'$ distribution of $\phi(0)$ is log
concave; by similar arguments to those above, we also have $w_\mu$
almost surely that $\nu' \prec \nu \prec \nu' + 1$ where $\nu$ and
$\nu'$ are the laws of $\phi(0)$ under $\mu$ and $\mu'$
respectively.  It follows that the tails of $\nu$ decay
exponentially; in particular, $\mu$ has a finite expectation at
every point in $\mathbb Z^d$, and also that $\mu (\phi_h(x))$ exists
for all $x \in \mathbb Z^d$ and has finite expectation.  By the
ergodic theorem, we have $\liminf_{n \rightarrow \infty}|\Lambda_n|
\sum_{x \in \Lambda_n} \phi_h(x) = \lim_{n\rightarrow \infty}
|\Lambda_n| \sum_{x \in \Lambda_n} \phi_h(x) = 0$, and the desired
equivalence follows. \qed \end{proof}

Applying Lemma \ref{offsetspectrum} gives the following:

\begin{cor} \label{offsetspectrumcor} The height offset spectrum $h(\mu)$ is a measure on $[0,1)$
which is ergodic with respect to the maps $x \mapsto x+u_i \pmod{1}$ for $1 \leq i \leq d$.  In
particular, if one of the components of $\mu$ is irrational, then $h(\mu)$ is uniformly
distributed.  Also, if $\mu$ is extremal---so that $h(\mu)$ is a point measure---then we must have
$u \in \tilde{\mathcal L}$, where $\tilde{\mathcal L}$ is the dual lattice of $\mathcal L$.
\end{cor}

\begin{thm} \label{HOVchar} If $\mu_1$ and $\mu_2$ are minimal gradient phases with the same slope
$u \in U_{\Phi}$ and the same height offset spectrum $\nu \in \mathcal P([0,1))$, then $\mu_1 =
\mu_2$. \end{thm}

\begin{proof} For any $\epsilon > 0$, we take $\mu_{\epsilon}$ to be $\mu_1 \times \mu_2$
conditioned on the event $A_{\epsilon}$ that the distance between $h(\phi_1)$ and $h(\phi_2)$ on
$[0,1)$ (viewed as a circle) is at most $\epsilon$; note that this event occurs with positive
probability. Since $\mu_{\epsilon}$ and $\mu_1 \times \mu_2$ conditioned on the complement of
$A^{\epsilon}$ are both $\mathcal L$-invariant gradient Gibbs measures, and since $SFE(\mu) =
2\sigma(u)$ is minimal, it follows from the affine property of $SFE$ that $SFE(\mu_{\epsilon}) = 2
\sigma(u)$. Note that the marginals of $\mu_{\epsilon}$ are $\mu_1$ and $\mu_2$.  Letting
$\epsilon$ tend to zero, by Theorem \ref{levelsetcompactness}, there is a limit point $\mu_0$ with
$SFE(\mu_0) = SFE(\mu)$ and at which $h(\phi_1) = h(\phi_2)$, $\mu_0$-almost surely.

As in the proof of Lemma \ref{avgHOV}, we note that $\mu_0$-almost surely (for appropriate choice
of additive constants), $\pi^{\phi_1} \prec \pi^{\phi_2} \prec \pi^{\phi_1} + 1$.  In particular,
this implies $h(\phi_1) \leq h(\phi_2) \leq h(\phi_1) + 1$; since $h(\phi_1)$ and $h(\phi_2)$
agree, $\mu_0$-almost surely, modulo one, we have either $h(\phi_1) = h(\phi_2)$ or $h(\phi_2) = h
(\phi_1)+1$.  Assume without loss of generality that the former is the case.  Since $\pi^{\phi_1}
\prec \pi^{\phi_2} \prec \pi^{\phi_1}+1$, an application of the ergodic theorem implies that
$\pi^{\phi_1} = \pi^{\phi_2}$ (otherwise, $h(\phi_1) \not = h(\phi_2)$ modulo one). \qed
\end{proof}

Theorem \ref{HOVchar} and Corollary \ref{offsetspectrumcor} imply the following: \begin{cor} If $u
\in U_{\Phi}$ and one of the components of $u$ is irrational, then the minimal gradient phase
$\mu_u$ of slope $u$ is unique. \end{cor}

We also have: \begin{cor}If $u \in U_{\Phi}$, all of the components
of $u$ are rational, then the smallest positive rational number
obtained as $(u,x)$, for $x \in \mathcal L$, has the form $1/n$ for
some $n\in \mathbb Z$, and each minimal gradient phase $\mu_u$ of
slope $u$ has height offset spectrum given by the uniform measure on
$\{c, c+1/n, c+2/n, \ldots, c+(n-1)/n \}$ for some $c \in [0, 1/n)$.
\end{cor} \begin{proof} The maps $g_x$, giving translation of
$[0,1)$ by $(u,x)$ modulo $1$, are elements of the group of all
rotations of $[0,1)$.  The map $x \rightarrow g_x$ is a homomorphism
from the additive group $\mathcal L$ into this abelian group.  Since
$u$ is rational, its image is a finite subgroup of the set of
rotations; letting $n$ be the order of this group, the result
follows from Corollary \ref{offsetspectrumcor} and the fact that
every measure on $[0,1)$ which is ergodic under translations by this
group is given by a uniform measure on $\{c, c+1/n, c+2/n, \ldots,
c+(n-1)/n \}$ for some $c \in [0, 1/n)$. \qed \end{proof}

Finally, we would like to describe precisely the way in which
$\mu_u$ decomposes into extremal measures---one extremal measure for
each ``height offset'' value modulo $1$.  We do this first for the
irrational case.  In this lemma, we say a function $f:\Omega \mapsto
\mathbb R$ is said to be {\it increasing} if for each $\phi_1,
\phi_2 \in \overline{\Omega}$ with $\phi_1 \leq \phi_2$, we have
$f(\phi_1) \leq f(\phi_2)$.  We say $f$ is {\it decreasing} if $-f$
is increasing.  We say an event $A \in \overline{\mathcal F}$ is
increasing (decreasing) if $1_A$ is increasing (decreasing).

\begin{thm} \label{offsetdecomposition} Suppose $u \in U_{\Phi}$, one of the components of $u$ is
irrational, and $\mu_u$ is the unique smooth minimal gradient phase
of slope $u$. Then there exists a unique family $\mu_{u,a}$ of
extremal Gibbs measures (one for each $a \in \mathbb R$) on
$(\Omega, \mathcal F)$ with all of the following properties:
\begin{enumerate}
\item {\bf Height offset property:} For each $a \in\mathbb R$, we have that for $\mu_{u,a}$ almost all $\phi$,
$h(\phi) = a$.
\item {\bf Stochastic domination property:} $\mu_{u,a_1}$ stochastically dominates $\mu_{u,a_2}$ whenever $a_1 \geq a_2$.
\item {\bf Vertical translational symmetry:} For each $b \in \mathbb Z$, we have $b+\mu_{u,a} = \mu_{u,a+b}$.
\item {\bf Decomposition property:} The restriction of $\int_{c}^{c+1} \mu_{u,a} da$ to $\mathcal F^{\tau}$ is equal
to $\mu_u$ for every $c \in \mathbb R$.
\item {\bf Right-continuity:} For each increasing event $A \subset \Omega$, $\mu_{u,a}(A)$ is increasing and
right-continuous in $a$.
\item {\bf Extremality:} For all $a \in \mathbb R$, $\mu_{u,a}$ is extremal.
\item {\bf Horizontal translational symmetry:} For each $x \in \mathcal L$ and $a \in \mathbb R$, we have
$\theta_x \mu_{u,a}=\mu_{u,a+(u,x)}$. \end{enumerate} \end{thm}

\begin{proof} First, we will prove uniqueness by showing that there is at most one definition of
the $\mu_{u,a}$ which satisfies all of the above properties.  Fix a
$c \in \mathbb R$ and extend $\mu_u$ to $(\Omega, \mathcal F)$ in
such a way that $h$ is $\mu_u$ almost surely in $[c,c+1)$.  By the
decomposition property and the height offset property, we can write
$\mu_u = \int_{a=c}^{c+1}\mu_{u,a}da$.

Now, for any $c_1,c_2$ with $c \leq c_1 < c_2 < c+1$, we write $\mu_{u,(c_1,c_2)}$ for the measure
$\mu_u$ conditioned on the positive-probability event $h \in (c_1,c_2)$; since height offset modulo
one is $\mathcal F^{\tau}$-measurable, the decomposition and height offset properties also imply
that $\mu_{u,(c_1,c_2)} = \int_{a = c_1}^{c_2} \mu_{u,a}da$.  The right continuity property implies
that for each $A$ and $a \in [c, c+1)$, we have $\mu_{u,a}(A) = \lim_{b \rightarrow 0}
\mu_{u,(a,a+b)}(A)$ (where $b\rightarrow 0$ from the right).  Since the increasing events generate
the $\sigma$-algebra $\mathcal F$, any two measures which agree on increasing events must agree on
all measurable events in $\mathcal F$; thus, if there exists a measure $\mu_{u,a}$ for which
$\mu_{u,a}(A) = \lim_{b \rightarrow 0} \mu_{u,(a,a+b)}(A)$ for all increasing $A$, that measure is
unique.

But we still have to prove that such a measure in fact exists.
First, we claim that this limit exists for every increasing set $A$
and for every $a \in [c,c+1)$. To see this, first observe that
$\mu_{u,(a_1,a_2)} \prec \mu_{u, (a_3, a_4)}$ whenever $0 \leq a_1 <
a_2 \leq a_3 < a_4 \leq 1$. To show this, it is enough to note from
Lemma \ref{clustercoupling} and Theorem \ref{smallerSFEaverage} that
$\mu_{u,(a_1,a_2)} \otimes \mu_{u,(a_3,a_4)} \otimes \pi$ almost
surely, $B^+ - B^- = 1$.  This implies that for extremal measures
$(\nu_1, \nu_2)$ chosen from the extremal decompositions of these
measures, we almost surely have $\nu_1 \prec \nu_2 +\alpha \prec
\nu_1+1$ for some value of $\alpha \in \mathbb Z$; but since the
height almost surely satisfies $h(\nu_1) < h(\nu_2) < h(\nu_1)+1$,
we may conclude that $\alpha=0$.  In fact, we can also note that
$\mu_{u,(a_1,a_2)} \prec \mu_{u, (a_1, a_3)}$ whenever $0 \leq a_1 <
a_2 \leq a_3 < 1$; this follows from the fact that
$\mu_{u,(a_1,a_3)}$ is a weighted average of $\mu_{u,(a_1,a_2)}$ and
a measure ---namely, $\mu_{u,(a_2,a_3)}$---which dominates
$\mu_{u,(a_1,a_2)}$.

This implies that for every increasing event $A \in \mathcal F$, $\mu_{u,(a,a+b)}(A)$ is a
decreasing function of $b$, and hence has a limit as $b$ tends to zero.  We would like to extend
this convergence to all measurable sets $A$. To this end, first, the reader may easily check that
the set of measures $\mu$ for which $\mu_u - 1 \prec \mu \prec \mu_u +1$ is sequentially compact in
the topology of local convergence.  Thus, for some subsequence of values of $b$ tending to zero,
the limit $\lim_{b \rightarrow 0} \mu_{u,(a,a+b)}$ exists as a measure on $(\Omega, \mathcal F)$,
which we denote by $\mu_{u,a}$.  Since the value of $\mu_{u,a}$ on increasing sets is independent
of the subsequence---when the subsequence is chosen so that $\mu_{u,a}$ is in fact a measure---then
the value of $\mu_{u,a}$ on all sets is independent of the subsequence (since increases sets
generate $\mathcal F$).

We can take the above limit as a definition for $\mu_{u,a}$ for each
$a \in [c,c+1)$; the extension to all $a$ is determined by the
vertical translation property: for $a \in \mathbb Z$, we have
$\mu_{u,a+a} = \mu_{u,a} + a$.  Now we must verify the list of
properties given above.  It is not hard to see that the above
definition is independent of the choice of $c$; the decomposition
property follows immediately.  Next, observe that $\mu_{u,(a-b,a)}
\prec \mu_{u,a} \prec \mu_{u,(a,a+b)}$ for all $a$ and $b > 0$;
letting $b$ tend to zero, the fact that $h = a$ for $\mu_{u,a}$
almost all $\phi$ follows from the definition of $h$.  The
stochastic domination property follows from the fact that $\mu_{u,a}
\prec \mu_{u,(a,b)} \prec \mu_{u,b}$ whenever $a < b$.

Now, from the stochastic domination property, it is clear that for every increasing cylinder set
$A$, $\mu_{u,a}(A)$ is increasing in $a$.  The decomposition property and height offset property
imply that $\mu_{u,(a,a+b)}(A)$ is the average of $\mu_{u,a'}(A)$ over $a \in (a,a+b)$; since
$\mu_{u,a}(A)$ is the limit of these values as $b$ tends to zero, it follows that $\mu_{u,a}(A)$ is
right-continuous in $A$.

We still need to verify extremality and the horizontal translational symmetry.  We will first check
that these properties hold for Lebesgue almost all $a$ and then use continuity arguments to extend
then to all $a$.

To see almost-sure extremality, let $\mu^b$ be the measure on triplets obtained as follows. To
sample from $\mu^b$, first choose $a$ uniformly in $[0,1)$ and then sample $\phi_1$ and $\phi_2$
independently from $\mu_{u,(a,a+b)}$.  Let $\mu$ be the limit of these measures as $b$ tends to
zero.  As in the proof of Theorem \ref{HOVchar}, it is not hard to see that this limit has minimal
specific free energy and that $B^+=B^-$ almost surely.  Using the limit definition of $\mu_{u,a}$,
it is also not hard to see that the following is an equivalent definition of $\mu$: to sample from
$\mu$, first choose $a$ uniformly from $[0,1)$ and then sample $(\phi_1, \phi_2, r)$ from
$\mu_{u,a} \otimes \mu_{u,a} \otimes \pi$.

Next, if the $\mu_{u,a}$ were not extremal for almost all $a$, then for some $\epsilon, \delta >
0$, there would be an $\epsilon$ fraction of $a$ values in $[0,1)$ for which the probability that
two extremal measures independently sampled from $w_{\mu_{u,a}}$ are different is at least
$\delta$.  But in this case, by Lemma \ref{clustercoupling}, we would have to have $B^+ \not = B^-$
with probability at least $\delta \epsilon$, a contradiction.  The horizontal translation symmetry
argument is the same, except that in this case, to sample from $\mu$, we first choose $a$ uniformly
in $[0,1)$, then choose $\phi_1$ from $\mu_{u,a}$ and $\phi_2$ from the measure $\theta_x
\mu_{u,a-(u,x)}$ (which has the same height almost surely as $\mu_{u,a}$ by the definition of
height offset variables).

Now, suppose that $\mu_{u,a}$ is not extremal; then it can be written as $p \nu_1 + (1-p) \nu_2$
for some $0 \leq p \leq 1$ and Gibbs measures $\nu_1$ and $\nu_2$ which differ on at least on
increasing event: without loss of generality, say $A$ is increasing in $\Omega$ and $\nu_1(A) <
\nu_2(A)$.  Now, write $f(\phi) = \pi^{\phi}(A)$---we can think of, $f(\phi)$ as describing the
probability of $A$ in the extremal measure from which $\phi$ was chosen.  For an decreasing
sequence of values $a_i$ lower limit is $a$, we have $\mu_{u,a_i}$ extremal, which implies that
$\pi^{\phi}(A)$ is $\mu_{u,a_i}$-almost surely constant for each $i$.  Thus, $\pi^{\phi}(A)$ is
$\mu_{u,a_i}$ almost surely equal to $\mu_{u,a_i}(A)$ for each $i$.  By right continuity, we know
that $\mu_{u,a}(A) = \lim_{i \rightarrow \infty} \mu_{u,a_i}$.  And this is in turn equal to
$\mu_{u,a}(f)$. Since $\mu_{u,a} \prec \mu_{u,a_i}$ for each $i$, the law of $f(\phi)$ when $\phi$
chosen from $\mu_{u,a_i}$ dominates the law of $f(\phi)$ when $\phi$ is chosen from $\mu_{u,a}$.
This implies that $f(\phi)$ (when $f$ is sampled from $\mu_{u,a}$) is dominated by a sequence of
constant random variables whose values converge to $\mu_{u,a}(f)$; this implies that, for
$\mu_{u,a}$ almost all $\phi$, we have $f(\phi) \leq \mu_{u,a}(f)$, which implies that $f$ is
$\mu_{u,a}$-almost surely constant, a contradiction.

The horizontal translational symmetry argument is simpler; we observe that from the almost-sure
invariance that whenever $b < c < b+1$, we have $$\int_b^c \mu_{u,a}da = \mu_{u,(b,c)} = $$
$$\theta_x \mu_{u,(b-(u,x), c-(u,x))} = \int_{b-(u,x)}^{c-(u,x)} \theta_x \mu_{u,a}da .$$ Taking
limits as $c$ approaches $b$ from above gives the result. \qed
\end{proof}

The rational case of Theorem \ref{offsetdecomposition} is straightforward and the proof is similar:
\begin{thm} \label{rationaloffsetdecomposition} Suppose $u \in U_{\Phi}$, all of the components of
$u$ are rational, and $\mu_u$ is a smooth minimal gradient phase of
slope $u$ with height spectrum given by uniform measure on $\{c,
c+1/n, \ldots, c+(n-1)/n \}$. Then there exists a unique family
$\mu_{u,a}$ of extremal Gibbs measures on $(\Omega, \mathcal F)$
(one for each $a \in c + \frac{1}{n} \mathbb Z $) with all of the
following properties: \begin{enumerate}
\item {\bf Height offset property:} For each $a \in c + \frac{1}{n} \mathbb Z$, we have that for $\mu_{u,a}$
almost all $\phi$, $h(\phi) = a$.
\item {\bf Stochastic domination property:} $\mu_{u,a_1}$ stochastically dominates $\mu_{u,a_2}$ whenever
$a_1 \geq a_2$.
\item {\bf Vertical translational symmetry:} For each $b \in \mathbb Z$, we have $b+\mu_{u,a} = \mu_{u,a+b}$.
\item {\bf Decomposition property:} The restriction of $\sum_{i=0}^{n-1} \mu_{u,a + \frac{i}{n}}$ to $\mathcal F^{\tau}$
is equal to $\mu_u$ for every $a \in c + \frac{1}{n} \mathbb Z$.
\item {\bf Extremality:} For all $a \in c + \frac{1}{n} \mathbb Z$, $\mu_{u,a}$ is extremal.
\item {\bf Horizontal translational symmetry:} For each $x \in \mathcal L$ and $a \in \mathbb R$, we have
$\theta_x \mu_{u,a}=\mu_{u,a+(u,x)}$. \end{enumerate} \end{thm}

\section{Example of slope-$u$ gradient phase multiplicity} \label{multexsection}

Here, we give an example of an LSAP $\Phi$ and a slope $u \in
U_{\Phi}$ for which the gradient phase of slope $u$ is {\it not}
unique (and a sketch of the proof that it is not unique).  For any
$x \in \mathbb Z^3$, write
$$\epsilon(x) = \begin{cases} 0 & \text{Either zero or one component of $x$ is odd.} \\
1 & \text{ Either two or three components of $x$ are odd.}\\ \end{cases}$$

Now, consider the LSAP $\Phi$ defined as follows.  When $\epsilon(x)=\epsilon(y)$, we have:
$$V_{x,y}(\eta) = \begin{cases} 0 & \eta = 0 \\ C & \eta = \pm 1 \\ \infty & \text{otherwise} \\
\end{cases}$$

When $\epsilon(x) \not = \epsilon(y)$, write
$$V_{x,y}(\eta) = \begin{cases} 0 & \eta \in \{ 0, \epsilon(y) - \epsilon(x) \} \\ \infty &
\text{otherwise} \\ \end{cases}$$ In order to observe the symmetries
of this potential better, we replace $\phi(x)$ with $\phi(x) -
\frac{\epsilon(x)}{2}$.  Now, $\phi(x)$ assumes values not in
$\mathbb Z$, necessarily, but in $\mathbb Z +
\frac{\epsilon(x)}{2}$; modify $\Phi$ accordingly.  In the modified
system, when $\epsilon(x)= \epsilon(y)$, then $\Phi_{x,y}(\phi)$ is
$0$ if $\phi(x)=\phi(y)$, $C$ if $|\phi(x) - \phi(y)|=1$ and
infinity otherwise.  But when $\epsilon(x) \not = \epsilon(y)$, we
have $\Phi_{x,y}(\phi) = 0$ if $|\phi(x) - \phi(y)| = \frac{1}{2}$
and $\infty$ otherwise.  The reader may check that $U_{\Phi}$ is the
interior of a symmetric polyhedron with the zero slope in its
interior.  Note that here, $\Phi$ restricted to the set of $x$ for
which $\epsilon(x) = 0$ (respectively, $\epsilon(x) = 1)$ is the
Ising potential on that set; the only difference is that $\phi(x)$
is allowed to assume values in $\mathbb Z$ (respectively,
$\frac{1}{2} + \mathbb Z$), instead of merely values in $\{ 0, 1
\}$, as in the Ising model.

Now, take $\mathcal L = 2 \mathbb Z^3$.  This potential has two
$\mathcal L$-invariant {\it ground states} (as defined in Definition
6.18 of \cite{G}): we have $\phi^+(x) = \frac{\epsilon(x)}{2}$ for
all $x \in \mathbb Z^d$, and $\phi^-(x) = - \frac{\epsilon(x)}{2}$
for all $x \in \mathbb Z^d$. (Each one is defined only up to an
additive constant.)  A standard argument due to Peierls (see, e.g.
Theorem 6.9, Theorem 18.25 of \cite{G}, and the surrounding
discussion) implies that for $C$ sufficiently large, there will be a
slope zero minimal gradient phase which is a small perturbation of
each of these ground states; in particular, there is more than one
slope zero minimal gradient phase.

The rough essence of these arguments is as follows: first, $SFE(\mu)
= 0$, when $\mu$ is one of the two ground states.  Thus, any $\mu$
for which $SFE(\mu)$ is minimal must satisfy $SFE(\mu) \leq 0$.
Recall that we can represent $SFE(\mu)$ as minus the entropy of
$\mu$ (i.e., $SFE_{\Phi_0}(\mu)$ where $\Phi_0$ is the potential
which is identically zero) plus the expected ``energy per site'' of
$\mu$, which we will write by $\mu(\Phi)$. The former value is
clearly at most $\log 2$ (since the number of finite-energy
configurations on an $n \times n $ box is at most $2^{n^2}$), which
implies that if $\mu$ is minimal, then we must have $\mu(\Phi) \leq
\log 2$.  This implies that the fraction of vertices whose heights
differ from those of a neighbor by $1$ is bounded above by
$\frac{\log 2}{C}$.  Define a {\bf ground state cluster} to be a
maximal connected set of vertices on which $\phi$ is equal (up to an
additive constant) to one of the two ground states. Next, one
samples a configuration $\phi$ on an $n \times n$ torus and uses
entropy considerations to argue that the probability that there
exists a cluster of size larger than $c$ containing the origin tends
to zero exponentially in $c$.  Thus, a typical configuration
contains one large ground state cluster with small ``islands''
spread throughout.  Taking a weak limit as $n$ tends to infinity,
one obtains a smooth shift-invariant measure which is equal to one
of the ground states on an infinite cluster ground state cluster
with small finite islands.  By symmetry, there exists such a measure
for each of the two ground states.  Again, the reader should consult
\cite{G} for full details.  (The proof of Theorem 6.9 in \cite{G} is
easier to read than the more general proof in Chapter 18 and
contains all of the ideas needed for the above example.)

Similar constructions to these give gradient phase multiplicity in
dimensions higher than than three.  In fact, if we relax the
condition that $\Phi$ be a nearest neighbor potential---allowing
$\Phi_{x,y}$ to be nonzero whenever either $|x-y|=1$ ($x$ and $y$
are adjacent) or $x-y \in (\pm 1, \pm 1)$ ($x$ and $y$ are
``diagonally adjacent'')---then we can construct a similar example
when $d=2$. In this case, define $\epsilon(x)$ to be $0$ or $1$ as
the sum of the coordinates of $x$ is respectively even and odd.
Exactly as in the previous example, when $\epsilon(x)= \epsilon(y)$,
then we take $\Phi_{x,y}(\phi)$ to be $0$ if $\phi(x)=\phi(y)$, $C$
if $|\phi(x) - \phi(y)|=1$ and infinity otherwise.  And when
$\epsilon(x) \not = \epsilon(y)$, we have $\Phi_{x,y}(\phi) = 0$ if
$|\phi(x) - \phi(y)| = \frac{1}{2}$ and $\infty$ otherwise.  Take
$\mathcal L = 2 \mathbb Z^2$.  As in the previous example, there are
two distinct $\mathcal L$-invariant ground states, and a Peierls
argument implies that at sufficiently low temperature, there exist
at least two distinct minimal $\mathcal L$-ergodic gradient phases,
each of which is a small perturbation of one of these ground states.

Given the simplicity of the above examples, it is perhaps surprising that---as we will see in
Chapter \ref{discretegibbschapter}---minimal gradient phase multiplicity {\it never} occurs for $u
\in U_{\Phi}$ when $d=2$ and $\Phi$ is a simply attractive potential.

\chapter{Discrete, two-dimensional gradient phases} \label{discretegibbschapter}
\section{Height offsets and main result} Throughout this chapter, we assume $m=1$ and $d=2$, $E = \mathbb Z$, and
$\Phi$ is simply attractive.  When this is the case, we can
completely classify the minimal gradient phases.  In the previous
chapter, we proved, for general $d \geq 1$, that when $E = \mathbb
R$, there is a unique minimal gradient phase of each slope $u \in
U_{\Phi}$, and this phase is extremal.  When $E = \mathbb Z$, we
found, for $d \geq 1$, that if there exists a rough minimal gradient
phase of slope $u \in U_{\Phi}$, then it is the unique minimal
gradient phase of slope $u$.  We also found, when $E = \mathbb Z$
and $d \geq 1$, that each smooth minimal gradient phase of slope $u
\in U_{\Phi}$ is completely determined by its ``height offset
spectrum''---the measure on $[0,1)$ given by $h(\mu)$ modulo $1$.
The main purpose of this chapter is to show that, in two dimensional
systems, this spectrum is trivial and the minimal gradient phase of
slope $u$ is unique:

\begin{thm} \label{2dspectrum} Suppose that $d=2$, $E = \mathbb Z$, and $\Phi$ is simply attractive.
Then for every $u \in U_{\Phi}$, there exists a unique minimal
gradient phase $\mu_u$ of slope $u$. This $\mu_u$ is extremal.  In
particular, if $\mu_u$ is a smooth phase, then its height offset
spectrum is a point mass. \end{thm}

By Lemma \ref{offsetspectrum}, the height offset spectrum of $\mu_u$ can only be a point mass if
$\mu \in \tilde{\mathcal L}$.  This implies the following corollary:

\begin{cor} Each $\mu_u$ described in Theorem \ref{2dspectrum} is a rough phase unless $u \in
\tilde{\mathcal L}$. \end{cor}

We refer to the slopes in $\tilde{\mathcal L}$ for which $\mu_u$ is
smooth as {\it smooth slopes}. In \cite{KOS}, a work by this author
and two other authors, we show how to determine---for a class of
Lipschitz simply attractive potentials $\Phi$ based on perfect
matchings of periodic bipartite graphs---precisely {\it which} of
the possible smooth slopes are smooth. Depending on $\Phi$, none,
all, or some nonempty proper subset of the vertices of
$\tilde{\mathcal L} \cap U_{\Phi}$ will be smooth slopes.  (We also
prove, in that context, a direct correspondence between the smooth
slopes and the non-strictly-convex regions---a.k.a. {\it
facets}---of certain surface-tension-minimizing surfaces subject to
boundary and volume constraints.)  In general, we do not know how to
determine explicitly {\it which} of the slopes in $\tilde{\mathcal
L} \cap U_{\Phi}$ are smooth for any other families of two
dimensional simply attractive models.

\section{FKG inequality}
Let $\mu$ be a probability density on $\Omega$ that is a finite combination of point measures
(i.e., measures supported on a single $\phi \in \Omega$).  We say that $\mu$ {\it has the $MTP_2$
(multivariate total positivity) property} if $\mu(\phi_1) \mu(\phi_2) \leq \mu(\max ( \phi_1,
\phi_2 )) \mu (\min (\phi_1, \phi_2))$ for all $\phi_1, \phi_2: \Lambda \mapsto E$.

As in the previous chapter, we say function $f:\Omega \mapsto
\mathbb R$ is said to be {\it increasing} if for each $\phi_1,
\phi_2 \in \overline{\Omega}$ with $\phi_1 \leq \phi_2$, we have
$f(\phi_1) \leq f(\phi_2)$.  We say $f$ is {\it decreasing} if $-f$
is increasing.  We say an event $A \in \overline{\mathcal F}$ is
increasing (decreasing) if $1_A$ is increasing (decreasing).  The
{\it FKG inequality} states that whenever $\mu$ has the $MTP_2$
property, any two bounded, increasing functions $f$ and $g$ on
$\Omega$ are non-negatively correlated; i.e., $\mu(fg) \geq \mu(f)
\mu(g)$.  (See the original paper by Fortuin, Kastelyn, and Ginibre
\cite{FKG}.)  In particular, this implies that
any two increasing events (or any two decreasing events) are
non-negatively correlated.  We say a general measure $\mu$ on
$(\Omega, \mathcal F)$ {\it satisfies the FKG inequality} if each
pair of increasing events in $\mathcal F$ is non-negatively
correlated.  (This implies the analogous statement about general
bounded functions $f$ and $g$, as is easily seen by approximating
$f$ and $g$ with step functions.)

We say a potential $\Phi$ is {\it submodular} if for every $\Lambda
\subset \subset Z^d$, $\Phi_{\Lambda}$ has the property that
$$\Phi_{\Lambda}(\phi_1) + \Phi_{\Lambda}(\phi_2) \geq
\Phi_{\Lambda} (\min(\phi_1,\phi_2)) +
\Phi_{\Lambda}(\max(\phi_1,\phi_2)).$$

In particular, every simply attractive potential is submodular.  Since $\Phi$ is submodular and
$Z_{\Lambda}(\phi)$ is finite, then it is clear that $\gamma^{\Phi}_{\Lambda}(\cdot, \phi)$ has the
$MTP_2$ property and hence satisfies the FKG inequality.

\begin{lem} \label{extremalFKG} If $\Phi$ is simply attractive and $\mu$ is an extremal Gibbs
measure, then $\mu$ satisfies the FKG inequality. \end{lem}
\begin{proof}  Let $A$ and $B$ be increasing events in $\mathcal F$.
By the reverse martingale theorem and the tail triviality of $\mu$,
we have $\gamma_{\Lambda_n}(A| \phi) \rightarrow \mu(A)$ for
$\mu$-almost all $\phi \in \Omega$; the same is true of increasing
events $B$ and $A \cap B$.  Since each $\gamma_{\Lambda_n}(\cdot|
\phi)$ satisfies the FKG inequality, the result follows.  \qed
\end{proof}

\begin{lem} \label{muFKG} Let $\mu_u$ be a smooth minimal gradient phase with slope $u \in
U_{\Phi}$. Extend $\mu_u$ to $(\Omega, \mathcal F)$ in such a way that $h(\phi)$ is $\mu_u$-almost
surely in $[0,1)$. Then $\mu_u$ satisfies the FKG inequality. \end{lem}

\begin{proof}  Let $A$ and $B$ be increasing events in $\mathcal F$.  By Theorems
\ref{offsetdecomposition} and \ref{rationaloffsetdecomposition}, we
can write $\mu_u = \int_0^1 \mu_{u,a}da$ (if $u$ has an irrational
component) or $\mu_u = \sum_{i=0}^{n-1} \mu_{u,c + i/n}$ (for some
$n$, if $u$ is rational).  The same theorems imply that
$\mu_{u,a}(A)$, $\mu_{u,a}(B)$, and $\mu_{u,a}(A\cap B)$ are all
increasing functions in $a$. Moreover, Lemma \ref{muFKG} implies
that $\mu_{u,a}(A \cap B) \geq \mu_{u,a}(A) \mu_{u,a}(B)$ for each
$a$.

In the irrational case, we have: $$\mu_u(A \cap B) = \int_0^1
\mu_{u,a}(A \cap B) da \geq \int_0^1 \mu_{u,a}(A) \mu_{u,a}(B).$$
Applying the FKG inequality to the increasing (in $a$) functions
$f(a) = \mu_{u,a}(A)$ and $g(a) = \mu_{u,a}(B)$ and the uniform
measure on $a \in [0,1)$, we can then say $$\int_0^1 \mu_{u,a}(A)
\mu_{u,a}(B) \geq \int_0^1 \mu_{u,a}(A) \int_0^1 \mu_{u,a}(B) =
\mu_u(A) \mu_u(B).$$ A similar argument applies in the rational
case; in this case, we replace uniform measure on $[0,1)$ with
uniform measure on $\{c,c+1/n, \ldots, c+(n-1)/n \}$. \qed
\end{proof}

We present one more straightforward fact about the FKG inequality.

\begin{lem} \label{FKGprod} If $\mu$ satisfies the FKG inequality and $\nu$ satisfies the FKG
inequality, then $\mu \otimes \nu$ satisfies the FKG inequality. \end{lem}

\begin{proof} We aim to show $$\int \int 1_A(x,y) 1_B(x,y) \nu(dy) \mu(dx) \geq \int \int 1_A(x,y)
\nu(dy) \mu(dx) \int \int 1_B(x,y) \nu(dy) \mu(dx),$$ where $A$ and
$B$ are increasing events and the integrals are over the spaces on
which $\mu$ and $\nu$ are defined.  The functions $f_A(x) = \int
1_A(x,y)\nu(dy)$, $f_B(x) = \int 1_B(x,y)\nu(dy)$, and $f_{A \cap
B}(x) = \int 1_{A \cap B}(x,y)\nu(dy)$ are clearly increasing
functions of $x$; moreover, $f_A(x) f_B(x) \leq f_{A \cap B}(x)$
pointwise by the FKG inequality for $\nu$. Combining this with the
FKG inequality for $\mu$, we have

\begin{eqnarray*} \int \int 1_{A \cap B}(x,y) \nu(dy) \mu(dx) & = & \int f_{A \cap B} \mu(dx)\\
 & \geq & \int f_A(x) f_B(x) \mu(dx) \\
 & \geq & \int f_A(x) \mu(dx) \int f_B(x) \mu(dx) \\
 & = & \int \int 1_A(x,y) \nu(dy) \mu(dx) \int \int 1_B(x,y) \nu(dy)
\mu(dx). \end{eqnarray*} \qed \end{proof}

\section{Reduction to statement about $\{0,1\}^{\mathbb Z^2}$ measures}
In this section, we let $\Omega_{\Gamma}$ be the space of subsets
$\Gamma$ of $\mathbb Z^2$; we may think of elements of
$\Omega_{\Gamma}$ as functions $\phi: \mathbb Z^2 \mapsto \{0,1\}$,
where $\phi(x) = 1$ if and only if $x \in \Gamma$.  Let $\mathcal
F_{\Gamma}$ be the usual product $\sigma$-algebra on
$\Omega_{\Gamma}$.  We will derive the main result of this chapter,
Theorem \ref{2dspectrum}, as a consequence of the following theorem:

\begin{thm} \label{nofunnymeasure} There exists no measure $\rho$ on $(\Omega_{\Gamma}, \mathcal
F_{\Gamma})$ which possesses all of the following properties: \begin{enumerate}
\item $\rho$ satisfies the FKG inequality---i.e., there exist no two increasing events $A$ and $B$ of $\Gamma$ for
which $\rho(A)\rho(B) < \rho(A \cap B)$.
\item With positive $\rho$ probability one, one of the following events occurs (the second occurring with positive
probability) \begin{enumerate} \item $\Gamma = \emptyset$ \item
$\Gamma$ and $\mathbb Z^2 \backslash \Gamma$ are both infinite
connected subsets of $\mathbb Z^2$.  (Equivalently, if $\Gamma$ is
treated in the dual sense as a subset of lattice squares of $\mathbb
Z^2$, then the boundary between $\Gamma$ and its complement consists
of a single infinite, non-self-intersecting path $P_{\Gamma}$.)
\end{enumerate}
\item The random infinite boundary path $P_{\Gamma}$ described in the previous item---conditioned on such a path
existing---has a law that is $\mathcal L$-invariant.  \end{enumerate} \end{thm}

(The reader may verify the equivalence asserted in the second statement.)

We will prove Theorem \ref{nofunnymeasure} in the next section.  In
this section, we prove that Theorem \ref{nofunnymeasure} implies
Theorem \ref{2dspectrum}.  To do this, we show that if either $\mu$
is a non-extremal minimal gradient phase of slope $u \in U_{\Phi}$
{\it or} $\mu_1$ and $\mu_2$ are distinct extremal minimal gradient
phases of slope $u \in U_{\phi}$, then we can construct a measure
$\rho$ which violates Theorem \ref{nofunnymeasure}.

In the former case, write $\overline{\mu} = \mu \otimes \mu \otimes \pi \otimes \pi$, where $\pi$
is the measure on $\Sigma$ as defined as in Section \ref{tripletsection}; here, we will view $\mu$
as being extended to $(\Omega, \mathcal F)$ in such a way that $h(\phi)$ is $\mu$-almost surely in
$[0,1)$.   By Lemma \ref{muFKG}, $\mu$ defined on $(\Omega, \mathcal F)$ in this way satisfies the
FKG inequality.  Now, given $(\phi_1, \phi_2, r_1, r_2)$
sampled from $\overline{\mu}$, we define $r$ on an edge $e=(x,y)$ by  $$r(e) = \begin{cases} r_1(e) & \text{$\phi_2(x) > \phi_1(x)$ and $\phi_2(y) > \phi_1(y)$} \\
r_2(e) & \text{otherwise} \\ \end{cases}$$  Clearly, $(\phi_1,
\phi_2, r)$ has the same law as $\mu \otimes \mu \otimes \pi$.  Our
reason for introducing $r_1$ and $r_2$ is that certain random sets
(described below) are monotone in these $r_1$ and $r_2$ --- so that
the FKG inequality applies --- even though they are not monotone in
$r$.

Namely, we claim that the indicator functions of $T^+_0$ --- defined
in terms of $(\phi_1, \phi_2, r)$, as in Section
\ref{moredefinitionssection} --- and $\mathbb Z^2 \backslash T^-_1$
are both increasing functions of the four-tuple $(\phi_2, - \phi_1,
r_2, -r_1)$.  (Recall that $T^+_0$ is an infinite set on which
$\phi_2>\phi_1$ and $T^-_1$ as a infinite set on which $\phi_2 \leq
\phi_1$.)

To see this, suppose that $(\phi_1', \phi_2', r_1', r_2')$ are such
that $$(\phi_2', - \phi_1', r_2', -r_1') \geq (\phi_2, - \phi_1, r_2,
-r_1).$$  Let $A_+$ be the set of points on which
$\phi_2(x)>\phi_1(x)$ and $A_-$ the complement---i.e., the set of
points on which $\phi_2(x) < \phi_1(x) + 1$; define $A_+'$ and
$A_-'$ analogously.  By definition, $T^+_0$ is the set of points in
those infinite components of $S_0^c$ (the complement of $S_0$) which
are contained in $A^+$. (Recall that any edge connecting a vertex in
$A^+$ to a vertex in $A^-$ necessarily belongs to $S_0$.)  Clearly,
$A_+'$ is a superset of $A_+$.  Now, we would like to show that if
$e=(x,y)$, with $x,y \in A^+$, and $e \in S_0^c$, then we also have
$e \in (S_0^c)'$. This will imply that $(T^+_0)'$ is indeed a
superset of $T^+_0$. (Here $(T^+_0)'$ and $(S_0^c)'$ are defined in
the obvious way, using $(\phi_1', \phi_2', r_1', r_2')$ instead of
$(\phi_1, \phi_2, r_1, r_2)$.)  Note first that since $-r_1 \geq
-r_1'$, we have $r'(e) \leq r(e)$.  Second, the amount of energy
required for a swap is
$$\[V_{x,y} (\phi_2(y) - \phi_1(x)) + V_{x,y}(\phi_1(y) - \phi_2(x))\] -$$ $$\[V_{x,y} (\phi_2(y) -
\phi_2(x)) + V_{x,y}(\phi_1(y) - \phi_1(x))\].$$

Since $V_{x,y}$ is convex and $\phi_2(x) > \phi_1(x)$, it is clear
that increasing $\phi_2(y)$ can only increase the value of the above
expression.  Similar observations show that increasing $\phi_2(x)$,
decreasing $\phi_1(y)$, and decreasing $\phi_1(x)$ can also only
increase the value of the above expression.  It follows that if $e
\in S_0^c$, we must also have $e \in (S_0^c)'$; and thus $T^+_0$ is
an increasing function of the four-tuple as claimed.  A similar
argument shows that the indicator of $T^-_1$ is a decreasing
function of the same four-tuple, and hence $\mathbb Z^2 \backslash
T^-_1$ is an increasing function.

Recall by Theorem \ref{offsetdecomposition} and Theorem
\ref{rationaloffsetdecomposition} that $\mu$ can be written $\mu =
\int_0^1 \mu_{u,a}da$ (if $u$ has an irrational component) or $\mu =
\sum_{i=0}^{n-1} \mu_{u,c + i/n}$ (for some $n$, if $u$ is
rational).  Thus, sampling $(\phi_1, \phi_2)$ from $\mu \otimes \mu$
is equivalent to first independently sampling $a_1$ and $a_2$ from
uniform measure on either $[0,1)$ or $\{c, c+1/n, \ldots,
c+(n-1)/n\}$, and then sampling $(\phi_1,\phi_2)$ from $(\mu_{u,a_1}
\otimes \mu_{u,a_2})$.  In either case, there is a $\mu \otimes \mu$
positive probability that $a_2 > a_1$, and thus $h(\phi_2) >
h(\phi_1)$. By Lemma \ref{offsetdom}, we have $\pi^{\phi_1} \prec
\pi^{\phi_2} \prec \pi^{\phi_1}+1$ in this case; by Lemma
\ref{clustercoupling}, this implies that $B^+ = 1$ and $B_-=0$---and
thus, both $T^+_0$ and $T^-_1$ are non-empty (conditioned on $a_2 >
a_1$). By Lemma \ref{onecluster}, each of these sets is
$\overline{\mu}$-almost surely a single infinite cluster; these
clusters are clearly disjoint.

Define $\tilde{T}^+_0$ to be the union of $T^+_-$ and all finite
components of its complement. Clearly, $\tilde{T}^+_0$ is also
increasing in the four-tuple, and both it an its complement are
connected. The same is true of $\mathbb Z^2 \backslash
\tilde{T}^-_1$, defined similarly.  When $\tilde{T}^+_0$ and
$\tilde{T}^-_1$ are disjoint, it is not hard to see that
$\tilde{T}^-_1$ and $\tilde{T}^+_0$ are also disjoint; in
particular, they are both infinite connected sets with infinite
connected complements.

An important observation is that $\tilde{T}^-_1$ and $\tilde{T}^+_0$ are gradient measurable
functions of the four-tuple: that is, given $(\phi_1,\phi_2, r_1, r_2)$ with $\phi_1$ and $\phi_2$
defined only up to additive constant, we can $\mu\otimes\mu\otimes\pi \otimes \pi$-almost surely
determine the additive constants that make $0 \leq  h(\phi_1) \leq 1)$ for $i \in \{1,2\}$, and
this determines $\tilde{T}^-_1$ and $\tilde{T}^+_0$.

A natural question to ask now is this: if we translate the four-tuple $(\phi_1,\phi_2, r_1, r_2)$
(with $\phi_1$ and $\phi_2$ defined only up to additive constant) by some $x \in \mathcal L$, does
this have the affect of translating $\tilde{T}^-_1$ and $\tilde{T}^+_0$ by $x$?  The answer is
almost yes.  Only the values of $h(\phi_1)$ and $h(\phi_2)$ modulo one are defined in the gradient
$\sigma$-algebra; when additive constants are chosen so both values are in $[0,1)$, write $\alpha =
h(\phi_2) - h(\phi_1)$; assume without loss of generality that $\alpha > 0$.  In any case, we will
have $-1 < \alpha < 1$.

Now, translating the four-tuple by $x$ has the affect of changing
$h(\phi_1)$ and $h(\phi_2)$ by $(u,x)$; now, if we add integer
constants to make $h(\phi_1) + (u,x)$ and $h(\phi_2) + (u,x)$ lie in
$[0,1)$, we will not necessarily find $h(\phi_2) - h(\phi_1) =
\alpha$; if adding $(u,x)$ does not change the integer parts of
$h(\phi_1)$ and $h(\phi_2)$ by the same amounts, we may find instead
that $h(\phi_2) - h(\phi_1) = \alpha - 1$. But regardless of $x$, we
will have either $h(\phi_2) - h(\phi_1) = \alpha$ or $h(\phi_2) -
h(\phi_1) = \alpha - 1$.  In the first case, translating the
four-tuple (defined up to additive constant) by $x$ does indeed have
the effect of translating $\tilde{T}^{-1}$ and $\tilde{T}^+_0$ by
$x$.  In the second case, it has the effect of translating
$\tilde{T}^{-1}$ and $\tilde{T}^+_0$ by $x$ and swapping the two
sets. (The reader may easily check that adding $1$ to $\phi_1$ and
subsequently the labels of $\phi_1$ and $\phi_2$ does indeed have
this effect.)

Thus, translating $(\phi_1,\phi_2,r_1,r_2)$ by $x$ also translates the unordered pair of sets
$(\tilde{T}^{-1}, \tilde{T}^+_0)$ by $x$.  Since the law of this four-tuple is $\mathcal
L$-invariant on the gradient $\sigma$-algebra, this implies that the law of the unordered pair
$(\tilde{T}^{-1}, \tilde{T}^+_0)$ is $\mathcal L$-invariant as well.

Now, we are ready to define our random set $\Gamma$.  Choose $\chi$
uniformly from the two-element set $\{0,1\}$, and choose the
four-tuple $(\phi_1, \phi_2, r_1, r_2)$ from $\overline{\mu}$.  We
now write:

$$\Gamma = \begin{cases} \tilde{T}^+_0 & h(\phi_2) > h(\phi_1), \chi = 0 \\
\mathbb Z^2 \backslash \tilde{T}^-_1 & h(\phi_2) > h(\phi_1), \chi = 1 \\
\emptyset & \text{otherwise} \\ \end{cases}$$

By the discussion above, the boundary between $\Gamma$ and its
complement has an $\mathcal L$-invariant law. Also, since
$\tilde{T}^+_0 \subset \mathbb Z^2 \backslash \tilde{T}^-_1$, the
indicator function of $\Gamma$ is an increasing function of
$(\phi_1, -\phi_2, r_1, -r_2)$ and $\chi$.  Each of the five
independent components satisfies the $FKG$ inequality, and it
follows from Lemma $\ref{FKGprod}$ that increasing events in this
five-tuple are non-negatively correlated. In particular, any two
increasing functions of $\Gamma$ are non-negatively correlated; in
other words, $\Gamma$ satisfies the FKG inequality.  Thus, the
random set $\Gamma$---which we produced using a non-extremal minimal
gradient phase $\mu_u$ of slope $u \in U_{\Phi}$---is in violation
of Theorem \ref{nofunnymeasure}. We conclude that no non-extremal
minimal gradient phase $\mu_u$ of slope $u \in U_{\Phi}$ can exist.

Now, recall that we also promised to used Theorem \ref{nofunnymeasure} to rule out the existence
distinct minimal gradient phases $\mu_1$ and $\mu_2$ of slope $u \in U_{\phi}$.  The above argument
implies that both such measures must be extremal.  Assume such measures exist and, without loss of
generality, $h(\mu_2) > h(\mu_1)$.  In this case, we can take $\overline{\mu} = \mu_1 \otimes \mu_2
\otimes \pi \otimes \pi$ and simply take $\Gamma = \tilde{T}^+_0$ as defined above.  As before, it
is clear that this set satisfies the FKG inequality; and in this case, $\Gamma = \tilde{T}^+_0$
does have a $\mathcal L$-invariant law.  Almost surely, it is connected and has a connected
complement and hence its existence contradicts Theorem \ref{2dspectrum}.  We have now proved that
Theorem \ref{nofunnymeasure} implies Theorem \ref{2dspectrum}.

\section{Proof of statement about $\{0,1\}^{\mathbb Z^2}$ measures}

\subsection{Definitions and overview of proof} \label{funnyoverview}
In this section, we prove Theorem \ref{nofunnymeasure}. Our proof is in some ways similar to the
original Russo-Seymore-Welsh proof of the non-existence of infinite clusters in critical
two-dimensional Bernoulli percolation (see, e.g., \cite{Gr} for this proof and many relevant
references).  However, their lemma and their proof relied heavily on $\mu$ having reflection
symmetries as well as translational symmetries. Because we are not assuming any reflection
symmetries, our construction will require a little bit more machinery than the analogous
construction in \cite{Gr}, including some topological results.

Suppose that $\rho$ is as described in Theorem \ref{nofunnymeasure}.
Let $A_{\infty}$ be the event that $\Gamma$ and $\mathbb Z^2
\backslash \Gamma$ are infinite connected sets; by assumption,
$\rho(A_{\infty}) > 0$.  Let $A_{\emptyset}$ be the event that
$\Gamma = \emptyset$.  Also by assumption, $\rho(A_{\infty} \cup
A_{\emptyset}) = 1$. Now, we take $\epsilon$ to be an extremely
small fixed constant: $\epsilon = 10^{-10000}
\rho(A_{\infty})^{10000}$ will comfortably suffice for all of our
arguments. Let $A_k$ be the event that $A_{\infty}$ occurs and that
both $\Gamma \cap \Lambda_k$ and $(\mathbb Z^2 \backslash A_k) \cap
\Lambda_k$ are non-empty.  (In this section, we assume for
convenience that $\Lambda_k$ is shifted to be centered at the
origin---i.e., $\Lambda_k = [\lfloor -k/2 \rfloor , \lfloor k/2-1
\rfloor]^2 \subset \mathbb Z^2$.)  Choose $k$ large enough so that
$\rho(A_k) > \rho(A_{\infty}) - \epsilon$.

Let $B_k$ be the event that there exists a path---consisting entirely of elements in $\Gamma
\backslash \Lambda_k$---which encircles the set $\Lambda_k$.  Clearly, $B_k$ is disjoint from both
$A_k$ and $A_{\emptyset}$, and hence $\rho(B_k) < \epsilon$.  We will use topological arguments and
the FKG inequality to prove that $\rho(B_k) \geq \epsilon$---hence proving Theorem
\ref{nofunnymeasure} by contradiction.

\begin{figure} \begin{center}
\leavevmode
 \epsfbox[0 0 450 250]{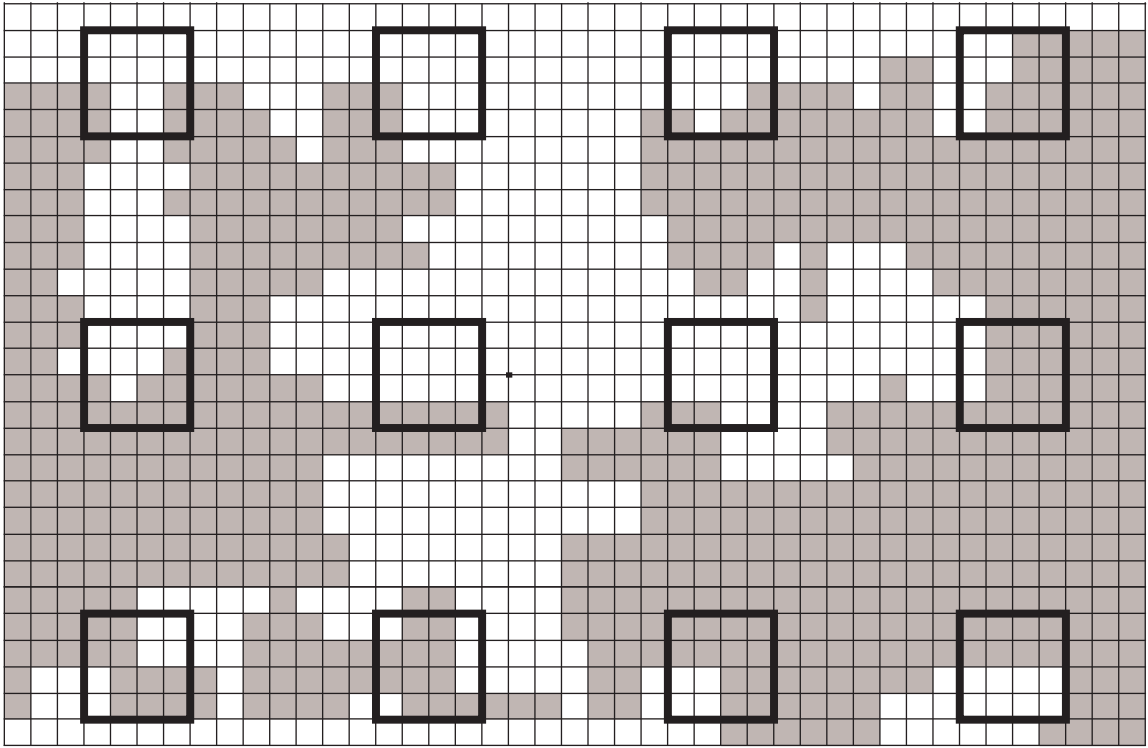} \end{center} \caption[Possible
illustration of $\Gamma$ and $\Delta$ in dual perspective] {Possible
illustration of $\Gamma$ and $\Delta$ in dual perspective when $k =
4$ and $n=11$: the shaded squares are elements of $\Gamma$, squares
inside $4 \times 4$ boxes are elements of $\Delta$, vertices on
boundaries of these boxes are members of $\tilde{\Delta}$.  The
infinite path $P_{\Gamma}$ is the boundary between $\Gamma$ and its
complement.}  \label{clusterfigure} \end{figure}

Fix an integer $n > 2k$.  Denote by $\Delta(v)$ the ``shifted box''
$nv + \Lambda_k$ and write $\Delta = \cup_{v \in \mathbb Z^2}
\Delta(v)$.  Write $\overline{\Delta}(v) = nv + \Lambda_{k+1}$. Let
$\tilde{\Delta}(v)$ be the outer band of square faces around
$\Delta(v)$---i.e., the set of square faces of $\mathbb Z^d$ that
are incident to at least one vertex of $\Delta(v)$ and at least one
vertex of $\overline{\Delta}(v) \backslash \Delta(v)$.  Also, take
$\tilde{\Delta} = \prod_{v \in\mathbb Z^d} \tilde{\Delta}(v)$.  See
Figure \ref{clusterfigure} for a dual version (i.e., squares
depicted as vertices and vice versa) of this picture.

We can think of the path $P_{\Gamma}$---which is a sequence of
square faces of $\mathbb Z^2$ (or equivalently, a sequence of
vertices in the dual graph)---as being oriented in such a way that
$\Gamma$ lies on its left.  Denote by $A(v)$ the event that
$P_{\Gamma}$ hits $\tilde{\Delta}_v$, and in between the first and
last times $P_{\Gamma}$ hits $\tilde{\Delta}(v)$, $P_{\Gamma}$ hits
no square which is fewer steps away from a $\tilde{\Delta}(w)$, with
$w \not = v$, than it is from $\tilde{\Delta}(v)$.  For example,
$A(v)$ occurs whenever $\Delta(v)$ is one of the $4 \times 4$ boxes
in Figure \ref{clusterfigure} {\it except} when $\Delta(v)$ is the
top left box.

Clearly, by choosing $n$ large enough so that the probability that
$P_{\Gamma}$ takes more than $n-k$ steps between visits to
$\tilde{\Delta}(v)$ is small, we can make the probability of
$A_{\infty} \backslash A(v)$ arbitrarily small.

We will henceforth assume that in additional to satisfying $n > 2k$,
we also have that $n \mathbb Z^2 \subset \mathcal L$ and that $n$ is
large enough so that $\rho(A(0)) > \rho(A_{\infty}) - \epsilon$. The
$\mathcal L$-invariance of the law of the path $P_{\Gamma}$ implies
that $\rho(A(v))
> \rho_k(A_{\infty}) - \epsilon$ for any $v \in \mathbb Z^2$; or equivalently, $ \rho(A_{\infty}
\backslash A(v)) < \epsilon$.  The values of $n$, $k$, and
$\epsilon$ will remain fixed throughout the proof.

Given the event $A(v)$, we define $v_+$ and $v_-$ to be such that
$\tilde{\Delta}(v_-)$ is the last band that the path $P_{\Gamma}$
hits before the first time it hits $\tilde{\Delta}(v)$. Similarly,
$\tilde{\Delta}(v_+)$ is the first band that the path $P_{\Gamma}$
hits after the last time it hits $\tilde{\Delta}(v)$.

Let $C(v,w)$ be the event that some vertex incident to
$\overline{\Delta}(v)$ and some vertex incident to
$\overline{\Delta}(w)$ are in the same connected component of
$\Gamma \backslash \Delta$.  In other words, $C(v,w)$ is the event
that there is a path in $\Gamma$ that connects vertices incident to
$\overline{\Delta}(v)$ and $\overline{\Delta}(w)$ without passing
through any other box $\Delta(x)$, $x \not \in \{v,w\}$.  Given that
$A(w)$ occurs, it is easy to see that $C(w,w_+)$ (respectively,
$C(w,w_-)$) must occur as well---we see this by taking a path
comprised of vertices that lie immediately to one side of the
portion of $P_{\Gamma}$ that connects $\tilde{\Delta}(w)$ and
$\tilde{\Delta}(w_+)$ (respectively, $\tilde{\Delta}(w_-)$).

Given a non-self-intersecting path $p = v_1, \ldots, v_r$ in $\Gamma
\backslash \Delta$ connecting $\overline{\Delta}(v)$ and
$\overline{\Delta}(w)$, we define a continuous,
non-self-intersecting path $\overline{p}$ from $nv$ to $nw$ by
connecting the dots in the sequence $nv, v_1, \ldots, v_r, nw$ with
straight line segments.  Since $v_1 \in \overline{\Delta}(v)
\backslash \Delta(v)$ and $v_r \in \overline{\Delta}(w) \backslash
\Delta(w)$ are initial and final points of $p$, it is not hard to
see the the initial and final segments of $\overline{p}$ do not
intersect any of the other segments.  Given any continuous path $q$
from $nv$ to $nw$ which is contained (except for its two endpoints)
in $\mathbb R^2 \backslash n \mathbb Z^2$, we denote by $C_q(v,w)$
the event that there exists a path $p$ in $\Gamma \backslash
\Delta$, connecting $\overline{\Delta}(v)$ and
$\overline{\Delta}(w)$ for which $\overline{p}$ is homotopically
equivalent to $q$ (i.e., there exists a continuous $Q :[0,1]^2
\rightarrow \mathbb R^2$ with $Q(t,0) = \overline{p}(t)$, $Q(t,1) =
q(t)$, $q(0,t) = v$, $q(0,t) = w$, for all $t \in [0,1]$ and $Q(t,u)
\in \mathbb R^2 \backslash \mathbb Z^2$ for all $(t,u) \in (0,1)
\times [0,1]$).  Denote by $C^+_q(v,w)$ the event that $w = v_+$ and
the path $P_{\Gamma}$ from $\overline{\Delta}(v)$ to
$\overline{\Delta}(w)$ is homotopic to $q$. Define $C^-_q(v,w)$
analogously using $v_+$ instead of $v_-$. Clearly, $C^+_q(v,w)$ and
$C^-_q(v,w)$ are both contained in the event $C_q(v,w)$.

Let $B(v)$ be the event that there exists a cycle in $\Gamma \backslash \Delta$ which disconnects
$\Delta(v)$ from infinity.  We will ultimately prove Theorem \ref{nofunnymeasure} by showing,
first, that certain combinations of events of the form $C_q(v,w)$ together imply, for topological
reasons, the event $B(0)$; we then bound the probabilities of these combinations of events using
the FKG inequality and the shift-invariance of the law of $P_{\Gamma}$.  This will enable us to
prove that $\rho(B(0)) > \epsilon$, a contradiction.

Our next step is to review some basic facts about the topology of the countably punctured plane.

\subsection{Homotopy classes of paths in countably punctured plane}

It is well known that the homotopy group of $\mathbb R^2$ minus a
discrete set of points is given by the free group generated by those
points. (See, e.g., Section 3.5 and Exercise 3.4 of \cite{Ma}.)
Instead of dealing with $\mathbb R^2 \backslash \mathbb Z^2$,
however, it will be convenient for us to deal with the {\it closed
countably punctured plane} defined as follows.  First, consider an
closed annulus of inner radius $\epsilon$ and outer radius
$\epsilon/2$. We can map this onto bijectively onto the space $D' =
D_{\epsilon} \backslash 0 \cup 0 \times S^1$, where $D_{\epsilon}$
is the closed disc of radius $\epsilon$ and $S^1$ is the unit
circle, by sending a point $(r, \theta)$ (defined in polar
coordinates) to $(2r-r_0, \theta)$, if $r \neq 0$, and to $0 \times
\theta$ if $r = 0$. We endow $D'$ with the topology that makes $D'$
and the annulus homeomorphic (when the annulus is endowed with the
standard topology induced by the Euclidean metric).

We define the {\it closed countably punctured plane} $W = \mathbb
R^2 \backslash \mathbb Z^2 \cup \mathbb Z^2 \times S^1$ analogously.
It is homeomorphic to $\mathbb R^2 \backslash [\mathbb Z^d + D]$
(where $D$ is a closed disc of any radius less than $1/2$), and it
is not hard to see that it has the same homotopy group as $\mathbb
R^2 \backslash \mathbb Z^2$. Intuitively, the closed countably
punctured plane is obtained by first poking a hole in $\mathbb R^2$
at each lattice point and then inserting and infinitesimal rivet at
that point; an advantage of that this space has over $\mathbb R^2
\backslash \mathbb Z^2$ is that (as we will see later---see Figure
\ref{taut}, in the next section) every homotopy class has a minimal
length representative.

We will now describe the homotopy group of $W$ more explicitly.  Let
$a = (a_1,a_2)$ be an arbitrary point in $\mathbb R^2 \backslash
\mathbb Z^2$ with irrational coordinates.  For each $x$ in $\mathbb
Z^2$, let $r_x$ be the closed line segment from $a$ to the point $x$
(including its limit point $x \times \arg(a-x)$).  Let $u_x$ be the
portion of the same ray (from $a$ through $x$ to infinity) which
lies between $x$ and $\infty$, together with its limit point $x
\times \arg(x-a)$.

We can describe the homotopy classes of paths in $\mathbb R^2
\backslash \mathbb Z^2$ which start and end at $a$ in the following
way. Let $\overline{x}$ be the homotopy class of a path which
follows $r_x$ from $a$ towards $x$, then makes a counterclockwise
loop around $x$, and then returns to $a$ along $r_x$. Every homotopy
class can be uniquely represented by a reduced word in elements of
the form $\overline{x}$, for $x \in \mathbb Z^2$ (i.e., a
finite-length word in the elements $\overline{x}$ and
$\overline{x}^{-1}$ in which no element $\overline{x}$ appears next
to its inverse $\overline{x}^{-1}$).

Now, let $p:[0,1] \mapsto \mathbb R^2$ be any cycle in $\mathbb R^2
\backslash \mathbb Z^2$ which starts and ends at $a$ and which has
only finitely many intersections with the rays $u_x, x \in \mathbb
Z^2$. Given $p$, we can generate the word corresponding to its
homotopy class as follows. Let $t$ vary between $0$ to $1$. Each
time $p(t)$ crosses a ray $u_x$ in a counterclockwise direction, add
$\overline{x}$ to the end of the word; each time it crosses a ray
$u_x$ in a clockwise direction, add $\overline{x}^{-1}$.  The reader
may easily verify (e.g., using induction on word length) that the
fundamental group element produced in this way describes the
homotopy class of $p$.

Through the remainder of our discussion, if $x \in \mathbb Z^2$, we
will treat $x$ as an element of $W$ by using $x$ as a shorthand for
$x \times 0$.  For any $x,y \in \mathbb Z^2$, let $P_{x,y}$ be the
set of continuous paths from $x$ to $y$ in $W$.  Let $r_x'$ denote
the linear segment from the $a$ to $x \times \arg(a-x)$, followed by
a counter-clockwise arc from $x \times \arg(a-x)$ to $x \times 0$.
Note that the map $p \mapsto (r'_x)^{-1} p r'_y$ (using the standard
concatenation definition of path multiplication; see, e.g., Chapter
2 of \cite{Ma}) induces a one-to-one correspondence between homotopy
classes of paths $p$ from $a$ to itself and homotopy classes of
paths from $x$ to $y$. (The inverse of this map is given by $p
\mapsto r'_x p (r'_y)^{-1}$.)

We can each view $p \in P_{x,y}$ as a function from $[0,1]$ to $W$,
defined up to a monotonic, continuous reparametrization of $[0,1]$;
by slight abuse of notation, we will sometimes also use $p$ to
denote the subset of $W$ contained in the image of $[0,1]$ under
this function.  If $p \in P_{x,y}$ for some $x$ and $y$, let $P_p$
be the set of paths in $P_{x,y}$ which are homotopically equivalent
to $p$ in $W$.  Let $P$ be the union of $P_{x,y}$ over all disjoint
pairs $x,y \in \mathbb Z^2$.

\subsection{Minimal length paths}
For which sets of paths $p_1,\ldots, p_k$ and points $x \in \mathbb
Z^2$ is it the case that $x$ lies in a bounded component of $\mathbb
R^2 \backslash \cup_{i=1}^k q_i$ for {\it every} $(q_1, \ldots, q_k)
\in \prod_{i=1}^k P_{p_i}$?  Roughly speaking, the answer is that
this is the case whenever the ``taut'' or ``minimal length'' paths
$$(\tilde{p}_1, \tilde{p}_2, \ldots, \tilde{p}_k) \in \prod_{i=1}^k
P_{p_i}$$ in these homotopy classes disconnect $x$ from infinity in
a certain ``strong'' sense.  The purpose of this subsection is to
make this statement precise.  Some of the results in this section
are related to algorithms in the computer science literature for
finding minimal length paths of given homotopy classes in regions
with polygonal obstacles (see, e.g., \cite{HSn} for examples and
additional references) and for testing equivalency between homotopy
classes by computing the unique minimal length representatives of
those classes.  However, we have been unable to find exactly what we
need in the literature, so we will give our own proofs of some of
the basic facts (such as the existence of paths of minimal length)
in our context.

A {\it line segment} in $W$ is an open line segment in $\mathbb R^2
\backslash \mathbb Z^2$ together with its limit points (which may be
in $\mathbb Z^2 \times S^1$).  An {\it arc} in $W$ is a closed path
in $x \times S^1$ (for some $x \in \mathbb Z^2$) which moves either
strictly clockwise or strictly counterclockwise around $S^1$.  A
{\it piecewise linear} path $p$ in $W$ is a path formed by
concatenating finitely many line segments and arcs (where no two
arcs appear sequentially in the concatenation); we will also assume
that $p$ is parametrized in such a way that it is not constant on
any interval of $[0,1]$.  It also suits our purposes to assume that
the endpoints of $p$ are points of $\mathbb Z^2$.

When $p$ is piecewise linear, we write $p'_+$ for the right
derivative and $p'_-$ for the left derivative of $p$ (so that
$p'_+(t) = p'_-(t)$ whenever $p$ is differentiable at $t$). If, in
some subinterval interval of $[0,1]$, we have $p(t) =
(x,\theta(t))$, with $x \in \mathbb Z^2$, then we write $p'(t) =
\frac{\partial}{\partial t} \theta(t) (-\sin\theta, \cos\theta)$; we
think of $p'(t)$ as a vector pointing ``around the infinitesimal
circle'' in the direction that $p$ is moving; we define $p'_+(t)$
and $p'_-(t)$ on such subintervals accordingly. If $p(0) = (x_0,
\theta_0)$, then we write $p'_-(0) = (-\cos \theta_0, -\sin
\theta_0)$, and if $p(1) = (x_1,\theta_1)$, we write $p'_+(1) =
(-\cos \theta_1, -\sin \theta_1)$.  (Informally, we think of $p$ as
``emerging from inside the hole at $x_0$'' when $t=0$ and ``turning
inwards into the hole at $y_0$'' when $t=1$.)  Now, whenever $p$ is
piecewise linear, $p'_+$ and $p'_-$ are defined throughout the
interval $[0,1]$. Unless otherwise stated, we will always assume
that a piecewise linear $p$ is parametrized according to Euclidean
length/arc length (so that $|p'_-|$ and $|p'_+|$ are both constant).
Also, we will assume that $p$ has no $U$-turns (i.e., points $t$ at
which $p'_-(t) = p'_+(t)$.)

If $p$ is piecewise linear, a {\it free corner} of $p$ is a point $t
\in [0,1]$ for which $p(t) \in \mathbb R^2 \backslash \mathbb Z^2$
and at which the path $p$ changes directions; we refer to $p(t)$ as
the {\it position} of the corner.  We also refer to each connected
component of $p^{-1}[\mathbb Z^2 \times S^1]$ as a {\it loop corner}
of $p$; the {\it length} of a loop corner is the length of the
corresponding arc; the {\it position} is the corresponding point $x
\in \mathbb Z^2$.  We say $p$ is {\it taut} if it contains no free
corners and the length of every loop corner is at least $\pi$. The
{\it length} of a piecewise linear path in $W$ is the sum of the
Euclidean lengths of its line segments.  We now prove the following:
\begin{lem} \label{tautunique} In every homotopy class $P_p \subset
P_{x,y}$ there exists exactly one taut path $\tilde{p}$.  This
$\tilde{p}$ has minimal length among all paths in $P_p$. \end{lem}
\begin{proof} Although this lemma may seem obvious intuitively, it
will take us a fair amount of space to prove it.  Suppose that $q_1$
and $q_2$ are two distinct homotopically equivalent taut paths in
$P_p \subset P_{x,y}$. Choose $a$ to be a generic point with
irrational coordinates that lies vertically below the entire paths
$q_1$ and $q_2$---i.e., the second coordinate of $a$ is less than
the smallest value $x_2$ which occurs in a point $(x_1,x_2)$ in
either of these paths; in particular, this allows us to assume that
whenever $x \in \mathbb Z^2$ and $\overline{x}$ occurs in the word
corresponding to $P_p$, $x$ is higher than (i.e., has higher second
coordinate than) the point $a$ in the plane $\mathbb R^2$.

\begin{figure} \begin{center} \leavevmode \epsfbox[0 0 360 200]{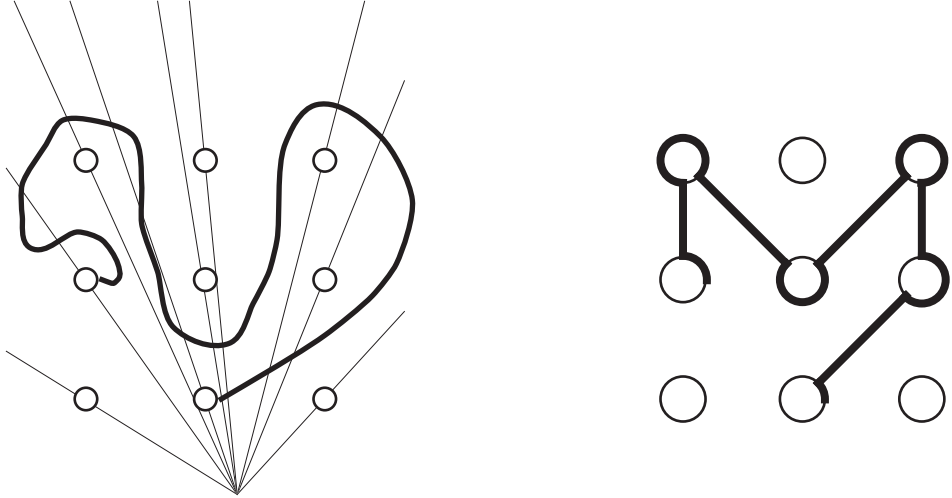} \end{center} \caption[A
path $p$ and the taut version $\tilde{p}$.] {A path $p$ (left) and the taut version $\tilde{p}$
(right). The infinitesimal ``rivets'' are of $W$ are shown as small circles on the left, slightly
larger circles on the right.  Rays $r(x)$ and $u(x)$ are depicted on the left.}  \label{taut}
\end{figure}

Now, for $q_1$, we can form a word as follows.  Let $X$ be the finite set of points $x$ for which
one of the paths $q_i$ either crosses $u_x$ at some point or contains a loop corner at position
$x$.  Order the points $x_1, x_2, \ldots, x_m$ in order of the arguments of $(x-a)$.  Let $t$ vary
between $0$ to $1$. Each time $q_1(t)$ crosses a ray $u_x$, with $x \in X$, in a counterclockwise
direction, add $\overline{x}$ to the end of the word; each time it crosses such a ray $u_x$ in a
clockwise direction, add $\overline{x}^{-1}$. Similarly, each time it cross $r_x$, with $x \in X$,
counterclockwise, add $\hat{x}$ to the end of the word; each time it crosses such an $r_x$
clockwise, add $\hat{x}^{-1}$ to the end of the word.  Denote this word by $w_1$ and the
analogously defined word for $q_2$ by $w_2$.  (Since $a$ is generic, and each of $q_1$ and $q_2$
contains only finitely many arcs and line segments, it is not hard to see that that the set of $t$
at which one of the $q_i(t)$ crosses a given $u_x$ or $r_x$ is finite, and that at each such $t$,
the path crosses the ray transversely---either clockwise or counterclockwise---so that the above
construction is well defined.)

The rays $r_{x_i} \cup u_{x_i}$ separate $W$ into wedge shaped open pieces; for $1 \leq i \leq
m-1$, denote by $W_i$ the piece between $r_{x_i} \cup u_{x_i}$ and $r_{x_{i+1}} \cup u_{x_{i+1}}$.
Denote by $W_0 = W_m$ the complement of the closures of the $W_i$ for $1 \leq i \leq m-1$.  Note
that we write down a symbol $\overline{x}_i$ or $\hat{x}_i$ each time $q_1$ passes from $W_i$ to
$W_{i+1}$ and $\overline{x}_i^{-1}$ or $\hat{x}_i^{-1}$ each time $q_i$ passes from $W_{i+1}$ to
$W_i$.  Now, we observe that if $q_1$ is taut, no two symbols of the form $\overline{x}_i$ and
$\overline{x}_i^{-1}$ (or, similarly, $\hat{x}_i$ and $\hat{x}_i^{-1}$) can occur in sequence.  If
they did, then there would be some $t_1$ and $t_2$ with $q_1((t_1,t_2))$ contained in $W_{i}$ while
$q_1(t_1)$ and $q_1(t_2)$ both belong to the same member of $\{r_{x_i}, u_{x_i}\}$. And it is easy
to see that at some point in this interval $(t_1,t_2)$ (e.g., the first point at which the argument
of $q_1$ achieves its maximum) $q_1$ must have either a loop corner at $x_{i+1}$ (with length less
than $\pi$) or a free corner. A similar argument implies that the opposite sequences
$\hat{x}_i^{-1}\hat{x}_i$ and $\overline{x}_i^{-1}\overline{x}_i$ cannot occur in either $w_1$ or
$w_2$.

Since symbols occur when $q_i$ passes from one wedge to another, it follows that whenever a
positive symbol ($\overline{x}_i$ or $\hat{x}_i$) and an inverse of a symbol (of form
$\overline{x}_i^{-1}$ or $\hat{x}_i^{-1}$) occur next to each other in the word, then they have the
same index $i$ and one is a ``hat'' and one is a ``bar'' symbol.  Also, it is easy to see that if
two positive symbols occur successively in $w_1$, and the index of the first is $i$, then the index
of the second is $i+1$.  If two negative symbols occur successively and the index of the first is
$i$, then the index of the second is $i-1$.

We know that the elements of the form $\overline{x}_i$ and $\overline{x}^{-1}$ that appear in the
word are determined by the homotopy class of $p$.  From the above paragraph, it is clear that
$\overline{x}_i$ and its inverse $\overline{x}_i^{-1}$ cannot appear in the word with only ``hat''
elements separating them; otherwise, a hat symbol would have to occur next to its inverse.  It
follows that, if the ``hat'' symbols are omitted, then the remaining ``bar'' symbols give a reduced
form expression of the fundamental group word of $p$.  Also, in between a pair of ``bar'' symbols,
the hat symbol sequence is completely determined by the rules of the above paragraph. Thus, $w_1 =
w_2$.

Now, let $f_1(i)$ be points where $q_1$ first intersects the ray
that corresponds to the $i$th element of the word.  Since $q_1$ is
taut, each $f_1(i)$ is either a point $x$ in $\mathbb Z^2$---at
which a loop corner of arc length at least $\pi$ occurs---or at a
point in the interior of $r_x$ or $u_x$.  Define $f_2$ accordingly.
It is not hard to see that each pair $f_1(i)$ and $f_1(i+1)$ must be
connected by a straight line segment and/or arc of $q_1$ (since the
pair is connected by a piecewise linear path which intersects no ray
transversely and has no corners); it follows that if $f_1$ = $f_2$,
then $q_1$ and $q_2$ must be equal.  In fact, if $f_1(i) \in \mathbb
Z^2$, then its argument is determined by the word ordering; thus,
$f_1(i)$ is determined by the value $g_1(i) = |f_1(i) - a|$.  So it
is enough for us to prove that $g_1(i) = g_2(i)$ for all $i$.
Suppose otherwise, and let $i$ be a value for which $g_1(i) -
g_2(i)$ is maximal.  Suppose that $f_1(a)$ and $f_2(a)$ lie along
the ray through $x_j$. If neither lies at the point $x_j$, then a
corner cannot occur at either one of them, and a simple geometric
argument shows that at least one of the values $g_1(i-1) - g_2(i-1)$
and $g_1(i+1) - g_2(i+1)$ is greater than $g_1(i) - g_2(i)$. On the
other hand, suppose that only $f_1(i)$ lies at the point $x_j$ and
$f_2(i)$ does not. Assume for simplicity that $f_1(i)$ and $f_2(i)$
lie on $u_x$ (the case when they lie on $r_x$ is similar). Then we
see, first, since $f_2(i)$ is part of a straight line segment of
$q_2$ from the ray through $x_{j-1}$ to the ray through $x_{j+1}$,
the values $f_2(i+1)$ and $f_2(i-1)$ will lie on this pair of rays,
as will $f_1(i+1)$ and $f_1(i-1)$ (since the words are equivalent).
It follows that the corner that occurs in $q_1$ at $f_1(i)$ has
angle between $\pi$ and $2\pi$ (if it were greater than $2\pi$, then
the path would have to intersect $r_x$ before proceeding one of the
other rays). Again, a simple geometric argument implies at least one
of the values $g_1(i-1) - g_2(i-1)$ and $g_1(i+1) - g_2(i+1)$ is
greater than $g_1(i) - g_2(i)$.  We have now proved that the class
$P_p \subset P_{x,y}$ contains at most one taut path $\tilde{p}$.

It remains to prove that $P_{x,y}$ contains {\it at least} one taut path $\tilde{p}$ and that this
path has minimal length. Given a word $w$, the above arguments determine the order in which any
taut path in the corresponding homotopy class from $x$ to $y$ would have to intersect $r_{x_i}$ and
$u_{x_i}$ (and the direction---clockwise or counterclockwise) for each such value). Let $P_w$ be
the set of paths from $x$ to $y$ which indeed intersect the rays in the given order; as seen above,
such paths are completely determined by the function $g_i$.  Since length is a continuous function
of the $g_i$, it is clear that $P_w$ contains an element $\tilde p$ for which the length is
minimal. It is not hard to see that if $\tilde p$ failed to be taut, we could decrease its length
by moving one of the $g_i$.

How do we see that $\tilde{p}$ has minimal length among all paths? Arguments similar to those above
imply that when looking for minimal length paths, we may restrict our attention to paths $p$ which
induce the same word $w$ as $\tilde{p}$ (as described above). (If this is not the case, then a
portion of the path will exit and enter one of the wedges $W_i$ along the same ray; and thus the
path can be shortened by pulling that portion taut.)  It is also not hard to see that we may assume
$p$ is piecewise linear (since otherwise, by ``straightening'' segments and arcs of $p$, we could
produce a piecewise linear $p' \in P_p$ with length less than or equal to that of $p$). The
arguments above then imply that $\tilde{p}$ has minimal length among paths of this form; hence it
has minimal length over all all paths in $P_p$. \qed \end{proof}

\begin{figure} \begin{center} \leavevmode \epsfbox[0 0 220 200]{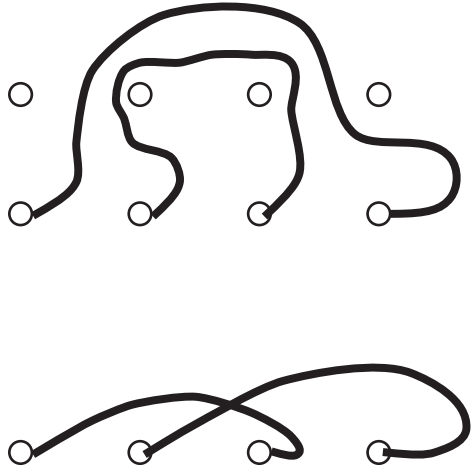} \end{center}
\caption[Essential and inessential crossings] {The taut versions of the upper two paths have points
of intersection but do not have an essential crossing; the taut versions of the lower two paths
have an essential crossing.} \label{essential} \end{figure}

If $r$ is a path, each of whose endpoints is an endpoint of either $q$ or $p$, then we say that $p$
and $q$ have a {\it crossing of type $r$} if there is a path $r' \in P_r$ for which the image of
$r'$ in $W$ lies in the union of the images of $p$ and $q$.   Now define a metric on the set of
paths in a particular homotopy class: $\delta(p,q) = \inf \[ \sup_{t \in [0,1]} |p'(t) - q'(t)|
\]$, where $|x|$ denotes the Euclidean norm of $x$ and the infimum runs over all parametrizations
$p'$ and $q'$ of the paths $p$ and $q$.   We say two paths $p_1$ and $p_2$ are {\it equivalent} if
$\delta(p_1,p_2) = 0$; so $\delta$ is actually a metric on equivalence classes, not paths.  Denote
by $B_\gamma(p)$ the ball of radius $\gamma$ about $p$ in this metric. We say that $p$ and $q$ have
an {\it essential crossing of type $r$} if for some sufficiently small $\gamma$, $p'$ and $q'$ have
a crossing of type $r$ whenever $(p',q') \in \[P_p \cap B_{\gamma}(p), P_q \cap B_{\gamma}(q)\]$.
See Figure \ref{essential}.  The usefulness of these concepts stems from the following lemma.

\begin{lem}  \label{essentialcrossing} If $p_1$ and $p_2$ are piecewise linear paths and
$\tilde{p}_1$ and $\tilde{p}_2$ have an essential crossing of type
$r$, then $p_1$ and $p_2$ also have an essential crossing of type
$r$.  \end{lem}

\begin{proof} From the definition of essential crossings, it is clear that the set $C_r$ of path
pairs $(q_1,q_2) \in P_{p_1} \times P_{p_2}$ with no essential
crossings of type $r$ is closed with respect to the product topology
generated by $\delta$ on $P_{p_1} \times P_{p_2}$.  It is also not
hard to see that the set of elements in $P_{p_1} \times P_{p_2}$
with combined length less than $L$, for some $L \in \mathbb R$, is
compact; in particular, the length function is lower
semi-continuous.  It follows that the length function achieves its
minimum over $C_r$ on some pair: we may assume this pair is $(p_1,
p_2)$.

Now, we claim that $p_1$ and $p_2$ are both taut.  Suppose
otherwise, and that without loss of generality $p_1$ is not taut.
Then there either exists an $s \in (0,1)$ for which $p_1(s)=x \in
\mathbb R^2 \backslash \mathbb Z^2$ and $p_1$ fails to be linear on
a neighborhood of $s$, or there exists an $s \in [0,1]$ for which
$p_1(s)=x \in \mathbb Z^2 \times S^1$ and the angle of the arc at
$p_1$ is less than $2\pi$.

Now, let $D$ be a small closed disc centered at $x$; we may assume
that $D$ contains no elements of $\mathbb Z^2$ (except $x$ if $x
\in\mathbb Z^2$) and that its radius is generic---so that $p_1$ and
$p_2$ each intersect the boundary of $D$ at only finitely many
points.

Let $q_1$ and $q_2$ be obtained from $p_1$ and $p_2$ by ``pulling taut'' the portions of $p_1$ and
$p_2$ inside $D$, i.e., replacing them with the minimal length piecewise linear paths with the same
endpoints on $D$ and in the same homotopy classes.  Clearly, $q_1$ and $q_2$ are shorter than $p_1$
and $p_2$.  We will be done if we can show that $q_1$ and $q_2$ (or some $q_1'$ and $q_2'$ whose
length can be made arbitrarily close to the length of $q_1$ and $q_2$) have no essential crossing
of type $r$.

\begin{figure} \begin{center} \leavevmode \epsfbox[0 0 420 200]{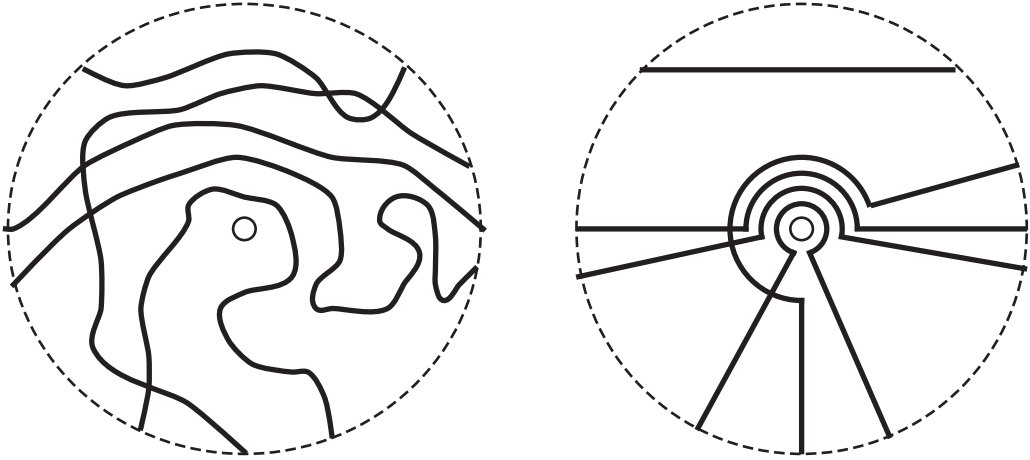} \end{center}
\caption[Pulling paths almost taut] {Five paths pulled ``almost taut.''  Each pair of paths with a
point of intersection on the right side must have an analogous point of intersection on the left
side.} \label{pulling} \end{figure}

Now, by assumption, there exist $p_1'$ and $p_2'$ arbitrarily close
to $(p_1,p_2)$ which have no crossing of type $r$.  Given such
$p_1'$ and $p_2'$, we claim that we can ``almost pull taut'' the
portions of $p_1'$ and $p_2'$ inside of $D$ in such a way that no
crossing occurs of type $r$ occurs.  First let us deal with the case
that $x$ is not a point in $\mathbb Z^2$; in this case, we simply
pull $p_1'$ and $p_2'$ taut in $D$ to produce $q_1'$ and $q_2'$.
Now observe that a pair of segments in $q_1'$ and $q_2'$ will
intersect one another if and only if the positions at which the
segments exit $D$ ordered in such a way that the endpoints of one
segment divide the circle into two pieces, one containing each of
the endpoints of the other segment (by adjusting $p_1'$ and $p_2'$
slightly if necessary, we may assume that their points of
intersection with the boundary $\partial D$ occur at distinct
locations).

The famous {\it Jordan curve theorem} states that any continuous simple closed curve in the plane
separates the plane into two disjoint regions, the inside and the outside.  A simple corollary is
that if distinct points $a,b,c,d$ are in cyclic order around a disc, and $p$ is a continuous path
in the disc from $a$ to $c$, and $q$ is a continuous path from $b$ to $d$, then $p$ and $q$ must
intersect.  It is not hard to see that if $q_1'$ and $q_2'$ have a crossing of type $r$, $p_1'$ and
$p_2'$ will have a crossing of the same type.  It follows by assumption that $q_1'$ and $q_2'$ have
no crossing of type $r$.  Finally, it is not hard to see (by choosing $p_1'$ and $p_2'$ close
enough to $p_1$ and $p_2$) that we can arrange for $q_1'$ and $q_2'$ to be arbitrarily close, in
the $\delta$ metric, to $q_1$ and $q_2$; hence, $q_1$ and $q_2$ have no essential crossing of type
$r$.

The above argument shows that if $p_1$ and $p_2$ are minimal length
paths with no essential $r$ crossing, then $p_1$ and $p_2$ must be
taut in a neighborhood of any point $x \in \mathbb R^2 \backslash
\mathbb Z^2$.   If $x \in \mathbb Z^2$, and either of $p_1$ or $p_2$
fails to be taut at $x$ (i.e., has an arc whose length is less than
$\pi$) then we can apply a similar argument.  In this case, to pull
$p_1'$ and $p_2'$ ``taut'' we replace the segments of these paths
passing through $D$ with the minimal length paths in the same
homotopy class; a path of this type will either be a straight line
segment connecting two points on $D$ or a straight line segment from
the boundary of $D$ followed by an arc of length at least $\pi$ and
another straight line segment out to the outer boundary of $D$.

Let $A$ be the set of all the arcs that occur in the paths of this form produced from the segments
of $p_i'$ and $q_i'$.   We can partially order these arcs via inclusion; it is well-known that
every partial ordering has an extension to a total ordering, and hence, we can replace each
infinitesimal arc with an arc with the same angle and small positive radius --- and choose the
radii in such a way that they are decreasing with respect to the total ordering.  See Figure
\ref{pulling}.  Adjusting $p_1'$ and $p_2'$ slightly if necessary, we may assume that no two of
these arcs share an endpoint.

It is clear that the paths defined in this way can be made to have
arbitrarily small radius. Two such paths will cross one another if
and only if the corresponding arcs $a_1$ and $a_2$ have the property
that $a_1$ is neither a subset of $a_2$ nor a subset of its
complement.  Another corollary of the Jordan Curve Theorem is that
if $\theta_1 <  \theta_2 < \theta_3 < \theta_4$ are angles and $D$
has radius $R$, and $p$ is a path in $D \backslash \{x \}$ from $(R,
\theta_1)$ to $(R, \theta_3)$ (where the angles, elements of
$\mathbb R$, may be viewed as determining points in the universal
cover of the annulus---so that the amount of winding of $p$ is given
by $\theta_3 - \theta_1$) and $q$ is an analogously defined path
from $(R, \theta_2)$ to $(R, \theta_4)$, then $p$ and $q$ must cross
(and in fact, must have a crossing of the same type).   As before,
it is not hard to see that we can arrange for $q_1'$ and $q_2'$ to
be arbitrarily close to $q_1$ and $q_2$; hence $q_1$ and $q_2$ have
no essential crossing of type $r$. \qed \end{proof}

We say that a collection of piecewise linear paths $p_1, \ldots, p_k$ {\it separates $x$ from $y$}
if there exists no path from $x$ to $y$ which does not cross at least one of the $p_i$.  We say
that $p_1, \ldots, p_k$ {\it essentially separate $x$ from $y$} if there exists no piecewise linear
path from $x$ to $y$ which does not have an essential crossing with at least one of the $p_i$.

We say that $p_1, \ldots, p_k$ {\it (essentially) bound $x$ away from infinity} if for all but
finitely many $y$, $p_1,\ldots, p_k$ (essentially) separates $x$ from $y$.

\begin{lem} \label{separation} If $\tilde{p}_1, \tilde{p}_2, \ldots, \tilde{p}_k$ essentially
separate $x$ from infinity, then $p_1, \ldots, p_k$ essentially separate $x$ from infinity.  In
particular, $p_1, \ldots, p_k$ separate $x$ from infinity. \end{lem}

\subsection{Completing proof of Theorem \ref{nofunnymeasure}} \label{completingnofunnymeasure}

Now, we return to the terminology of Section \ref{funnyoverview}.  We have already defined an
extremely small constant, $\epsilon$.  We will also need $\delta$ (a very small constant) and
$\gamma$ (a quite small constant) and $\beta$ (a small constant). The exact values are unimportant.
The following will comfortably suffice for our purposes:

$$\epsilon = 10^{-10000} \mu(A_{\infty})^{10000}$$

$$\delta = 10^{-1000} \mu(A_{\infty})^{1000}$$

$$\gamma = 10^{-100} \mu(A_{\infty})^{100}$$

$$\beta = 10^{-10} \mu(A_{\infty})^{10}$$

Recall that our aim is to prove that $\rho(B(0)) \geq \epsilon$,
thereby deriving a contradiction.

We call $w$ a {\it $\delta$-preferred direction} if the coordinates
of $w$ are relatively prime to one another and $\rho
(C_{p_v}(v,v+w)) \geq \delta$ for all $v \in \mathbb Z^2$ where
$p_v$ is a straight path connecting $v$ to $v+w$.

\begin{lem}  \label{twodelta} If there exist two distinct $\delta$-preferred directions $w_1$ and
$w_2$ (with $w_1 \not = - w_2$), then $\rho(B(0)) \geq \epsilon$. \end{lem}

\begin{proof} The FKG inequality implies that with probability at least $\delta^8$, the events
$C_{p_i}(a_i,a_{i+1})$ occur for each $1 \leq i \leq 8$, where the values $a_i$ are given by $v,
v+w, w,-v+w, -v, -v-w, -w, v-w, v$ (so $a_9 = a_1$), successively, and each $p_i$ is the straight
path from $a_i$ to $a_{i+1}$.  In this case, there are eight paths---call them $q_1, \ldots, q_i$
contained in $\Gamma \backslash \Delta$---with each $q_i$ connecting $\overline{\Delta}(a_i)$ to
$\overline{\Delta}(a_{i+1})$ in a way homotopically equivalent $p_i$.

We claim that if each of these events occurs and $A(a_i)$ occurs for
each $i$, then $B(0)$ must also occur. To see this, first, say a
vertex $b \overline{\Delta}(v) \cap \Gamma$ is {\it exposed} if
there exists a path in $\Gamma \backslash \Delta$ connecting $b$ to
some $b' \in \overline{\Delta(w)}$, for some $w \not = v$. Then
observe that if a given $A(v)$ occurs, then any two exposed points
in $\overline{\Delta}(v)$ can be connected by a ``short'' path in
$\overline{\Delta}(v)$---i.e., a path which contains only points
that are closer to $\Delta(v)$ than to any other $\Delta(w)$.  If $v
\not = 0$, then such a path is homotopically retractable in $\mathbb
R^2 \backslash \{(0,0)\}$ to a straight line (since $v \not = 0$).
By concatenating these short paths with the $q_i$'s, we produce a
cycle in $\overline{\Delta}(a_{i+1})$ in $\Gamma \backslash \Delta$
which is homotopically equivalent in $\mathbb R^2 \backslash
\{(0,0)\}$ to the concatenation of the $p_i$'s (which surrounds
$0$). It follows that $B(0)$ must occur.

Given the events $C_{p_i}(a_i,a_{i+1})$---which imply the event
$A_{\infty}$---the $A(a_i)$ fail to occur with $\rho$ probability at
most $8\epsilon$.  Thus, $\rho(B(0)) \geq \delta^8 - 8 \epsilon \geq
\epsilon$. \qed \end{proof}

\begin{lem} \label{onegamma} There exists at least one $\gamma$-preferred direction. \end{lem}

\begin{proof} We will assume that there is no $\gamma$-preferred direction, and attempt to derive a
contradiction. Now, let $u(v)$ be the angle of the direction of the
first line segment in the taut path with the same homotopy class as
the path $P_{\Gamma}$ assumes between $v$ and $v_+$ (this is
well-defined given the event $A(v)$).  Given two distinct angles
$(\theta_1, \theta_2)$, let $C_{(\theta_1, \theta_2)}(v)$ be the
event that $C_p(v,w)$ occurs for some path $p$ which, when pulled
taut, leaves $v$ in a direction which lies strictly in the interval
$(\theta_1,\theta_2)$ (on the counterclockwise side of $\theta_1$).

Now suppose that with probability at least $\beta$, $u(0)$ lies in $(0,\pi/2)$ and with probability
at least $\beta$ it lies in $(\pi/2,\pi)$; by shift-invariance of the law of $P_{\Gamma}$, this
implies the same is true of $u(v)$ for any $v$.  This also implies that for each $v$, the
increasing events $C_{(0,\pi/2)}(v)$ and $C_{(\pi/2,\pi)}(v)$ occur with probability at least
$\beta$.  Hence their union occurs with probability at least $\beta^2$.  Now, we claim that
$C_{0,\pi/2}((0,0))$ and $C_{\pi/2,\pi}((1,0))$ together imply $C_{p_0}((0,0),(1,1))$where $p_0$ is
the straight path from $(0,0)$ to $(1,0)$.  To see this, first make the geometric observation that
any taut path from $(0,0)$ with starting direction in $(0,\pi/2)$ and a taut path from $(1,0)$ with
starting direction in $(\pi/2,\pi)$ must intersect; in fact, the first linear segments of these two
paths must cross transversely at some point in the set $\{(x_1,x_2): 0 < x_1 < 1\}$, forming an
essential crossing of type $p_0$.  Lemma \ref{essentialcrossing} then implies that
$C_{p_0}((0,0),(1,1))$ follows from $C_{0,\pi/2}((0,0))$ and $C_{\pi/2,\pi}((1,0))$.  Thus, in this
case, $(1,0)$ is a $\beta^2>\gamma$-preferred direction.  Using similar arguments for $(0,1)$, we
conclude that no two consecutive quadrants can each contain $u(0)$ with probability $\beta$.

The same argument applies if we replace $(1,0)$ and $(0,1)$ with any
pair of generators for the lattice $\mathbb Z$ (i.e., any pair of
vectors $w_1$ and $w_2$, whose integer span is $\mathbb Z^2$, or
equivalently, any $w_1$ and $w_2$, each of which contains a
relatively prime pair of coordinates, for which the parallelogram
with sides determined by $w_1$ and $w_2$ has unit area).  In this
case, we cannot have $u(0)$ contained in two consecutive quadrants
of the form $(\arg(\pm v), \arg(\pm w))$ each with $\beta$
probability.  This also implies that at least one pair of opposite
quadrants has combined probability less than $2\beta$.

Now, suppose that $u(0)$ is equal to $\pi/2$ with probability $\beta$.  When this occurs, the taut
version $\tilde{p}$ of the section of $P_{\Gamma}$ between $(0,0)$ and $(0,0)_+$ will start in the
$(0,1)$ direction and pass the points $(0,1), (0,2), (0,3), \ldots$ (on either the left or right,
with arcs of length $\pi$) up to some vertex at which it either turns to the left or right or ends.
Assume without loss of generality that with probability at least $\beta/2$, $\tilde{p}$ does not
pass the first vertex on the right with angle $\pi$.  Conditioned on this event, let $i$ denote the
first $(0,i)$ in the sequence that $\tilde{p}$ does {\it not} pass on the left (with angle $\pi$ or
greater---see Figure \ref{leftpassing}). This $i$ is a random variable; let $i_0$ be its median
value (conditioned on $u(0)$ and on $\tilde{p}$ not passing the first vertex on the right with
angle $\pi$).  Now, let $C_{i_0}^+(v)$ be the event that there is some path $p$ in $\Gamma
\backslash \Delta$ between $\tilde{\Delta}(v)$ and some $\tilde{\Delta}(w)$, which, when pulled
taut starts out by moving in the $(0,1)$ direction (as before, {\it not} passing first vertex on
the right at angle $\pi$) and passes at {\it least} $i_0$ vertices on the left before it stops or
turns; let $C_{i_0}^-(v)$ be defined analogous except that the taut path passes at {\it most} $i_0$
vertices on the left before it stops or turns. Each of these events has probability at least
$\beta/4$. Now, we claim that $C_{i_0}^+((0,0))$ and $C_{i_0}^-(0,1)$ together imply
$C_{p_0}((0,0),(0,1))$, where $p_0$ is the straight path from $(0,0)$ to $(0,1)$. To see this, as
before, by Lemma \ref{essentialcrossing}, it is enough to let $p_1$ and $p_2$ be the paths whose
existence is guaranteed by the events $C_{i_0}^-((0,0))$ and $C_{i_0}^+(0,1)$ and to show that
$\tilde{p}_1$ and $\tilde{p}_2$ must have an essential crossing of type $p_0$. (We leave these
details to the reader.)  It follows that $(0,1)$ is a $(\beta/4)^2 > \gamma$-preferred direction.
We conclude that $u(0)$ cannot be equal to $(0,1)$ with probability $\beta$; a similar argument
implies that $u(0)$ cannot be equal to any single vector $v$ with probability $\beta$.

\begin{figure} \begin{center} \leavevmode \epsfbox[0 0 280 150]{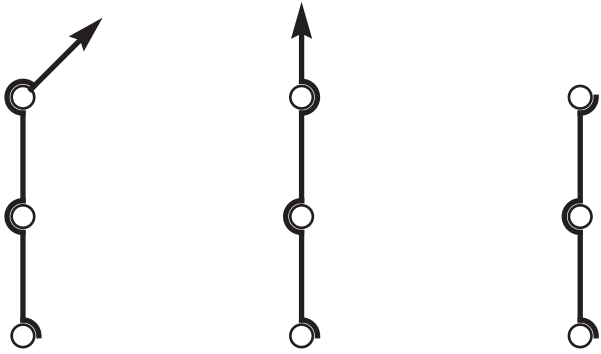} \end{center}
\caption[Passing vertices on the left.]  {If the bottom vertex is the origin, then the first path
shown passes two vertices (namely, $(0,1)$ and $(0,2)$ on the left).  The other two paths each pass
one vertex on the left.} \label{leftpassing} \end{figure}

Now, with probability $A_{\infty}$, $u(0)$ has some direction.  Thus, for any pair of generators
$(v,w)$, since each of $\pm v$ and $\pm w$ has probability less than $\beta$, and one pair of
opposite quadrants has combined probability less than $2 \beta$, the other pair of opposite
quadrants must have combined probability equal to at least $A_{\infty} - 6\beta$.  Assume without
loss of generality that $u(0)$ lies in the intervals $(\arg v, \arg w)$ and $(\arg (-v), \arg(-w))$
with probability at least $A_{\infty} - 6 \beta$.  Equivalently, $u'(0) \in (\arg v, \arg w)$ where
$u'(0)$ is $u(0)$ modulo $\pi$.

Now, we can replace $v,w$ with either the set of generators $(v,v+w)$ or $(v+w, w)$.  Since $u$ is
equal to $\pm (v+w)$ with probability at most $\beta$, we have that $u'(0)$ will lie in one of the
two intervals $(\arg v, \arg (v+w))$ and $(\arg (v+w), \arg w)$ with probability at least
$(A_{\infty} - 7 \beta)/2$.  Assume without loss of generality that this is the first interval.
Then for the set of generators $(v,v+w)$, the quadrant between $v$ and $v+w$ and its opposite must
be the highly probable ones; that is, we must have $u'(0) \in (\arg v,\arg (v+w))$ with probability
at least $A_{\infty} - 6 \beta$.

We can repeat this process, each time sending a generating pair
$(v,w)$ to either $(v,v+w)$ or $(v+w, w)$.  Now, assume (changing
bases if necessary to make this the case) that the initial pair of
vectors was $v = (1,0)$ and $w = (0,1)$.  Then with each
modification, the coordinates of $v$ and $w$ remain positive but the
sum of the coordinates increases.  It follows that the
parallelograms defined by $v$ and $w$ has its opposite corner grow
progressively longer (i.e., $v+w$ increases in norm) and skinnier
(since the parallelograms always have area one) with each iteration,
and at each step, the probability that $u(0)$ belongs to the narrow
angle defined by the parallelogram is at least $A_{\infty} - 6
\beta$.  However, this cannot be true for a nested sequence of
arbitrarily small angles, because $u(v)$ has no sufficiently large
point masses (i.e., it achieves no single value with probability
more than $\beta$), so this is a contradiction. \qed \end{proof}

To complete the proof, suppose that there exists one
$\gamma$-preferred direction $v$ and no other direction which is
$\delta$-preferred.  Assume for now that $v = (1,0)$.  Now, replace
$\Delta$ with $\Delta'$ obtained by removing every second row from
$\Delta$; that is, write $\Delta'((x,y)) =  n(x,2y) + \Lambda_k$ and
$\Delta' = \cup_{v \in \mathbb Z^2} \Delta(v)$.  Now, observe that,
if we redefine preferred directions as before but using $\Delta'$
instead of $\Delta$, then although $(0,1)$ may no longer be a
$\gamma$-preferred direction, it is still a $\gamma^2 -
\epsilon$-preferred direction (in particular, it is a
$\delta$-preferred direction). To see this, observe that if there
are paths homotopic to the straight ones in $\Gamma \backslash
\Delta$ from $\overline{\Delta}((0,0))$ to
$\overline{\Delta}((1,0))$ and from $\overline{\Delta}((1,0))$ to
$\overline{\Delta}((2,0))$, and $A((1,0))$ occurs, then there must
be a path from $\overline{\Delta}((0,0))$ to
$\overline{\Delta}((2,0))$ which is homotopic to the straight one.
It is not hard to see that Lemma \ref{twodelta} and Lemma
\ref{onegamma} still apply to the modified system produced by
replacing $\Delta$ with $\Delta'$.

Now, we repeat this process of switching $\Delta$ with $\Delta'$; if at some point $(1,0)$ ceased
to be a $\gamma$-preferred direction then, as we have observed, it would still be a
$\delta$-preferred direction.  By Lemma \ref{onegamma}, there would have to be another direction
which was $\gamma$-preferred, and by Lemma \ref{twodelta}, this would imply $\rho(B(0)) \geq
\epsilon$.

On the other hand, it is not hard to see that if we repeat this process for $m$ steps (so that
$\Delta'$ contains every $2^m$th column of $\Delta$), then as $m$ gets large, the probably that
$0_+$ lies on the vertical coordinate axis will tend to $A_{\infty}$.  By the argument used in the
proof of Lemma \ref{onegamma}, this implies that for $m$ large enough, $(1,0)$ will be a
$\delta$-preferred direction, and again, Lemma \ref{twodelta} will imply that $\rho(B(0)) \geq
\epsilon$.  If $v$ is not equal to $(1,0)$, we can change the basis for the integer lattice
indexing the boxes of $\Delta$ so that it is equal to $(1,0)$ and apply the same argument as the
one above to show that $\rho(B(0)) \geq \epsilon$.

This concludes our proof by contradiction of Theorem \ref{nofunnymeasure}.

\subsection{A corollary}

The following easy corollary of Theorem \ref{nofunnymeasure} may be of independent interest:

\begin{cor}\label{nofunnymeasuretwo} There exists no $\mathcal L$-invariant Gibbs measure $\rho$ on
the space of subsets of $\mathbb Z^2$ which possesses the following
properties: \begin{enumerate}
\item $\rho$ satisfies the FKG inequality.
\item With $\rho$ probability one, $\Gamma$ and $\Gamma^c$ each have a single infinite connected
component.\end{enumerate} \end{cor}

\begin{proof} Let $\Gamma'$ be the union of $\Gamma$ and all of the finite components of the
complement of $\Gamma$.  If $\Gamma$ is a random variable whose law satisfies the conditions of
Corollary \ref{nofunnymeasuretwo}, then $\Gamma'$ is a random variable whose law satisfies the
conditions of Theorem \ref{nofunnymeasure}. \qed \end{proof}

\section{Ergodic gradient Gibbs measures of slope $u \in \partial U_{\Phi}$} \label{2dboundary}
We have now described the $\mathcal L$-ergodic Gibbs measures with slopes in $U_{\Phi}$. The
ergodic Gibbs measures with slopes in $\partial U_{\Phi}$ are much easier to describe.  Let $X$ be
a single closed edge of the polygon $\partial U_{\Phi}$.

We say $x$ and $y$ lie in the same {\it frozen band} of $\mathbb
Z^2$ if for every $\mu \in U_{\Phi}$ with slope in $X$, the value
$\phi(y) - \phi(x)$ is almost surely constant.  By Lemma
\ref{tautpaths}, there exists an infinite sequence of parallel
``frozen bands'' (call them $b_i$) which are subsets of $\mathbb
Z^2$. We say a function $\phi$ is {\it $X$-direction taut} if the
height differences of $\phi$ on every such $b_i$ are precisely these
values.

Let $a_i$ be the first vertex on the horizontal coordinate axis of
$\mathbb Z^2$ which intersects $b_i$.  (If $a_i$ is not defined for
every $i$---because the bands $b_i$ are actually parallel to the
horizontal axis---then replace the horizonal coordinate axis with
the vertical axis in this definition.)  Let $\Delta_i$ be the set of
possible differences $\phi(a_i) - \phi(a_{i-1})$ for functions
$\phi$ which have the defined differences on the frozen bands. (Note
that $\Delta_i$ may be all of $\mathbb Z$.)  Let $k_X$ be the
smallest integer for which there exists a $v$ in $\mathcal L$ for
which $\theta_v b_0 = b_{k_X}$.

Let $\Gamma$ be the set of functions $f:\mathbb Z \mapsto \mathbb Z$
for which $f(i) \in \Delta_i$ for all $i \in \mathbb Z$.  It is now
easy to verify the following:

\begin{thm} The set of extremal Gibbs measures $\mu$ on $Z^2$ for which $\phi$ is $\mu$-almost
surely $X$-direction taut is in one-to-one correspondence with
elements of $\Gamma$.  The set of $\mathcal L$-ergodic Gibbs
measures $\mu$ on $Z^2$ with slope in $X$ is in one-to-one
correspondence with measures on $\Gamma$ with finite expectations,
and which are ergodic under translations by $k_X \mathbb Z$.
\end{thm}

\chapter{Open Problems} \label{openproblemschapter}
\section{Universality when $m=1$ and $d = 2$}
Several of the most intriguing questions about random surfaces arise in the case that $m=1$ and
$d=2$ and $\Phi$ is a simply attractive potential.

\subsection{Infinite differentiability of $\sigma$ away from roughening transition slopes} In the
case of periodically weighted domino tilings, the surface tension
$\sigma$ is infinitely differentiable away from the slopes in the
dual of $\mathcal L$ \cite{KOS}.  We conjecture that this is the
case for general discrete simply attractive potentials when $m=1$
and $d=2$.  We further conjecture that the smooth phases with slopes
in $U_{\Phi}$ occur at precisely those slopes at which $\sigma$ has
a cusp (as in the dimer model case \cite{KOS}).

\subsection{Central limit theorems: convergence to Gaussian free field} Kenyon recently
proved that with certain kinds of boundary conditions, random domino
tiling height functions have a scaling limit (when the lattice
spacing tends to zero) which is the ``massless free field," a
conformally invariant Gaussian process whose coefficients in the
eigenbasis of the Dirichlet Laplacian are independent Gaussians
\cite{Ke2}. Naddaf and Spencer proved a similar result for
Ginzburg-Landau $\nabla \phi$-interface models \cite{NS}.  We
conjecture that a similar result holds for general simply attractive
models a in rough phase.

\subsection{Level set scaling limits} If a height function $\phi$ defined on $\mathbb Z^2$ is
interpolated to a function $\overline{\phi}$ which is continuous and
piecewise linear on simplices, then the {\it level sets} $C_a$,
given by $\overline{\phi}^{-1}(a)$, for $a \in \mathbb R$, are
unions of disjoint cycles.  What do the typical ``large'' cycles
look like when $\Phi$ is simply attractive and $\phi$ is sampled
from a rough gradient phase?  The answer is given in \cite{SS} in
the simplest case of quadratic nearest neighbor potentials---in this
case, the ``scaling limit'' of the loops, as the mesh size gets
finer, is well defined, and the limiting loops look locally like a
variant of the Schramm Loewner evolution with parameter $\kappa =
4$.  We conjecture that this limit is universal---i.e., that the
level sets have the same limiting law for all simply attractive
potentials in a rough phase.

\subsection{Strong uniqueness}
We proved that when $m=1$ and $d=2$, $E = \mathbb Z^d$, and $\Phi$
is Lipschitz, simply attractive, and $\mathcal L$-ergodic, then the
$\mathcal L$-ergodic gradient Gibbs measure of slope $u \in
U_{\Phi}$ is unique.  However, we did not address the existence of
non-$\mathcal L$-ergodic gradient Gibbs measures of slope $u$.

To be precise, we say a non-$\mathcal L$-invariant gradient Gibbs measure $\mu \in \mathcal
P(\Omega, \mathcal F^{\tau})$ has {\it approximate slope $u$} if for any $\epsilon$, the
probability that $[\phi(x) - \phi(0)]- (u,x) \geq \epsilon n$ for some $x \in \Lambda_n$ tends to
zero as $n$ tends to $\infty$.

\begin{con} For each $u \in U_{\Phi}$, there exists only one gradient Gibbs measure with
approximate slope $u$. \end{con}

Such a result would be analogous to the two-dimensional Ising model result which says that there
exist no non-translation-invariant Gibbs measures (see \cite{GH} for a new proof of this fact and a
history of the problem); it is possible that techniques similar to those of \cite{GH} will be
useful in this context as well.

\section{Universality when $m=1$ and $d \geq 3$} \subsection{Central limit theorems:
convergence (after rescaling) to Gaussian free field} Using different scalings (which probably only
make sense when the spin space is continuous), Naddaff and Spencer extended their central limit
theorems to higher dimensions in the Ginzburg-Landau setting.  In what situations do similar
results hold for perturbed simply attractive potentials?

\subsection{Existence of rough phases}
We conjecture that there are no rough phases for simply attractive
potentials when $d \geq 3$, $m=1$, and $E = \mathbb R$; this is
known to be the case for Ginzburg-Landau models (see, e.g.,
\cite{FS} for more references) but it is not known in general.  We
also conjecture (although this may be riskier) that there are no
rough phases for any simply attractive potentials when $d \geq 3$
and $E = \mathbb Z$.

\section{Refinements of results proved here}
\subsection{Fully anisotropic potentials} The large deviations principle and variational principle
results in Chapters \ref{LDPempiricalmeasurechapter} and \ref{LDPchapter} were proved for isotropic
potentials and for Lipschitz potentials in the discrete setting.  We suspect that the large
deviations principle and variational principle theorems have analogs for perturbed {\it
anisotropic} (i.e., not necessarily isotropic) simply attractive potentials.  In particular:

\begin{con} The variational principle (Theorem \ref{variationalprinciple}) applies
to all perturbed simply attractive potentials. \end{con}

For the purposes of deriving the strongest possible anisotropic
general large deviations principles, it may be useful to use the
anisotropic analogs of the Orlicz-Sobolev space results we used in
Chapter \ref{orliczsobolevchapter}---several of these analogous
results are proved in \cite{C3}.  It would also be nice to have a
proof of the variational principle that does not rely as heavily on
analysis as the proof we presented here for the perturbed isotropic
case.

\subsection{Measures of infinite specific free energy}
In our version of the variational principle, we showed that for
every ergodic $\mathcal L$-invariant gradient Gibbs measure $\mu$
with slope $u \in U_{\Phi}$, $SFE(\mu) = \sigma(\mu)$ whenever
$SFE(\mu)$ is finite.  Say an $\mathcal L$-ergodic gradient Gibbs
measure $\mu$ on $(\Omega, \mathcal F^{\tau})$ is non-trivial if
$\mu$-almost surely, $Z_{\Lambda}(\phi) < \infty$ for all $\Lambda
\subset \subset \mathbb Z^d$.

\begin{con} For every $u \in U_{\phi}$, for every perturbed simply attractive potential, there
exists no non-trivial (in the sense described above) gradient Gibbs measure $\mu$ for which $S(\mu)
= u$ and $SFE(\mu) = \infty$. \end{con}

\subsection{More general domains $D$} There is a range of weaker Orlicz-Sobolev theorems that apply
when weaker regularity conditions are placed on $D$ (see, e.g., Remark 3.12 of \cite{C2}); indeed,
much of the literature on Orlicz-Sobolev spaces (see, e.g., the reference text \cite{MP}) is
focused on extending Orlicz-Sobolev bounds and embedding theorems to domains with strange
boundaries.  It is probably possible to use these more general results to prove weaker large
deviations principles for random surfaces on appropriate mesh approximations of these more general
domains.

\section{Substantially different potentials}
\subsection{Strongly repulsive lattice particles} We proved large deviations principles for classes of
perturbed simply attractive potentials $\Phi$. For what other kinds
of nearest-neighbor potentials $\Phi$ do similar large deviations
principles apply? Consider the case the $m=d=3$; in this case, we
might think of $\phi$ as describing the spatial position of atoms in
a three-dimensional solid lattice.  To take into account repulsive
forces between atoms, let $\Phi = \Phi^1 + \Phi^2$ where $\Phi^1$ is
simply attractive and $\Phi^2$ is ``strongly repulsive'' potential
given by $\Phi^2_{x,y}(\eta) = V(\eta)$ for {\it every} pair $x,y
\in \mathbb Z^d$, where $V:\mathbb R \mapsto \mathbb R$ is symmetric
and satisfies $\lim_{|\eta| \rightarrow 0} V(\eta) = \infty$ and
$\lim_{|\eta| \rightarrow \infty} |\eta|^d V(\eta) = 0$.  In this
case, as before, we will define $\sigma(u)$ to be the minimal
specific free energy among {\it ergodic} Gibbs measures $u$ of a
given slope (here $u$ is a $3 \times 3$ matrix). It is not hard to
see that $\sigma$ is symmetric and $\sigma(u) = \infty$ if and only
if $u=0$ and that $\sigma(u)$ tends to infinity as the determinant
of $u$ tends to $0$; this latter fact is a frequently imposed
condition in the study of continuous variational problems and
partial differential equations.  (See, for example, Section 9.2 of
\cite{RR}.)

Clearly, the surface tension will not be convex for models of this
form.  But are the gradient phases unique?  If not, is it possible
to classify them or to prove a large deviations principle similar to
the one proved in this text? Under some conditions, we might expect
the gradient phases to be random perturbations of periodic lattice
packings, so this problem may be related to sphere packing problems.

\subsection{Large deviations for more general tiling problems} Domino tilings are in one-to-one
correspondence with the finite-energy height functions $\phi:\mathbb
Z^2 \mapsto \mathbb Z$ for an appropriate potential $\Phi$.  Many
other classes of tilings (e.g., ribbon tilings, tilings by $1 \times
a$ and $b \times 1$ blocks, etc.) are in one-to-one correspondence
with the finite-energy height functions $\Phi: \mathbb Z^2 \mapsto
\mathbb Z^m$, where $m>1$ and $\Phi$ is an appropriately chosen
convex nearest-neighbor gradient potential.  (See \cite{S2} or
\cite{CL}.) However, none of the results in this text (variational
principle, large deviations principle, Gibbs measure classification,
etc.) is known for any non-trivial tiling problem in which the
height function space has dimension $m > 1$.

\subsection{Non-nearest-neighbor interactions}
As observed in Section \ref{multexsection}, the gradient phase
uniqueness arguments of Chapter \ref{discretegibbschapter} fail to
hold if we consider convex pair potentials which are {\it not}
nearest neighbor potentials.  What, in fact, can be said about the
gradient phases in the non-nearest neighbor case?  When $d=2$ and $E
= \mathbb Z$, is there a convex, finite range gradient potential for
which the smooth gradient phases have arbitrarily many slopes (not
merely slopes in the dual of $\mathcal L$)?  Is there a convex,
infinite range potential for which there is a smooth gradient phase
of every rational slope?

\subsection{General submodular potentials}
It is natural to wonder whether the cluster-swapping arguments used
in this text really only apply for pair potentials.

We say a potential $\Phi$ is {\it submodular} if for every $\Lambda \subset \subset Z^d$,
$\Phi_{\Lambda}$ has the property that $\Phi_{\Lambda}(\phi_1) + \Phi_{\Lambda}(\phi_2) \geq
\Phi_{\Lambda} (\min(\phi_1,\phi_2)) + \Phi_{\Lambda}(\max(\phi_1,\phi_2))$.

When $\Phi$ is a gradient pair potential, the property of submodularity is equivalent to the
convexity of the potential functions $V_{x,y}$.  Is there some variant of the cluster swapping
argument that applies to general finite-range, submodular potentials?

\appendix

\chapter{$SFE$ and the lexicographic past} \label{SFEdecompositionchapter} We present here an
alternative proof of Theorem \ref{SFEisanexpectation} which also
leads to a representation of the specific free energy in terms of
the entropy of $\mu$ on one period of $\mathcal L$, conditioned on
the lexicographic past. This approach is analogous to the approach
used in Chapter 15 of \cite{G} for non-gradient measures.  We
restate the theorem here for convenience.  \begin{thm} The function
$\alpha:\Omega \mapsto \mathbb R$, defined by $\alpha(\phi) =
SFE(\pi_{\mathcal L}^{\phi})$, is $\mathcal T \cap \mathcal
I_{\mathcal L}$-measurable, is bounded below, and satisfies
$$SFE(\mu) = \mu(\alpha) = \mu(SFE(\pi_{\mathcal L}^*))$$ for all
$\mu \in \mathcal P_{\mathcal L}(\Omega, \mathcal F)$. \end{thm}

\begin{cor} If $\mu$ can be written $$\mu = \int_{\ex\mathcal P_{\mathcal L}(\Omega, \mathcal F)}
\nu w_{\mu} (d\nu)$$ then $$SFE(\mu) = \int_{\text{ex}\mathcal
P_{\mathcal L}(\Omega, \mathcal F^{\tau})} SFE(\nu) w_{\mu} (d\nu) =
w_{\mu}(SFE).$$ \end{cor}

We would like to prove this by citing an analogous result proved for
non-gradient measures.  First, we need some notation. Whenever $\mu
\in \mathcal P(\Omega, \mathcal F)$ and $\Delta \subset Z^d$, write
$\mu\lambda_{\Delta} \in \mathcal P(\Omega, \mathcal F)$ to mean the
independent product of $\mu_{\mathbb Z^d \backslash \Delta}$ (to
determine $\phi(x)$ when $x \not \in \Delta$) and
$\lambda^{|\Delta|}$ (to determine $\phi(x)$ for $x \in \Delta$).
(This is a finite measure if $\lambda$ is a finite measure.)
Sometimes we will replace $\lambda^{\Delta}$ with
$f\lambda^{\Delta}$ where $f$ is an $\mathcal F_{\Delta}$-measurable
function, but the definition is the same.  Next, given $x,y \in
\mathbb Z^d$, write $x \prec y$ if $x$ precedes $y$ in the
lexicographic ordering of $\mathbb Z^d$ (i.e., $x_j < y_j$ where $j
= \inf \{ i | x_i \neq y_i \}$). For each $y \in \mathbb Z^d$, write
$\Gamma(y) = \{ y \} \cup \{ x \in \mathbb Z^d | x \prec y \}$ and
$\Gamma^*(y) = \{ x \in \mathbb Z^d | x \prec y \}$.  Now,
Proposition 15.16 of \cite{G} states the following.  (Since we
intend apply this result to an ``alternate'' measure space
$(\overline{\Omega}, \overline{\mathcal F})$, and not the space
$(\Omega, \mathcal F)$ we have been using throughout this text, we
will use bars over these and other variables in the theorem
statement to avoid confusion with our analogously defined ``global''
variables.)

\begin{lem} \label{specificentropylexicographic} Let $\overline{\mu}$ be a shift-invariant measure
on $(\overline{\Omega}, \overline{\mathcal F})$ where
$\overline{\Omega} = \overline{E}^{\mathbb Z^d}$, $(\overline{E},
\overline{\mathcal{E}})$ is a compact standard Borel space with a
finite underlying measure $\overline{\lambda}$, and
$\overline{\mathcal F}, \overline{\mathcal I}, \overline{\mathcal
T}$ are the corresponding Borel product, shift-invariant, and tail
$\sigma$-algebras.  Define $h(\overline{\mu}) = \lim_{n \rightarrow
\infty} |\Lambda_n|^{-1} \mathcal H(\overline{\mu} , \overline{\mu}
\overline{\lambda}_{|\Lambda_n|})$. Then $h$ is well-defined and
$\overline{\mathcal T} \cap \overline{\mathcal I}$-measurable.
Moreover, for each $y \in \mathbb Z^d$, $$h(\overline{\mu}) =
\mathcal H_{\mathcal
\overline{F}_{\Gamma(y)}}(\overline{\mu},\overline{\mu}\overline{\lambda}_{\{y\}}).$$
\end{lem}

In the setting of Lemma \ref{specificentropylexicographic}, $h(\overline{\mu})$ is called the {\it
specific entropy} of $\overline{\mu}$.  In less formal terms, the lemma states that
$h(\overline{\mu})$ is equal to the $\overline{\mu}$-expected conditional entropy of the random
variable $\phi(y)$ (with respect to $\overline{\lambda}$) {\it given} the values $\phi(x)$ for $\{
x| x \prec y \}$.  The next lemma, Proposition 15.20 of \cite{G}, is analogous to Theorem
\ref{SFEisanexpectation}.

\begin{lem} \label{specificentropyisanexpectation} If, in the setting of Lemma
\ref{specificentropylexicographic}, $\overline{\mu}$ has a representation $$\overline{\mu} =
\int_{\ex\mathcal P_{\mathcal L}(\overline{\Omega}, \overline{\mathcal F})} \overline{\nu}
w_{\overline{\mu}} (d\overline{\nu})$$ where $\overline{\mathcal F}$ is the Borel $\sigma$-algebra
of the product topology on $\overline{\Omega}$, then
 $$h(\overline{\mu}) = \overline{\mu}(h(\pi^*_{\mathcal L})) = \int_{\ex \mathcal P_{\mathcal L}(\overline{\Omega},
 \overline{\mathcal F}^{\tau})} h(\overline{\nu}) w_{\overline{\mu}}(d\overline{\nu}) = w_{\overline{\mu}}(h).$$ \end{lem}

We will now derive Theorem \ref{SFEisanexpectation} from these
lemmas.  Our first step will be to use $\mu$ to construct a new
measure space and a related measure $\overline{\mu}$ to which the
above lemmas apply.  Let $\Lambda$ be a fundamental domain of
$\mathcal L$. Though our argument applies to any sublattice
$\mathcal L$ of $\mathbb Z^d$, we will assume, for notational
simplicity, that $\mathcal L = k \mathbb Z^d$ for some integer $k
\geq 1$, so that $\Lambda = [0,k-1]^d$.  Let $e_1,\ldots, e_d$ be
the standard basis vectors for $\mathbb Z^d$.  Define $\overline E =
E^{|\Lambda|}$.  Given $\phi \in \Omega $, we define
$\overline{\phi} \in \overline{\Omega}$ by writing
$$\overline{\phi}(x) = (\phi(kx + a_1) - \phi(kx + a_0), \ldots,
\phi(kx + a_{|\Lambda|}) - \phi(kx + a_0))$$ where $\{a_i \}$ is an
enumeration of the points in $\Lambda \cup \{-e_d \}$, with $a_0 =
-e_d$. Note that if $k=1$, then $\overline{\phi}$ is simply a
discrete derivative of $\phi$ in the $e_d$ direction. We define a
measure $\overline{\lambda} = \exp \[H^o_{\Lambda \cup \{-e_d
\}}\]\lambda^{|\Lambda|}$. If $\Phi$ is a perturbed simply
attractive model, then $\overline{\lambda}$ is easily seen to be a
finite measure; we can add $\infty$ to $\overline{E}$ to make it a
compact standard Borel space (by the definition of \cite{G}) with
finite underlying measure.  For later use, for $1 \leq j \leq d$,
define $\overline{\lambda}^j$ analogously to $\overline{\lambda}$ by
replacing $e_d$ with $e_j$. Finally, let $\overline{\mu}$ be the law
on $(\overline{\Omega}, \overline{\mathcal F})$ induced by $\mu$ and
the map $\phi \mapsto \overline{\phi}$.  Now $\overline{\mu}$ and
$(\overline{\Omega}, \overline{\mathcal F})$ satisfy the conditions
in Lemma \ref{specificentropylexicographic} and Lemma
\ref{SFEisanexpectation}. Throughout this section, if $\Delta_1
\subset \mathbb Z^d$ and $\Delta_2 \subset \mathbb Z^d$, we use the
notation $\Delta_1 + \Delta_2 = \{x+y|x \in \Delta_1, y \in \Delta_2
\}$.  We define a modified potential as follows: $$\Phi'_{\Delta} =
\begin{cases} 0 & \text{$\Delta \subset \[\Lambda \cup \{-e_d \}\] +
x$ for some $x \in \mathcal
L$} \\
\Phi_{\Delta} & \text{ otherwise } \end{cases}$$

We now claim that $$SFE(\mu) = h(\overline{\mu}) + \mu(\Phi'),$$
where $\mu(\Phi')$ is the {\it specific energy} of $\mu$ with
respect to $\Phi'$, defined by $$\mu(\Phi') = \lim_{n \rightarrow
\infty} |\Lambda_n|^{-1} \mu H^{\Phi'}_{\Lambda_n}.$$ Note that
since both $\mu$ and $\Phi'$ are $\mathcal L$-invariant and $\Phi'$
has finite range, it is also possible to write $\mu(\Phi')$ as the
$\mu$ expectation of a single positive cylinder function (see
Chapter 15 of \cite{G}).  We claim that if we can show that
$SFE(\mu) = h(\overline{\mu})+\mu(\Phi')$, then this fact, along
with Lemmas \ref{specificentropylexicographic} and
\ref{SFEisanexpectation}, will imply Theorem
\ref{SFEisanexpectation}. To see this, observe that by Lemma
\ref{decompositionergodic} and the definition of ergodic
compositions at the beginning of this chapter, it follows that
$\mu(F) = \mu(\pi_{\mathcal L}^*(F)) = \int \nu(F) w_{\mu}(d\nu)$
whenever $F$ is the indicator function of a cylinder event; since a
positive cylinder function can be written as a positive linear sum
of countably many indicator functions of cylinder events, it easily
follows that $\mu(\Phi') = \mu(\pi_{\mathcal L}^*(\Phi')) = \int
\nu(\Phi') w_{\mu}(d\nu)$.

It also follows trivially from definitions that $\mu = \int \nu
w_{\mu}(d\nu)$ implies $\overline{\mu} = \int \overline{\nu}
w_{\mu}(d\nu)$.  Thus, $$SFE(\mu) = h(\overline{\mu}) + \mu(\Phi') =
\int h(\overline{\nu}) w_{\mu}(d\nu) + \int \nu(\Phi') w_{\mu}(d\nu)
= $$ $$\int SFE(\nu) w_{\mu}(d\nu) = w_{\mu}(SFE).$$  By Lemma
\ref{specificentropylexicographic}, $h(\overline{\nu})$ is $\mathcal
T \cap \mathcal I_{\mathcal L}$-measurable as a function of $\nu \in
\mathcal P_{\mathcal L}(\Omega, \mathcal F)$; in particular, it is
$\mathcal F$-measurable, and the same is trivially true of
$\nu(\Phi')$.  Thus by Lemma \ref{decompositionergodic} and the
definition of ergodic decompositions, we also have the rest of the
theorem statement: $$SFE(\mu) = \mu(SFE(\pi_{\mathcal L}^*)).$$ It
now remains only to prove that $SFE(\mu) = h(\overline{\mu}) +
\mu(\Phi')$.  We do this by checking two equalities, which we state
separately as lemmas.  First, define $\mu \overline{\lambda}_y^j \in
\mathcal P(\Omega, \mathcal F)$ as follows: sampling $\phi$ from
that measure is equivalent to first sampling $\phi$ from $\mu$ and
then re-sampling $\phi(x)$ for $x \in \Lambda + kj$ in such a way
that the values $\{\phi(kj + a_i) - \phi(kj - e_j)|1\leq i \leq
|\Lambda| \}$ obey the law of $\overline{\lambda}^j$ (where
$\overline{\lambda}^j$ is as defined above).

\begin{lem} For each $y \in \mathbb Z^d$, if $SFE(\mu)<\infty$, then $$SFE(\mu) = \mathcal H_{\mathcal
F^{\tau}_{k\Gamma(y)+\Lambda}}(\mu , \mu \overline{\lambda}^d_y)+
\mu(\Phi').$$ Moreover, if $SFE(\mu) = \infty$, then we still have
$$SFE(\mu) = \mathcal H_{\mathcal F^{\tau}_{k\Gamma(y)+\Lambda}}(\mu
, \mu \overline{\lambda}^j_y)+ \mu(\Phi')$$ for some $1 \leq j \leq
d$ (in which case we may relabel the coordinates axes so that
$j=d$). \end{lem}

\begin{proof} Let $\Phi'' = \Phi - \Phi'$.  Write $H^o_{\Delta} = H^1_{\Delta} + H^2_{\Delta}$
where the latter two terms are the components of $H^o_{\Delta}$
coming from $\Phi'$ and $\Phi''$ respectively.  If $\Delta \subset
\mathbb Z^d$, write $\overline{\Delta} = k\Delta + \Lambda$.  It
follows from the definition of free energy that:
$$|\overline{\Lambda}_n|^{-1} FE_{\overline{\Lambda}_n}(\mu) =
|\overline{\Lambda}_n|^{-1} \mathcal H_{\mathcal
F_{\overline{\Lambda}_n}^{\tau}}(\mu,\exp\[-H^o_{\overline{\Lambda}_n}\]
\lambda^{|\overline{\Lambda}_n|-1}) = $$
$$|\overline{\Lambda}_n|^{-1} \mathcal H_{\mathcal
F_{\overline{\Lambda}_n}^{\tau}}(\mu,\exp\[-H^2_{\overline{\Lambda}_n}\]
\lambda^{|\overline{\Lambda}_n|-1}) + |\overline{\Lambda}_n|^{-1}
\mu (H^1_{\overline{\Lambda}_n}).$$ As $n$ tends to infinity, the
left hand side tends to $SFE(\mu)$ and the second term on the right
hand side tends to $\mu(\Phi')$.  So it is enough to check that:
$$\lim_{n\rightarrow \infty} |\overline{\Lambda}_n|^{-1} \mathcal H_{\mathcal
F_{\overline{\Lambda}_n}^{\tau}}(\mu,\exp\[-H^2_{\overline{\Lambda}_n}\]
\lambda^{|\overline{\Lambda}_n|-1}) = \mathcal H_{\mathcal
F^{\tau}_{\overline{\Gamma}(y)}}(\mu , \mu
\overline{\lambda}^d_y).$$

We now prove this fact using an argument similar to the one in the
proof of \ref{specificentropylexicographic} in \cite{G}.  Write
$\alpha_j = \mathcal H_{\mathcal F^{\tau}_{\Gamma(y)}}(\mu , \mu
\overline{\lambda}^j_y)$.  By Lemma \ref{conditionalentropy},
$\alpha_d$ is the expected conditional entropy of
$\overline{\phi}(y)$ (with respect to $\overline{\lambda}^j_y$) {\it
given} the differences $\phi(x_1) - \phi(x_2)$ for all $x_1, x_2 \in
\overline{\Gamma}(y)$.  We will now express $\mathcal H_{\mathcal
F_{\overline{\Lambda}_n}^{\tau}}(\mu,\exp\[-H^1_{\overline{\Lambda}_n}\]
\lambda^{|\overline{\Lambda}_n|-1})$ as a sum of $|\Lambda_n|$
expected conditional entropies, all of which (except for a boundary
set) are between zero and $\alpha_d$, and most of which are very
close to $\alpha = \alpha_d$.

Now, define $\alpha_{\Lambda_n}(y) = \mathcal H_{\mathcal
F^{\tau}_{\overline{\Lambda}_n \cap \overline{\Gamma}(y)}}(\mu, \mu
\overline{\lambda}_y^j)$ where $j = j_y = \sup \{ i | y - e_i \in
\Lambda_n \}$. (Separately define $\alpha_{\Lambda_n}(0) = \mathcal
H_{\mathcal F^{\tau}_{k\Gamma^*(0)}}(\mu, \mu
\lambda^{|\Lambda|-1})$.)  In words, we can think of each
$\alpha_{\Lambda_n}(y)$ as an expected conditional entropy of an
$n^d$-dimensional random variable $\psi(y)$ (given by the
$\phi(ky+i) - \phi(ky-e_{j_y})$) with respect to $\lambda_y^j$ {\it
given} the values of $\psi(x)$ for $x \prec y$ and $x \in
\Lambda_n$.  (Note that $\psi(0)$ is only an $n^d-1$ dimensional
variable.)  Let $\tilde{H}_{\overline{\Lambda}_n}$ be the sum of the
energies that give the Radon-Nikodym derivatives for these measures
$\overline{\lambda}_y^j$: that is, $H_{\overline{\Lambda}_n}) =
H_{\Lambda} + \sum_{y \in \Lambda_n, y \not = 0} H_{[y+\Lambda] \cup
\{j_y\}}$.  Now, repeated applications of Lemma
\ref{conditionalentropy} to this sequence of expected conditional
entropies yields
 $$\mathcal H_{\mathcal
F_{\Lambda}^{\tau}}(\mu,\exp\[-\tilde{H}_{\overline{\Lambda}_n}\]
\lambda^{|\overline{\Lambda}_n|-1}) = \sum_{y \in \Lambda_n}
\alpha_{\Lambda_n}(y)$$  Now, $\tilde{H}_{\overline{\Lambda}_n}$ and
$H^1_{\overline{\Lambda}_n}$ differ only by the inclusion and
exclusions of energy terms coming from the boundary of $\Lambda_n$.
Unless one of these expected energy terms is infinite (in which case
it is clear that $SFE(\mu) = \infty$ and $\alpha^j=\infty$ for some
$j$), $\mathcal L$ -invariance implies that
 $$\mathcal H_{\mathcal
F_{\Lambda}^{\tau}}(\mu,\exp\[-H^1_{\overline{\Lambda}_n}\]
\lambda^{|\overline{\Lambda}_n|-1}) = \sum_{y \in \Lambda_n}
\alpha_{\Lambda_n}(y) + o(n^d).$$

If $\alpha^j = \infty$ for some $j$, then we will assume that the
coordinate axes were labelled in such a way that $j=d$, in which
case the above expression still holds.  (We will satisfy the lemma
statement in this case by showing that $SFE(\mu) = \infty$.)  By
Lemma \ref{relativeentropylemma}, we have $\alpha_{\Lambda_n}(y)
\leq \alpha$ whenever both $y \in \Lambda$ and $y - e_d \in \Lambda$
(i.e., except for boundary terms). By Lemma
\ref{increasingalgebraentropy}, for any $\epsilon$, there exists a
finite subset $\Delta \subset \Gamma^*(0)$ for which $\mathcal
H_{\mathcal F^{\tau}_{k\Delta}}(\mu, \mu \overline{\lambda}_y^d)
\geq \alpha - \epsilon$.  Then Lemma \ref{relativeentropylemma}
implies that whenever $\alpha_{\Lambda_n}(y) \geq \alpha - \epsilon$
whenever $y + \Delta \subset \Lambda_n$.  Letting $n$ tend to
infinity, we conclude that $$\alpha - \epsilon \leq
\lim_{n\rightarrow \infty} |\overline{\Lambda}_n|^{-1} \mathcal
H_{\mathcal
F_{\overline{\Lambda}_n}^{\tau}}(\mu,\exp\[-H^2_{\overline{\Lambda}_n}\]
\lambda^{|\overline{\Lambda}_n|-1}) \leq \alpha$$ Since $\epsilon$
was arbitrary, the expression is equal to $\alpha$. \qed \end{proof}
\begin{lem}$h(\overline{\mu}) = \mathcal H_{\mathcal
F^{\tau}_{\overline{\Gamma}(y)}}(\mu,\mu \overline{\lambda}_y)$
\end{lem}

\begin{proof} By our construction and Lemma \ref{specificentropylexicographic}, $h(\overline{\mu})
= \mathcal H_{\mathcal A_y}(\mu,\mu \overline{\lambda}_y)$, where
$\mathcal A_y$ is the $\sigma$-algebra on $\Omega$ formed by pulling
back $\overline{\mathcal F}_{\Gamma(y)}$ via the map $\phi \mapsto
\overline{\phi}$; in other words, $\mathcal A_y$ is the smallest
$\sigma$-algebra in which $\phi(x_1) - \phi(x_2)$ is measurable for
every pair $(x_1,x_2)$ contained in $\Lambda \cup \{-e_d \} \cup
\Lambda + kj$ for some $y \prec 0$.  Define $f_{\phi}(x_1,x_2) =
\phi(x_1) - \phi(x_2)$ for all such pairs.  So now, we need only
show that $\mathcal H_{\mathcal A_y}(\mu,\mu \overline{\lambda}_y) =
\mathcal H_{\mathcal F^{\tau}_{\overline{\Gamma}(y)}}(\mu,\mu
\overline{\lambda}_y)$. The two $\sigma$-algebras are slightly
different; in the former, functions of $\phi(x_1) - \phi(x_2)$ are
measurable only when $x_1, x_2 \in \overline{\Gamma}(y)$ and $x_1$
and $x_2$ are in the same ``row''---i.e., $x_1 = y_1+z_1$ and $x_2 =
y_2+z_2$ where each $z_i$ is in $\Lambda$ and the $y_i \in k
[\Gamma(y)]$ are vectors that agree on all but the $d$th coordinate.
In particular, for each $y = (y_1,\ldots, y_d) \in \Gamma_0$, define
an ``average row distance'' function that is measurable with respect
to $\mathcal F^{\tau}_{\Gamma(y)}$ but not $\mathcal A_y$:
$$g_{\phi}(y) = \lim_{j \rightarrow \infty} j^{-1} \sum_{i=1}^j
\phi(y_1,y_2,\ldots, y_{d-1}, i) - \phi(0,0,\ldots, 0,i).$$ Since we
may assume that $\mu$ is shift invariant with finite slope
(otherwise it is clear that $SFE(\mu) = \infty$ and $\mu(\Phi') =
\infty$, so we must have $SFE(\mu) = h(\overline{\mu}) + \mu(\Phi')$
in that case), it follows from the ergodic theorem that $g_\phi$ is
well-defined $\mu$-almost surely.  Now, the lemma is equivalent to
the statement that the expected conditional entropy of
$\overline{\phi}(0)$ given $f_{\phi}$ is the same as the expected
conditional entropy given $f_{\phi}$ {\it and} $g_{\phi}$. It
suffices to show that given $f_{\phi}$, the regular conditional
probability distributions for the random variables
$\overline{\phi}(0)$ and $g_{\phi}$ are almost surely independent.

Suppose otherwise.  Then there would be, with positive probability,
some event $A$ depending only on $g_{\phi}$ and some $\epsilon > 0$
such that $|\mu(A |f_{\phi}, \overline{\phi}(0)) - \mu (A |
f_{\phi})| > \epsilon$ (where here $\mu(A|*)$ denotes the regular
conditional distribution---given $*$---integrated over $A$). By the
(one-dimensional) ergodic theorem, this statement would also have to
be true with positive probability for a positive fraction of the
shifted functions $\theta_{j e_d} \phi$ with $j \in \mathbb Z$; but
since the probability of $A$ with respect to an the increasing
sequence (in $j$) of subalgebras $\mathcal A_{y + je_d}$ is a
martingale, this contradicts the martingale convergence theorem.
\qed \end{proof}

\chapter{Summary of notations}

\begin{longtable}{p{1.1in}p{5in}}
\\ \multicolumn{2}{c}{\bf Vectors and lattices} \\
\hline

$d$ & dimension of the configuration lattice $\mathbb Z^d$ \\
$E$ & target space for random functions (usually $\mathbb R^m$ or $\mathbb Z^m$) \\
$m$ & dimension of $E$ \\
$e_1,\ldots, e_d$ & standard basis vectors in domain space $\mathbb Z^d$ \\
$e^1, \ldots, e^m$ & standard basis vectors in range space $E$ \\
$\mathcal E$ & $\sigma$-algebra on $E$ \\
$\lambda$ & underlying measure (usually counting or Lebesgue) on $(E, \mathcal E)$ \\
$\Lambda, \Delta$ & subsets of $\mathbb Z^d$ \\
$\Lambda \subset \subset \mathbb Z^d$ & $\Lambda \subset \mathbb
Z^d$ and $|\Lambda| < \infty$ \\
$\mathcal L$ & rank-$d$ sublattice of $\mathbb Z^d$ \\
$\Lambda_n$ & $[0,kn - 1]^d \subset \mathbb Z^d$ where $k$ is chosen
so that $k \mathbb Z^d \subset \mathcal L$ \\
$\Theta$ & group of translations of $\mathbb Z^d$ by elements of
$\mathcal L$ \\
$\tau$ & group of translations of $E$ \\
$\overline{\mathcal L}$ & dual lattice of $\mathcal L$ \\

\\ \multicolumn{2}{c}{\bf Configuration space and $\sigma$-algebras} \\
\hline

$\Omega$ & configuration space, set of functions from $\mathbb Z^d$ to $E$ \\
$\phi, \psi$ & elements of $\Omega$ \\
$\theta_x(\phi)$ & translation of $\phi \in \Omega$ by $x \in
\mathbb Z^d$, i.e., $(\theta_x
\phi)(y) = \phi(x+y)$ \\
$\mathcal F_\Lambda$ & smallest $\sigma$-algebra on $\Omega$ in which values on $\Lambda$ are $\mathcal E$-measurable\\
$\mathcal F$ & $\sigma$-algebra generated by $\mathcal F_\Lambda$, $\Lambda \subset \subset \mathbb Z^d$ \\
$\mathcal T_\Lambda$ & $\mathcal F_{\mathbb Z^d \backslash
\Lambda}$\\
$\mathcal T$ & tail $\sigma$-algebra, defined as $\cap_{\Lambda \subset \subset \mathbb Z^d} \mathcal T_\Lambda$ \\
$\mathcal F^\tau$ & $\sigma$-algebra of $\tau$-invariant elements of
$\mathcal F$ \\
$\mathcal F_\Lambda^\tau$ & $\mathcal F_\Lambda \cap \mathcal
F^\tau$\\
$\mathcal T_\Lambda^\tau$ & $\mathcal T_\Lambda \cap \mathcal F^\tau$ \\
$\mathcal T^\tau$ & $\mathcal T \cap \mathcal F^\tau$ \\




\\ \multicolumn{2}{c}{\bf Gibbs Potentials and Hamiltonians} \\
\hline

$\Phi, \Psi$ & Gibbs potentials; $\Phi = \{ \Phi_{\Lambda}:\Lambda
\subset \subset \mathbb Z^d \}$, where each $\Phi_\Lambda : \Omega
\rightarrow
\mathbb R \cup \{\infty\}$ is $\mathcal F_\Lambda$ measurable \\

$H^\Phi_\Lambda(\phi)$, $H_\Lambda(\phi)$ & Hamiltonian in
$\Lambda$, defined as
$\sum_{\Delta \cap \Lambda \not = \emptyset} \Phi_\Delta(\phi)$ \\
$H^o_\Lambda(\phi)$ & interior Hamiltonian in $\Lambda$, defined as
$\sum_{\Delta \subset \Lambda} \Phi_\Delta(\phi)$ \\
$Z^\Phi_{\Lambda}(\phi)$, $Z_{\Lambda}(\phi)$ & partition function
on $\Lambda$ with boundary condition $\phi$ on $\mathbb Z^d
\backslash \Lambda$, i.e.,
$\int \prod_{x \in \Lambda} d \phi(x) \exp\[-H_{\Lambda}(\phi)\]$ \\
$\gamma^\Phi_\Lambda$, $\gamma_\Lambda$ & transition kernel
corresponding to a Gibbs rerandomization on $\Lambda$, i.e.,
$\gamma^{\Phi}_{\Lambda}(A, \phi) = Z_{\Lambda}(\phi)^{-1} \int
\prod_{x \in \Lambda} d \phi(x) \exp\[-H_{\Lambda}(\phi)\] 1_A(\phi).$ \\
$Z^o_{\Lambda}$ & free boundary partition function on $\Lambda$ with
respect to $\Phi$, i.e., integral of $e^{-H^o_{\Lambda}}$ over
entire
space $E^{|\Lambda|-1}$ of functions (defined up to additive constant) on $\Lambda$\\
$W(\Lambda)$ & $-\text{log} Z^o_{\Lambda}$\\
$\mathcal D_{\mathcal L}$ & space of positive $\mathcal L \times
\tau$-invariant potentials for which $W(e)$ is finite for every edge
$e$ \\

\\ \multicolumn{2}{c}{\bf Nearest-neighbor Gibbs Potentials} \\
\hline
$V$, $V_{x,y}$ & nearest-neighbor difference potential \\
$\overline{V}$ & wedge-normalization of $V$, defined as $V(\eta) -
\log g(F(\eta))$ where $g(\eta) = 2-4|\eta - \frac{1}{2}|$ and
$F(\eta) = \frac{
\int_{-\infty}^{\eta}e^{-V(\eta)}d\eta}{\int_{-\infty}^{\infty}
e^{-V(\eta)}d\eta}$
\\
$\eta$ & nearest neighbor height difference (input to $V$) \\
SAP & simply attractive potential (a.k.a., convex nearest-neighbor,
periodic difference potential) \\
ISAP & isotropic simply attractive potential\\
LSAP & Lipschitz simply attractive potential\\
$\Phi_V$ & ISAP in which $V_{x,y} = V$ for all adjacent pairs $x,y\in \mathbb Z^d$\\

\\ \multicolumn{2}{c}{\bf Spaces of probability measures on configuration space} \\
\hline

$\mathcal P (\Omega, \mathcal F)$ & set of probability measures
on $(\Omega, \mathcal F)$ \\
$\mathcal G(\Omega, \mathcal F), \mathcal G$ & set of Gibbs measures
on $(\Omega, \mathcal F)$, i.e., measures $\mu$ such that for all
$\Lambda \subset \subset Z^d$, $0 < Z_\Lambda(\phi) < \infty$
$\mu$-a.s.\ and $\mu \gamma_\Lambda = \mu$
\\
$\mathcal P (\Omega, \mathcal F^\tau)$ & set of probability measures
on $(\Omega, \mathcal F^\tau)$ \\
$\mathcal P_{\mathcal L} (\Omega, \mathcal F^\tau)$ & set of
$\mathcal L$-invariant probability measures on $(\Omega, \mathcal
F^\tau)$ \\
$\mathcal G(\Omega, \mathcal F^\tau), \mathcal G^\tau$ & set of
gradient Gibbs measures on $(\Omega, \mathcal F^\tau)$, i.e.,
measures $\mu$ such that for all $\Lambda \subset \subset Z^d$, $0 <
Z_\Lambda(\phi) < \infty$ $\mu$-a.s.\ and $\mu \gamma_\Lambda = \mu$
\\
$\mathcal P (\Omega, \mathcal F^\tau)$ & set of probability measures
on $(\Omega, \mathcal F^\tau)$ \\
$\mathcal G_{\mathcal L}(\Omega, \mathcal F^\tau)$ & set of
$\mathcal L$-invariant gradient Gibbs measures on $(\Omega, \mathcal F^\tau)$ \\
$\ex\mathcal P_{\mathcal L} (\Omega, \mathcal F^\tau)$ & set of
$\mathcal L$-ergodic probability measures on $(\Omega, \mathcal
F^\tau)$ \\
$\ex\mathcal G(\Omega, \mathcal F^{\tau})$ & set of extremal
gradient Gibbs
measures on $(\Omega, \mathcal F^\tau)$ \\
$\ex\mathcal G_{\mathcal L}(\Omega, \mathcal F^{\tau})$ & set of
$\mathcal L$-ergodic gradient Gibbs measures on $(\Omega, \mathcal F^\tau)$ \\
$\mu, \nu$ & measures, usually elements of $\mathcal P(\Omega, \mathcal F^\tau)$ \\


\\ \multicolumn{2}{c}{\bf Relative entropy, free energy, and specific free energy} \\
\hline

$\mathcal H(\mu, \nu)$ & relative entropy of $\mu$ with
respect to $\nu$ \\
$\mathcal H_\mathcal A(\mu, \nu)$ & relative entropy of $\mu$ with
respect to $\nu$ on sub $\sigma$-algebra $\mathcal A$ \\

$\mu_\Lambda$ & restriction of $\mu \in \mathcal P_{\mathcal
L}(\Omega, \mathcal F^{\tau})$ to $\mathcal F_{\Lambda}^\tau$
\\
$FE_{\Lambda}(\mu)$ & free energy of $\mu \in \mathcal P_{\mathcal
L}(\Omega, \mathcal F^{\tau})$ in $\Lambda$, i.e., $\mathcal H
(\mu_{\Lambda},e^{-H^o_{\Lambda}}\lambda^{|\Lambda|-1})$\\

$SFE_\Lambda(\mu)$ & specific free energy
of $\mu$ in $\Lambda$, i.e., $|\Lambda|^{-1} FE_{\Lambda}(\mu)$\\

$SFE(\mu)$ & specific free energy of $\mu$, i.e., $\lim_{n
\rightarrow
\infty} SFE_{\Lambda_n}(\mu)$\\


\\ \multicolumn{2}{c}{\bf Slopes and surface tension} \\
\hline

$S(\mu)$ & slope of $\mu$ (where $\mu \in \mathcal P_{\mathcal
L}(\Omega, \mathcal F^{\tau})$) \\
$u$ & slope variable (a linear function from $\mathbb R^d$ to $\mathbb R^m$) \\
$\sigma^\Phi, \sigma$ & surface tension, $\sigma(u) =
\inf \{SFE^\Phi(\mu): \mu \in \mathcal P_{\mathcal L} (\Omega, \mathcal F^\tau), S(\mu) = u \}$ \\
$P(\Phi)$ & pressure of $\Phi$, $P(\Phi) = \inf \{SFE^\Phi(\mu):\mu
\in \mathcal P_{\mathcal L} (\Omega, \mathcal F^\tau)\} =
\inf_{u \in \mathbb R^{m \times d}}\sigma(u)$ \\
$U_\Phi$ & interior of set of slopes $u$ with $\sigma(u) < \infty$
\\


\\ \multicolumn{2}{c}{\bf Extremal and ergodic decompositions} \\
\hline

$e_A$ & evaluation map $\mu \rightarrow \mu_A$ \\
$e(\chi)$ & smallest $\sigma$ algebra on a subset $\chi$ of
$\mathcal P_{\mathcal L} (\Omega, \mathcal F)$ or $\mathcal
P_{\mathcal L} (\Omega, \mathcal F^\tau)$ that makes $e_A$ measurable for each $A \in \mathcal F$ \\
$w_\mu$ & extremal decomposition of $\mu \in \mathcal G^\tau$, a
measure on $(\ex\mathcal G^\tau, e(\ex\mathcal G^\tau))$ {\em or}
ergodic decomposition of $\mu \in \mathcal P_{\mathcal L} (\Omega,
\mathcal F^\tau)$, a measure on $$(\ex \mathcal P_{\mathcal L}
(\Omega, \mathcal F^\tau),
e(\ex \mathcal P_{\mathcal L} (\Omega, \mathcal F^\tau)))$$ \\

$\pi^\phi$ & limit of Gibbs rerandomizations $\gamma_{\Lambda_n}(\cdot | \phi)$ of $\phi$ on $\Lambda_n$ \\
$\pi^\phi_\mathcal L$ & shift-averaged limit of Gibbs rerandomizations of $\phi$ \\

\\ \multicolumn{2}{c}{\bf Topologies on probability spaces} \\
\hline

$\mathcal P(X, \mathcal X)$ & space of probability measures on a
measure space $(X, \mathcal X)$ \\
$\tau$-topology & smallest topology on $\mathcal P(X, \mathcal X)$
in which $\nu \mapsto \nu(A)$ is continuous for every $A \in \mathcal X$ \\

weak topology & smallest topology on $\mathcal P(X, \mathcal X)$ in
which $\nu \mapsto \nu(f)$ is continuous for every bounded
continuous function $f$ on $X$ \\

$\mathcal A$ & topology of local convergence on $\mathcal P(\Omega,
\mathcal F^\tau)$, i.e., smallest topology in which the maps $\mu
\mapsto \mu(f)$ are continuous for every bounded function $f:\Omega
\rightarrow \mathbb R$ that is $\mathcal
F_{\Lambda}^{\tau}$-measurable for some $\Lambda \subset \subset
\mathbb Z^d$
 \\

$\mathcal B$ & basis for $\mathcal A$ given by set of finite
intersections of sets of the form $\{\mu: \mu(F) < \epsilon \}$,
where $F: \Omega \rightarrow \mathbb R$ is bounded and $\mathcal
F^\tau_\Lambda$-measurable for some
$\Lambda \subset \subset \mathbb Z^d$.\\

\\ \multicolumn{2}{c}{\bf Lattice approximations to continuous domains} \\
\hline

$D_n$ & subset of $\mathbb Z^d$ that approximates $n D$ (e.g., $n D \cap \mathbb Z^d$) \\
$\hat D_n$ & simplex domain derived from $D_n$ \\
$\tilde D_n$ & $\frac{1}{n} \hat D_n$ \\
$\phi_n$ & a function from $D_n$ to $E$ \\
$\hat \phi_n$ & piecewise linear interpolation of $\phi_n$ to
$\hat{D}_n$. \\
$\tilde{\phi}_n$ & rescaling of $\hat \phi_n$ to $\tilde{D}_n$ given
by $\tilde{\phi}(\eta) = \frac{1}{n} \hat{\phi}(n\eta)$ \\

\\ \multicolumn{2}{c}{\bf Young functions} \\
\hline $A$ & Young function, i.e., a convex, even function $A:
\mathbb R \mapsto \mathbb R^+\cup \{\infty\}$ for which $A(0)=0$,
$A$ is finite on some open interval, and $A$ is not
identically zero \\
$A_d$ & Sobolev conjugate of $A$ in $d$ dimensions\\
$A^*$ & sub-conjugate of $A$, i.e., any Young function increasing
essentially more slowly near infinity than $A_d$ \\
$A(v), v \in \mathbb R^d$ & $\sum_{i=1}^d A(v_i)$\\

\\ \multicolumn{2}{c}{\bf Orlicz-Sobolev spaces} \\
\hline
$D$ & a domain in $\mathbb R^d$, usually a member of $\mathbb G(\frac{d-1}{d})$ (defined below) \\
$|D|$ & Lebesgue measure of $D$ \\
$\alpha$ & a multi-index $\alpha = (\alpha_1, \ldots, \alpha_d)$
with $0 \leq \alpha_i \in \mathbb Z^d$ \\
$D^{\alpha}$ & distributional $\alpha$ derivative \\

$||f||_{A,D}$, $||f||_A$ &
$\inf \{ k| \int_D A(\frac{f(\eta)}{k})d\eta \leq 1 \}$\\
$L^A(D)$, $L^A$ & Orlicz space $\{f : ||f||_{A,D} < \infty \}$\\

$W^{j,A}(D)$ & Orlicz-Sobolev space $\{ f \in L^A(D): D^{\alpha}f
\in L^A(D) \text { for } 0 \leq |\alpha| \leq j \}.$ \\

$||\nabla f||_{A,D}$ & $\inf \{ k| \int_D A\[\frac{ \nabla f(\eta)
}{k}\] d\eta \leq 1 \}.$\\

$P(E;D)$ & perimeter of $E$ relative to $D$, i.e., total variation
over $D$ of the gradient of the characteristic function of $E$ \\
$\mathbb G(z)$ & set of bounded domains $\{ D \subset \mathbb R^n
\}$ for which there exists a constant $C$ such that $[\min \{|E|,
|D-E| \}]^{z} \leq CP(E;D)$ for all Lebesgue measurable subsets
$E$ of $D$. \\

$L^A_0(D)$ & $L^{A}(D)\cup_{n=1}^{\infty} L^{A}(\tilde{D}_n)$\\

\\ \multicolumn{2}{c}{\bf Empirical measures and large deviations} \\
\hline

$L_n(\phi)$ & $|\Lambda_n \cap \mathcal L|^{-1} \sum_{x \in
\Lambda_n \cap \mathcal L} \delta_{\theta_x \phi}$, called the
empirical measure of $\phi$ on $\Lambda_n$, a member of $\mathcal
P(\Omega, \mathcal F)$ \\
$B_n$ & $\{ \phi : L_n(\phi) \in B \}$ where $B \in \mathcal B$\\

$C^u_n$ & $\{\phi : |[\phi(x) - \phi(x_0)] - [\phi_u(x) -
\phi_u(x_0)]| \leq \epsilon \text{ for all } x \in \partial
\Lambda_n \backslash \{x_0\}$, where $\epsilon$ is fixed independently of $n$ and $\phi_u$
is plane of slope $u$\\

$PBL^u_B(\mu)$ & $\limsup_{n \rightarrow \infty} - |\Lambda_n|^{-1}
\log
\[ \int 1_{C^u_n \cap B_n} e^{-H^o_{\Lambda}(\phi)} \prod_{x \in
\Lambda_n \backslash \{x_0\} }d\phi(x) \]$\\

$PBL^u(\mu)$&$\sup_{B \ni \mu, B \in \mathcal B} PBL^u_B(\mu)$\\
$PBL(\mu)$&$PBL^{S(\mu)}(\mu)$\\
$FBL_B(\mu)$ & $\liminf_{n \rightarrow \infty} -|\Lambda_n|^{-1}
\log
\[ \int 1_{B_n} e^{-H^o_{\Lambda}(\phi)} \prod_{x \in \Lambda_n
\backslash \{x_0\} }d\phi(x) \]$\\
$FBL(\mu)$ & $\sup_{B \ni \mu, B \in \mathcal B} FBL_B(\mu)$\\

$R_{\phi_n, n}$ &  $\int_{D} \delta_{(x, \theta_{\lfloor nx \rfloor}
\phi_n)}dx$, called the empirical profile measure of $\phi_n$ and
defined as an element of $\mathcal P(D \times \Omega)$, with
$\sigma$-algebra understood to be
Lebesgue measure times $\mathcal F^\tau$ \\
$\mu_n$ & Gibbs measure $\mu_n$ on $(\Omega, \mathcal F_{D_n}^\tau)$
defined by Gibbs potential $H^o_{D_n}$\\
$\rho_n$ & measure on $\mathcal P(D \times \Omega) \times
L^{\overline V^*}_0$ induced by $\mu_n$
and the map $\phi_n \rightarrow (R_{\phi_n, n}, \tilde \phi_n)$ \\
$S(\mu(D',\cdot))$ & slope of the probability measure $\mu(D',
\cdot)/\mu(D' \times \Omega)$, for $\mu \in \mathcal P(D \times \Omega)$ \\
$I(\mu,f)$ & rate function of large deviations principle satisfied
by $\rho_n$, given by $$\begin{cases}
        SFE (\mu(D, \cdot)) - P(\Phi)  &
        \text{$\mu (\cdot, \Omega)$ is Lebesgue measure on $D$ and} \\ & \text{$\mu (D', \cdot)$ is
$\mathcal L$-invariant when $|D'|>0$ and}\\
& S(\mu(x,\cdot)) = \nabla f(x)\\
        \infty & \text{otherwise} \\
        \end{cases} $$ \\
$\mathcal X$ & the topology on $\mathcal P(D \times \Omega) \times
L^{\overline{V}^*}_0(D)$ for the large deviations principle, given
by the product of (on the first coordinate) the smallest topology in
which $\mu\rightarrow \mu(D' \times f)$ is measurable for all
rectangular subsets $D'$ of $D$ and bounded cylinder functions $f$
and (on the second coordinate) the $L^{\overline{V}^*}_0(D)$
topology \\
\\ \multicolumn{2}{c}{\bf Cluster swapping and triplets} \\
\hline

$\mathbb E^d$ & set of edges of the lattice $\mathbb Z^d$  \\
$\Sigma$ & $[0,\infty)^{\mathbb E^d}$ \\
$\overline{\Omega}$ & $\Omega \times \Omega \times \Sigma$\\
$\overline{\mathcal F}^{\tau}$ & $\sigma$-algebra generated by
$\mathcal F^{\tau} \times \mathcal F^{\tau}$ times the product
topology on $\Sigma$ \\
$\Phi_{\Lambda} (\phi_1, \phi_2, r)$ & $\Phi_{\Lambda}(\phi_1) +
\Phi_{\Lambda}(\phi_2) + \sum_e r(e)$\\
$R(\phi_1, \phi_2, r)$ & cluster swapping map defined on $\Omega \times \Omega \times \Sigma$\\

$h$ & height offset variable for $\mu \in \mathcal P_{\mathcal
L}(\Omega, \mathcal F^{\tau})$, i.e., a function tail-measurable,
$\mu$-a.s.\ finite function $h:\Omega \mapsto \mathbb R \cup \{
\infty \}$ such that $h(\phi+c) = h(\phi)+c$ for all $\phi \in
\Omega, c \in E$ and $\mu$-a.s.\ $h(\phi) = h(\theta_v \phi) +
(u,v)$ when $v \in \mathcal L$ and $u= S(\pi^\phi_\mathcal L)$\\

$h(\mu)$ & height offset spectrum of $\mu \in \mathcal P_{\mathcal
L}(\Omega, \mathcal F^{\tau})$, i.e., the law of $h(\phi)$ modulo
one, if $\phi$ is chosen from $\mu$, viewed as a measure on $[0,1)$
\\

$\mathcal S_c$ & set of all edges for which $(\phi_1 + c, \phi_2,
r)$ is swappable\\

$T^+_c$ & set of vertices $v$ in infinite clusters of $\mathcal S_c$
complement for which $\phi_2(v) > \phi_1(v) + c$ throughout the
cluster \\
$T^-_c$ & set of vertices $v$ in infinite clusters of $\mathcal S_c$
complement
for which $\phi_2(v) < \phi_1(v) + c$ throughout the cluster \\

$B^+$ & $\inf \{c: \text{ $T^+_c$ is empty} \}$\\
$B^-$ & $\sup \{c: \text{ $T^+_-$ is empty} \}$\\

$\mu_u$ & minimal gradient phase of slope $u$ (uniquely defined, under some conditions)\\

$\mu_{u,a}$ & extremal Gibbs measures, such that $\mu_{u,a}$ a.s.\
$h(\phi) = a$ and $\mu_u$ is the weighted average of $\mu_{u,a}$
where $a$ is chosen from $h(\mu)$ \\

\\ \multicolumn{2}{c}{\bf Random subsets of $\mathbb Z^2$ (Chapter \ref{discretegibbschapter} only)} \\
\hline

$\Omega_{\Gamma}$ & set of all of subsets $\Gamma$ of $\mathbb Z^2$
\\
$\mathcal F_{\Gamma}$ & product $\sigma$-algebra on
$\Omega_{\Gamma}$\\

$P_{\Gamma}$ & a single infinite non-self-intersecting path forming
the boundary of $\Gamma$ (when such a path exists) \\

$A_{\infty}$ & the event that $\Gamma$ and $\mathbb Z^2 \backslash
\Gamma$ are infinite connected sets \\

$A_{\emptyset}$ & the event that $\Gamma = \emptyset$\\

$A_k$ & event that $A_{\infty}$ occurs and that both $\Gamma \cap
\Lambda_k$ and $(\mathbb Z^2 \backslash A_k) \cap \Lambda_k$ are
non-empty \\

$\Lambda_k$ & assumed in this section to be shifted to be centered
at the origin---i.e., $\Lambda_k = [\lfloor -k/2 \rfloor , \lfloor
k/2-1 \rfloor]^2 \subset \mathbb Z^2$ \\

$B_k$  & event that there exists a path---consisting entirely of
elements in $\Gamma \backslash \Lambda_k$---which encircles the set
$\Lambda_k$ \\

$\Delta(v)$ & ``shifted box'' $nv + \Lambda_k$ \\
$\Delta$ & $\cup_{v \in \mathbb Z^2} \Delta(v)$ \\
$\overline{\Delta}(v)$ & $nv + \Lambda_{k+1}$ \\

$\tilde{\Delta}(v)$ & outer band of square faces around
$\Delta(v)$---i.e., the set of square faces of $\mathbb Z^d$ that
are incident to at least one vertex of $\Delta(v)$ and at least one
vertex of $\overline{\Delta}(v) \backslash \Delta(v)$ \\
$\tilde{\Delta}$ & $\prod_{v \in\mathbb Z^d} \tilde{\Delta}(v)$\\

$A(v)$ & event that $P_{\Gamma}$ hits $\tilde{\Delta}_v$, and in
between the first and last times $P_{\Gamma}$ hits
$\tilde{\Delta}(v)$, $P_{\Gamma}$ hits no square which is fewer
steps away from a $\tilde{\Delta}(w)$, with $w \not = v$, than it is
from $\tilde{\Delta}(v)$ \\

$v_+, v_-$ & given the event $A(v)$, vertices such that
$\tilde{\Delta}(v_-)$ is the last band that the path $P_{\Gamma}$
hits before the first time it hits $\tilde{\Delta}(v)$, and
$\tilde{\Delta}(v_+)$ is the first band that the path $P_{\Gamma}$
hits after the last time it hits $\tilde{\Delta}(v)$ \\

$C(v,w)$ & event that some vertex incident to $\overline{\Delta}(v)$
and some vertex incident to $\overline{\Delta}(w)$ are in the same
connected component of $\Gamma \backslash \Delta$ \\

$\overline{p}$ & continuous interpolation of discrete path $p$\\

$C_q(v,w)$ & event that there exists a path $p$ in $\Gamma
\backslash \Delta$, connecting $\overline{\Delta}(v)$ and
$\overline{\Delta}(w)$ for which $\overline{p}$ is homotopically
equivalent to $q$ \\
$C^+_q(v,w)$ & event that $w = v_+$ and the path $P_{\Gamma}$ from
$\overline{\Delta}(v)$ to $\overline{\Delta}(w)$ is homotopic to $q$
\\
$C^-_q(v,w)$ & defined analogously using $v_+$ instead of $v_-$ \\

$B(v)$ & event that there exists a cycle in $\Gamma \backslash
\Delta$ which disconnects $\Delta(v)$ from infinity \\

$W$ & closed countably punctured plane $\mathbb R^2 \backslash
\mathbb Z^2 \cup \mathbb Z^2 \times S^1$, a space homeomorphic to
$\mathbb R^2 \backslash [\mathbb Z^2 + D$ where $D$
is any disc of radius less than $1/2$ \\

$a = (a_1,a_2)$ & fixed point in $\mathbb R^2 \backslash \mathbb
Z^2$ with irrational coordinates\\
$r_x$ & closed line segment from $a$ to the point $x$ (including
its limit point $x \times \arg(a-x)$) \\
$u_x$ & portion of same ray (from $a$ through $x$ to infinity) which
lies between $x$ and $\infty$, together with its limit point $x
\times \arg(x-a)$ \\

$\overline{x}$ & homotopy class of a path which follows $r_x$ from
$a$ towards $x$, then makes a counterclockwise loop around $x$, and
then returns to $a$ along $r_x$ \\

$P_{x,y}$ & set of continuous paths from $x$ to $y$ in $W$ \\

$r_x'$ & linear segment from the $a$ to $x \times \arg(a-x)$,
followed by a counter-clockwise arc from $x \times \arg(a-x)$ to $x
\times 0$ \\

$p'_+$, $p'_-$ & right, left derivative of $p$ \\

$P_p \subset P_{x,y}$ & homotopy class \\
$\tilde{p}$ & taut version of path $p$\\
$\epsilon$ & $10^{-10000} \mu(A_{\infty})^{10000}$\\
$\delta$ & $10^{-1000} \mu(A_{\infty})^{1000}$\\
$\gamma$ & $10^{-100} \mu(A_{\infty})^{100}$\\
$\beta$ & $10^{-10} \mu(A_{\infty})^{10}$

\end{longtable}

    \bibliographystyle{abbrv}
    \bibliography{randomsurface}

\def\cprime{$'$}
\begin{thebibliography}{10}

\bibitem{A}
R.~A. Adams.
\newblock {\em Sobolev spaces}.
\newblock Academic Press [A subsidiary of Harcourt Brace Jovanovich,
  Publishers], New York-London, 1975.
\newblock Pure and Applied Mathematics, Vol. 65.

\bibitem{AMO}
R.~K. Ahuja, T.~L. Magnanti, and J.~B. Orlin.
\newblock {\em Network flows}.
\newblock Prentice Hall Inc., Englewood Cliffs, NJ, 1993.
\newblock Theory, algorithms, and applications.

\bibitem{BC}
T.~Baker and L.~Chayes.
\newblock On the unicity of discontinuous transitions in the two-dimensional
  {P}otts and {A}shkin-{T}eller models.
\newblock {\em J. Statist. Phys.}, 93(1-2):1--15, 1998.

\bibitem{Ba}
R.~J. Baxter.
\newblock {\em Exactly solved models in statistical mechanics}.
\newblock Academic Press Inc. [Harcourt Brace Jovanovich Publishers], London,
  1982.

\bibitem{BD}
G.~Ben~Arous and J.-D. Deuschel.
\newblock The construction of the {$(d+1)$}-dimensional {G}aussian droplet.
\newblock {\em Comm. Math. Phys.}, 179(2):467--488, 1996.

\bibitem{BG}
D.~Bertacchi and G.~Giacomin.
\newblock Enhanced interface repulsion from quenched hard-wall randomness.
\newblock {\em Probab. Theory Related Fields}, 124(4):487--516, 2002.

\bibitem{BGV}
T.~Bodineau, G.~Giacomin, and Y.~Velenik.
\newblock On entropic reduction of fluctuations.
\newblock {\em J. Statist. Phys.}, 102(5-6):1439--1445, 2001.

\bibitem{B}
E.~Bolthausen.
\newblock Random walk representations and entropic repulsion for gradient
  models.
\newblock In {\em Infinite dimensional stochastic analysis (Amsterdam, 1999)},
  volume~52 of {\em Verh. Afd. Natuurkd. 1. Reeks. K. Ned. Akad. Wet.}, pages
  55--83. R. Neth. Acad. Arts Sci., Amsterdam, 2000.

\bibitem{BDG}
E.~Bolthausen, J.-D. Deuschel, and G.~Giacomin.
\newblock Entropic repulsion and the maximum of the two-dimensional harmonic
  crystal.
\newblock {\em Ann. Probab.}, 29(4):1670--1692, 2001.

\bibitem{BEF}
J.~Bricmont, A.~El~Mellouki, and J.~Fr{\"o}hlich.
\newblock Random surfaces in statistical mechanics: roughening, rounding,
  wetting,{$\ldots\,$}.
\newblock {\em J. Statist. Phys.}, 42(5-6):743--798, 1986.

\bibitem{BK}
R.~M. Burton and M.~Keane.
\newblock Density and uniqueness in percolation.
\newblock {\em Comm. Math. Phys.}, 121(3):501--505, 1989.

\bibitem{Ch}
L.~Chayes.
\newblock Percolation and ferromagnetism on {${\bf Z}\sp 2$}: the {$q$}-state
  {P}otts cases.
\newblock {\em Stochastic Process. Appl.}, 65(2):209--216, 1996.

\bibitem{CM}
L.~Chayes, D.~McKellar, and B.~Winn.
\newblock Percolation and {G}ibbs states multiplicity for ferromagnetic
  {A}shkin-{T}eller models on {$\mathbb Z\sp 2$}.
\newblock {\em J. Phys. A}, 31(45):9055--9063, 1998.

\bibitem{C}
A.~Cianchi.
\newblock Continuity properties of functions from {O}rlicz-{S}obolev spaces and
  embedding theorems.
\newblock {\em Ann. Scuola Norm. Sup. Pisa Cl. Sci. (4)}, 23(3):575--608, 1996.

\bibitem{C1}
A.~Cianchi.
\newblock Boundedness of solutions to variational problems under general growth
  conditions.
\newblock {\em Comm. Partial Differential Equations}, 22(9-10):1629--1646,
  1997.

\bibitem{C2}
A.~Cianchi.
\newblock Some results in the theory of {O}rlicz spaces and applications to
  variational problems.
\newblock In {\em Nonlinear analysis, function spaces and applications, Vol. 6
  (Prague, 1998)}, pages 50--92. Acad. Sci. Czech Repub., Prague, 1999.

\bibitem{C3}
A.~Cianchi.
\newblock A fully anisotropic {S}obolev inequality.
\newblock {\em Pacific J. Math.}, 196(2):283--295, 2000.

\bibitem{CEP}
H.~Cohn, N.~Elkies, and J.~Propp.
\newblock Local statistics for random domino tilings of the {A}ztec diamond.
\newblock {\em Duke Math. J.}, 85(1):117--166, 1996.

\bibitem{CKP}
H.~Cohn, R.~Kenyon, and J.~Propp.
\newblock A variational principle for domino tilings.
\newblock {\em J. Amer. Math. Soc.}, 14(2):297--346 (electronic), 2001.

\bibitem{CL}
J.~H. Conway and J.~C. Lagarias.
\newblock Tiling with polyominoes and combinatorial group theory.
\newblock {\em J. Combin. Theory Ser. A}, 53(2):183--208, 1990.

\bibitem{D}
A.~J. de~Souza.
\newblock An extension operator for {O}rlicz-{S}obolev spaces.
\newblock {\em An. Acad. Brasil. Ci\^enc.}, 53(1):9--12, 1981.

\bibitem{DZ}
A.~Dembo and O.~Zeitouni.
\newblock {\em Large deviations techniques and applications}, volume~38 of {\em
  Applications of Mathematics (New York)}.
\newblock Springer-Verlag, New York, second edition, 1998.

\bibitem{DG2}
J.-D. Deuschel and G.~Giacomin.
\newblock Entropic repulsion for the free field: pathwise characterization in
  {$d\geq3$}.
\newblock {\em Comm. Math. Phys.}, 206(2):447--462, 1999.

\bibitem{DG}
J.-D. Deuschel and G.~Giacomin.
\newblock Entropic repulsion for massless fields.
\newblock {\em Stochastic Process. Appl.}, 89(2):333--354, 2000.

\bibitem{DGI}
J.-D. Deuschel, G.~Giacomin, and D.~Ioffe.
\newblock Large deviations and concentration properties for {$\nabla\phi$}
  interface models.
\newblock {\em Probab. Theory Related Fields}, 117(1):49--111, 2000.

\bibitem{DS}
J.-D. Deuschel and D.~W. Stroock.
\newblock {\em Large deviations}, volume 137 of {\em Pure and Applied
  Mathematics}.
\newblock Academic Press Inc., Boston, MA, 1989.

\bibitem{DKS}
R.~Dobrushin, R.~Koteck{\'y}, and S.~Shlosman.
\newblock {\em Wulff construction}, volume 104 of {\em Translations of
  Mathematical Monographs}.
\newblock American Mathematical Society, Providence, RI, 1992.
\newblock A global shape from local interaction, Translated from the Russian by
  the authors.

\bibitem{DV}
M.~D. Donsker and S.~R.~S. Varadhan.
\newblock Large deviations from a hydrodynamic scaling limit.
\newblock {\em Comm. Pure Appl. Math.}, 42(3):243--270, 1989.

\bibitem{DFJ}
B.~Durhuus, J.~Fr{\"o}hlich, and T.~J{\'o}nsson.
\newblock Self-avoiding and planar random surfaces on the lattice.
\newblock {\em Nuclear Phys. B}, 225(2, FS 9):185--203, 1983.

\bibitem{ES}
R.~G. Edwards and A.~D. Sokal.
\newblock Generalization of the {F}ortuin-{K}asteleyn-{S}wendsen-{W}ang
  representation and {M}onte {C}arlo algorithm.
\newblock {\em Phys. Rev. D (3)}, 38(6):2009--2012, 1988.

\bibitem{FK}
C.~M. Fortuin and P.~W. Kasteleyn.
\newblock On the random-cluster model. {I}. {I}ntroduction and relation to
  other models.
\newblock {\em Physica}, 57:536--564, 1972.

\bibitem{FKG}
C.~M. Fortuin, P.~W. Kasteleyn, and J.~Ginibre.
\newblock Correlation inequalities on some partially ordered sets.
\newblock {\em Comm. Math. Phys.}, 22:89--103, 1971.

\bibitem{FW}
M.~I. Freidlin and A.~D. Wentzell.
\newblock {\em Random perturbations of dynamical systems}, volume 260 of {\em
  Grundlehren der Mathematischen Wissenschaften [Fundamental Principles of
  Mathematical Sciences]}.
\newblock Springer-Verlag, New York, second edition, 1998.
\newblock Translated from the 1979 Russian original by Joseph Sz\"ucs.

\bibitem{FS1}
J.~Fr{\"o}hlich, C.-E. Pfister, and T.~Spencer.
\newblock On the statistical mechanics of surfaces.
\newblock In {\em Stochastic processes in quantum theory and statistical
  physics (Marseille, 1981)}, volume 173 of {\em Lecture Notes in Phys.}, pages
  169--199. Springer, Berlin, 1982.

\bibitem{FS2}
J.~Fr{\"o}hlich and T.~Spencer.
\newblock The {B}ere\v zinski\u\i -{K}osterlitz-{T}houless transition
  (energy-entropy arguments and renormalization in defect gases).
\newblock In {\em Scaling and self-similarity in physics (Bures-sur-Yvette,
  1981/1982)}, volume~7 of {\em Progr. Phys.}, pages 29--138. Birkh\"auser
  Boston, Boston, MA, 1983.

\bibitem{FZ}
J.~Fr{\"o}hlich and B.~Zegarli{\'n}ski.
\newblock The phase transition in the discrete {G}aussian chain with {$1/r\sp
  2$} interaction energy.
\newblock {\em J. Statist. Phys.}, 63(3-4):455--485, 1991.

\bibitem{F}
T.~Funaki.
\newblock Recent results on the {G}inzburg-{L}andau {$\nabla\phi$} interface
  model.
\newblock In {\em Hydrodynamic limits and related topics (Toronto, ON, 1998)},
  volume~27 of {\em Fields Inst. Commun.}, pages 71--81. Amer. Math. Soc.,
  Providence, RI, 2000.

\bibitem{FN}
T.~Funaki and T.~Nishikawa.
\newblock Large deviations for the {G}inzburg-{L}andau {$\nabla\phi$} interface
  model.
\newblock {\em Probab. Theory Related Fields}, 120(4):535--568, 2001.

\bibitem{FO}
T.~Funaki and S.~Olla.
\newblock Fluctuations for {$\nabla\phi$} interface model on a wall.
\newblock {\em Stochastic Process. Appl.}, 94(1):1--27, 2001.

\bibitem{FS}
T.~Funaki and H.~Spohn.
\newblock Motion by mean curvature from the {G}inzburg-{L}andau {$\nabla \phi$}
  interface model.
\newblock {\em Comm. Math. Phys.}, 185(1):1--36, 1997.

\bibitem{GKR}
A.~Gandolfi, M.~Keane, and L.~Russo.
\newblock On the uniqueness of the infinite occupied cluster in dependent
  two-dimensional site percolation.
\newblock {\em Ann. Probab.}, 16(3):1147--1157, 1988.

\bibitem{Gar}
R.~J. Gardner.
\newblock The {B}runn-{M}inkowski inequality.
\newblock {\em Bull. Amer. Math. Soc. (N.S.)}, 39(3):355--405 (electronic),
  2002.

\bibitem{G}
H.-O. Georgii.
\newblock {\em Gibbs measures and phase transitions}, volume~9 of {\em de
  Gruyter Studies in Mathematics}.
\newblock Walter de Gruyter \& Co., Berlin, 1988.

\bibitem{GH}
H.-O. Georgii and Y.~Higuchi.
\newblock Percolation and number of phases in the two-dimensional {I}sing
  model.
\newblock {\em J. Math. Phys.}, 41(3):1153--1169, 2000.
\newblock Probabilistic techniques in equilibrium and nonequilibrium
  statistical physics.

\bibitem{Gi4}
G.~Giacomin.
\newblock Anharmonic lattices, random walks and random interfaces.
\newblock {\em Recent research developments in statistical physics}, I:97--118,
  2000.

\bibitem{Gi}
G.~Giacomin.
\newblock Aspects of statistical mechanics of random surfaces.
\newblock {\em Preliminary (and partial) notes of the lectures given at IHP},
  2001.

\bibitem{Gi3}
G.~Giacomin.
\newblock Limit theorems for random interface models of {G}inzburg-{L}andau
  {$\nabla\phi$} type.
\newblock In {\em Stochastic partial differential equations and applications
  (Trento, 2002)}, volume 227 of {\em Lecture Notes in Pure and Appl. Math.},
  pages 235--253. Dekker, New York, 2002.

\bibitem{Gi2}
G.~Giacomin.
\newblock On stochastic domination in the {B}rascamp-{L}ieb framework.
\newblock {\em Math. Proc. Cambridge Philos. Soc.}, 134(3):507--514, 2003.

\bibitem{GOS}
G.~Giacomin, S.~Olla, and H.~Spohn.
\newblock Equilibrium fluctuations for {$\nabla\phi$} interface model.
\newblock {\em Ann. Probab.}, 29(3):1138--1172, 2001.

\bibitem{Go}
H.~H. Goldstine.
\newblock {\em A history of the calculus of variations from the 17th through
  the 19th century}, volume~5 of {\em Studies in the History of Mathematics and
  Physical Sciences}.
\newblock Springer-Verlag, New York, 1980.

\bibitem{Gr}
G.~Grimmett.
\newblock {\em Percolation}, volume 321 of {\em Grundlehren der Mathematischen
  Wissenschaften [Fundamental Principles of Mathematical Sciences]}.
\newblock Springer-Verlag, Berlin, second edition, 1999.

\bibitem{GPV}
M.~Z. Guo, G.~C. Papanicolaou, and S.~R.~S. Varadhan.
\newblock Nonlinear diffusion limit for a system with nearest neighbor
  interactions.
\newblock {\em Comm. Math. Phys.}, 118(1):31--59, 1988.

\bibitem{HS}
B.~Helffer and J.~Sj{\"o}strand.
\newblock On the correlation for {K}ac-like models in the convex case.
\newblock {\em J. Statist. Phys.}, 74(1-2):349--409, 1994.

\bibitem{HSn}
J.~Hershberger and J.~Snoeyink.
\newblock Computing minimum length paths of a given homotopy class.
\newblock {\em Comput. Geom.}, 4(2):63--97, 1994.

\bibitem{H}
M.~Huber.
\newblock A bounding chain for {S}wendsen-{W}ang.
\newblock {\em Random Structures Algorithms}, 22(1):43--59, 2003.

\bibitem{I}
D.~Ioffe.
\newblock Exact large deviation bounds up to {$T\sb c$} for the {I}sing model
  in two dimensions.
\newblock {\em Probab. Theory Related Fields}, 102(3):313--330, 1995.

\bibitem{K}
P.~Kastelyn.
\newblock The statistics of dimers on a lattice, i. the number of dimer
  arrangements on a quadratic lattice.
\newblock {\em Physica}, 27:1209--1225, 1965.

\bibitem{Kem}
J.~H.~B. Kemperman.
\newblock On the {FKG}-inequality for measures on a partially ordered space.
\newblock {\em Nederl. Akad. Wetensch. Proc. Ser. A 80=Indag. Math.},
  39(4):313--331, 1977.

\bibitem{Ke4}
R.~Kenyon.
\newblock Local statistics of lattice dimers.
\newblock {\em Ann. Inst. H. Poincar\'e Probab. Statist.}, 33(5):591--618,
  1997.

\bibitem{Ke}
R.~Kenyon.
\newblock Conformal invariance of domino tiling.
\newblock {\em Ann. Probab.}, 28(2):759--795, 2000.

\bibitem{Ke5}
R.~Kenyon.
\newblock Long-range properties of spanning trees.
\newblock {\em J. Math. Phys.}, 41(3):1338--1363, 2000.
\newblock Probabilistic techniques in equilibrium and nonequilibrium
  statistical physics.

\bibitem{Ke3}
R.~Kenyon.
\newblock The planar dimer model with boundary: a survey.
\newblock In {\em Directions in mathematical quasicrystals}, volume~13 of {\em
  CRM Monogr. Ser.}, pages 307--328. Amer. Math. Soc., Providence, RI, 2000.

\bibitem{Ke2}
R.~Kenyon.
\newblock Dominos and the {G}aussian free field.
\newblock {\em Ann. Probab.}, 29(3):1128--1137, 2001.

\bibitem{KOS}
R.~Kenyon, A.~Okounkov, and S.~Sheffield.
\newblock {Dimers and Amoebae}.

\bibitem{KW}
R.~W. Kenyon and D.~B. Wilson.
\newblock Critical resonance in the non-intersecting lattice path model.
\newblock {\em Probab. Theory Related Fields}, 130(3):289--318, 2004.

\bibitem{LM2}
J.~L. Lebowitz and A.~Martin-L{\"o}f.
\newblock On the uniqueness of the equilibrium state for {I}sing spin systems.
\newblock {\em Comm. Math. Phys.}, 25:276--282, 1972.

\bibitem{LM}
J.~L. Lebowitz and A.~E. Mazel.
\newblock On the uniqueness of {G}ibbs states in the {P}irogov-{S}inai theory.
\newblock {\em Comm. Math. Phys.}, 189(2):311--321, 1997.

\bibitem{LRS}
M.~Luby, D.~Randall, and A.~Sinclair.
\newblock Markov chain algorithms for planar lattice structures.
\newblock {\em SIAM J. Comput.}, 31(1):167--192 (electronic), 2001.

\bibitem{Ma}
W.~S. Massey.
\newblock {\em A basic course in algebraic topology}, volume 127 of {\em
  Graduate Texts in Mathematics}.
\newblock Springer-Verlag, New York, 1991.

\bibitem{M}
V.~G. Maz'ja.
\newblock {\em Sobolev spaces}.
\newblock Springer Series in Soviet Mathematics. Springer-Verlag, Berlin, 1985.
\newblock Translated from the Russian by T. O. Shaposhnikova.

\bibitem{MP}
V.~G. Maz{\cprime}ya and S.~V. Poborchi.
\newblock {\em Differentiable functions on bad domains}.
\newblock World Scientific Publishing Co. Inc., River Edge, NJ, 1997.

\bibitem{MMR}
A.~Messager, S.~Miracle-Sol{\'e}, and J.~Ruiz.
\newblock Convexity properties of the surface tension and equilibrium crystals.
\newblock {\em J. Statist. Phys.}, 67(3-4):449--470, 1992.

\bibitem{NS}
A.~Naddaf and T.~Spencer.
\newblock On homogenization and scaling limit of some gradient perturbations of
  a massless free field.
\newblock {\em Comm. Math. Phys.}, 183(1):55--84, 1997.

\bibitem{P}
A.~Pisztora.
\newblock Surface order large deviations for {I}sing, {P}otts and percolation
  models.
\newblock {\em Probab. Theory Related Fields}, 104(4):427--466, 1996.

\bibitem{Pr}
J.~Propp.
\newblock Generating random elements of finite distributive lattices.
\newblock {\em Electron. J. Combin.}, 4(2):Research Paper 15, approx. 12 pp.
  (electronic), 1997.
\newblock The Wilf Festschrift (Philadelphia, PA, 1996).

\bibitem{PW}
J.~G. Propp and D.~B. Wilson.
\newblock Exact sampling with coupled {M}arkov chains and applications to
  statistical mechanics.
\newblock In {\em Proceedings of the Seventh International Conference on Random
  Structures and Algorithms (Atlanta, GA, 1995)}, volume~9, pages 223--252,
  1996.

\bibitem{RR}
M.~Renardy and R.~C. Rogers.
\newblock {\em An introduction to partial differential equations}, volume~13 of
  {\em Texts in Applied Mathematics}.
\newblock Springer-Verlag, New York, second edition, 2004.

\bibitem{Ro}
J.~Rougemont.
\newblock Evaporation of droplets in the two-dimensional {G}inzburg-{L}andau
  equation.
\newblock {\em Phys. D}, 140(3-4):267--282, 2000.

\bibitem{R}
W.~Rudin.
\newblock {\em Real and complex analysis}.
\newblock McGraw-Hill Book Co., New York, third edition, 1987.

\bibitem{S2}
S.~Sheffield.
\newblock Ribbon tilings and multidimensional height functions.
\newblock {\em Trans. Amer. Math. Soc.}, 354(12):4789--4813 (electronic), 2002.

\bibitem{S}
S.~Sheffield.
\newblock {Gaussian free fields for mathematicians}.
\newblock {\em arXiv:math.PR/0312099}, 2004.

\bibitem{SS}
S.~Sheffield and O.~Schramm.
\newblock {Contour lines of the 2D Gaussian free field}.
\newblock {\em In preparation}, 2006.

\bibitem{Sp}
H.~Spohn.
\newblock Interface motion in models with stochastic dynamics.
\newblock {\em J. Statist. Phys.}, 71(5-6):1081--1132, 1993.

\bibitem{SW}
R.~H. Swendsen and J.-S. Wang.
\newblock Nonuniversal critical dynamics in {M}onte {C}arlo simulations.
\newblock {\em Physical Review Letters}, 58:86--88, 1987.

\bibitem{W}
G.~Wulff.
\newblock Zur frage der geschwindigkeit des wachsthums und der aufl\"osung der
  kristallfl\"achen.
\newblock {\em Zeit. Krystall. Min.}, 34:449--530, 1901.

\end{thebibliography}
\nocite{*}
\end{document}